\documentclass[11pt, 
]{amsart}
\nonstopmode
\usepackage{graphicx}
\usepackage{amssymb,amsmath,amsthm,amscd,accents}
\usepackage{mathrsfs}
\usepackage{lscape}
\usepackage{enumerate}
\usepackage{upgreek}
\usepackage[usenames,dvipsnames]{color}
\usepackage[colorlinks=true, pdfstartview=FitV,
 linkcolor=blue,citecolor=blue,urlcolor=blue]{hyperref}
\usepackage{tikz}
\usetikzlibrary{arrows.meta,arrows}
\usepackage{bm}
\usepackage{marginnote}
\usepackage{verbatim}
\usepackage{comment}
\usepackage{tcolorbox}

\usepackage[all]{xy}

\usepackage[shortlabels]{enumitem} 

\allowdisplaybreaks[3]

\newcommand{\nc}{\newcommand}
\numberwithin{equation}{section}

\nc{\hs}{\hspace*}
\nc{\ms}{\mspace}

\nc{\qR}[1]{\ttq_{\mspace{-2mu}\raisebox{-.8ex}{${\scriptstyle{#1}}$}}}

\theoremstyle{plain}
\newtheorem{lemma}{Lemma}[section]
\newtheorem{proposition}[lemma]{Proposition}
\newtheorem{theorem}[lemma]{Theorem}
\newtheorem{algorithm}[lemma]{Algorithm}

\newtheorem{corollary}[lemma]{Corollary}
\newtheorem{conjecture}{Conjecture}
\newtheorem{convention}{Convention}

\newtheorem{Conj}{Conjecture}

\theoremstyle{definition}
\newtheorem{remark}[lemma]{Remark}
\newtheorem{example}[lemma]{Example}

\newtheorem{definition}[lemma]{Definition}

\newcommand{\FM}[1]{F\bl #1 \br }
\newcommand{\FMq}[1]{F_q\bl #1 \br }
 
\newcommand{\Bq}{\mathbf{B}_q}
\newcommand{\At}{\mathbf{A}_t}
\newcommand{\Zq}{\Z[q^{\pm \frac{1}{2}}]}

\newcommand{\leN}{ \preceq_{\mspace{-2mu}\raisebox{-.5ex}{\scalebox{.6}{${\rm N}$}}} }
\newcommand{\lN}{ \prec_{\mspace{-2mu}\raisebox{-.5ex}{\scalebox{.6}{${\rm N}$}}} }

\renewcommand{\le}{\leqslant}
\renewcommand{\ge}{\geqslant}
\renewcommand{\preceq}{\preccurlyeq}

\newcommand{\Ker}{{\operatorname{Ker}}}
\newcommand{\st}{\mathop{\mbox{\normalsize$*$}}\limits}
\newcommand{\sta}{\mathop{\mbox{\normalsize$\star$}}\limits}

\newcommand{\seteq}{\mathbin{:=}}

\newcommand{\soplus}{\mathop{\mbox{\normalsize$\bigoplus$}}\limits}
\newcommand{\ev}{{\operatorname{ev}}}

\newcommand{\lxi}{{  {}^\xi  } \hspace{-.2ex}}
\newcommand{\lhs}{{  {}^{s}  }  \hspace{-.2ex} }
\newcommand{\lhbs}{{  {}^{(s)}  }  \hspace{-.2ex} }
\newcommand{\lhz}{{  {}^{0}  }  \hspace{-.2ex} }
\newcommand{\lh}[1]{{  {}^{#1}  }  \hspace{-.2ex} }
\newcommand{\luxi}{{  {}^\uxi  }}

\newcommand{\chit}{\chi_{\mspace{-2mu}\raisebox{-.5ex}{${\scriptstyle{t}}$}}}

\newcommand{\lcalX}{ {}^L \hspace{-.3ex}\calX}
\newcommand{\lScalMp}{ {}^\sfS \hspace{-.5ex}\calM_+}

\newcommand{\uip}{ u_{i,p}  }
\newcommand{\vip}{ v_{i,p}  }

\newcommand{\Irr}{{\rm Irr}}

\newcommand{\K}{\sfK}
\newcommand{\ex}{\mathrm{ex}}
\newcommand{\fr}{\mathrm{fr}}
\newcommand{\Kex}{\K_\ex}
\newcommand{\Kfr}{\K_\fr}

\newcommand{\seed}{\scrS}

\newcommand{\g}{\mathfrak{g}}

\newcommand{\n}{\mathfrak{n}}

\newcommand{\C}{\mathbb{C}}
\newcommand{\Q}{\mathbb{Q}}
\newcommand{\Z}{\mathbb{Z}\ms{1mu}}

\newcommand{\al}{{\ms{1mu}\alpha}}

\newcommand{\la}{\lambda}
\newcommand{\be}{{\ms{1mu}\beta}}
\newcommand{\ga}{\gamma}

\newcommand{\La}{\Lambda}

\newcommand{\Dynkin}{\triangle}
\newcommand{\ssq}[1]{  {}^{#1} \hspace{-.4ex} \overset{\gets}{\para} }
\newcommand{\rssq}[1]{  {}^{#1} \hspace{-.4ex} \overset{\to}{\para} }

\newcommand{\Dynkinf}{{\triangle \hspace{-1.37ex} {\raisebox{.2ex}{{\scalebox{.53}{${\mathbf f}$}}}} \hspace{1.3mm}  }}
\newcommand{\sDynkinf}{{\triangle \hspace{-.92ex} \raisebox{.18ex}{{\scalebox{.38}{${\mathbf f}$}}} \hspace{1mm} }}

\newcommand{\bDynkinf}{\bDynkin \hspace{-2.0ex} \raisebox{.2ex}{{\scalebox{.53}{  {\color{white} ${\mathbf f}$}}}} \ }

\newcommand{\bDynkin}{\mathop{\mathbin{\mbox{\large $\blacktriangle$}}}}

\newcommand{\lb}{(\hspace{-0.3ex}(}
\newcommand{\rb}{)\hspace{-0.3ex})}
\newcommand{\im}{\imath}
\newcommand{\jm}{\jmath}

\newcommand*\reced[1]{ \fontsize{8}{8}\selectfont \tikz[baseline=(char.base)]{
 \node[shape=rectangle,draw,inner sep=1.4pt] (char) {#1};} \fontsize{12}{12}\selectfont }

\newcommand*\circled[1]{ \fontsize{6}{6}\selectfont \tikz[baseline=(char.base)]{
  \node[shape=circle,draw,inner sep=0.4pt] (char) {#1};} \fontsize{12}{12}\selectfont }

\newcommand*\bcircled[1]{\fontsize{6}{6}\selectfont \tikz[baseline=(char.base)]{
    \node[shape=circle, fill= black, draw=black, text=white,  inner sep=0.4pt] (char) {#1};} \fontsize{12}{12}\selectfont}

\newcommand*\rcircled[1]{\fontsize{6}{6}\selectfont \tikz[baseline=(char.base)]{
    \node[shape=circle, fill= gray, draw=black, text=white, inner sep=0.4pt] (char) {#1};} \fontsize{12}{12}\selectfont}

\newcommand{\bii}{ \textbf{\textit{i}}}
\newcommand{\bij}{ \textbf{\textit{j}}}

\newcommand{\biQ}{ \textbf{\textit{Q}}}

\newcommand{\tA}{\widetilde{A}}
\newcommand{\tvee}{{\widetilde{\vee}}}
\newcommand{\tB}{\widetilde{B}}
\newcommand{\teta}{\widetilde{\eta}}

\newcommand{\tX}{\widetilde{X}}
\newcommand{\tY}{\widetilde{Y}}

\newcommand{\tS}{\widetilde{S}}
\newcommand{\ts}{\widetilde{s}}
\newcommand{\tbfm}{\widetilde{\bfm}}

\newcommand{\tpara}{ \overset{\gets}{\para}}
\newcommand{\opara}{ \overset{\to}{\para}}
\newcommand{\tw}{{\widetilde{w}}}
\newcommand{\tm}{\widetilde{m}}

\newcommand{\tuB}{\ms{1mu}\widetilde{\usfB}}
\newcommand{\tfb}{\widetilde{\sfb}}

\newcommand{\tp}{\widetilde{p}}

\newcommand{\oxi}{\overline{\xi}}
\newcommand{\oi}{{\overline{\im}}}
\newcommand{\ochi}{\overline{\chi}}
\newcommand{\osigma}{\overline{\sigma}}

\newcommand{\osfB}{\overline{\mathsf{B}}}
\newcommand{\usfB}{\underline{\mathsf{B}}}

\newcommand{\usfC}{\underline{\mathsf{C}}}
\newcommand{\uw}{{\underline{w}}}
\newcommand{\ubfm}{{\underline{\bfm}}}
\newcommand{\usfm}{{\underline{\sfm}}}
\newcommand{\ucalN}{\ms{1mu}\underline{\calN}}
\newcommand{\ucalM}{\ms{1mu}\underline{\calM}}

\newcommand{\utm}{\underline{\tm}}

\newcommand{\um}{\underline{m}}
\newcommand{\umn}[1]{\underline{m_{#1}}}

\newcommand{\uxi}{{\underline{\xi}}}

\newcommand{\frakU}{\mathfrak{U}}
\newcommand{\fraku}{\mathfrak{u}}
\newcommand{\frakv}{\mathfrak{v}}
\newcommand{\frakK}{\mathfrak{K}}

\newcommand{\frakB}{\mathfrak{B}}

\newcommand{\frakS}{\mathfrak{S}}
\newcommand{\frakF}{\mathfrak{F}}

\newcommand{\frakz}{\mathfrak{z}}
\newcommand{\sfC}{\mathsf{C}}
\newcommand{\sfS}{\mathsf{S}}
\newcommand{\sfJ}{\mathsf{J}}
\newcommand{\sfn}{\mathsf{n}}
\newcommand{\sfc}{\mathsf{c}}

\newcommand{\sfE}{\mathsf{E}}
\newcommand{\sfF}{\mathsf{F}}
\newcommand{\sfL}{\mathsf{L}}

\newcommand{\sfD}{\mathsf{D}}
\newcommand{\sfP}{\mathsf{P}}
\newcommand{\sfQ}{\mathsf{Q}}
\newcommand{\sfW}{\mathsf{W}}

\newcommand{\sfh}{\mathsf{h}}

\newcommand{\sfT}{\mathsf{T}}

\newcommand{\sfN}{\mathsf{N}}
\newcommand{\sfX}{\mathsf{X}}
\newcommand{\sfY}{\mathsf{Y}}
\newcommand{\sfb}{\mathsf{b}}
\newcommand{\sfg}{\mathsf{g}}
\newcommand{\sfd}{\mathsf{d}}
\newcommand{\sfm}{\mathsf{m}}

\newcommand{\sfK}{\mathsf{K}}
\newcommand{\sfk}{\Bbbk}

\newcommand{\bbK}{\mathbb{K}}

\newcommand{\bfY}{\mathbf{Y}}
\newcommand{\bfB}{\mathbf{B}}
\newcommand{\bfF}{\mathbf{F}}
\newcommand{\bfE}{\mathbf{E}}
\newcommand{\bfL}{\mathbf{L}}

\newcommand{\bfg}{\mathbf{g}}
\newcommand{\bfn}{\mathbf{n}}
\newcommand{\bfN}{\mathbf{N}}
\newcommand{\bfc}{\mathbf{c}}

\newcommand{\bfm}{\mathbf{m}}
\newcommand{\bfa}{\mathbf{a}}
\newcommand{\bfb}{\mathbf{b}}
\newcommand{\bcdot}{{\boldsymbol{\cdot}}}
\newcommand{\bplus}{{\boldsymbol{+}}}
\newcommand{\bminus}{{\boldsymbol{-}}}
\newcommand{\bmu}{\boldsymbol{\mu}}

\newcommand{\bmup}{\bmu_\bplus}
\newcommand{\bmum}{\bmu_\bminus}

\newcommand{\Cz}{{\scrC^0_\g}}

\newcommand{\Cbz}{{\scrC^0_\bfg}}
\newcommand{\Cxi}{{ \scrC^\xi}}
\newcommand{\Cuxi}{{ \scrC^\uxi}}

\newcommand{\Yxi}{{ \calY^\xi}}

\newcommand{\Mxi}{{{}^\xi\calM}}

\newcommand{\tDynkinxi}{{{}^\xi\tDynkin}}

\newcommand{\calC}{\mathcal{C}}
\newcommand{\calQ}{\mathcal{Q}}

\newcommand{\calA}{\mathcal{A}}
\newcommand{\calS}{\mathcal{S}}

\newcommand{\calK}{\mathcal{K}}
\newcommand{\calR}{\mathcal{R}}
\newcommand{\calI}{\mathcal{I}}

\newcommand{\calX}{\mathcal{X}}
\newcommand{\calY}{\mathcal{Y}}
\newcommand{\calM}{\mathcal{M}}
\newcommand{\calN}{\mathcal{N}}

\newcommand{\calT}{\mathcal{T}}
\newcommand{\calW}{\mathcal{W}}
\newcommand{\calL}{\mathcal{L}}
\newcommand{\calJ}{\mathcal{J}}

\newcommand{\scrQ}{\mathscr{Q}}
\newcommand{\scrC}{\mathscr{C}}
\newcommand{\scrA}{\mathscr{A}}

\newcommand{\scrS}{\mathscr{S}}
\newcommand{\scrP}{\mathscr{P}}
\newcommand{\scrF}{\mathscr{F}}

\newcommand{\ttt}{\mathtt{t}}

\newcommand{\ttq}{{\ms{1mu}\mathtt{q}\ms{1mu}}}
\newcommand{\To}[1][{\hspace{2ex}}]{\xrightarrow{\,#1\,}}

\newlength{\mylength}
\newcommand*{\para}{%
  \rlap{\rotatebox{-30}{\rule[.05ex]{.4pt}{.77em}}}%
  \kern.04em%
  \rlap{\kern.36em\raisebox{0.649519052835em}{\rule{.6em}{.4pt}}}%
  \rule{.6em}{.4pt}\kern-.04em%
  \rotatebox{-30}{\rule[.05ex]{.4pt}{.77em}}}

\newcommand{\hDynkin}{{\overline{\Dynkin}}}
\newcommand{\hDynkinf}{{\overline{\Dynkinf}}}

\newcommand{\tDynkin}{{\widetilde{\Dynkin}}}
\newcommand{\tDynkinf}{{\widetilde{\Dynkinf}}}
\newcommand{\tbDynkinf}{{\widetilde{\bDynkinf}}}

\newcommand{\wl}{\sfP}
\newcommand{\rl}{\sfQ}

\newcommand{\weyl}{\sfW}
\newcommand{\lan}{\langle}
\newcommand{\ran}{\rangle}

\newcommand{\ee}{\end{enumerate}}
\newcommand{\bitem}{\begin{itemize}}
\newcommand{\eitem}{\end{itemize}}
\newcommand{\ben}{\begin{enumerate}[{\rm (1)}]}
\newcommand{\bnum}{\begin{enumerate}[{\rm (i)}]}
\newcommand{\bNum}{\begin{enumerate}[{\rm (I)}]}
\newcommand{\bnump}{\begin{enumerate}[{\rm (i)$'$}]}
\newcommand{\bna}{\begin{enumerate}[{\rm (a)}]}
\newcommand{\bnA}{\begin{enumerate}[{\rm (A)}]}
\newcommand{\bc}{\begin{cases}}
\newcommand{\ec}{\end{cases}}

\newenvironment{myequation}
{\relax\setlength{\arraycolsep}{1pt}\begin{eqnarray}} {\end{eqnarray}}
\newenvironment{myequationn}
{\relax\setlength{\arraycolsep}{1pt}\begin{eqnarray*}} {\end{eqnarray*}}

\nc{\eq}{\begin{myequation}}
\nc{\eneq}{\end{myequation}}
\nc{\eqn}{\begin{myequationn}}
\nc{\eneqn}{\end{myequationn}}
\nc{\cl}{\colon}
\nc{\ake}[1][1ex]{\rule[-#1]{0ex}{1ex}}
\nc{\akew}[1][1ex]{\rule[-1ex]{#1}{0ex}}
\nc{\akeu}[1][1ex]{\rule[#1]{0ex}{1ex}}
\nc{\id}{\mathrm{id}}
\nc{\bl}{\bigl(}
\nc{\br}{\bigr)}
\nc{\qt}[1]{\quad\text{#1}}
\nc{\qtq}[1][{and}]{\quad\text{{#1}}\quad}
\nc{\eqs}[1]{\underset{\raisebox{.4ex}[.7ex][0ex]{$\scriptstyle{#1}$}}{=}}
\nc{\snoi}{\smallskip \noindent}
\nc{\mnoi}{\medskip \noindent}
\nc{\Mat}{\mathrm{Mat}}
\nc{\ol}{\overline}
\nc{\ul}{\underline}
\nc{\ang}[1]{\boldsymbol{\langle}{#1}\boldsymbol{\rangle}}
\nc{\rang}[1]{\boldsymbol{\langle}{#1}\boldsymbol{\rangle} \hspace{-.6ex} \boldsymbol{\rangle}}
\nc{\ba}{\begin{array}}
\nc{\ea}{\end{array}}
\nc{\noi}{\noindent}
\nc{\evq}{ \ev_{q=1}}

\nc{\yi}{x_{i}}
\nc{\yj}{x_{j}}
\nc{\yk}{x_{k}}

\nc{\yjm}{x_{j,m}}
\nc{\yjmp}{x_{j,m+1}}
\nc{\yjmm}{x_{j,m-1}}
\nc{\yim}{x_{i,m}}
\nc{\yip}{x_{i,p}}
\nc{\yipp}{x_{i,p+1}}
\nc{\yipm}{x_{i,p-1}}
\nc{\yjp}{x_{j,p}}
\nc{\yjpp}{x_{j,p+1}}
\nc{\yjpm}{x_{j,p-1}}
\nc{\yimp}{x_{i,m+1}}
\nc{\yimm}{x_{i,m-1}}
\nc{\ykm}{x_{k,m}}
\nc{\ynm}{x_{n,m}}
\nc{\ynmp}{x_{n,m+1}}
\nc{\ynmpp}{x_{n,m+2}}
\nc{\ynnm}{x_{n-1,m}}
\nc{\ynnmp}{x_{n-1,m+1}}
\nc{\ynnmpp}{x_{n-1,m+2}}
\nc{\ykmp}{x_{k,m+1}}
\nc{\ykmm}{x_{k,m-1}}
\nc{\ykp}{x_{k,p}}
\nc{\ykpp}{x_{k,p+1}}
\nc{\ykpm}{x_{k,p-1}}

\nc{\seq}[1]{ \boldsymbol{(} {#1} \boldsymbol{)}   }

\title[Quantum virtual Grothendieck rings and quantum cluster algebras]{Quantization of virtual Grothendieck rings and their structure including quantum cluster algebras}

\author[I.-S. Jang]{Il-Seung Jang$^\ddagger$}
\address[I.-S. Jang]{$^\ddagger$Department of Mathematics, Incheon National University, Incheon 22012, Korea}
\email{ilseungjang@inu.ac.kr}
\urladdr{https://sites.google.com/view/isjang/home/}
\thanks{$^\ddagger$ I.-S. Jang was supported by Incheon National University Research Grant in 2023.}

\author[K.-H. Lee]{Kyu-Hwan Lee$^{\star}$}
\thanks{$^{\star}$ K.-H. Lee was partially supported by a grant from the Simons Foundation (\#712100).}
\address[K.-H. Lee]{$^{\star}$Department of Mathematics, University of Connecticut, Storrs, CT 06269, U.S.A.}
\email{khlee@math.uconn.edu}
\urladdr{https://www.math.uconn.edu/~khlee/}

\author[S.-j.~Oh]{Se-jin Oh$^{\dagger}$}
\address[S.-j.~Oh]{$^{\dagger}$Department of Mathematics, Sungkyunkwan University, Suwon, South Korea}
\email{sejin092@gmail.com}
\urladdr{https://sites.google.com/site/mathsejinoh/}
\thanks{$^{\dagger}$ S.-j.\ Oh was supported by the Ministry of Education of the Republic of Korea and the National Research Foundation of Korea (NRF-2022R1A2C1004045).}

\keywords{quantum affine algebra, t-quantized Cartan matrix, virtual Grothendieck ring, quantization, Kazhdan--Lusztig theory, cluster algebra, positivity}
\subjclass[2020]{13F60, 17B37,17B10, 17B67, 18N25}

\date{\today}

\begin{document}

\begin{abstract}
The quantum Grothendieck ring of a certain category of finite-dimensional modules over a quantum loop algebra associated with a complex finite-dimensional simple Lie algebra $\g$ has a quantum cluster algebra structure of skew-symmetric type. Partly motivated by a search of a ring corresponding to a quantum cluster algebra of {\em skew-symmetrizable} type, the quantum {\em virtual} Grothendieck ring, denoted by $\frakK_q(\g)$, is recently introduced by Kashiwara--Oh \cite{KO23} as a subring of the quantum torus based on the $(q,t)$-Cartan matrix specialized at $q=1$. In this paper, we prove that $\frakK_q(\g)$ indeed has a quantum cluster algebra structure of skew-symmetrizable type. This task essentially involves constructing distinguished bases of $\frakK_q(\g)$ that will be used to make cluster variables and generalizing the quantum $T$-system associated with Kirillov--Reshetikhin modules to establish a quantum exchange relation of cluster variables. Furthermore, these distinguished bases naturally fit into the paradigm of Kazhdan--Lusztig theory and our study of these bases leads to some conjectures on quantum positivity and $q$-commutativity.
\end{abstract}

\setcounter{tocdepth}{1}

\maketitle
\tableofcontents

\section{Introduction}

\subsection{Background} \label{subsec:background}
Let $\sfC=(\sfc_{i,j})$ be a Cartan matrix of finite type, and let $\g$ be the finite-dimensional simple Lie algebra over $\mathbb{C}$ associated with $\sfC$. Since its inception as trigonometric solutions to the quantum Yang--Baxter equation \cite{Dr, Ji}, the quantum loop algebra $U_q(\calL\g)$ of $\g$ has been one of the central objects in representation theory and mathematical physics, and various algebraic and geometric approaches have been taken to study the finite-dimensional modules over $U_q(\calL\g)$. Moreover, for the last 15 years or so, as categorification became one of the major trends in representation theory and cluster algebra structures were discovered ubiquitously, the category $\scrC_{\g}$ of finite-dimensional $U_q(\calL\g)$-modules became a focal point of research where these new ideas and methods could be applied fruitfully,
since the  quantum  Grothendieck ring of
 $\scrC_\g$ provides a categorification of a (quantum) cluster algebra and  generalizes the Kazhdan--Lusztig(KL) theory.

To be more precise,
the quantum cluster algebra $\mathcal A$, introduced by Berenstein--Fomin--Zelevinsky (BFZ) in \cite{BZ05,FZ02}, is a non-commutative $\Z[q^{\pm 1/2}]$-algebra contained in the quantum torus $\Z[ \tX_k^{\pm 1} | k \in \K]$ which is equipped with a distinguished set of generators ({\it quantum cluster variables}) grouped into subsets ({\it quantum clusters}), where $\K$ is an index set. Each cluster is defined inductively by a sequence of certain combinatorial algebraic operations ({\it mutations}) from an initial cluster. 
Since then, numerous connections and applications have been discovered in  various fields of mathematics.

It is well-known that the  quantum  cluster algebra was introduced in an attempt to create an algebraic framework for the dual-canonical/upper-global basis $\mathbf B^*$ \cite{L90, K91, K93} of the quantum group $U_q(\g)$. Indeed, it is shown in \cite{GLS13, GY17} that the unipotent quantum coordinate algebra
$A_q(\mathfrak{n})$ of $U_q(\g)$, which is the graded dual of the half of $U_q(\g)$, has a quantum cluster algebra structure, and intensive research has been performed to understand the structure in relation with $\mathbf B^*$ (see \cite{Kas18} for a survey). In these efforts, it turned out that categorification provides powerful methods \cite{KL1,KL2,R08,KKKO18}.

When $\g$ is of simply-laced type with its set of positive roots denoted by $\Phi^+_\g$, we can consider the path algebra $\C Q$ of the Dynkin quiver $Q$ associated with $\g$ and obtain the Auslander--Reiten (AR) quiver $\Gamma_Q$ of $\C Q$. In turn, $\Gamma_Q$ can be understood as a
{heart} of the AR-quiver $\widehat{\Dynkin}$ of the
derived category $D^b({\rm Rep}(\C Q))$, called the repetition quiver.
%
  
In \cite{HL15}, which culminates preceding works \cite{Rin90, L90, VV03, Nak04, Toen, H04, HL10},
Hernandez and Leclerc defined the heart subcategory $\scrC^Q_\g$ of $\scrC_\g$ by using $\Gamma_Q$, 
and proved that the \emph{quantum Grothendieck ring} $\calK_\ttt(\scrC^Q_\g)$ of $\scrC^Q_\g$ is isomorphic to $A_q(\mathfrak{n})$ and that the isomorphism sends the basis of $\calK_\ttt(\scrC^Q_\g)$ consisting of the elements corresponding to simple objects in $\scrC^Q_\g$ to $\mathbf B^*$ of $A_q(\mathfrak{n})$ (cf.~\cite{Nak11}).


To extend the results of \cite{HL15,HL16} to non-simply-laced types, the Q-datum $\scrQ$ is introduced in \cite{FO21} as a generalization of the Dynkin quivers of types $ADE$. Through the Q-datum for any finite type, the (combinatorial) AR-quiver $\Gamma_\scrQ$, the repetition quiver $\widehat{\Dynkin}^\sigma$, and the heart subcategory $\scrC^\scrQ_\g$ of $\scrC_{\g}$ are naturally defined, where $\sigma$ is the Dynkin diagram automorphism of simply-laced $\bfg$ whose orbits produce the Dynkin diagram of $\g$.
One could possibly expect that $\calK_\ttt(\scrC^\scrQ_\g)$ would be isomorphic to $A_q(\mathfrak{n})$ of $U_q(\mathfrak g)$ when $\g$ is of non-simply-laced type, generalizing the result in types $ADE$ to all types. However, further studies \cite{KO19,OS19,HO19,FHOO} show that the quantum Grothendieck ring $\calK_\ttt(\scrC^\scrQ_\g)$ is actually isomorphic to $A_q(\mathsf{\bfn})$ of $U_q(\bfg)$ associated with $\bfg$ of simply-laced type. Hence the structure of $\calK_\ttt(\scrC^\scrQ_\g)$ is intrinsically relevant to the counterpart of simply-laced type, and the quantum cluster algebra
structure associated with $\calK_\ttt(\scrC^\scrQ_\g)$ is still of skew-symmetric type.

\subsection{Overview of this paper}
Since there are quantum cluster algebras of skew-symmetrizable type, a natural question arises:

\vskip -1.2em

\begin{eqnarray*} &&
\parbox{85ex}{
 {\em Can we extend $\calK_\ttt(\scrC^Q_\g)$ $($or $\calK_\ttt(\scrC_\g))$ in such a way to have a quantum cluster algebra structure of skew-symmetrizable type?}
}
\end{eqnarray*}

\noindent
Partly motivated by this question, Kashiwara and Oh introduced the {\em quantum virtual Grothendieck ring} $\frakK_q(\g)$ inside the quantum torus $\calX_q(\g)$ with respect to the $(q,t)$-Cartan matrix specialized at $q=1$ in a recent paper \cite{KO23}.
Pursuing the direction further, in this paper, we prove that $\frakK_q(\g)$ indeed has a quantum cluster algebra structure of skew-symmetrizable type.  In a subsequent paper, our result will be utilized to fully answer the above question and to genuinely extend the results of \cite{HL15} in the sense that $A_q(\mathfrak{n})$ is involved even for $\g$ of non-simply-laced type. 
We remark that 
the evaluation of $\frakK_q(\g)$ at $q=1$ coincides with the \emph{folded $t$-character ring} (Remark~\ref{rem: evaluation of Kqg}),
denote by $\calK^-(\g)$, 
which is introduced by Frenkel--Hernandez--Reshetikhin in \cite{FHR21} to explore a (conjectural) quantum integrable model corresponding to what is called the folded Bethe Ansatz equation (see Remark~\ref{rem: relation with FHR}).

Though we do not yet have an actual category that will replace $\scrC_{\g}$ for our purpose (cf.~\cite[Remark 3.2, Remark 5.1]{FHR21}), we can still utilize an algebraic characterization of $\calK_\ttt(\Cz)$ as the intersection of the kernels of screening operators in $\calY_{\ttt}(\g)$, where $\Cz$ is the skeleton subcategory of $\scrC_{\g}$ and $\calY_{\ttt}(\g)$ is the quantum torus with respect to the $(q,t)$-Cartan matrix specialized at $t=1$.

In order to give a quantum cluster algebra structure on $\frakK_q(\g)$ in this paper, we need to construct  quantum cluster variables and exchange relations for mutations. The former requires constructing distinguished bases for $\frakK_q(\g)$ and the latter amounts to generalizing the quantum $T$-system associated with Kirillov--Reshetikhin (KR) modules as explained briefly below.

We establish three bases of $\frakK_q(\g)$, denoted by $\sfF_q, \sfE_q$, and $\sfL_q$ respectively. The basis $\sfF_q$ is constructed by a generalization of Frenkel--Mukhin (FM) algorithm \cite{FM01}, which plays a crucial role in 
studying $\frakK_q(\g)$. Furthermore, it induces two other important bases $\sfE_q$ and $\sfL_q$ of $\frakK_q(\g)$. 
Then we take a $q$-commuting subset of $\sfF_q$ consisting of \emph{KR-polynomials} (Definition~\ref{def: KR-type}) as the quantum cluster of initial seed and develop a {\em quantum folded $T$-system} among the KR-polynomials in $\sfF_q$ to serve as the quantum exchange relation. After making compatible pairs available for our use (cf.~\cite{KO23}), we establish a quantum cluster algebra structure on a subalgebra and extend it to $\frakK_q(\g)$.

It is worthwhile to remark that when $\g$ is simply-laced, the basis $\sfL_q$ (resp.~$\sfE_q$) comes from simple (resp.~standard) modules in $\scrC^0_{\g}$, and the entries of the transition matrix between $\sfL_q$ and $\sfE_q$ are understood as analogues of the KL-polynomials.  
Thus our construction of $\sfL_q$ and $\sfE_q$ for all the finite types extends the KL-theory for $\scrC_\g^0$. 
Moreover, we have conjectures related (i) to positivity on  KR-polynomials in $\sfF_q$ and  \emph{real} elements  in $\sfL_q$, and (ii) to BFZ-expectation that every quantum cluster monomial is an element in the canonical basis  (see Conjecture~\ref{conj: real positive} below).


Throughout this paper, the interplay between $\g$ and its simply-laced type counterpart $\bfg$ and the Dynkin diagram automorphism $\sigma$ (cf.~\eqref{eq: D auto} and \eqref{eq: pair of gg}) provides important viewpoints leading to natural definitions. 
However, we emphasize that none of our main constructions, including  bases $\sfF_q, \sfE_q$, and $\sfL_q$, is obtained merely from  combining objects in each orbit of $\sigma$. That is, none of our results is a consequence of simple folding. Rather, there seem to exist quite intriguing features of non-simply-laced type objects at the quantum level.

In the following subsections, we will review known results in Sections \ref{subsec_intro:qca} and \ref{subsec_intro:qcr and qca} with some details,  and present our results more rigorously in Section \ref{subsec:main}, and mention our future work in Section \ref{subsec:fw}.

\subsection{Quantum Grothendieck ring and quantum loop analogue of KL-theory} \label{subsec_intro:qca}
From the study for $q$-deformation of $\calW$-algebras, the $q$-character\footnote{In the main body of this paper, we sometimes call  it \emph{$t$-character} by replacing the role of $q$ by $t$.} theory
for $\scrC^0_\g$ was invented by Frenkel--Reshetikhin  \cite{FR99} and further developed by Frenkel--Mukhin \cite{FM01}, which says that the (non-quantum) Grothendieck ring $K(\scrC^0_\g)$ of $\scrC^0_\g$ is isomorphic to the commutative ring generated by the $q$-characters of fundamental modules $L(Y_{i,p})$ under the Chari--Pressley's classification \cite{CP95, CP}. 
For simply-laced type $\bfg$, Nakajima \cite{Nak04} and Varagnolo--Vasserot \cite{VV03} constructed a non-commutative $t$-deformation of $K(\scrC^0_\bfg)$ in a quantum torus $\calY_\ttt(\bfg)$, denoted by $\calK_\ttt(\scrC^0_\bfg)$, based on a geometrical point of view.
Since the specialization of $\calK_\ttt(\scrC^0_\bfg)$ at $\ttt=1$ recovers $K(\scrC^0_\bfg)$, we call $\calK_\ttt(\scrC^0_\bfg)$ the quantum Grothendieck ring associated with $\scrC^0_\bfg$.

In particular, Nakajima established a KL-type algorithm to describe the composition multiplicity $P_{m,m'}$ of a simple module $L(m')$ inside a standard module $E(m)$ through equations in $K(\scrC^0_\bfg)$:
Denoting by $\calM_+$ the parametrizing set of simple modules in $\scrC^0_\bfg$, we have
$$
[E(m)] = [L(m)] + \sum_{m'\in \calM_+; \; m' \lN m} P_{m,m'} \, [L(m')].
$$
 It is proved by Nakajima \cite{Nak03,Nak04} that the multiplicity $P_{m,m'}$ is equal to the specialization at $t=1$ of a polynomial $P_{m,m'}(t)$ with non-negative coefficients, which can be understood as a quantum loop analogue of KL-polynomial.

One step further, each $q$-character of simple module $L(m)$ (resp.~standard module $E(m)$) allows a $t$-deformation in $\calK_\ttt(\scrC^0_\bfg)$, denoted by $L_t(m)$ (resp.~$M_t(m)$), whose coefficients in $\Z[t^{\pm1/2}]$ are non-negative. Its specialization at $t=1$ recovers the $q$-character of $L(m)$ (resp. $M(m)$) and the transition map between
$\bfL_\ttt =\{  L_t(m)\}$ and $\bfE_\ttt =\{  E_t(m) \}$ in $\calK_\ttt(\scrC^0_\bfg)$ satisfies the following equation:
\begin{align} \label{eq: KL paradigm}
 E_t(m) =  L_t(m) + \sum_{m' \in \calM_+; \, m' \lN m} P_{m,m'}(t)\,  L_t(m')  \quad \text{where $P_{m,m'}(t) \in t \Z_{\ge0}[t]$}.
\end{align}
We call $\bfL_\ttt$ the \emph{canonical basis}  and $\bfE_t$ the \emph{standard basis} of $\calK_\ttt(\scrC^0_\bfg)$, respectively (see Remark~\ref{rem: KL theory} also). In what follows, {\em positivity} generally means that polynomials of interest have non-negative coefficients as is the case with $P_{m,m'}(t) \in t \Z_{\ge0}[t]$.
We remark that, in these developments, 
 the geometry of quiver varieties plays an essential role.

Despite the absence of fully developed theory of quiver varieties for general type $\g$, Hernandez \cite{H03,H04} constructed a conjectural KL-theory for $\scrC_\g^0$ in a purely algebraic way. Let us explain this more precisely. Using the \emph{quantum Cartan matrix} $\sfC(q)$, Hernandez constructed the quantum torus $\calY_\ttt(\g)$
and defined $\calK_\ttt(\scrC^0_\g)$ to be the intersection of the kernels of the $t$-deformed screening operators $S_{i,\ttt}$'s on $\calY_\ttt(\g)$. 
Then he constructed a basis $\bfF_\ttt = \{  F_t(m) \}$ by deforming the FM-algorithm and proved the positivity of $F_t(Y_{i,p})= L_t(Y_{i,p})$.
Then the basis $\bfF_\ttt$ induces two other bases $\bfE_\ttt=\{E_t(m)\}$ and $\bfL_\ttt=\{L_t(m)\}$ satisfying \eqref{eq: KL paradigm} that enable us to establish a conjectural KL-theory, expecting the positivity of analogues of KL-polynomials and $L_t(m)$'s.

Recently, large parts of the conjectures for non-simply-laced $g$ are proved by Fujita--Hernandez--Oh--Oya through so-called \emph{propagation of positivity}. Let $\bfg$ be an unfolding of $g$ as follows:
$$
(g,\bfg) = (B_n,A_{2n-1}),  \  \  (C_n,D_{n+1}), \  \ (F_4,E_{6}), \  \  (G_2,D_{4}).
$$
Then it is proved in \cite{HL15,FHOO} that
$$ \text{
 $\bbK_\ttt(\scrC^0_\bfg)$ and $\bbK_\ttt(\scrC_g^0)$ have the same presentation,}$$ where $\bbK_\ttt(\scrC^0_\g) \seteq \Q(q^{1/2})\otimes_{\Z[q^{\pm 1/2}]}  \calK_\ttt(\scrC_\g^0)$. Hence the ring $\bbK_\ttt(\scrC^0_g)$ can be interpreted  as the boson-extension of $A_q(\bfn)$ of the simply-laced $\bfg$. Then the KL-theory and positivity are established for type $B_n$ using the quantum Schur--Weyl duality functor \cite{KKK18,KKO} between $\scrC_{A_{2n-1}}^0$ and $\scrC_{B_{n}}^0$, and similar conjectures for $CFG$-types are mostly resolved in \cite{FHOO,FHOO2} using the quantum Schur--Weyl duality functor \cite{KKK15,KO19,OS19} for these types and
the \emph{degrees} (also called \emph{$g$-vectors}) of (quantum) cluster algebra theory. As indicated above, the presentation of $\bbK_\ttt(\scrC^0_g)$ is of simply-laced type even for non-simply-laced $g$.

\subsection{ Quantum cluster algebra structure of skew-symmetric type on $\calK_\ttt(\Cz)$} \label{subsec_intro:qcr and qca} In the seminal paper \cite{HL16}, Hernandez--Leclerc proved that $K(\scrC_\g^-)$ for a subcategory $\scrC_\g^-$ of $\scrC_\g^0$ has a cluster algebra structure of skew-symmetric type for any $\g$ of finite type.
To show the cluster algebra structure, they employed the T-system among Kirillov--Reshetikhin (KR) modules proved by Nakajima \cite{Nak03} for simply-laced types and by Hernandez \cite{H06} for non-simply-laced types.
Then the result of \cite{HL16} is extended to $\calK_\ttt(\scrC_\g^0)$ in \cite{B21,HO19,FHOO,KKOP1,KKOP2,FHOO2} to obtain quantum cluster algebras of skew-symmetric type. 
Some important features of these works can be summarized as follows: 
\bna
\item The extension to whole category $\scrC_\g^0$ in~\cite{KKOP1,KKOP2} involves a categorical language.
\item The main idea of the extension to quantum cluster algebra in \cite{B21,HO19,FHOO} is the quantization of T-system among KR modules.
\item The monoidal categorification result in~\cite{KKOP2}  tells us that
every quantum cluster monomial of $\calK_\ttt(\scrC_\g^0)$ corresponds to an element of $\bfL_\ttt$. This gives an affirmative answer to the BFZ-conjecture \cite{FZ02} on $\bfB^*$ and the quantum cluster monomials.
\item  As every KR-polynomial $F_t(m)$ appears as a quantum cluster variable of $\calK_\ttt(\scrC_\g^0)$, it is proved in~\cite{KKOP2,FHOO2} that $F_t(m)  =  L_t(m)$ for any KR-module $L(m)$.
\ee
Here we remark that the result of \cite{KKOP2} is for $K(\scrC_\g^0)$ and extended to $\calK_\ttt(\scrC_\g^0)$ in \cite{FHOO2}.

\subsection{Main results of this paper} \label{subsec:main} In this paper, we initiate a study of $\frakK_q(\g)$ in the perspective of Sections \ref{subsec_intro:qca} and \ref{subsec_intro:qcr and qca}. Due to lack of a representation theory corresponding to $\frakK_q(\g)$, we approach the ring $\frakK_q(\g)$ by analyzing its construction in~\cite{KO23} and by exploiting \ref{it:qvgr1} and \ref{it:qvgr2}, where 
\bNum
\item \label{it:qvgr1} $\frakK_q (\g)$ is a $q$-deformation of the commutative ring $\mathcal{K}^{-} (\g)$, which
is the specialization of the refined ring $\overline{\calK}_{\ttq,\ttt,\upalpha}(\g)$ of interpolating $(\ttq, \ttt)$-characters in \cite{FHR21} at $\ttq = 1$ and $\upalpha = d$,
\item \label{it:qvgr2} $\mathcal{K}^{-}(\bfg) \simeq K(\scrC_\bfg^0)$   if $\bfg$ is of simply-laced type,
\ee
(see Section~\ref{subsec: VGR} and \cite[Introduction]{KO23}). 
Here $\upalpha$ is a factor to interpolate several characters (see \cite[Remark 6.2(1)]{FHR21}) and $d$ is the lacing number of $\g$.
In particular, if $\sfg$ is of non-simply-laced type, there exist a simply-laced $\bfg$ containing $\sfg$ as a non-trivial Lie subalgebra (e.g.~see \cite[Proposition 7.9]{Kac} with \eqref{eq: pair of gg})
and a surjective homomorphism
\begin{align}\label{eq: naive folding}
\mathcal{K}^{-} (\bfg) \twoheadrightarrow \mathcal{K}^{-}(\sfg) \simeq \frakK(\sfg),
\end{align}
which is induced from the folding of generators of $\mathcal{K}^{-}(\bfg)  \simeq K(\scrC_\bfg^0)$. 
\smallskip

The main results of this paper can be summarized into two statements:
\bnA
\item \label{it: goalA} we construct bases $\sfF_q$, $\sfE_q$, and $\sfL_q$ of $\frakK_q(\g)$, which play similar roles of
$\bfF_\ttt$, $\bfE_\ttt$, and $\bfL_\ttt$,
\item \label{it: goalB}  we establish \emph{skew-symmetrizable} quantum cluster algebra structures on subrings of $\frakK_q(\g)$ (including itself) using the bases in \ref{it: goalA}.
\ee
Here we emphasize that our results can \emph{not} be obtained from the folding in~\eqref{eq: naive folding},
 as we do \emph{not} have a surjective homomorphism
$A_q(\bfn)  \twoheadrightarrow A_q(\sfn)$
from the canonical surjection $\C[\bfN] \twoheadrightarrow \C[\sfN]$, where $\C[N]$ denotes the unipotent coordinate ring of $N$
of $\g$.

\subsubsection{\em \bf Construction of bases and KL-paradigm for $\frakK_q(\g)$ } 
Let $\usfC(t)$ be the $(q,t)$-Cartan matrix specialized at $q=1$, which is called \emph{$t$-quantized Cartan matrix}.
To construct the basis $\sfF_q$ of $\frakK_q(\g)$, we apply a $q$-deformed version of FM-algorithm with respect to $\usfC(t)$.
However, there is no guarantee that the algorithm terminates in finite steps. To avoid this problem, we prove that the monomials (not including coefficients) of $F_q(X_{i,p})$ $((i,p) \in \widehat{\Dynkin}^\sfg_0)$ in $\sfF_q$ is obtained from those of the $q$-character of $L(Y_{\im,p})$ of type $\bfg$ via~\eqref{eq: naive folding} for $(\im,p) \in \widehat{\Dynkin}^\bfg_0$. 
Furthermore, we prove that a similar phenomenon occurs for a KR-polynomial $F_q(m^{(i)}[p,s])$ (Proposition \ref{prop:folding lemma}). This result implies that the outputs of the algorithm are indeed contained in
$\frakK_q(\g)$ and form a basis $\sfF_q$. 
The basis $\sfF_q$ nicely characterizes an element in $\frakK_q(\g)$ since each element in $\sfF_q$ has a unique dominant monomial (Theorem~\ref{thm: F_q}). Here we emphasize once more that general elements in  $\sfF_q$ are not susceptible of similar manipulations based on \eqref{eq: naive folding} even in the specialization at $q=1$ (Example~\ref{ex: Fq neq osigma Ft in general}), and determining the $\Z[q^{\pm 1/2}]$-coefficients of $F_q(m^{(i)}[p,s])$ is a completely different problem even for a KR-polynomial $F_q(m^{(i)}[p,s])$.

We investigate properties of the KR-polynomials in $\sfF_q$ in detail,
since they will be used as the quantum cluster variables of $\frakK_q(\g)$ (Propositions~\ref{prop:properties of Fq}
and~\ref{prop: range of folded KR}).
By applying the framework in \cite{H04}, we construct the standard basis $\sfE_q=\{ E_q(m) \}$ and the canonical basis $\sfL_q= \{ L_q(m) \}$ fitting into the paradigm of Kazhdan--Lusztig theory:
$$
E_q(m)  = L_q(m) + \sum_{m' \in \calM;\, m' \lN m} P_{m,m'}(q) \, L_q(m') \quad \text{where $P_{m,m'}(q) \in q\Z[q]$}.
$$

\subsubsection{\em \bf Quantum cluster algebra structure of skew-symmetrizable type on $\frakK_q(\g)$}
Based on the construction of bases for $\frakK_q(\g)$, we show quantum cluster algebra structures on subrings of $\frakK_q(\g)$ as the first task in the second part of this paper.

In \cite{KO23}, Kashiwara and Oh constructed
a compatible pair $(\Lambda,\tB)$ arising from the isomorphism between the subtorus $\calX_{q,Q}(\g)$ of $\calX_q(\g)$
and the torus containing $A_q(\n)$, in which the exchange matrix $\tB$
is skew-symmetrizable. Here $Q=(\Dynkin,\xi)$ is a Dynkin quiver of type $\g$. Interpreting entries in $\Lambda$ as pairing of KR-monomials (Theorem \ref{thm: compatible pairs}), we form
an initial quantum cluster consisting of certain KR-polynomials $F_q(m)$ for each Dynkin quiver $Q=(\Dynkin,\xi)$ and its corresponding subring $\frakK_{q,\xi}(\g)$.

As a quantum cluster should consist of mutually $q$-commutative elements, we prove that the family of $F_q(m)$ in the initial cluster are mutually $q$-commutative,  using the truncation homomorphism (Proposition \ref{prop: inj}) and the properties 
of KR-polynomials.
By investigating $q$-commuting conditions (Lemmas~\ref{lem: m(i;p,s)},
~\ref{lem: m[p,s] m(p,s)}, and~\ref{lem: last term in T-system}) and multiplicative structure among KR-polynomials $F_q(m)$,
we obtain \emph{the quantum folded T-systems} among KR-polynomials $F_q(m)$ (Theorem \ref{thm: quantum folded}):
\begin{align*}
\FMq{ \um^{(i)}[p,s) } \hspace{-.3ex} * \hspace{-.3ex}  \FMq{\um^{(i)}(p,s]} \hspace{-.3ex} = \hspace{-.3ex} q^{\al(i,k)}  \FMq{\um^{(i)}(p,s)} \hspace{-.3ex} * \hspace{-.3ex} \FMq{\um^{(i)}[p,s]} \hspace{-.3ex} + \hspace{-.3ex} q^{\ga(i,k)}  \hspace{-2.3ex} \prod_{j; \; d(i,j)=1} \hspace{-2.3ex}  \FMq{\um^{(j)}(p,s)}^{-\sfc_{j,i}}.
\end{align*}
Then we prove that $\frakK_{q,\xi}(\g)$ has a quantum cluster algebra structure of skew-symmetrizable type (Theorem~\ref{them:main2}) by using the quantum folded T-systems as mutation relations and applying special sequences of mutations.
In the proof, we adopt the setup of \cite{HL16,B21} and use the \emph{valued quivers} (Section~\ref{subsec: Bi-deco}) (equivalent to exchange matrices) for the sequences of mutations.
As applications, we obtain a quantum cluster algorithm to compute KR-polynomials $F_q(m)$ (Proposition~\ref{prop: m times mu}) and a sufficient condition for $q$-commutativity of certain pairs of  KR-polynomials $F_q(m)$ (Theorem~\ref{thm: folded KR-commuting condition}).

As the second task, we extend the result on $\frakK_{q,\xi}(\g)$  to the whole ring $\frakK_q(\g)$.  For this purpose,
we construct a new quantum seed, whose valued quiver  is a ``sink-source" quiver reflecting features of $\g$ and whose initial quantum cluster consists of certain KR-polynomials $F_q(m)$. Here the $q$-commutativity of the initial quantum cluster follows from Theorem~\ref{thm: folded KR-commuting condition}.
Finally, we prove that $\frakK_{q}(\g)$ has a quantum cluster algebra structure of skew-symmetrizable type by establishing (a)
a mutation equivalence between the valued quiver of $\frakK_{q,\xi}(\g)$ and that of $\frakK_{q}(\g)$,
and finding out (b) special sequences of mutations that yield every KR-polynomial $F_q(m)$ as a cluster variable. 

Since every KR-polynomial $F_q(m)$ appears as a cluster variable and every quantum cluster monomial is expected to be a canonical basis element and \emph{real}, we have the following conjecture:
\begin{Conj}  \label{conj: real positive} 
\bna
\item Every quantum cluster monomial  of $\frakK_q(\g)$ is contained in $\sfL_q$.
\item  For every KR-polynomial $F_q(m)$, we have $F_q(m)=L_q(m)$   and  $F_q(m)$ has non-negative coefficients.
\item  If $L_q(m)$ is real, that is, for any $k\in \Z_{\ge1}$, there exists $t \in \Z$ such that $L_q(m)^k =q^t L_q(m^k)$, then it  has non-negative coefficients.
\ee
\end{Conj}

\noindent
Also, we have two more conjectures on the $q$-commutativity of KR-polynomials $F_q(m)$ in Conjectures~\ref{conj: condition for q-comm} and ~\ref{conjecture 3}, which can be understood as natural generalizations of the results in~\cite{OS19} and~\cite{KKOP2,FHOO2}, respectively.

\subsection{Future work} \label{subsec:fw} 
In a forthcoming paper \cite{JLO2},  we study the heart subring $\frakK_{q,Q}(\g)$ of $\frakK_q(\g)$ 
in terms of a generalization $Q$ of the Dynkin quiver to non-simply-laced type, where the AR-quiver $\Gamma_Q$ and the repetition quiver $\widehat{\Dynkin}$ are defined for $\g$ of any finite type including BCFG.
Since it is shown in this paper that $\frakK_q(\g)$ has a quantum cluster algebra structure (of skew-symmetrizable type), as it is with $\calK_\ttt(\scrC_\g^0)$ in \cite{HL15,FHOO}, it will be shown that each  heart subring $\frakK_{q,Q}(\g)$ is isomorphic to $A_{\Z[q^{\pm 1/2}]}(\n)$ via a certain isomorphism $\Uppsi_Q$ and that the \emph{normalized} dual-canonical/upper-global basis of $A_{\Z[q^{\pm 1/2}]}(\n)$ corresponds to the subset $\sfL_{q,Q} \seteq \sfL_q \cap \frakK_{q,Q}(\g)$ under $\Uppsi_Q$.  
This justifies the name of $\sfL_q$, the \emph{canonical basis}.
Here we would like to make an emphasis on the difference between the known result and our new result when $g$ is non-simply-laced:
in the previous $\calK_\ttt(\scrC_g^Q)$-case,
the corresponding $A_q(\bfn)$ is of simply-laced type $\bfg$,
while in the new $\frakK_{q,Q}(g)$-case, the type of $A_q(n)$ is the same as that of $g$.
Based on some investigation of the heart subrings, we will also clarify the
presentation of $$\bbK_q(\g) \seteq \Q(q^{1/2}) \otimes_{\Z[q^{\pm 1/2}]} \frakK_q(\g),$$ which says that $\bbK_q(\g)$ can be understood as a boson-extension of $A_q(\n)$, as $\bbK_\ttt(\scrC^0_\bfg)$ is for $A_q(\bfn)$
of simply-laced type $\bfg$. Then we will show that the automorphisms of $\frakK_q(\g)$, arising from the reflections on Dynkin quivers $Q$ and the isomorphisms $\Uppsi_Q$, preserve the canonical basis $\sfL_q$ of $\frakK_q(\g)$ and induce a braid group action on $\frakK_q(\g)$.

\subsection*{Convention} Throughout this paper, we use the following convention.
\begin{enumerate}[\,\,$\bullet$]
\item For a statement $\mathtt{P}$, we set $\delta(\mathtt{P})$ to be $1$ or $0$
depending on whether $\mathtt{P}$ is true or not.
As a special case, we use the notation $\delta_{i,j} \seteq \delta(i=j)$ (Kronecker's delta).
\item For $k,l \in \Z$ and $s \in \Z_{\ge1}$, we write $k \equiv_s l$ if $s$ divides $k-l$ and $k \not\equiv_s l$, otherwise.
\item For a monoidal abelian category $\calC$, we denote its Grothendieck ring by $K(\calC)$.
The class of an object $X \in \calC$ is denoted by $[X] \in K(\calC)$.
\item  A monomial in a Laurent polynomial ring $\Z[x_j^{\pm 1} \ | \   j \in J ]$ is said to be \emph{dominant} (resp. \emph{anti-dominant}) if it is a product of non-negative (resp. non-positive) powers of $x_i$'s.
\item For elements $\{r_{j} \}_{j \in J}$ in a ring $(R,\star)$,  parameterized by a totally ordered set $J=\{ \cdots < j_{-1}<j_0 <j_1 < \cdots \}$, we write
 $$  \sta_{j \in J}^{\to}  r_j \seteq \cdots \star r_{j_{-1}} \star r_{j_0} \star r_{j_1} \star \cdots. $$
\item For integers $a,b \in \Z$, we set
\begin{align*}
&[a,b] \seteq \{ x\in \Z \; | \;  a \le x \le b\} &(a,b] \seteq  \{ x\in \Z \; | \;  a <x \le b\} \\
&[a,b) \seteq \{ x\in \Z \; | \;  a \le x < b\} &(a,b) \seteq  \{ x\in \Z \; | \;  a < x < b\}
\end{align*}
We refer to subsets of these forms as \emph{intervals}.
\item Let $X=\{ x_j \ | \ j \in J \}$ be a parameterized by an index set $J$. Then for $j \in J$ and a subset $\mathcal{J} \subset J$, we set
$$   (X)_j \seteq x_j   \quad \text{ and } \quad  (X)_{\calJ} \seteq \{x_j \ | \ j \in \calJ \}.$$
\end{enumerate}

\section{Preliminaries} \label{sec:preliminaries}

\subsection{Cartan datum} Let $\g$ be a Kac--Moody algebra of a symmetrizable type. We denote its Cartan matrix by $\sfC=(\sfc_{i,j})_{i,j \in I}$,
Dynkin diagram\footnote{Our convention is a variation of the Coxeter--Dynkin diagram in the sense that we connect vertices with single edges only. See the examples for the finite types. We will call them {\em Dynkin diagrams} for simplicity.}
by $\Dynkin$, weight lattice by $\wl$, set of simple roots by $\Pi=\{ \al_i \ | \ i \in I \}$ and set of simple coroots
by $\Pi^\vee = \{ h_i \ | \  i \in I \}$.    

Let $\sfD = {\rm diag}(d_i  \in \Z_{\ge 1}  \;  | \; i \in I )$ denote a diagonal matrix such that  
$$\text{$\osfB=\sfD\sfC$ and $\usfB=\sfC\sfD^{-1}$ become symmetric.}$$
We take $\sfD$ and the scalar product $(  \cdot, \cdot )$ on $\wl$ such that
\begin{align}\label{eq: sfD}
(\al_i,\al_j) = d_i \sfc_{i,j}= d_j \sfc_{j,i} \in \Z \ \text{ and }  \ (\al_i,\al_i) \in 2\Z_{\ge 1} \quad \text{ for all $i \in I$}.
\end{align}
We also denote by $\Phi_{\pm}$ the set of positive (resp. negative) roots of $\g$.
For each $i \in I$, we choose $\varpi_i \in \wl$ such that $\lan h_i,\varpi_j \ran=\delta_{i,j}$ $(j \in I)$. The free abelian group $\rl \seteq \soplus_{i \in I}\Z\al_i$ is called \emph{the root lattice}.

Throughout this paper, we use the following convention of finite Dynkin diagrams:
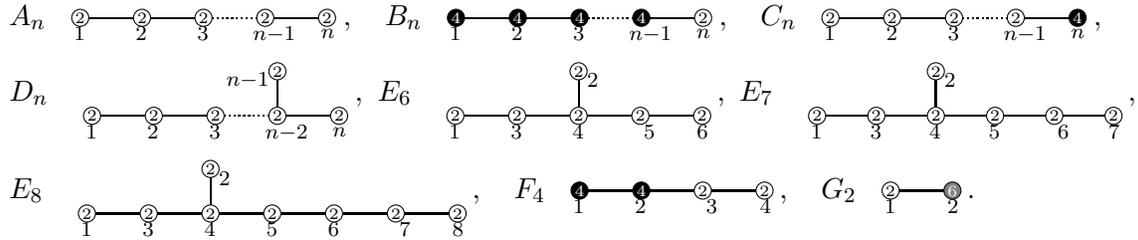
\begin{align*}
&A_n  \ \   \xymatrix@R=0.5ex@C=3.5ex{ *{\circled{2}}<3pt> \ar@{-}[r]_<{ 1 \ \  } & *{\circled{2}}<3pt> \ar@{-}[r]_<{ 2 \ \  } & *{\circled{2}}<3pt> \ar@{.}[r]_<{ 3  \ \  }
&*{\circled{2}}<3pt> \ar@{-}[r]_>{ \ \ \raisebox{-0.16cm}{\scriptsize $n$}} &*{\circled{2}}<3pt> \ar@{}[l]^>{   n-1} }, \quad
 B_n  \ \   \xymatrix@R=0.5ex@C=3.5ex{ *{\bcircled{4}}<3pt> \ar@{-}[r]_<{ 1 \ \  } & *{\bcircled{4}}<3pt> \ar@{-}[r]_<{ 2 \ \  } & *{\bcircled{4}}<3pt> \ar@{.}[r]_<{ 3  \ \  }
&*{\bcircled{4}}<3pt> \ar@{-}[r]_>{ \ \ \raisebox{-0.16cm}{\scriptsize $n$}} &*{\circled{2}}<3pt> \ar@{}[l]^>{   n-1} }, \quad
C_n  \ \   \xymatrix@R=0.5ex@C=3.5ex{ *{\circled{2}}<3pt> \ar@{-}[r]_<{ 1 \ \  } & *{\circled{2}}<3pt> \ar@{-}[r]_<{ 2 \ \  } & *{\circled{2}}<3pt> \ar@{.}[r]_<{ 3  \ \  }
&*{\circled{2}}<3pt> \ar@{-}[r]_>{ \ \ \raisebox{-0.16cm}{\scriptsize $n$}} &*{\bcircled{4}}<3pt> \ar@{}[l]^>{   n-1} }, \  \allowdisplaybreaks\\
&D_n  \ \  \raisebox{1em}{ \xymatrix@R=2ex@C=3.5ex{  &&& *{\circled{2}}<3pt> \ar@{-}[d]_<{ n-1  }  \\
*{\circled{2}}<3pt> \ar@{-}[r]_<{ 1 \ \  } & *{\circled{2}}<3pt> \ar@{-}[r]_<{ 2 \ \  } & *{\circled{2}}<3pt> \ar@{.}[r]_<{ 3  \ \  }
&*{\circled{2}}<3pt> \ar@{-}[r]_>{ \ \ \raisebox{-0.16cm}{\scriptsize $n$}} &*{\circled{2}}<3pt> \ar@{}[l]^>{   n-2} }},  \
E_{6}  \ \  \raisebox{1em}{   \xymatrix@R=2ex@C=3.5ex{  &&  *{\circled{2}}<3pt> \ar@{-}[d]^<{ 2\ \ }   \\
*{\circled{2}}<3pt>  \ar@{-}[r]_<{ 1 \ \  } & *{\circled{2}}<3pt> \ar@{-}[r]_<{ 3 \ \  } & *{\circled{2}}<3pt> \ar@{-}[r]_<{ 4 \ \  }
&*{\circled{2}}<3pt> \ar@{-}[r]_>{ \ \ 6} &*{\circled{2}}<3pt> \ar@{}[l]^>{   5} } },  \
E_{7}  \ \  \raisebox{1em}{    \xymatrix@R=2ex@C=3.3ex{  &&  *{\circled{2}}<3pt> \ar@{-}[d]^<{ 2\ \ }   \\
 *{\circled{2}}<3pt>  \ar@{-}[r]_<{ 1 \ \  } & *{\circled{2}}<3pt>  \ar@{-}[r]_<{ 3 \ \  } & *{\circled{2}}<3pt> \ar@{-}[r]_<{ 4 \ \  } & *{\circled{2}}<3pt> \ar@{-}[r]_<{ 5 \ \  }
&*{\circled{2}}<3pt> \ar@{-}[r]_>{ \ \ 7} &*{\circled{2}}<3pt> \ar@{}[l]^>{   6} } },  \allowdisplaybreaks \\
&  E_{8}  \ \  \raisebox{1em}{    \xymatrix@R=2ex@C=3.5ex{ &&  *{\circled{2}}<3pt> \ar@{-}[d]^<{ 2\ \ }   \\
 *{\circled{2}}<3pt>  \ar@{-}[r]_<{ 1 \ \  } & *{\circled{2}}<3pt>  \ar@{-}[r]_<{ 3 \ \  } & *{\circled{2}}<3pt> \ar@{-}[r]_<{ 4 \ \  } & *{\circled{2}}<3pt> \ar@{-}[r]_<{ 5 \ \  }   & *{\circled{2}}<3pt> \ar@{-}[r]_<{ 6 \ \  }
&*{\circled{2}}<3pt> \ar@{-}[r]_>{ \ \ 8} &*{\circled{2}}<3pt> \ar@{}[l]^>{   7} } }, \quad
F_{4}   \ \   \xymatrix@R=0.5ex@C=3.5ex{    *{\bcircled{4}}<3pt> \ar@{-}[r]_<{ 1 \ \  } & *{\bcircled{4}}<3pt> \ar@{-}[r]_<{ 2 \ \  }
&*{\circled{2}}<3pt> \ar@{-}[r]_>{ \ \ 4  } &*{\circled{2}}<3pt> \ar@{}[l]^>{   3} }, \quad
 G_2 \ \   \xymatrix@R=0.5ex@C=3.5ex{  *{\circled{2}}<3pt> \ar@{-}[r]_<{ 1 \ \  }  & *{\rcircled{6}}<3pt> \ar@{-}[l]^<{ \ \ 2  } }.
\end{align*}
Here $ \xymatrix@R=0.5ex@C=3.5ex{    *{\circled{t}}<3pt>}_k$ means that $(\al_k,\al_k)=t$.  
For $i,j \in I$, we denote by $d(i,j)$ the number of edges between $i$ and $j$ in $\Dynkin$ (whenever it is well-defined).
For example, in the finite $B_n$-case,
$d(n,n-1)=d(n-1,n)=1$ and $d(n,n-2)=d(n-2,n)=2$,
and in the finite $D_n$-case,
$d(n,n-1)=d(n-1,n)=2$ and $d(n,n-2)=d(n-2,n)=1$.

We denote by $\Dynkin_0$ the set of vertices and $\Dynkin_1$  the set of edges.
\emph{Throughout this paper, we consider only connected Dynkin diagrams.}
We sometimes use $\bDynkin$ for non-simply-laced  types to distinguish them from those of simply-laced types, and use $\Dynkinf$ for finite types
and,   when an emphasis is needed,    $\bDynkinf$ for finite non-simply-laced types.
For each $\Dynkinf$,  our convention amounts to taking   $$\sfD \seteq {\rm diag}( (\al_i,\al_i)/2 \ | \  i \in \Dynkinf_0 ) \text{ such that } \min( (\al_i,\al_i)/2 ) = 1.$$

The Weyl group $\weyl$ of $\g$ is generated by the reflections $s_i$ $(i \in I)$ acting on $\wl$ by
$$   s_i(\la) = \la - \lan \la,h_i \ran \al_i \qquad (\la \in \wl, \; i \in I).$$
A \emph{Coxeter element} of $\weyl$ is a product of the form $s_{i_1}\cdots s_{i_{|I|}}$ such that $\{ i_k \ | \ 1 \le k \le |I| \} =I$.
All Coxeter elements are conjugate in $\weyl$ when $\Dynkin$ is a tree \cite{Col,Hu}, and their common order in $\weyl$ is finite when $\weyl$ is finite \cite{Sp09}, in which case the order is called the \emph{Coxeter number} and denoted by $\sfh$.

A  bijection $\sigma$ from $\Dynkin_0$ to itself is said to be a \emph{Dynkin diagram automorphism}
if $\lan h_i,\al_j \ran= \lan h_{\sigma(i)},\al_{\sigma(j)} \ran$ for all $i,j\in \Dynkin_0$.
\emph{Throughout this paper, we assume that  Dynkin diagram automorphisms $\sigma$ satisfy the following condition:}
\begin{align} \label{eq: auto cond}
\text{ \emph{there is no $i \in \Dynkin_0$ such that $d(i,\sigma(i))=1$}}.
\end{align}
The condition in~\eqref{eq: auto cond} is referred to as an \emph{admissibility} (see~\cite[\S 12.1.1]{L94}).

For each Dynkin diagram $\Dynkinf$ of finite type $A_{2n-1}$, $D_n$ or $E_6$, there exists a unique non-identity Dynkin diagram automorphism $\vee$ of order $2$ (except $D_4$-type, in which case, there are three automorphism of order $2$ and
two non-identity automorphisms $\widetilde{\vee}$ and $\widetilde{\vee}^2$ of order $3$) satisfying the condition in~\eqref{eq: auto cond}.
\begin{equation} \label{eq: D auto}
\begin{aligned}
&A_{2n-1}  \ \   \xymatrix@R=0.5ex@C=5ex{ *{\circled{2}}<3pt> \ar@{-}[r]_<{ 1 \ \  } \ar@/^1.5pc/@{<-->}[rrrr]|{\vee}
& *{\circled{2}}<3pt> \ar@{-}[r]_<{ 2 \ \  } \ar@/^1pc/@{<-->}[rr]
& *{\circled{2}}<3pt> \ar@{.}[r]
&*{\circled{2}}<3pt> \ar@{-}[r]_>{ \ \ 2n-1} &*{\circled{2}}<3pt> \ar@{}[l]^>{   2n-2} }, \quad
&&D_n  \ \  \raisebox{1em}{ \xymatrix@R=2ex@C=5ex{  &&& *{\circled{2}}<3pt> \ar@{-}[d]_<{ n-1  }  \ar@/^.6pc/@{<-->}[dr]|{\vee}  \\
*{\circled{2}}<3pt> \ar@{-}[r]_<{ 1 \ \  } & *{\circled{2}}<3pt> \ar@{-}[r]_<{ 2 \ \  } & *{\circled{2}}<3pt> \ar@{.}[r]
&*{\circled{2}}<3pt> \ar@{-}[r]_>{ \ \ n} &*{\circled{2}}<3pt> \ar@{}[l]^>{   n-2} }},  \allowdisplaybreaks\\
& E_{6}  \ \  \raisebox{1em}{   \xymatrix@R=2ex@C=3.5ex{  &&  *{\circled{2}}<3pt> \ar@{-}[d]^<{ 2\ \ }   \\
*{\circled{2}}<3pt>  \ar@{-}[r]_<{ 1 \ \  } \ar@/_1.5pc/@{<-->}[rrrr]|{\vee}
& *{\circled{2}}<3pt> \ar@{-}[r]_<{ 3 \ \  } \ar@/_1pc/@{<-->}[rr]
& *{\circled{2}}<3pt> \ar@{-}[r]_<{ 4 \ \  }
&*{\circled{2}}<3pt> \ar@{-}[r]_>{ \ \ 6} &*{\circled{2}}<3pt> \ar@{}[l]^>{   5} } },
&&D_4  \ \  \raisebox{1em}{ \xymatrix@R=2ex@C=5ex{  & *{\circled{2}}<3pt> \ar@{-}[d]_<{ 3  }  \ar@/^.6pc/@{<--}[dr]|{\tvee}  \ar@/_.6pc/@{-->}[dl]|{\tvee} \\
*{\circled{2}}<3pt> \ar@{-}[r]_<{ 1 \ \  }  \ar@/_1.5pc/@{-->}[rr]|{\tvee} &*{\circled{2}}<3pt> \ar@{-}[r]_>{ \ \ 4} &*{\circled{2}}<3pt> \ar@{}[l]^>{  2} }},  \\
\end{aligned}
\end{equation}

For a  Lie algebra $\bfg$ of simply-laced finite type associated to $\Dynkinf$ and a Dynkin diagram automorphism $\sigma (\ne {\rm id})$ on $\Dynkinf$,
we denote by $\sfg$ the Lie subalgebra of $\bfg$ such that it is non-simply-laced type \cite[Proposition 7.9]{Kac} and obtained via $\sigma$:
\begin{align}\label{eq: pair of gg}
 (\sfg \; | \; (\bfg,\sigma)): \   (C_n \; | \; (A_{2n-1},\vee)), \ (B_n\; | \; (D_{n+1},\vee)), \ (F_4\; | \; (E_{6},\vee)) , \ ( G_2 \; | \; (D_4,\widetilde{\vee})).
\end{align}

\noindent
Note that there exists a natural surjective map from $I^\bfg$ to $I^\sfg$ sending $I^\bfg \ni \im \mapsto \ol{\im} \in I^\sfg$, where $\ol{\im}$ is an index in $I^\sfg$
which can be also understood as the orbit of $i$ under $\sigma$.

\subsection{Dynkin quiver} A {\it Dynkin quiver} $Q=(\Dynkin, \xi)$ of $\Dynkin$ is an oriented graph, whose underlying graph is $\Dynkin$, together with a function $\xi: \Dynkin_0 \to \Z$, called a height function of $Q$, which satisfies the following condition:
\begin{equation} \label{eq: condition of height function}
\xi_i = \xi_j+1 \qquad \text{ if } d(i,j)=1 \text{ and } i \to j \text{ in } Q.
\end{equation}

\begin{remark}
We emphasize here that \emph{not} every Dynkin diagram $\Dynkin$ has a Dynkin quiver. For instance,
 if $\Dynkin$ is of affine type $A_{2n}^{(1)}$, there is no Dynkin quiver associated with $\Dynkin$.  
 Thus, when we mention a Dynkin quiver $Q=(\Dynkin,\xi)$, it implies that
$\Dynkin$ has one (see also \cite[\S 14.1]{L94}).  
\end{remark}

Note that, since $\Dynkin$ is connected,
\begin{align} \label{eq: differ int}
\text{height functions of $Q$ differ by integers. }
\end{align}
Conversely, to a Dynkin diagram $\Dynkin$ and a function $\xi:\Dynkin \to \Z$
satisfying $|\xi_i-\xi_j|=1$ for $i,j\in I$ with $d(i,j)=1$, we can define an orientation on $\Dynkin$ to obtain a Dynkin quiver in an obvious way. Thus it is enough to specify a pair $(\Dynkin,\xi)$ of a Dynkin diagram and a height function to present a Dynkin quiver.

For a Dynkin quiver $Q=(\Dynkin,\xi)$, we call $i \in\Dynkin_0$  a \emph{source} (resp. \emph{sink}) of $Q$ (or $\xi$) if $\xi_i >\xi_j$ (resp. $\xi_i <\xi_j$) for all $j \in \Dynkin_0$ with $d(i,j)=1$.
For a Dynkin quiver $Q=(\Dynkin,\xi)$ and its source $i$, we denote by $s_iQ$ the Dynkin quiver $(\Dynkin,s_i\xi)$, where $s_i\xi$ is the height function defined as follows:
\begin{align}\label{eq: si xi}
(s_i\xi)_j = \xi_j -2 \times \delta_{i,j}.
\end{align}
We call the operation from $Q$  to $s_iQ$ the \emph{reflection of $Q$ at a source $i$ of $Q$}. Note that for Dynkin quivers $Q=(\Dynkin,\xi)$ and $Q'=(\Dynkin,\xi')$ with $\xi_i \equiv_2 \xi'_i$ for all $i \in \Dynkin_0$, there exists a sequence of reflections,
including $s_i^{-1}$ at sink, from $Q$ to $Q'$.

For a reduced expression $\uw=s_{i_1} \cdots s_{i_l}$ of $w \in \weyl$ or a sequence $\tw=(i_1,\ldots, i_l)_{i_1, \dots, i_l \in \Dynkin_0}$
of indices,  we say that $\uw$ (or $\tw$) is \emph{adapted to $Q =(\Dynkin,\xi)$} if
$$  \text{ $i_k$ is a source of $s_{i_{k-1}}s_{i_{k-2}} \cdots s_{i_1} Q$ for all $1 \le k \le l$.}  $$

\noindent
For a Dynkin quiver $Q=(\Dynkin,\xi)$, let $s_{i_1}\cdots s_{i_{n}}$ be a $Q$-adapted reduced expression of a Coxeter element. Then the height function $\xi'$ of the Dynkin quiver
$s_{i_{n}} \cdots s_{i_1}Q$ is given by
\begin{align} \label{eq: -2}
 \xi'_i = \xi_i - 2  \quad \text{for any $i \in \Dynkin_0$.}
\end{align}

Note that, for $\g$ of finite type, we can obtain a Dynkin quiver $Q=(\Dynkinf,\xi)$ of the same type by assigning orientations to edges in $\Dynkinf$, and  there exists the Coxeter element $\tau_Q$ all of whose reduced expressions are adapted to $Q$.
Note that, in finite type, there exists a unique element $w_0$ in $\sfW$  whose length is the largest. Also the element $w_0$ induces an involution $^* : I \to I$
given by $w_0(\al_i) = -\al_{i^*}$.

\begin{convention} \label{conv: height}
Throughout this paper, we take a height function $\xi$ on a finite Dynkin quiver $\Dynkinf$ such that $\xi_1 \equiv_2 0$.
\end{convention}

\smallskip

Let $Q=(\Dynkin,\xi)$ be a Dynkin quiver and $\sigma$ be a
non-trivial Dynkin diagram automorphism of $\Dynkin$ satisfying~\eqref{eq: auto cond}. We call a
Dynkin quiver $Q$ \emph{$\sigma$-fixed} if $\xi_i =
\xi_{\sigma^k(i)}$ for $0 \le k < |\sigma|$.
For a $\sigma$-fixed Dynkin quiver $Q=(\Dynkinf^\bfg,\xi)$ of finite simply-laced  type $\bfg$ and the pair
$(\bfg,\sfg)$ obtained via $\sigma$ in~\eqref{eq: pair of gg}, we obtain a Dynkin quiver $\ol{Q}=(\Dynkinf^\sfg,\oxi)$ of non-simply-laced type $\sfg$
by defining
$\oxi_\oi = \xi_\im$ for all $\im \in I^\bfg$.

\subsection{$t$-quantized Cartan matrix} \label{subsec:t-Cartan}
For an indeterminate $x$ and  integers $k \ge  l \ge 0$, we set
\begin{align*}
[k]_x \seteq \dfrac{x^k-x^{-k}}{x - x^{-1}}, \quad [k]_x! \seteq \prod_{u=1}^k [u]_x  \ \ \text{ and } \ \  \left[ \begin{matrix} k \\ l \end{matrix}\right]_x \seteq \dfrac{[k]_x!}{[k-l]_x![l]_x!}.
\end{align*}
For an indeterminate $q$ and $i \in I$, we set $q_i = q^{d_i}$ where $\sfD = {\rm diag}(d_i \in \Z_{\ge 1} \ | \ i \in I)$ satisfies \eqref{eq: sfD}.
For a given Cartan matrix $\sfC$, we set $\calI = (\calI_{i,j})_{i,j\in I}$ the \emph{adjacent matrix} of $\sfC$ by 
$\calI_{i,j}= -\delta(i \ne j)\sfc_{i,j}$.

In \cite{FR98}, the $(q,t)$-deformation of Cartan matrix $\sfC(q,t) =(\sfc_{i,j}(q,t))_{i,j \in I}$ is introduced, where
\begin{align*}
\sfc_{i,j}(q,t) \seteq (q_it^{-1}+q_i^{-1}t)\delta_{i,j} - [\calI_{i,j}]_q.
\end{align*}
Then we have two kinds of specializations of $\sfC(q,t)$, one of which is $\sfC(q)\seteq \sfC(q,1)$, called the \emph{quantum Cartan matrix}, and the another is $\usfC(t)\seteq\sfC(1,t)$, called the \emph{$t$-quantized Cartan matrix}.

Throughout this paper, we mainly consider the following symmetric matrix
\begin{align} \label{eq: usfB}
\usfB(t) \seteq \usfC(t)\sfD^{-1}.
\end{align}
Note that $\usfB(t)\vert_{t=1}=\usfB \in {\rm GL}_{|I|}(\Q)$. We regard $\usfB(t)$ as an element of ${\rm GL}_{|I|}(\Q(t))$ and denote its inverse
by $\tuB(t)=(\tuB_{i,j}(t))_{i,j\in I}$ provided it exists. Let
\begin{equation} \label{eq: inverse of tBt}
\begin{split}
\tuB_{i,j}(t) =\sum_{u\in\Z} \tfb_{i,j}(u)t^u
\end{split}
\end{equation}
be the Laurent expansion of $\tuB_{i,j}(t)$ at $t=0$. Note that $\tuB_{i,j}(t) = \tuB_{j,i}(t)$ for all $i,j \in I$.
The closed formulae of $\usfB(t)$ and $\tuB_{i,j}(t)$ for all finite types can be found in \cite{KO23B, KO23} (see also references therein).

\begin{lemma}[\cite{HL15,FM22,KO23}] \label{lem: b range}
Let $\tuB(t)$ be associated with $\sfC$ of finite type. Then, for any $i,j \in I$ and $u \in \Z$, we have
\ben
\item $\tfb_{i,j}(u)=0$ if $u \le d(i,j)$ or $d(i,j) \equiv_2 u$,
\item  $\tfb_{i,j}(d(i,j)+1)=\max(d_i,d_j)$.
\ee
\end{lemma}

For a Dynkin quiver $Q$, we choose a subset $\tDynkin_0$ of $\Dynkin_0 \times \Z$ as follows:
\begin{align*}
\tDynkin_0 \seteq \{ (i,p) \in I \times \Z \; | \; p- \xi_i \in 2\Z\}.
\end{align*}
By Convention~\ref{conv: height}, $\tDynkin_0$ does not depend on the choice of $Q$.
For $i,j \in \tDynkinf_0$, we define an {\em even function} $\teta_{i,j}:\Z \to \Z$ as follows:
\begin{align}\label{eq: teta}
\teta_{i,j}(u) = \tfb_{i,j}(u) + \tfb_{i,j}(-u)  \qquad \text{ for } u \in \Z.
\end{align}

\begin{lemma} [{\cite{B21,KO23}}]\label{lem: teta}
We have
$$   \teta_{i,j}(u-1)+ \teta_{i,j}(u+1) + \sum_{k; \; d(k,j)=1} \lan h_k,\al_j \ran \teta_{i,k}(u) = \delta_{u,1}\delta_{i,j} \times 2d_i.$$
\end{lemma}

\subsection{Valued quiver} \label{subsec: Bi-deco}
Let $\sfK$ be a (possibly infinite) countable index set with a decomposition $\sfK = \sfK_\ex \sqcup \sfK_\fr$. We call $\sfK_\ex $ the set of \emph{exchangeable indices} and $\sfK_\fr $ the set of \emph{frozen indices}.

We call an integer-valued $\sfK \times \Kex$ matrix $\tB=(b_{i,j})_{i \in \sfK, j \in \Kex}$  an \emph{exchange matrix} if it satisfies the following properties:
\begin{eqnarray} &&
\parbox{85ex}{
\bna
\item For each $j \in \sfK_\ex$, there exist finitely many $i \in \sfK$ such that $b_{i,j} \ne 0$.
\item \label{it: skew sym} Its principal part $B \seteq (b_{i,j})_{i,j \in \sfK_\ex}$ is \emph{skew-symmetrizable}; i.e., there exists a sequence $S=(\mathsf{t}_i \ |  \ i \in \sfK_\ex, \mathsf{t}_i  \in \Z_{\ge 1} )$ such that
$\mathsf{t}_i b_{i,j} = - \mathsf{t}_j b_{j,i}$ for all $i,j \in \sfK_\ex$.
\ee
}\label{eq: ex matrix}
\end{eqnarray}

For an exchange matrix $\tB$, we associate a \emph{valued  quiver} $\calQ_{\tB}$
whose set of vertices is $\sfK$ and  arrows between vertices are assigned by the following rules:
\begin{align}     \label{eq: bivalued}
\bc
\xymatrix@R=0.5ex@C=6ex{ *{\bullet}<3pt> \ar@{->}[r]^{{\ulcorner}a,b {\lrcorner} }_<{k}  &*{\bullet}<3pt> \ar@{-}[l]^<{l}}
& \text{ if  $l , k \in \Kex$, $l \ne k$, $b_{kl}=a \ge 0$ and $b_{lk}=b \le 0$}, \\
\xymatrix@R=0.5ex@C=6ex{ *{\circ}<3pt> \ar@{->}[r]^{{\ulcorner}a,0 {\lrcorner} }_<{k}  &*{\bullet}<3pt> \ar@{-}[l]^<{l}} \ \ \text{(resp.}  \  \xymatrix@R=0.5ex@C=6ex{ *{\circ}<3pt> \ar@{<-}[r]^{{\ulcorner}0,b {\lrcorner} }_<{k}  &*{\bullet}<3pt> \ar@{-}[l]^<{l}} \text{)}
& \text{ if  $l \in \Kex$, $k \in \Kfr$ and $b_{kl}=a \ge 0$ (resp. $b_{kl}=b \le 0$)}.
\ec
\end{align}
Here we do not draw an arrow between $k$ and $l$ if $b_{kl} = 0$ (and $b_{lk} = 0$ when $l, k \in \Kex$).
Note that $\circ$ denotes a vertex in $\Kfr$, and
We call $\ulcorner a,b  \lrcorner$ the \emph{value} of an arrow.

\begin{convention} \label{conv: bi-deco}
For some special values $\ulcorner a, b \lrcorner$, we will use the following scheme to draw a valued quiver for convenience: For $l , k \in \Kex$ $l \ne k$,
\ben
\item if $b_{kl}=1$ and $b_{lk}=-b<0$, use $\xymatrix@R=0.5ex@C=6ex{ *{\bullet}<3pt> \ar@{->}[r]|{<b}_<{k}  &*{\bullet}<3pt> \ar@{}[l]^<{l}} $,
\item if $b_{kl}=2$ and $b_{lk}=-b<0$, use $\xymatrix@R=0.5ex@C=6ex{ *{\bullet}<3pt> \ar@{=>}[r]|{<b}_<{k}  &*{\bullet}<3pt> \ar@{}[l]^<{l}} $,
\item if $b_{kl}=3$ and $b_{lk}=-b<0$, use $\xymatrix@R=0.5ex@C=6ex{ *{\bullet}<3pt> \ar@{=>}[r]|{<b}_<{k}  &*{\bullet}<3pt> \ar@{-}[l]|{ \; \; \; \; }^<{l}} $,
\item we usually skip $< \hspace{-.7ex} 1$ in an arrow $($when $\ulcorner a ,-1 \lrcorner$ and $1 \le a\le 3)$ for notational simplicity,
\ee
and for $l\in \Kex$ and $k \in \Kfr$,
\ben
\item[{\rm (5)}] if $b_{kl}= 1$ $($resp. $b_{kl}= -1)$, use $\xymatrix@R=0.5ex@C=6ex{ *{\circ}<3pt> \ar@{->}[r]_<{k}  &*{\bullet}<3pt> \ar@{-}[l]^<{l}} $ $($resp. $\xymatrix@R=0.5ex@C=6ex{ *{\circ}<3pt> \ar@{}[r]_<{k}  &  *{\bullet}<3pt> \ar@{->}[l]^<{l}} )$,
\item[{\rm (6)}] if $b_{kl}=2$ $($resp. $b_{kl}= -2)$, use $\xymatrix@R=0.5ex@C=6ex{ *{\circ}<3pt> \ar@{=>}[r]_<{k}  &*{\bullet}<3pt> \ar@{}[l]^<{l}} $ $($resp. $\xymatrix@R=0.5ex@C=6ex{ *{\circ}<3pt> \ar@{}[r]_<{k}  &*{\bullet}<3pt> \ar@{=>}[l]^<{l}} )$,
\item[{\rm (7)}] if $b_{kl}=3$ $($resp. $b_{kl}= -3)$, use $\xymatrix@R=0.5ex@C=6ex{ *{\circ}<3pt> \ar@{=>}[r]_<{k}  &*{\bullet}<3pt> \ar@{-}[l]^<{l}} $ $($resp.  $\xymatrix@R=0.5ex@C=6ex{ *{\circ}<3pt> \ar@{-}[r]_<{k}  &*{\bullet}<3pt> \ar@{=>}[l]^<{l}} )$.
\ee
\end{convention}

\emph{Throughout this paper, we always apply Convention~\ref{conv: bi-deco}.}

\begin{definition} \label{def: tB associated with Dynkin}
Let $\Dynkin$ be a Dynkin diagram. We set  $\tDynkin_0 \times \tDynkin_0$-matrix $\tB_{\tDynkin_0}$ whose entries $b_{(i,p),(j,s)} $ are defined as follows:
\begin{align} \label{eq: ex tDynkin}
b_{(i,p),(j,s)} =
\bc
(-1)^{\delta(s>p)}\sfc_{i,j}  &\text{ if  $|p-s|=1$ and  $i \ne j$},\\
(-1)^{\delta(s>p)}   & \text{ if  $|p-s|=2$ and  $i = j$},  \\
0 & \text{otherwise}.
\ec
\end{align}
\end{definition}
\noindent
Note that $\tB_{\tDynkin_0}$ satisfies ~\eqref{eq: ex matrix} with a sequence
 $S \seteq ( s_{i,p} \; | \;  s_{i,p} =d_i )$ and without frozen vertices.
We denote by
$\tDynkin$ the valued quiver associated  to $\tB_{\tDynkin_0}$.

\smallskip

We call the arrows $(i,p) \gets (i,p+2)$ in $\tDynkin$ the \emph{horizontal arrows}
and  the arrows between $(i,p)$ and $(j,p+1)$ for $d(i,j)=1$ the \emph{vertical\,\footnote{Visually, they are slant.} arrows}.

\begin{convention} \label{conv: bi-deco2}
We use dashed arrows
$\xymatrix@R=0.5ex@C=4ex{  \ar@{<.}[r]   & } $ for horizontal arrows in $\tDynkin$
to distinguish them with vertical arrows in $\tDynkin$.
\end{convention}

\begin{example}  \label{ex: valued quivers}
Under Conventions~\ref{conv: bi-deco} and \ref{conv: bi-deco2}, 
%
when $\bDynkinf$ is of finite type $B_{3}$,
the valued  quiver $\tbDynkinf$ is depicted as
$$
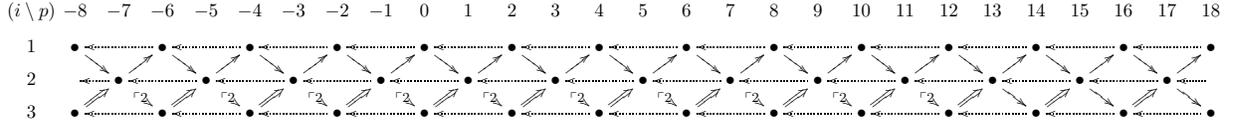
\raisebox{3mm}{
\scalebox{0.65}{\xymatrix@!C=0.5mm@R=2mm{
(i\setminus p) & -8 & -7 & -6 &-5&-4 &-3& -2 &-1& 0 & 1& 2 & 3& 4&  5
& 6 & 7 & 8 & 9 & 10 & 11 & 12 & 13 & 14 & 15 & 16 & 17& 18 \\
1&\bullet \ar@{->}[dr] && \bullet \ar@{->}[dr]\ar@{.>}[ll] &&\bullet\ar@{->}[dr]\ar@{.>}[ll]
&&\bullet \ar@{->}[dr]\ar@{.>}[ll] && \bullet \ar@{->}[dr]\ar@{.>}[ll] &&\bullet \ar@{->}[dr]\ar@{.>}[ll] &&  \bullet \ar@{->}[dr]\ar@{.>}[ll]
&&\bullet \ar@{->}[dr]\ar@{.>}[ll] && \bullet \ar@{->}[dr]\ar@{.>}[ll] &&\bullet \ar@{->}[dr]\ar@{.>}[ll]  && \bullet\ar@{->}[dr]\ar@{.>}[ll] &&
\bullet\ar@{->}[dr]\ar@{.>}[ll] && \bullet\ar@{->}[dr]\ar@{.>}[ll]  && \bullet\ar@{.>}[ll]\\
2&&\bullet \ar@{->}[dr]|{\ulcorner 2}\ar@{->}[ur] \ar@{.>}[l]&& \bullet \ar@{->}[dr]|{\ulcorner 2}\ar@{->}[ur] \ar@{.>}[ll]&&\bullet \ar@{->}[dr]|{\ulcorner 2}\ar@{->}[ur]\ar@{.>}[ll]
&& \bullet \ar@{->}[dr]|{\ulcorner 2}\ar@{->}[ur]\ar@{.>}[ll]&& \bullet\ar@{->}[dr]|{\ulcorner 2} \ar@{->}[ur]\ar@{.>}[ll]&& \bullet \ar@{->}[dr]|{\ulcorner 2}\ar@{->}[ur]\ar@{.>}[ll]&&\bullet \ar@{->}[dr]|{\ulcorner 2}\ar@{->}[ur]\ar@{.>}[ll]&
&\bullet \ar@{->}[dr]|{\ulcorner 2}\ar@{->}[ur]\ar@{.>}[ll]&&\bullet\ar@{->}[dr]|{\ulcorner 2} \ar@{->}[ur]\ar@{.>}[ll]&& \bullet \ar@{->}[dr]|{\ulcorner 2}\ar@{->}[ur]\ar@{.>}[ll]
&&\bullet \ar@{->}[dr]\ar@{->}[ur]\ar@{.>}[ll]&& \bullet \ar@{->}[dr] \ar@{->}[ur]\ar@{.>}[ll]&&\bullet \ar@{->}[dr]\ar@{->}[ur]\ar@{.>}[ll] & \ar@{.>}[l]\\
3&\bullet  \ar@{=>}[ur]&& \bullet  \ar@{=>}[ur] \ar@{.>}[ll]&&\bullet \ar@{=>}[ur]\ar@{.>}[ll]
&&\bullet  \ar@{=>}[ur] \ar@{.>}[ll]&& \bullet \ar@{=>}[ur]\ar@{.>}[ll] &&\bullet  \ar@{=>}[ur]\ar@{.>}[ll]&&  \bullet  \ar@{=>}[ur]\ar@{.>}[ll]
&&\bullet  \ar@{=>}[ur]\ar@{.>}[ll] && \bullet  \ar@{=>}[ur]\ar@{.>}[ll]&&\bullet  \ar@{=>}[ur] \ar@{.>}[ll]&& \bullet \ar@{=>}[ur]\ar@{.>}[ll]&&
\bullet \ar@{=>}[ur]\ar@{.>}[ll] && \bullet \ar@{=>}[ur]\ar@{.>}[ll]  && \bullet \ar@{.>}[ll] }}}
$$
\end{example}

\begin{remark}
The valued quivers for simply-laced finite types coincide with the infinite quivers in \cite[Section 2.1.3]{HL16} where the infinite quivers are denoted by $\Gamma$.
\end{remark}

\begin{definition}  [{cf. \cite[Definition 5.5]{FHOO}}] \label{def: convex subset} \hfill
\ben
\item \label{def: hDynkin} We denote by $\hDynkin$ the quiver obtained from $\tDynkin$ by removing all horizontal arrows. We call $\hDynkin$ the \emph{valued repetition quiver} of $\Dynkin$\footnote{When we replace valued arrows with usual arrows, it is the usual repetition quiver $\widehat{\Dynkin}$ (see \cite{KO23} for non-simply-laced types).}.
\item \label{def: convexity} A subset $\calR \subset \tDynkin_0=\hDynkin_0$ is said to be \emph{convex} if it satisfies the following condition: For any oriented path $(x_1 \to x_2 \to \cdots \to x_l)$ consisting of  (vertical) arrows
in $\hDynkin$, we have $\{ x_1,x_2,\ldots,x_l\} \subset \calR$ if and only if $\{ x_1,x_l\} \subset \calR$.
\item We say that a convex subset $\calR \subset \tDynkin_0$ has a \emph{upper bound} if there exists $\max(p \ | \ (i,p) \in \calR )$ for each $i \in \Dynkin_0$.
\item For a convex subset $\calR \subset \tDynkin_0$, we set $\calR_\fr \seteq \{ (i,p) \; | \; p= \min(k \in \Z \; | \; (i,k) \in \calR ) \}$
and  $\calR_\ex \seteq \calR \setminus \calR_\fr$. We denote by ${}^\calR\tDynkin$ the valued quiver
associated to ${}^\calR\tB \seteq (b_{(i,p),(j,s)})_{(i,p) \in \calR, (j,s) \in \calR_\ex}$.
\item For a height function $\xi$ on $\Dynkin$, let $\lxi\tB \seteq (b_{(i,p),(j,s)})_{(i,p), (j,s) \in \lxi\tDynkin_0}$ and denote by $\lxi\tDynkin$ the valued quiver associated to $\lxi\tB$, where
\begin{align*}
\lxi\tDynkin_0 \seteq \{ (i,p) \in \tDynkin_0 \ | \  p \le \xi_i \}.
\end{align*}
Note that $\lxi\tDynkinf_0$ is a convex subset of $\tDynkinf$ for any height function $\xi$ on $\Dynkinf$.
\ee
\end{definition}

\section{\texorpdfstring{$t$}--characters of quantum loop algebra and virtual Grothendieck rings} \label{sec:t-chars and vGrothendieck}

\noindent
In this section, we first review the important properties of $t$-characters of finite-dimensional representations over quantum loop algebra briefly
(see~\cite{FR99, FM01, Nak03, H04, H06} for more details).
Then we recall the  virtual Grothendieck ring  $\frakK(\g)$
for any finite type $\g$ (see \cite{FHR21,KO23} for non-simply-laced types). 

\subsection{Quantum loop algebras} 
Let $t$ be an indeterminate.
We denote by $\sfk \seteq \ol{\Q(t)}$  the algebraic closure of the field $\Q(t)$ inside $\bigcup_{m \in \Z_{\ge 0}}
\ol{\Q}\lb t^{1/m} \rb $.  Let $\bfg$ be a complex finite-dimensional simple Lie algebra of simply-laced type. Note that, in this case, we can identify $\sfC(q)$ with $\usfC(t)$ by exchanging $q$ with $t$.

\smallskip
\emph{Throughout this paper, we use \textbf{bold symbols} to emphasize that those symbols are of simply-laced finite types. We also use $\im,\jm$ for indices in $I^\bfg$ for the same purpose.}
\smallskip

We denote by $U_{t}(\calL\bfg)$ the quantum loop algebra associated to $\bfg$, which
is the $\sfk$-algebra given by the set of infinite generators, called the Drinfeld generators,
subject to certain relations \cite{Dr88, Beck94}.
The quantum loop algebra $U_{t}(\calL\bfg)$ is a quotient of the corresponding
(untwisted) quantum affine algebra $U_t'(\widehat{\bfg})$ and hence has a Hopf algebra structure.

\medskip

\subsection{Finite dimensional modules and their $t$-characters} \label{subsec: review for fd mod and t-char}
We denote by $\scrC_\bfg$ the category of
finite-dimensional $U_{t}(\calL\bfg)$-modules of type $\mathbf{1}$.
The category $\scrC_\bfg$ is a $\sfk$-linear rigid non-braided monoidal category.
We say that $V$ and $W$ \emph{commute} if $V \otimes W \simeq W \otimes V$
as $U_{t}(\calL\bfg)$-modules.  We denote by $K(\scrC)$ the Grothendieck ring of $\scrC_\bfg$. Note that the set of simple objects in
$K(\scrC_\bfg)$ are parameterized by the set $(1 + z\sfk[z])^{I^\bfg}$ of $I^\bfg$-tuples of monic polynomials, which is called \emph{Drinfeld polynomials.}

In this paper, we usually consider the \emph{skeleton} subcategory $\Cbz$ of $\scrC_\bfg$. The subcategory $\Cbz$ contains every \emph{prime} simple module in $\scrC_\bfg$ up to \emph{parameter shifts}. To explain $\Cbz$, we need to consider
the Laurent polynomial $\calY$ generated by the set of variables $\{ Y^{\pm 1}_{\im,p} \}_{(\im,p) \in \widetilde{\sDynkinf}_0}$. Let us denote by  $\calM$ (resp. $\calM_+$ and $\calM_-$) the set of all  monomials (resp. dominant monomials and anti-dominant monomials)
of $\calY$.   For a monomial $\bfm$ in $\calY$, we write
\begin{align}\label{eq: uip}
\bfm = \prod_{(\im,p) \in \widetilde{\sDynkinf}_0} Y_{\im,p}^{u_{\im,p}(\bfm)}
\ \ \text{ and } \ \
\bfm_- = \prod_{(\im,p) \in \widetilde{\sDynkinf}_0} Y_{\im,p}^{-u_{\im,p}(\bfm)}
\end{align}
with $u_{\im,p}(\bfm) \in \Z$. For each $\bfm \in \calM_+$, we denote by $L(\bfm)$ the simple module in $\scrC$ whose Drinfeld polynomial is $\big(\prod_p (1-q^p)^{u_{\im,p}(\bfm)} \big)_{\im \in I^\bfg}$. Then the subcategory $\Cbz$ can be characterized by the Serre subcategory
of $\scrC_\bfg$ generated by $\{ L(\bfm) \ | \  \bfm \in \calM_+ \}$. Note that $\Cbz$ is a monoidal rigid subcategory of $\scrC_\bfg$.
In \cite{FR99}, Frenkel-Reshetikhin proved that there exists an injective ring homomorphism
$$\chit:K(\Cbz) \to \calY,$$
called \emph{the $t$-character homomorphism}\footnote{  It is usually called the $q$-character homomorphism in the literature.  }.
The existence of $\chit$ tells us that the Grothendieck ring $K(\scrC_\bfg)$ is commutative, even though $\scrC_\bfg$ is not braided.

For an interval $[a,b] \subset \Z$, $\im \in I^\bfg$, $k \in \Z_{\ge 1}$ and $(\im,p) \in \tDynkinf_0$, we set dominant monomials
\begin{equation} \label{eq:dominant monomials for KR module}
\bfm^{(\im)}[a,b]  \seteq \displaystyle\prod_{ (\im,s) \in \widetilde{\sDynkinf}_0; \ s \in [a,b]} Y_{\im,s} \quad \text{ and }\quad
\bfm^{(\im)}_{k,p} \seteq \displaystyle \prod_{s=0}^{k-1} Y_{\im,p+2s},
\end{equation}
and $\bfm^{(\im)}(a,b]$, $\bfm^{(\im)}[a,b)$, and $\bfm^{(\im)}(a,b)$ are defined similarly.

The simple module $L(\bfm^{(\im)}[p,s])$ $(p \le s)$ is called a \emph{Kirillov--Reshetikhin\,$($KR$)$ module}. When $p=s$ and $(i,p) \in \tDynkinf_0$, we call $L(Y_{\im,p})$ a \emph{fundamental module}.
Note that the Grothendieck ring $K(\Cbz )$ is a polynomial ring in the isomorphism classes of the fundamental modules $L(Y_{\im,p})$ \cite{FR99}.

\smallskip

For $\im \in I^\bfg, a \in \sfk^\times$, we set
$$
A_{\im,a} \seteq Y_{\im,at^{-1}} Y_{\im,at}   \prod_{\jm \colon d(\im,\jm)=1} Y_{\jm,a}^{-1}
=Y_{\im,at^{-1}} Y_{\im,at}
  \prod_{\jm \ne \im} Y_{\jm,a}^{\bfc_{\jm,\im}}.
$$

Note that there is an ordering $\leN$ on the set of monomials, called the \emph{Nakajima order}, defined as follows:
\begin{equation} \label{eq: Nakajima order}
\text{$\bfm \leN \bfm'$ if and only if
$\bfm^{-1}\bfm'$ is a product of elements in $\{A_{\im,a} \mid \im \in I^\bfg, a \in \sfk^\times \}$.}
\end{equation}

\begin{theorem} [\cite{FR98,FM01}]
For each dominant monomial $\bfm$, the monomials appearing in $\chit(L(\bfm))-\bfm$
are strictly less that $\bfm$ with respect to $\leN$.
\end{theorem}

The $t$-characters of KR-modules satisfies a system of functional equations called $T$-systems:

\begin{theorem}[{\cite[Theorem 1.1]{Nak03}}] $($See also \cite[Theorem 3.4]{H06}.$)$
For each $(i,p),(i,s) \in \tDynkinf^\bfg_0$ with $p \le  s$, we have
\fontsize{10}{10}
\begin{equation}\label{eq: T-sys}
\begin{aligned}
\chit \bl L(\bfm^{(\im)}[p,s) )\br \chit \bl L(\bfm^{(\im)}(p,s]) \br & = \chit\bl L(\bfm^{(\im)}[p,s]) \br \chit\bl L(\bfm^{(\im)}(p,s))
+ \prod_{\jm: \; d(\im,\jm)=1}  \chit\bl L(\bfm^{(\jm)}(p,s)). 
\end{aligned}
\end{equation}
\fontsize{11}{11}
\end{theorem}


Let $\xi$ be a height function on $\Dynkinf^\bfg$.
We  denote by $ \Mxi_+$ the set of all dominant monomials
in the variables $Y_{\im,p}$'s for $(\im,p) \in \lxi\tDynkinf_0$.

\begin{definition}
We define the subcategory $\scrC^\xi_\bfg$ as the Serre subcategory of $\scrC_\bfg$ such that
$\Irr \; \scrC^\xi_\bfg = \{ L(\bfm) \ |  \ \bfm \in \Mxi_+\}$.
\end{definition}

Since $\lxi\tDynkinf_0$ is a convex subset of $\tDynkinf_0$, we have
the following proposition:

\begin{proposition}
The category $\scrC^\xi_\bfg$ is a monoidal  subcategory of $\scrC_\bfg$.
\end{proposition}

\begin{proof}
This assertion follows from the same argument of the proof of \cite[Proposition 3.10]{HL16}.
\end{proof}

\subsection{Truncation}
We denote by
$\Yxi$ the Laurent polynomial ring generated by $Y_{\im,p}$'s for $(\im,p) \in \lxi\tDynkinf_0$.
We define a linear map $ ( \cdot )_{\le \xi}: \calY \to \Yxi$ by sending
the monomials which contain some $Y_{\im,p}$ with $(\im,p) \not\in \lxi\tDynkinf$ to zero and by keeping all the other terms.

\begin{proposition}
For a height function $\xi$, the $\Z$-linear map $(  \cdot  )_{\le \xi}: K(\scrC^\xi_\bfg) \to \Yxi$ given by
\begin{align*}
[V] \mapsto \lxi \chi_{\mspace{-2mu}\raisebox{-.5ex}{${\scriptstyle{t}}$}}(V)  \seteq (( \cdot )_{\le \xi}  \circ \chit )(V)
\end{align*}
gives an injective ring homomorphism $K(\scrC^\xi_\bfg) \hookrightarrow \Yxi$.
\end{proposition}

\begin{proof}
We can prove the assertion in the same way as in the proof of \cite[Proposition 6.1]{HL10}.
\end{proof}

\subsection{(Virtual) Grothendieck rings} \label{subsec: VGR}
Recall that when $\bfg$ is of simply-laced finite type, the $t$-character homomorphism $\chit$ is an injection from $K(\scrC^0_\bfg)$ into $\calY^\bfg$.
Thus we can identify $K(\scrC^0_\bfg)$ with
\begin{align*}
\frakK(\bfg) \seteq \chit \bl K(\scrC^0_\bfg)\br.
\end{align*}
We call $\frakK(\bfg)$ the Grothendieck ring of type $\bfg$ also.

\begin{proposition}[{\cite[Corollary 5.7]{FM01}}]
When $\bfg$ is of simply-laced type, we have
$$ \frakK(\bfg) = \bigcap_{\im \in I^\bfg} \left(\Z[Y^{\pm 1}_{\jm,l} \ |  \  (\jm,l) \in \tDynkinf^\bfg_0, \jm \ne \im]
\otimes \Z[Y_{\im,l}(1+A_{\im,l+1}^{-1})  \ |  \  (\im,l)  \in \tDynkinf^\bfg_0] \right) \subsetneq \calY^\bfg.$$
\end{proposition}

\medskip

Now we move on to non-simply-laced finite types. For $\sfg$ associated with $ (\bfg,\sigma)$ in~\eqref{eq: pair of gg}, we consider the Laurent polynomial ring
defined as follows: We first set
\begin{align*}
\calY^\sfg \seteq \Z[X_{i,p}^{\pm 1} \ | \  (i,p) \in \tDynkinf_0^\sfg].
\end{align*}
Then there exists a surjective ring homomorphism
\begin{align} \label{eq: osigma}
 \ol{\sigma}: \calY^\bfg \To \calY^\sfg \qquad \text{sending} \qquad   Y_{\sigma^k(\im),p} \longmapsto X_{\ol{\im},p}
\end{align}
for any $(\im,p) \in \tDynkinf^\bfg_0$ and  $0 \le k < |\sigma|$ (see Convention~\ref{conv: height}). Finally, we set
\begin{align*}
\frakK(\sfg) \seteq \ol{\sigma}(\frakK(\bfg))
\end{align*}
and call it the \emph{virtual Grothendieck ring} of type $\sfg$.
We call $\osigma(L(\bfm))$  the \emph{folded $t$-character} of $L(\bfm)$.

\medskip

Now we would like to unify the expression for $\frakK(\g)$ for  \emph{any} finite type $\g$ by replacing variables $Y_{i,p}$'s with $X_{i,p}$'s.
Let $\calX^\g$ be the Laurent polynomial ring
$\Z[X_{i,p}^{\pm 1} \ | \ (i,p)  \in \tDynkinf^\g_0]$. For $(i,p+1) \in \tDynkinf_0^\g$, we set
\begin{align} \label{eq: Bip}
\begin{split}
B_{i,p} \seteq X_{i,p-1}X_{i,p+1}  \prod_{j \colon d(i,j)=1} X_{j,p}^{\sfc_{j,i}}.
\end{split}
\end{align}

\begin{definition} \cite[\S 3.4]{FHR21} We define the commutative ring
\begin{align}\label{eq: VGR}
\begin{split}
 \frakK(\g) = \bigcap_{i \in I_\g} \left(\Z[X^{\pm 1}_{j,l} \ |  \  (j,l) \in \tDynkinf^\g_0, j \ne i]
\otimes \Z[X_{i,l}(1+B_{i,l+1}^{-1})  \ |  \  (i,l)  \in \tDynkinf^\g_0] \right) \subsetneq \calX^\g.
\end{split}
\end{align}
\end{definition}

\begin{remark} \label{rmk: sometimes}
Even though, we unify the expression for $\frakK(\g)$ by using $X_{i,p}$, $\calX$ and $B_{i,p}$, we sometimes use $Y_{\im,p}$, $\calY$ and $A_{\im,p}$ to emphasize that they are associated with $\bfg$ of simply-laced finite type.
\end{remark}

\begin{theorem} [{\cite[Proposition 3.3, Theorem 4.3]{FHR21}}] \label{thm: FHR} \hfill
\ben
\item Every element of $\frakK(\g)$ is characterized by the multiplicities of the dominant monomials contained in it.
\item For each $m \in \calM_+$, there is a unique element $F(m)$ of $\frakK(\g)$ such that
$m$ is the unique dominant monomial of $F(m)$ with its coefficient $1$. Therefore we have a
basis $\{ F(m) \ | \  m \in \calM_+^\g \} $ of $\frakK(\g)$
parameterized by dominant monomials $m$.
\item \label{it: surjection} For each pair $(\bfg,\sfg)$ obtained via $\sigma$, the map $\osigma$ induces a surjective ring homomorphism from $\frakK(\bfg)$ to  $\frakK(\sfg)$.
\ee
\end{theorem}

An $\calX$-monomial $\bfm$  is said to be \emph{right-negative} if the factors $X_{j,l}$ appearing in $m$, for which $l$ is maximal, have negative powers.

\begin{corollary} \label{cor: some F to F}
  For each pair $(\bfg,\sfg)$ obtained via $\sigma$ and $\bfm \in \calM_+^\bfg$, assume that
\begin{align} \label{eq: right negative-condition}
\text{every monomial in $F(\bfm)-\bfm$ is right-negative. }
\end{align}
Then $\osigma( F(\bfm) ) =  \FM{\osigma(\bfm) } \in \frakK(\sfg)$.
\end{corollary}

\begin{proof}
By Theorem~\ref{thm: FHR}~\ref{it: surjection} and ~\eqref{eq: right negative-condition}, $\osigma( F(\bfm))$ is an element in $\frakK(\sfg)$ containing the unique dominant monomial $\osigma(\bfm)$.
Thus our assertion follows.
\end{proof}

\begin{example} \label{ex: Fq neq osigma Ft in general}
For finite $A_5$-type, $\osigma(F(Y_{4,-2}Y_{2,0}))$ does not coincide with
$F(X_{2,-2}X_{2,0})$ of finite type $C_3$, since $\osigma(F(Y_{4,-2}Y_{2,0}))$ does not satisfies~\eqref{eq: right negative-condition}. More precisely,
$F(Y_{4,-2}Y_{2,0})$ contains $Y_{3,-1}Y_{5,-1}Y_{2,0}Y_{4,0}^{-1}$. On the other hand,  $\osigma(F(Y_{2,-2}Y_{2,0})) = F(X_{2,-2}X_{2,0})$.
\end{example}

Note that if $\bfm,\bfm' \in \calM^\bfg$ with $\bfm \leN \bfm'$, then we have
\begin{align}  \label{eq: preserving Nakajima ordering}
\osigma(\bfm) \leN \osigma(\bfm') \in \calM^\sfg.
\end{align}

It is proved in~ \cite{FM01,H06} that, for $\bfm^{(\im)}[p,s] \in \calM_+^\bfg$, $F(\bfm^{(\im)}[p,s])$ satisfies the condition in~\eqref{eq: right negative-condition}  and
\begin{align*}
F(\bfm^{(\im)}[p,s]) = \chit(L(\bfm^{(\im)}[p,s])).
\end{align*}
Thus we have
\begin{align} \label{eq: F to F KR}
\osigma(F(\bfm^{(\im)}[p,s])) = F(m^{(i)}[p,s])
\end{align}
and~\eqref{eq: T-sys} is changed into the following form:
For any finite type $\g$ and $(i,p),(i,s) \in \tDynkinf_0$ with $p \le  s$, we have
\begin{align}\label{eq: v T-sys}
\FM{ m^{(i)}[p,s) }\FM{m^{(i)}(p,s]} = \FM{m^{(i)}[p,s]}\FM{m^{(i)}(p,s)} + \hspace{-1.5ex} \prod_{j; \; d(i,j)=1}  \hspace{-1.5ex}  \FM{m^{(j)}(p,s)}^{-\sfc_{j,i}}.
\end{align}
We call~\eqref{eq: v T-sys} the \emph{folded $T$-systems}.

\begin{definition} \label{def: truncated virtual Grothendieck ring}  
\mbox{}
\ben
\item For a height function $\xi$ on $\Dynkinf^\bfg$ of simply-laced finite type, we set $$\lxi\frakK(\bfg) \seteq \lxi\chit(K(\Cxi)).$$
\item
For a height function $\xi$ on $\Dynkinf^\sfg$ of non-simply-laced finite type, we set
\begin{align*}
\lxi\frakK(\sfg) \seteq \osigma\bl \luxi\frakK(\bfg) \br,
\end{align*}
where $\uxi$ is the $\sigma$-fixed height function on $\Dynkinf^\bfg$  such that
$$\uxi_{\sigma^k(\im)} = \xi_{\oi} \qquad \text{ for any } 0 \le k < |\sigma| \text{ and } \im \in \sigma^{-1}(\oi).$$
\ee
\noindent
We call $\lxi\frakK(\sfg)$ the \emph{truncated virtual Grothendieck ring}
and  $\lxi\ochi_t(\bfm)$  the \emph{truncated folded $t$-character} of $L(\bfm)$ with respect to $\xi$, defined as below:
$$\xymatrix@R=0.5ex@C=8ex{K(\Cuxi^\bfg)  \ar@/^1.5pc/[rr]^{\lxi\ochi_t} \ar[r]_{\luxi\chit}& \luxi\frakK(\bfg) \ar[r]_{\osigma} &  \lxi\frakK(\sfg) } $$
\end{definition}

\begin{remark} \label{rem: relation with FHR} 
Let $G$ be a simply-connected complex Lie group associated with $\sfg$ of \emph{non-simply-laced type}.
In \cite{FHR21}, the authors formulate (conjectural) \emph{folded integrable models} of $\sfg$ corresponding to \emph{folded Bethe Ansatz equations}.
Then $\frakK(\sfg)$, denoted by $\mathcal{K}_t^-(\sfg)$ in \cite{FHR21}\footnote{In our introduction, we use $\mathcal{K}^-(\sfg)$ instead.}, plays the role of describing the spectra of the transfer-matrix $t_V(z,u)$ with a finite-dimensional $U_{t}(\calL\bfg)$-module $V$ in the folded integrable model, as in the role of  $\frakK(\bfg) \simeq K(\scrC_\bfg^0)$  in the integrable models for simply-laced types (cf.~\cite{FH15, FHR21} for more details).
We remark that our main interest is to study the structure of the \emph{quantization of $\frakK(\sfg)$} introduced independently in \cite{KO23} with other motivations related to \emph{canonical basis} and \emph{quantum cluster algebra structure}. 
In contrast, the authors of \cite{FHR21} mainly focus on a study of the folded integrable models associated with $\sfg$.
It would be interesting to find connections between our results and those in \cite{FHR21}.
\end{remark}

\section{Quantization} \label{sec:quantization}

\noindent
In this section, we quantize the Laurent polynomial ring $\calX$ with the resulting ring denoted by $\calX_q$, via the inverse matrix $\tuB(t)$ of \eqref{eq: usfB} associated with $\usfC(t)$ following \cite{KO23} (see also \cite{FHR21}), and define its subalgebra $\frakK_q(\g)$ that is regarded as a quantization of $\frakK(\g)$.

\subsection{Quantum torus}
Let $q$ be an indeterminate.
Let us recall that $\tfb_{i,j}(u)$ ($u \in \Z$) in \eqref{eq: inverse of tBt} and the even function $\teta_{i,j}: \Z \to \Z$ defined in~\eqref{eq: teta}.
\begin{definition}[{\cite{Nak04,VV03,H04,KO23}}] \label{def: calX_q}
Let $(\calX_q,*)$ be the $\Z[q^{\pm \frac{1}{2}}]$-algebra with the
generators $\{ \tX_{i,p}^{\pm 1} \; | \;  (i,p) \in \tDynkinf_0 \}$
with the defining relations
$\tX_{i,p}  * \tX_{i,p}^{-1}=\tX_{i,p}^{-1} *   \tX_{i,p}=1$  and 
$\tX_{i, p}  *  \tX_{j, s} = q^{\ucalN(i,p;j,s)}\tX_{j,s}  *  \tX_{i,p}$,
%
where $(i,p)$, $(j,s) \in \tDynkinf_0$ and
\begin{equation}\label{eq: ucalN}
\begin{aligned}
\ucalN(i,p;j,s) 
\seteq
\tfb_{i,j}(p-s-1)-\tfb_{i,j}(s-p-1)-\tfb_{i,j}(p-s+1)+\tfb_{i,j}(s-p+1).
\end{aligned}
\end{equation}
We call $\calX_q$  the \emph{quantum torus associated with $\usfC(t)$} (see Definition~\ref{def: quantum torus} below).
\end{definition}

\begin{remark}
For simply-laced finite types,
the quantum torus $\calX_q$ was already defined in \cite{Nak04, VV03, H04}, whereas for non-simply-laced finite types, it is introduced in \cite{KO23} very recently.
\end{remark}

Note that since $\tuB(t)$ is symmetric,
\begin{align*}
\ucalN(i,p;j,s) = \ucalN(j,p;i,s)=-\ucalN(i,s;j,p)=-\ucalN(j,s;i,p),
\end{align*}
and it follows from Lemma \ref{lem: b range} that
\begin{align}\label{eq: pos}
\ucalN(i,p;j,s) = \tfb_{i,j}(p-s-1)-\tfb_{i,j}(p-s+1)   \qquad  \text{ if } p >s.
\end{align}
Moreover, for $p \in \Z$ and $i ,j \in \Dynkinf_0$ such that $ (i,p),(j,p) \in \tDynkinf_0$,  Lemma~\ref{lem: b range} tells that
\begin{align} \label{eq: commute}
  \tX_{i,p} * \tX_{j,p} = \tX_{j,p} * \tX_{i,p}.
\end{align}

By specializing $q$ at $1$, the quantum torus $\calX_q$ recovers the commutative
Laurent polynomial ring $\calX$, while $\calX_q$ is non-commutative;
i.e., there exists a $\Z$-algebra homomorphism $\evq: \calX_q \to
\calX$ given by
$q^{\frac{1}{2}} \mapsto 1$ and $\tX_{i,p} \mapsto X_{i,p}$.

We say that $\tm \in \calX_q$ is a {\em $\calX_q$-monomial} if it is a product of the generators $\tX_{i,p}^{\pm1}$ and $q^{\pm \frac{1}{2}}$.
For a $\calX_q$-monomial $\tm \in \calX_q$, we set $u_{i,p}(\tm) \seteq u_{i,p}\bl \evq(\tm) \br$ (see~\eqref{eq: uip}).
An $\calX_q$-monomial $\tm$ is said to be \emph{right-negative} if  $\ev_{q=1}(\tm)$ is right-negative.
Note that a product of right negative $\calX$-monomials (resp. $\calX_q$-monomials) is right negative.
A $\calX_q$-monomial $\tm$ is called \emph{dominant} if
$\evq(\tm)$ is dominant.  Moreover, for $\calX_q$-monomials $\tm,\tm'$ in $\calX_q$, we define
$$   \tm \leN \tm'  \quad\text{ if and only if } \quad\evq(\tm) \leN \evq(\tm').$$

For $i \in \Dynkinf_0$, we call $\calX$-monomial $m$ (resp.~$\calX_q$-monomial $\tm$) {\it $i$-dominant} if  $u_{i,p}(m) \ge 0$ (resp.~$u_{i,p}(\tm) \ge 0$) for all   $p$ such that $(i,p) \in \tDynkinf_0$.
For $J \subset \Dynkinf_0$, we call $\calX$-monomial $m$ (resp.~$\calX_q$-monomial $\tm$) {\it $J$-dominant} if  $m$ (resp.~$\tm$) is $j$-dominant for all $j \in J$.
For monomials $\tm,\tm'$ in $\calX_q$, we define
\begin{align} \label{eq: ucalN for two monomials}
\ucalN(\tm,\tm') \seteq \sum_{(i,p),(j,s) \in\tDynkin_0} u_{i,p}(\tm)u_{j,s}(\tm')\ucalN(i,p;j,s).
\end{align}

There exists the $\Z$-algebra anti-involution $\ol{( \cdot )}$ on $\calX_q$ (\cite{H04,KO23}) given by
\begin{align} \label{eq: bar involution}
q^{\frac{1}{2}} \mapsto  q^{-\frac{1}{2}}, \qquad \tX_{i,p}  \mapsto q_{i}\tX_{i,p}.
\end{align}
Thus, for any $\calX_q$-monomial $\tm \in \calX_q$, there exists a unique $r \in \frac{1}{2}\Z$ such that
$q^{r}\tm$ is $\overline{( \cdot )}$-invariant.
A monomial of this form is called \emph{bar-invariant} and denoted by $\utm$.
For an example,
\begin{align*}
\text{$ \mathsf{X}_{i,p} \seteq q^{\frac{d_i}{2}}\tX_{i,p}$ is bar-invariant. }
\end{align*}
More generally, for a family $\left( u_{i,p} \ \left| \  (i,p) \in \tDynkinf_0
\right.\right)$ of integers with finitely many non-zero components, the expression
\begin{equation}   \label{eq: inv mono1}
\begin{split}
q^{\frac{1}{2} \sum_{(i,p)<(j,s)}   u_{i,p}u_{j,s} \ucalN(j,s;i,p) }  \st_{ (i,p) \in \widetilde{\sDynkinf}_0  }^{\to} \mathsf{X}_{i,p} ^{u_{i,p}}
\end{split}
\end{equation}
does not depend on the choice of an ordering on $\tDynkinf_0$ and is bar-invariant.

\begin{remark} \label{rem: bar-invariant} 
Note that the relations in Definition~\ref{def: calX_q} do not change when we replace $\tX_{i,p}$ with $\sfX_{i,p}$, and
$\utm$ depends only on $\evq(\tm)$. Therefore, for every monomial $m$ in $\calX$, we denote by $\um$ the bar-invariant monomial in $\calX_q$
corresponding to $m$.
Also the notation $Y_{i,p}$ of $(\bfY_{t},*)$  in \cite[Section 3]{HL15} corresponds to $\sfX_{i,p}$, the bar-invariant monomial, in this paper.
\end{remark}

For $(i,p)\in \tDynkinf_0$, we set
\begin{equation} \label{eq: variable tB}
 \tB_{i,p}  \seteq \underline{B_{i,p}} \ \in \calX_q.
\end{equation}

\begin{definition}
Let $\Bq^{-}$ be the $\Z[q^{\pm 1/2}]$-subalgebra of $\calX_q$ generated by $\tB_{i,p}^{-1}$'s for $(i, p) \in I \times \Z$.
For $k \in \Z_{\ge 1}$, we denote by $\Bq^{-k}$ the $\Z[q^{\pm 1/2}]$-span of the monomials  $\st^\to_{1 \le s \le k}\tB_{i_s,p_s}^{-1}$. 
\end{definition}

For bar-invariant $\calX_q$-monomials $\umn{1}$ and $\umn{2}$, we set $\umn{1} \cdot \umn{2} \seteq \ul{m_1m_2}$, and for
$\um_k$ $(k \in \Z_{\ge 1})$, we set
\begin{equation} \label{eq: dot product}
\prod_{ k }    \umn{k}   \seteq  \ul{\prod_{k} m_k}.
\end{equation}

\begin{definition}  [{cf.\ \cite[Definition 5.5]{FHOO}}] \label{def: quantum tori XqQ}
For a  subset $\sfS \subset  \tDynkinf_0$, we   denote by  ${}^\sfS \hspace{-.3ex} \calX_{q}$
the   quantum subtorus  of $\calX_{q}$ generated by $\tX^{\pm1}_{i,p}$ for $(i,p) \in \sfS \subset \tDynkinf_0$.
In particular, for a height function $\xi$ on $\Dynkinf$, we denote by $\lxi \calX_q$ the  quantum subtorus generated by $\tX^{\pm 1}_{i,p}$ for $(i,p) \in \lxi\tDynkinf_0$.
\end{definition}

\begin{proposition} [{\cite[Proposition 5.7]{KO23}}] \label{prop: YA com}
For $i,j\in I$ and $p,s,t,u \in \Z$ with $(i,p), (j,s+1), (i,t+1), (j,u+1) \in \tDynkinf_0$, we have
$$
\tX_{i,p} *\tB_{j,s}^{-1} = q^{\,\be(i,p;j,s)}\,\tB_{j,s}^{-1}*\tX_{i,p}
\quad \text{and} \quad
\tB^{-1}_{i,t}*\tB_{j,u}^{-1} = q^{\,\al(i,t;j,u)}\,\tB_{j,u}^{-1}*\tB_{i,t}^{-1}.
$$
Here,
\eq
\be(i,p;j,s) &&= \delta_{i,j}(-\delta_{p-s,1}+\delta_{p-s,-1})(\al_i,\al_i),
\label{eq;beta} \allowdisplaybreaks \\ [1ex]
\al(i,t;j,u)  &&=
\bc
\pm(\al_i,\al_i)&\text{if $(i,t)=(j,u\pm2)$,}\\
\pm2(\al_i,\al_j)&\text{if $d(i,j)=1$ and $t=u\pm1$,}\\
0&\text{otherwise.}\ec\label{eq;alpha}
\eneq
\end{proposition}

\subsection{Quantization $\frakK_q(\g)$ of $\frakK(\g)$} \label{subsec:quantization}
We briefly recall the construction of $\frakK_q(\g)$, defined in \cite{Nak04,VV03,H04,KO23},  by mainly following the argument in \cite{H03,H04}.
%
For each $i \in I$, we define the free $\calX_q$-left module
\begin{equation} \label{eq:free Xq-module}
\lcalX_{i,q}   \seteq \soplus_{ r: \; (i,r) \in \widetilde{\sDynkinf}_0 }    \calX_q  \bcdot \ts_{i,r}
\end{equation}
whose basis elements are denoted by $\ts_{i,r}$.
We also regard $\lcalX_{i,q}$ as a $\calX_q$-bimodule by defining right $\calX_q$-module action $\bcdot$ as follows:
\begin{align} \label{eq: right action on LXiq}
	 \ts_{i,r} \!\bcdot \tm = q_i^{-2u_{i,r}(\tm)}  \tm \bcdot \ts_{i,r},
\end{align}
where $\tm$ is an $\calX_q$-monomial (see Remark \ref{rem: some remark for screenings}, cf.~\cite[Lemma 4.6]{H04}).
Let $\calX_{i,q}$ be the quotient of $\lcalX_{i,q}$ by the $\calX_q$-submodule generated by the elements
\begin{equation} \label{eq:generators of submodules}
\tB_{i,r+1} \; \ts_{i,r} - q_i \ts_{i,r+2} \quad \text{ for }  (i,r) \in \widetilde{\Dynkinf}_0 .
\end{equation}
By following arguments in  \cite[Proposotion 4.8]{H04} and \cite[Lemma 4.3.1]{B21}, we have the following lemma:

\begin{lemma}   \label{lem:Xiq is free}
For each $l$ with $(i,l) \in \tDynkinf_0$,
the $\calX_q$-left module $\calX_{i,q}$ is free over any $\{ \ts_{i,r_0} \}$, where $(i,r_0) \in \tDynkinf_0$.
\end{lemma}

For all $i \in I$, we define
\begin{equation} \label{eq: ith q-screening operator}
S_{i,q} : \xymatrix@R=0.5ex@C=8ex{\calX_q  \ar@/^1pc/[rr]  \ar[r]_{\tS_{i,q}} &  \lcalX_{i,q}  \ar@{->>}[r]  &   \calX_{i,q}},
\end{equation}
where each map is defined as follows (recall \eqref{eq:free Xq-module} for definition of $\lcalX_{i,q}$):
\bna
\item the map  $\tS_{i,q}$    is defined by
\begin{align*}
	\tS_{i,q}(\tm) =
	\frac{1}{q_i^{-2}-1} \sum_{ r: \; (i,r) \in \widetilde{\sDynkinf}_0 }
	\left[ \ts_{i,r},\, \tm \right]
\end{align*}
for an $\calX_q$-monomial $\tm$, where $\lcalX_{i,q}$ is regarded as the $\calX_q$-bimodule,
\smallskip
\item the map from $\lcalX_{i,q}$ to $\calX_{i,q}$, denoted by an double-headed arrow, is the surjective map sending an element of $\lcalX_{i,q}$ to its image in $\calX_{i,q}$ (recall \eqref{eq:generators of submodules}).
\ee

\noindent
By direct computation, we have the following:
 
\begin{proposition} \label{prop: Siq is derivation}
The map $S_{i,q}$ is a $\Z[q^{\pm \frac{1}{2}}]$-linear map and derivation with respect to $*$, that is,
\begin{align} \label{eq: derivation}
	S_{i,q}(\tm_1 * \tm_2) =
	\tm_1 \bcdot S_{i,q}(\tm_2) + S_{i,q}(\tm_1) \bcdot \tm_2,
\end{align}
where the $\bcdot$ indicates the $\calX_q$-bimodule actions of $\calX_{i,q}$ induced from $\lcalX_{i,q}$.
\end{proposition}

\begin{definition} \label{def: virtual quantum G ring}
For $i\in\Dynkinf_0$, we denote by $\frakK_{i,q}(\g)$
the $\Z[q^{\pm \frac{1}{2}}]$-subalgebra of $\calX_q$ generated by
\begin{equation*}
\tX_{i,l} * (1+q_i^{-1}\tB_{i,l+1}^{-1})       \ \ \text{ and } \ \ \tX_{j,s}^{\pm1}  \quad  \text{ for } j\in\Dynkinf_0\setminus\{ i \} \ \ \text{ and } \ \  (i,l),(j,s) \in \tDynkinf_0.
\end{equation*}
\end{definition}

\noindent
By using the same arguments as in~\cite{FM01,H03,H04}, we have
\begin{equation} \label{eq: characterization of Kiq}
	\frakK_{i,q}(\g) = {\rm Ker}(S_{i,q}).
\end{equation}
Therefore, we call $S_{i,q}$ \emph{the $i$-th $q$-screening operator} with respect to $\frakK_{i,q}(\g)$.

\begin{definition}   \cite{KO23}  \label{def: QGR}
We set
\begin{align*}
\frakK_q(\g) \seteq  \bigcap_{i \in I} \frakK_{i,q}(\g)
\end{align*}
and call it the \emph{quantum virtual Grothendieck ring associated to $\usfC(t)$}.
\end{definition}
\noindent

\begin{remark} \label{rem: some remark for screenings} 
%
Using the fact that $S_{i,q}$ is a $\Z[q^{\pm \frac{1}{2}}]$-linear derivation (or by its definition with \eqref{eq: right action on LXiq}), one can check that
	$S_{i,q}(\tX_{j,l}^{-1}) = -\delta_{i,j} \tX_{i,l}^{-1} \cdot  \ts_{i,l}$.
Then it follows from the definition of $S_{i,q}$, \eqref{eq: right action on LXiq} and \eqref{eq:generators of submodules} that
\begin{align*}
\begin{split}
	S_{i,q}(\tX_{i,l}^{-1} + q_i^{-1} \tX_{i,l}^{-1} *\tB_{i,l-1}) &=
	(-\tX_{i,l}^{-1}) \ts_{i,l} + (q_i^{-1} \tX_{i,l}^{-1} * \tB_{i,l-1} ) \ts_{i,l-2}
	= 0.
\end{split}
\end{align*}
In fact, $\frakK_{i,q}(\g)$ is realized as the $\Z[q^{\pm \frac{1}{2}}]$-subalgebra of $\calX_q$ generated by $\tX_{i,l}^{-1} +q_i^{-1}\tX_{i,l}^{-1} * \tB_{i,l-1}$  
 and $\tX_{j,s}^{\pm1}$ for $j\in\Dynkinf_0\setminus\{ i \}$ and $(i,l),(j,s) \in \tDynkinf_0$ (cf.~\eqref{eq: characterization of Kiq}).
\end{remark}

\begin{remark} \label{rem: evaluation of Kqg}
Since the following diagram commutes (cf.~\cite{H04})
\begin{equation} \label{eq: commutative diagram for finite types}
\begin{split}
	\xymatrixcolsep{4pc}\xymatrixrowsep{1pc}
	\xymatrix{
		\calX_q \ar@{->}[r]^{S_{i,q}} \ar@{->}[d]_{\evq} & \calX_{i,q} \ar@{->}[d]^{\evq} \\
		\calX \ar@{->}[r]_{S_i} & \calX_i
		}
\end{split}
\end{equation}
where $S_i$ is the $i$-th screening operator with respect to $\usfC(t)$, $\evq \bl   \frakK_q(\g) \br \subset \frakK(\g)$.
However, the opposite inclusion is not trivial (for non-simply-laced types).
We resolve this issue in the next section.
\end{remark}

\section{Bases of \texorpdfstring{$\frakK_q(\g)$ } aand Kazhdan--Lusztig analogues} \label{sec:bases and KL analogues}
\noindent
Let $(\bfg, \sfg)$ be a pair in \eqref{eq: pair of gg}.
It is known in \cite{Nak00, Nak01} (see also \cite{H04}) that the basis $\bfF_q$ of $\frakK_q(\bfg)$
with properties \eqref{eq:properties of Ft} below can be constructed algorithmically by using a deformed Frenkel--Mukhin (FM for short) algorithm (cf.~\cite{FM01}) with respect to $\sfC(q)$ (so-called $t$-algorithm \cite{H04}).
This basis enables us to construct other important bases of $\frakK_q(\bfg)$ (see  \eqref{eq:basis Et}, Theorem \ref{thm:basis Lt}). 
In the second part of this section, we will construct a basis $\sfF_q$ of $\frakK_q(\sfg)$ by a deformed FM-algorithm with respect to $\usfC(t)$, and verify that it has similar properties to \eqref{eq:properties of Ft} by following the framework in \cite{H04}.
Moreover, we also construct other bases $\sfE_q$ and $\sfL_q$ of $\frakK_q(\sfg)$ from the basis $\sfF_q$ in the spirit of \cite{Nak00, H04} in which they studied analogues of Kazhdan--Lusztig polynomials \cite{KL} (see Theorem \ref{thm: bar L_q(m)}, Remarks \ref{rem: KL theory} and  \ref{rem: new KL-type polynomials}).

\subsection{Bases of $\frakK_q(\bfg)$}

Note that $\sfC(q)$ coincides with $\usfC(t)$ for simply-laced finite types, when we replace $q$ with $t$. Thus,
\begin{eqnarray*} &&
\parbox{95ex}{
\emph{throughout this subsection, we switch the roles of $q$ and $t$. Also, we use $\At^-$ instead of $\Bq^-$.}
}
\end{eqnarray*}
This makes our notations more compatible with the literature where only simply-laced types are considered.

 In \cite{H04} (cf.~\cite{Nak00,Nak04}), 
the algorithm for constructing basis  $\bfF_t \seteq \{ F_t(\ubfm) \ |  \ \bfm \in \calM^\bfg_+\} $  was proposed, so called {\it $t$-algorithm}.
The  structure and  properties of the algorithm can be summarized as follows:
\eq &&  \label{eq:properties of Ft}
\parbox{90ex}{
\bna
\item For each dominant $\calY_t$-monomial $\tbfm$,
we construct an element $F_t(\tbfm)$ by adding monomials  $\tbfm' \in  \tbfm \At^{-k}$ in an inductive way as $k$ increases from $0$. In the process, the coefficient for each monomial is also determined in an inductive way.
\item \label{it: ter} If there appears a unique $\tbfm'$  with the smallest $k \in \Z_{\ge 1}$
satisfying
\bnum
\item $\tbfm'$ is anti-dominant and $\tbfm' \in \tbfm \At^{-k}$ is generated in the performing step,
\item any monomial generated in the previous step is contained in $\tbfm\At^{-s}$  $(0 \le s <k)$, not anti-dominant,
and  strictly larger than $\tbfm'$ with respect to $\lN$,
\ee
then, the coefficient of $\tbfm'$ is contained in $t^{\frac{1}{2}\Z}$.
Furthermore, the sum of all monomials with coefficients obtained from the steps so far, denoted by $F_t(\tbfm)$, is contained in the kernel of $S_{\im,t}$ for all $\im$. Hence
$F_t(\tbfm)$ is an element of $\frakK_t(\bfg)$ and the $t$-algorithm terminates. %
\ee

Furthermore, each $F_t(\tbfm)$ satisfies the following properties:
\ben
\item \label{it: finite step} $F_t(\tbfm) \in \frakK_t(\bfg) \cap \tbfm \At^{-}$.
\item $F_t(\tbfm) $ is bar-invariant if $\tbfm$ is bar-invariant.
\item \label{it: less than}  
Every monomial of $F_t(\tbfm) - \tbfm$ is   strictly less than $\tbfm$ with respect to $\lN$. 
\ee
}\eneq
%
The $t$-algorithm might progress infinitely many times.
In fact, $F_t(\tbfm)$ was constructed in a completion of $\frakK_t(\bfg)$ at first.
Interestingly, the property \ref{it: finite step} in \eqref{eq:properties of Ft} is guaranteed once we prove 
\begin{equation} \label{eq:finiteness of FtYip}
F_t( \tY_{\im,p}) \in \frakK_t(\bfg).
\end{equation}
More precisely, \eqref{eq:finiteness of FtYip} implies $\bfE_t := \left\{ E_t(\ubfm) \, | \, \ubfm \in \calM^\bfg_+ \right\} \subset \frakK_t(\bfg)$, where $E_t(\ubfm)$ is given in  \eqref{eq:basis Et}. 
Then it is known (e.g.~see the proof of \cite[Proposition 6.3]{H04} for more detail) that $\bfE_t$ has the unit-triangular property with $\bfF_t$, that is, $F_t(\ubfm)$ can be written as a linear combination of elements in $\bfE_t \subset \frakK_t(\bfg)$,
so the proof for \ref{it: finite step} in \eqref{eq:properties of Ft} is reduced to prove \eqref{eq:finiteness of FtYip}. 
Then \eqref{eq:finiteness of FtYip} is deduced from \cite{Nak00,Nak04}.
%

\begin{remark} \label{rem: reversed q-algorithm}
  Another characterization of $\frakK_{\im,t}(\bfg)$ in Remark \ref{rem: some remark for screenings} 
allows us to consider the lowest $\ell$-weight version of the $t$-algorithm,
that is, a $t$-deformation of {\it reversed} Frenkel--Mukhin algorithm which is an algorithm starting from the lowest $\ell$-weight monomial.
For instance, the formulas in \cite[Lemma 4.13]{H04} can be re-formulated in terms of anti-dominant monomial with $\tA_{\im,k}$'s.
The reversed algorithm seems to be already known to experts in the theory of $q$-characters (e.g. see \cite{FM01}, \cite{Mou}).

Let $\tbfm_{-}$ be an anti-dominant (bar-invariant) $\calY_t$-monomial. We denote by $F_t(\tbfm_{-})$ the unique element of $\frakK_q(\bfg)$ generated by the {\it reversed} $t$-algorithm (referred above) with respect to $\tbfm_{-}$. Then one can verify that $F_t(\tbfm_{-})$ satisfies similar properties to \eqref{eq:properties of Ft} after modifying notations and terminologies associated with $\tbfm_{-}$.
For example, the property~\ref{it: less than} in \eqref{eq:properties of Ft} associated with $\tbfm_{-}$ is restated as every monomial appearing in $F_t(\tbfm_-) - \tbfm_-$ is \ strictly {\it greater} than $\tbfm_-$ with respect to $\lN$. 
Throughout this section, we often refer to these properties.
\end{remark}

\begin{theorem} \cite[Theorem 3.1]{Nak03} \cite[Theorem 4.1, Lemma 4.4]{H06} \label{thm:monomials for KR modules in types ADE}
For $(\im,p),(\im,s) \in \tDynkinf_0$ with $p<s$, the element $F_t( \ubfm^{(\im)}[p,s] ) \in \frakK_t(\bfg) $ is of the form
\begin{align*}
F_t( \ubfm^{(\im)}[p,s] )   =   \ubfm^{(\im)}[p,s] *  ( 1+ \tA^{-1}_{\im,s+1}  * \chi),
\end{align*}
where $\ubfm^{(\im)}[p,s] \seteq \ul{\bfm^{(\im)}[p,s]}$ and
$\chi$ is a $($non-commutative$)$ $\Z[t^{\pm \frac{1}{2}}]$-polynomial in $\tA^{-1}_{\jm,k+1}$ $(\jm,k) \in \tDynkinf_0$. In particular, we have
\begin{align*}
F_t( \ubfm^{(\im)}[p,s] ) = F_t( \ubfm^{(\im^*)}_-[p+\sfh,s+\sfh] ),
\end{align*}
where $\ubfm^{(\im^*)}_-[p+\sfh,s+\sfh] \seteq (\ubfm^{(\im^*)}[p+\sfh,s+\sfh])_-$.
Furthermore, $F_t( \ubfm^{(\im)}[p,s] )$ satisfies the following properties:
\ben
\item $ F_t( \ubfm^{(\im)}[p,s] )$ has the unique $($anti-$)$dominant monomial $ \ubfm^{(\im)}[p,s] $ $($resp.  $\ubfm^{(\im^*)}_-[p+\sfh,s+\sfh])$.
\item Each $\calY_t$-monomial of $F_t(\ubfm^{(\im)}[p,s]) - \ubfm^{(\im)}[p,s]-\ubfm^{(\im^*)}_-[p+\sfh,s+\sfh]$ is a product of $\tY_{\jm,u}^{\pm1}$ with 
$p \le u \le s+\sfh$, 
having at least one of it factors from $p < u < s+\sfh$, 
and is right-negative.  
In particular, each $\calY_t$-monomial of $F_t(\underline{\tY_{\im,p}}) - \underline{\tY_{\im,p}}-\underline{\tY_{\im^*,p+\sfh}}$ is a product of $\tY_{\jm,u}^{\pm1}$ with $p < u < p+\sfh$. 
\item \label{it: commute} For $((\im,p),(\jm,p) \in \tDynkinf_0, \; \jm \ne \im)$,  $ F_t( \ubfm^{(\im)}[p,s] )$ and $F_t( \ubfm^{(\jm)}[p,s] )$   commute; i.e.,   $ F_t( \ubfm^{(\im)}[p,s] ) *  F_t( \ubfm^{(\jm)}[p,s] ) =   F_t( \ubfm^{(\jm)}[p,s] )*  F_t( \ubfm^{(\im)}[p,s] )$.
\ee
\end{theorem}

It is well known that, for $r \in 2\Z$ and $\im \in \Dynkinf_0$,
\begin{align} \label{eq: translation by r}
\sfT_r( F_t( \ubfm^{(\im)}[p,s] ) )  = F_t( \ubfm^{(\im)}[p+r,s+r] ),
\end{align}
where $\sfT_{r}$ is the $\Z[t^{\pm \frac{1}{2}}]$-algebra automorphism of $\calY_t$ sending $\tY_{i,p}$ to $\tY_{i,p+r}$.

\begin{theorem} \cite[Theorem 5.11]{H04}  \label{thm: F_t}
\bna
	\item For every dominant $($resp.~anti-dominant$)$ monomials $\tbfm \in \calY_t$, $F_t(\tbfm)$ is the
unique element in $\frakK_t(\bfg)$ such that $\tbfm$ is the unique
dominant $($resp.~anti-dominant$)$ monomial of $F_t(\tbfm)$. 
	
	\item  Every monomial appearing in $F_t(\tbfm) - \tbfm$ is  strictly less $($resp.~strictly greater$)$ than $\tbfm$ with respect to $\lN$.
	
	\item The set $\bfF_t \seteq \{   F_t(\ubfm) \ | \ \bfm \in \calM^\bfg_{+} \}$   forms a
bar-invariant $\Z[t^{\pm \frac{1}{2}}]$-basis of $\frakK_t(\bfg)$.
\ee
\end{theorem}


\begin{remark} \label{rem: F vs Ft in ADE}
We remark that an element in $\frakK_t(\bfg)$ is characterized by the multiplicities of its dominant monomials by Theorem~\ref{thm: F_t}.
Then it yields that ${\rm ev}_{t=1}(F_t(\tbfm)) = F({\rm ev}_{t=1}(\tbfm))$.
\end{remark}

\begin{example} \label{ex: F20 in D4}
We present $F_t(\tY_{2,0})$ of type $D_4$ (cf. \cite[Example 5.3.2]{Nak00}) by organizing the monomials appearing in $F_t(\tY_{2,0})$ as a directed graph $\Gamma(\tY_{2,0})$ such that $F_t(\tY_{2,0})$ is the sum of the monomials on the vertices of the directed graph, see \eqref{eq: D4 graph}.
Note that in this example, we write the $\calY_t$-monomials according to the order given by
  \begin{equation} \label{eq: order for monomials}
    (\im, p) < (\jm, s) \quad \Longleftrightarrow \quad (\,p < s\,) \,\, \text{or} \,\, (\,p = s \,\,\text{and}\,\, \im < \jm\,).
  \end{equation}

We use the convention of \cite{FR99,Nak00} for the directed oriented graph $\Gamma(\tY_{2,0})$: For monomials $\tbfm_1$ and $\tbfm_2$, we use an colored directed edge $f(t)\;\tbfm_1\overset{\im,k}{\longrightarrow} g(t)\;\tbfm_2$ if $
\ev_{t=1}(\tbfm_2) = \ev_{t=1}(\tbfm_1 \tA_{\im,k}^{-1})$, where $f(t), g(t) \in \Z[t^{\pm \frac{1}{2}}]$.
Then the directed colored graphs $\Gamma(\tY_{2,0})$ of $F_t(\tY_{2,0})$ is given as below:
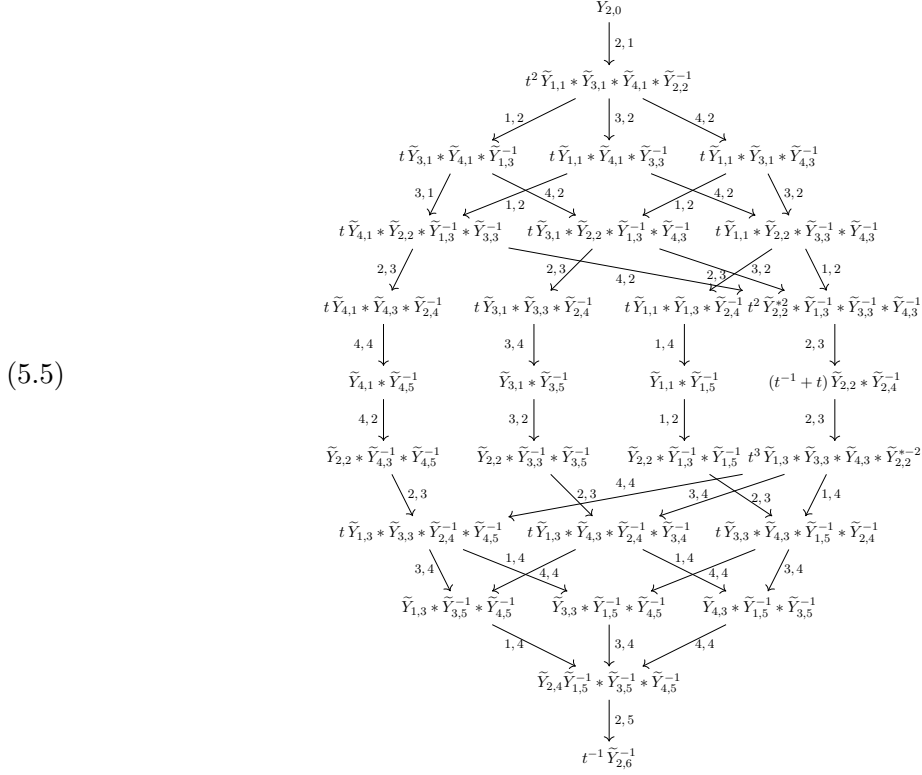
\begin{equation} \label{eq: D4 graph}
\scalebox{0.5}{$
    \begin{tikzpicture}[baseline=(current  bounding  box.center), every node/.style={scale=1.1}]
            %
            \node (0) at (0, 11) {$\tY_{2,0}$};
            \node (1) at (0, 9) {$t^2\, \tY_{1,1} * \tY_{3,1} *\tY_{4,1} *\tY_{2,2}^{-1}$};
            \node (2-1) at (-4, 7) {$t\, \tY_{3,1} *\tY_{4,1} *\tY_{1,3}^{-1}$};
            \node (2-2) at (0, 7) {$t\, \tY_{1,1} *\tY_{4,1} *\tY_{3,3}^{-1}$};
            \node (2-3) at (4, 7) {$t\, \tY_{1,1} *\tY_{3,1} *\tY_{4,3}^{-1}$};
            \node (3-1) at (-5, 5) {$t\, \tY_{4,1} *\tY_{2,2} *\tY_{1,3}^{-1} *\tY_{3,3}^{-1}$};
            \node (3-2) at (0, 5) {$t\, \tY_{3,1} *\tY_{2,2} *\tY_{1,3}^{-1} *\tY_{4,3}^{-1}$};
            \node (3-3) at (5, 5) {$t\, \tY_{1,1} *\tY_{2,2} *\tY_{3,3}^{-1} *\tY_{4,3}^{-1}$};
            \node (4-1) at (-6, 3) {$t\, \tY_{4,1} *\tY_{4,3} *\tY_{2,4}^{-1}$};
            \node (4-2) at (-2, 3) {$t\, \tY_{3,1} *\tY_{3,3} *\tY_{2,4}^{-1}$};
            \node (4-3) at (2, 3) {$t\, \tY_{1,1} *\tY_{1,3} *\tY_{2,4}^{-1}$};
            \node (4-4) at (6, 3) {$t^2\, \tY_{2,2}^{*2} *\tY_{1,3}^{-1} *\tY_{3,3}^{-1}* \tY_{4,3}^{-1}$};
            \node (5-1) at (-6, 1) {$\tY_{4,1}* \tY_{4,5}^{-1}$};
            \node (5-2) at (-2, 1) {$\tY_{3,1}* \tY_{3,5}^{-1}$};
            \node (5-3) at (2, 1) {$\tY_{1,1} *\tY_{1,5}^{-1}$};
            \node (5-4) at (6, 1) {$(t^{-1}+t)\, \tY_{2,2} *\tY_{2,4}^{-1}$};
            \node (6-1) at (-6, -1) {$\tY_{2,2} *\tY_{4,3}^{-1} *\tY_{4,5}^{-1}$};
            \node (6-2) at (-2, -1) {$\tY_{2,2} *\tY_{3,3}^{-1} *\tY_{3,5}^{-1}$};
            \node (6-3) at (2, -1) {$\tY_{2,2} *\tY_{1,3}^{-1} *\tY_{1,5}^{-1}$};
            \node (6-4) at (6, -1) {$t^{3}\, \tY_{1,3} *\tY_{3,3}* \tY_{4,3} *\tY_{2,2}^{*-2}$};
            \node (7-1) at (-5, -3) {$t\, \tY_{1,3} *\tY_{3,3}* \tY_{2,4}^{-1} *\tY_{4,5}^{-1}$};
            \node (7-2) at (0, -3) {$t\, \tY_{1,3} *\tY_{4,3} *\tY_{2,4}^{-1} *\tY_{3,4}^{-1}$};
            \node (7-3) at (5, -3) {$t\, \tY_{3,3} *\tY_{4,3} *\tY_{1,5}^{-1} *\tY_{2,4}^{-1}$};
            \node (8-1) at (-4, -5) {$\tY_{1,3} *\tY_{3,5}^{-1} *\tY_{4,5}^{-1}$};
            \node (8-2) at (0, -5) {$\tY_{3,3} *\tY_{1,5}^{-1} *\tY_{4,5}^{-1}$};
            \node (8-3) at (4, -5) {$\tY_{4,3} *\tY_{1,5}^{-1} *\tY_{3,5}^{-1}$};
            \node (9) at (0, -7) {$\tY_{2,4} \tY_{1,5}^{-1} *\tY_{3,5}^{-1} *\tY_{4,5}^{-1}$};
            \node (10) at (0, -9) {$t^{-1}\, \tY_{2,6}^{-1}$};
            \draw[-{>[scale=1.5]}] (0) -- (1) node[midway, right] {\footnotesize $2,1$};
            \draw[-{>[scale=1.5]}] (1) -- (2-1) node[midway, left] {\footnotesize $1,2$\,\,};
            \draw[-{>[scale=1.5]}] (1) -- (2-2) node[midway, right] {\footnotesize $3,2$};
            \draw[-{>[scale=1.5]}] (1) -- (2-3) node[midway, right] {\,\,\footnotesize  $4,2$};
            \draw[-{>[scale=1.5]}] (2-1) -- (3-1) node[midway, left] {\footnotesize $3,1$};
            \draw[-{>[scale=1.5]}] (2-1) -- (3-2) node[midway, right] {\,\,\footnotesize $4,2$};
            \draw[-{>[scale=1.5]}, draw] (2-2) -- (3-1) node[midway, below] {{\footnotesize $1,2$}};
            \draw[-{>[scale=1.5]}, draw] (2-2) -- (3-3) node[midway, right] {\,\,{\footnotesize $4,2$}};
            \draw[-{>[scale=1.5]}] (2-3) -- (3-2) node[midway, below] {\footnotesize $1,2$};
            \draw[-{>[scale=1.5]}] (2-3) -- (3-3) node[midway, right] {\footnotesize $3,2$};
            \draw[-{>[scale=1.5]}] (3-1) -- (4-1) node[midway, left] {\footnotesize $2,3$};
            \draw[-{>[scale=1.5]}] (3-2) -- (4-2) node[midway, left] {\footnotesize $2,3$};
            \draw[-{>[scale=1.5]}] (3-3) -- (4-3) node[midway, left] {\raisebox{-1.2pc}{\footnotesize $2,3$\,\,\,\,}};
            \draw[-{>[scale=1.5]}, draw] (3-1) -- (4-4) node[midway, below] {{\footnotesize $4,2$}};
            \draw[-{>[scale=1.5]}, draw] (3-2) -- (4-4) node[midway, right] {\,\,\,\quad \footnotesize $3,2$};
            \draw[-{>[scale=1.5]}, draw] (3-3) -- (4-4) node[midway, right] {{\footnotesize $1,2$}};
            \draw[-{>[scale=1.5]}] (4-1) -- (5-1) node[midway, left] {\footnotesize $4,4$\,\,};
            \draw[-{>[scale=1.5]}] (4-2) -- (5-2) node[midway, left] {\footnotesize $3,4$\,\,};
            \draw[-{>[scale=1.5]}] (4-3) -- (5-3) node[midway, left] {\footnotesize $1,4$\,\,};
            \draw[-{>[scale=1.5]}] (4-4) -- (5-4) node[midway, left] {\footnotesize $2,3$\,\,};
            \draw[-{>[scale=1.5]}] (5-1) -- (6-1) node[midway, left] {\,\,\footnotesize $4,2$};
            \draw[-{>[scale=1.5]}] (5-2) -- (6-2) node[midway, left] {\,\,\footnotesize $3,2$};
            \draw[-{>[scale=1.5]}] (5-3) -- (6-3) node[midway, left] {\,\,\footnotesize $1,2$};
            \draw[-{>[scale=1.5]}] (5-4) -- (6-4) node[midway, left] {\footnotesize $2,3$\,\,};
            \draw[-{>[scale=1.5]}] (6-1) -- (7-1) node[midway, right] {\footnotesize $2,3$};
            \draw[-{>[scale=1.5]}] (6-2) -- (7-2) node[midway, right] {\footnotesize $2,3$};
            \draw[-{>[scale=1.5]}] (6-3) -- (7-3) node[midway, right] {\raisebox{-0.9pc}{\,\,\footnotesize $2,3$}};
            \draw[-{>[scale=1.5]}, draw] (6-4) -- (7-1) node[midway, above] {{\footnotesize $4,4$}};
            \draw[-{>[scale=1.5]}, draw] (6-4) -- (7-2) node[midway, left] {{\footnotesize $3,4$}\,\,\,};
            \draw[-{>[scale=1.5]}, draw] (6-4) -- (7-3) node[midway, right] {{\footnotesize $1,4$}};
            \draw[-{>[scale=1.5]}] (7-1) -- (8-1) node[midway, left] {\footnotesize $3,4$};
            \draw[-{>[scale=1.5]}] (7-1) -- (8-2) node[midway, above] {\footnotesize $1,4$};
            \draw[-{>[scale=1.5]}, draw] (7-2) -- (8-1) node[midway, right] {\raisebox{-0.9pc}{\footnotesize $4,4$}};
            \draw[-{>[scale=1.5]}, draw] (7-2) -- (8-3) node[midway, above] {{\footnotesize $1,4$}};
            \draw[-{>[scale=1.5]}] (7-3) -- (8-2) node[midway, right] {\raisebox{-0.9pc}{\footnotesize $4,4$}};
            \draw[-{>[scale=1.5]}] (7-3) -- (8-3) node[midway, right] {\footnotesize $3,4$};
            \draw[-{>[scale=1.5]}] (8-1) -- (9) node[midway, left] {\footnotesize $1,4$\,\,};
            \draw[-{>[scale=1.5]}] (8-2) -- (9) node[midway, right] {\footnotesize $3,4$};
            \draw[-{>[scale=1.5]}] (8-3) -- (9) node[midway, right] {\,\,\footnotesize $4,4$};
            \draw[-{>[scale=1.5]}] (9) -- (10) node[midway, right] {\footnotesize $2, 5$};
    \end{tikzpicture}
    $}
\end{equation}
\end{example}

For a dominant monomial $\bfm \in \calM_+^\bfg$, we set
\begin{equation} \label{eq:basis Et}
E_t(\ubfm) \seteq  t^a
\left(  \st_{p\in\Z}^{\to} \left(  \st_{\im\in I^\bfg; (\im,p)\in\widetilde{\sDynkinf}_0} F_t(\tY_{\im,p})^{u_{\im,p}(\bfm)}  \right) \right),
\end{equation}
where $a$ is an element in $\frac{1}{2}\Z$ such that $\ubfm$ appears in $E_t(\ubfm)$ with the coefficient $1$.
Here
$ \st_{\im} F_t(\tY_{\im,p})^{u_{\im,p}(\bfm)} $ is well-defined by Theorem~\ref{thm:monomials for KR modules in types ADE}~\ref{it: commute}.
Note that $E_t(\ubfm)$ contains $\ubfm$  as its maximal monomial with respect to $\lN$. In particular, by Theorem~\ref{thm: F_t}, we have
\begin{align}\label{eq: uni t1}
E_t(\ubfm) = F_t(\ubfm) + \sum_{\bfm' \lN \bfm} C_{\bfm,\bfm'} F_t(\ubfm')
\end{align}
with $C_{\bfm,\bfm'} \in \Z[t^{\pm \frac{1}{2}}]$.
Note that the set $\bfE_t\seteq \{ E_t(\ubfm) \ | \ \bfm \in \calM_+ \}$ also forms a $\Z[t^{ \pm \frac{1}{2}}]$-basis since
\begin{align} \label{eq: finiteness}
\sharp \{  \bfm' \in \calM_+ \ | \  \bfm' \lN \bfm \} < \infty \quad \text{ for each $\bfm \in \calM_+$}.
\end{align}
We call $\bfE_t$ the \emph{standard basis} of $\frakK_t(\bfg)$.

Note that $Y_{\im,p}$ is a minimal element in $\calM_+$ with respect to the partial order $\leN$. Thus ~\eqref{eq: uni t1} tells that
$$  E_t(\ul{Y_{\im,p}}) = F_t(\ul{Y_{\im,p}}).$$

Using the bases $\bfF_t$ and $\bfE_t$, the third basis $\bfL_\ttt \seteq \{ L_t(\ubfm) \}$ of $\frakK_t(\bfg)$ have been constructed in an inductive way using $\leN$ such that
\begin{align} \label{eq: Et = Ft = Lt when level 1}
E_t(\ul{Y_{\im,p}}) = F_t(\ul{Y_{\im,p}})  = L_t(\ul{Y_{\im,p}})
\end{align}
and $L_t(\ubfm)$ for general $\bfm \in \calM_+$ is characterized as
in the following theorem.

\begin{theorem}\cite{Nak04} $($see also \cite{H04}$)$ \label{thm:basis Lt}
For a dominant monomial $\bfm \in \calM^\bfg_+$, there exists a unique element $L_t(\ubfm)$ in $\frakK_t(\bfg)$ such that
 $\overline{L_t(\ubfm)}=L_t(\ubfm)$ and
\begin{equation} \label{eq: uni t2}
\begin{split}
E_t(\ubfm) = L_t(\ubfm) + \displaystyle\sum_{\bfm' \lN \bfm }  P_{\bfm,\bfm'}(t) L_t(\ubfm') \quad \text{ with $P_{\bfm,\bfm'}(t) \in t\Z[t]$. }
\end{split}
\end{equation}
\end{theorem}
We call $\bfL_t$ the \emph{canonical basis} of $\frakK_t(\bfg)$.

\begin{remark} \label{rem: KL theory} 
In a highly influential paper \cite{KL}, Kazhdan and Lusztig conjectured a realization of the composition multiplicities of Verma modules for $\g$ in terms of a certain class of polynomials defined by Iwahori--Hecke algebras, so-called {\em Kazhdan--Lusztig polynomials} (KL polynomials, for short).
The \emph{Kazhdan--Lusztig conjecture} is that the specialization of the KL polynomials at $1$ coincides with the composition multiplicities of Verma modules. 
This is proved independently by Beilinson--Bernstein \cite{BB81, BB93} and Brylinski--Kashiwara \cite{BK81}.
 Moreover, it is shown in \cite{KL80} that the KL polynomials can be interpreted as the Poincar\'e polynomials for local intersection cohomology of Schubert varieties. This geometric interpretation gives the positivity of the KL polynomials.

 A similar story has been developed in the representation theory of quantum loop algebras.
In \cite{Nak00, Nak01, Nak04}, 
it is proved by Nakajima that the specialization of $P_{\bfm,\bfm'}(t)$ at $t=1$ gives the composition multiplicity of $L(\ubfm')$ in the standard module $E(\ubfm)$. Furthermore, $P_{\bfm,\bfm'}(t)$ coincides with the Poincar\'e polynomial of intersection cohomology of graded quiver varieties, which implies the positivity of $P_{\bfm,\bfm'}(t)$.
Consequently, the polynomials $P_{\bfm,\bfm'}(t)$'s may be viewed as analogs of KL polynomials.
It is worthwhile to remark that the recent development associated with $P_{\bfm,\bfm'}(t)$ in the direction of quantum loop algebras beyond ADE-types, see \cite{FHOO, FHOO2}.
\end{remark}

\begin{theorem} \cite{Nak04} \label{thm: Nak positivity} \hfill
\bna
\item \label{it:positivity of Ltm} For a dominant monomial $\bfm \in \calM^\bfg_+$, every monomial in $L_t(\ubfm)$ has a quantum positive coefficient;
that means, each coefficient of a monomial in $L_t( \ubfm)$ contained in $\Z_{\ge0}[t^{\pm \frac{1}{2}}]$.
In particular, we have $\ev_{t=1}(L_t(\ubfm))  = \chi_q(L(\ubfm)).$
\item For each monomial $\bfm^{(\im)}[p,s]$, we have $F_t(\ubfm^{(\im)}[p,s])=L_t(\ubfm^{(\im)}[p,s])$.
\item The coefficient $P_{\bfm,\bfm'}(t)$ in~\eqref{eq: uni t2} is actually contained in $t \Z_{\ge 0}[t]$.
\ee
\end{theorem}

\begin{remark} \label{rem: framework}
Let recapitulate the main points in this subsection.
From the $t$-algorithm, we obtain a basis  $\{ F_t(\ubfm) \ |  \ \bfm \in \calM^\bfg_+\} $ of $\frakK_t(\bfg)$. One crucial step is to prove that $F_{t}(\tY_{\im,p})$ is contained in $\frakK_t(\bfg)$.
Then it is proved in \cite{Nak04,H04} that there are frameworks for
constructing other two bases $\{ E_t(\ubfm) \ |  \ \bfm \in \calM^\bfg_+\} $ and $\{ L_t(\ubfm) \ |  \ \bfm \in \calM^\bfg_+\} $
of $\frakK_t(\bfg)$.
In particular,
the basis  $\{ L_t(\ubfm) \ |  \ \bfm \in \calM^\bfg_+\} $ is constructed using the other two bases through
the induction on $\calM_+$ via $\leN$, and there are uni-triangular transition maps  ~\eqref{eq: uni t1} and ~\eqref{eq: uni t2}
between the three bases.
\end{remark}

As $L_t(\ubfm)$ can be understood as a $t$-quantization of
$L(\ubfm)$ by Theorem~\ref{thm: Nak positivity}\,\ref{it:positivity of Ltm}, the $T$-system among KR modules is also $t$-quantized as follows:

\begin{theorem}\label{thm: Quantum T-system} \cite[Proposition 5.6]{HL15} $($see also \cite[Section 4]{Nak03}$)$
For $(\im,p),(\im,s) \in \tDynkinf_0$ with $p<s$, there exists an equation in $\frakK_t(\bfg)$:
\begin{align}\label{eq: quantum T-sys}
L_t \bl \ubfm^{(\im)}[p,s)  \br *  L_t \bl\ubfm^{(\im)}(p,s] \br = t^x L_t\bl \ubfm^{(\im)}[p,s] \br *  L_t\bl \ubfm^{(\im)}(p,s)\br + t^y \hspace{-1ex} \prod_{\jm; \; d(\im,\jm)=1}  \hspace{-1ex}  L_t\bl \ubfm^{(\jm)}(p,s) \br,
\end{align}
where $L_t\bl \ubfm^{(\jm)}(p,s))$ and $L_t\bl \ubfm^{(\jm')}(p,s))$ $(\jm,\jm' \in I)$ are pairwise commute and
$$ y = \dfrac{\tfb_{\im,\im}(2(s-p)+1)+\tfb_{\im,\im}(2(s-p)-1)}{2} \qquad \text{ and }  \qquad x = y- 1.$$
\end{theorem}

\subsection{Bases of $\frakK_q(\sfg)$} \label{subsec:bases of Kqg}
Assume that $\sfg$ is of {\it non-simply-laced finite type}.
Since $\sfC(q)$ can not be identified with $\usfC(t)$ anymore,
$$
\emph{we come back to the convention of the previous sections $($not the previous subsection$)$.}
$$

\noindent
Let $\frakK_q^\infty(\sfg)$ be a completion of $\frakK_q(\sfg)$ defined by the method in \cite[Section 5.2]{H04}.
By following the construction of $\{ F_t(\ubfm) \ |  \ \bfm \in \calM^\bfg_+\}$ in \cite{H03,H04}, we can establish an analog of the $t$-algorithm in \cite[Definition 5.19]{H04} on $\frakK_q^\infty(\sfg)$, called {\it $q$-algorithm} under the setting of Section \ref{subsec:quantization}.
Roughly speaking, the algorithm is given inductively by computing all possible quantized $i$-expansions (cf. \cite[Lemma 4.13]{H04}) while determining ``correct'' coefficients in $\Z[q^{\pm \frac{1}{2}}]$ of resulting monomials \cite[Definition 5.19]{H04} (cf. \cite[Section 5.5]{FM01}).

The $q$-algorithm employs the $t$-algorithm by considering the another quantization $\usfC(t)$ of the Cartan matrix $\sfC$.   Thus, in the computational view point, $\tX_{i, p}$ and  $q_i\tB_{i, p}$ in the $q$-algorithm play the roles of
$\tY_{i, p}$ and  $t\tA_{i, p}$ in the $t$-algorithm:
\begin{equation} \label{eq:replacing}
  \tY_{i,p} \,\,\longmapsto\,\, \tX_{i, p}, \qquad
  t\tA_{i,p} \,\,\longmapsto\,\, q_i\tB_{i, p}.
\end{equation}
We say that
\begin{enumerate}[(a)]
	\item the $q$-algorithm is {\it well-defined for step $s$} if the $k$-th coefficients defined as in \cite[Definition 5.19]{H04} with \eqref{eq:replacing} are verified for all $k \le s$,
	\item the $q$-algorithm {\it never fails} if it is well-defined for all steps.
\end{enumerate}

When the $q$-algorithm never fails, it yields, for each dominant monomial $\tm$ in $\calX_q$,
\begin{align} \label{eq: fqm in infinity}
F_q(\tm) \in \frakK_q^\infty(\sfg) = \bigcap_{i \in I} \frakK^\infty_{i,q}(\sfg),
\end{align}
containing $\tm$ as a unique dominant monomial in $F_q(\tm)$.  Here,  $\frakK_{i,q}^\infty(\sfg) $ is the completion of $\frakK_{i,q}(\sfg)$ satisfying
	$\frakK_{i,q}^\infty(\sfg) \cap \calX_q  = \frakK_{i,q}(\sfg) = \Ker (S_{i,q})$
(see Lemma~\ref{lem: characterization of KJq} below for more detail).
It is straightforward to check that the $q$-algorithm is well-defined and never fails
by following the framework of \cite[Section 5.3]{H04}.
Since the proof is quite parallel to \cite{H04},
the details are left to the reader.
As a result, we have the following.
\begin{proposition} \label{prop: basis Fq in Kqf_infty}
Let $\frakK_q^{\infty,\, f}(\sfg)$ be the $\Z[q^{\pm \frac{1}{2}}]$-submodule of $\frakK_q^\infty(\sfg)$ generated by elements in $\frakK_q^\infty(\sfg)$ with finitely many dominant monomials.
Then the set
\begin{equation*}
	\left\{\, F_q(\tm) \, | \, \text{$\tm$ is a dominant monomial in $\calX_q$} \,\right\}
\end{equation*}
is a $\Z[q^{\pm \frac{1}{2}}]$-basis of $\frakK_q^{\infty,\, f}(\sfg)$.
Indeed, $\frakK_q^{\infty,\, f}(\sfg)$ is a $\Z[q^{\pm \frac{1}{2}}]$-subalgebra of $\frakK_q^\infty(\sfg)$
\end{proposition}
\noindent
Instead of the rigorous proofs of the above results,
let us illustrate the $q$-algorithm in Example \ref{ex: F25 in G2} below by using results and notations in the parts after the example, frequently.

\begin{example} \label{ex: F25 in G2} 
We consider $F_{q}(\tX_{2,5})$ for type $G_2$.
Let us briefly summarize the procedure of $q$-algorithm. 
For $\calX$-monomials $m_1$ and $m_2$, we use an colored directed edge $m_1 \overset{i,\,p}{\longrightarrow} m_2$ if $m_2 = m_1 B_{i,p}^{-1}$.
For $\calX$-monomials $m$ and $m'$, we say that $m'$ is generated from $m$ if there exists a finite sequence $\left\{ (i_1, p_1), (i_2, p_2), \cdots, (i_\ell, p_\ell) \right\} \subset I \times \Z$ such that
\begin{equation*}
	\xymatrixcolsep{3.5pc}\xymatrixrowsep{0pc}
	\xymatrix{
	m = m_0\, \ar@{->}[r]^{\,\,\,i_1,\, p_1} & \,m_1\, \ar@{->}[r]^{i_2,\, p_2} & \cdots \ar@{->}[r]^{i_{\ell-1},\, p_{\ell-1}\quad\,\,\,} & \,m_\ell = m',
	}
\end{equation*}
where $m_k$ is a $\calX$-monomial of $E_J(m_{k-1})$ defined in \eqref{eq: EJm} for $1 \le k \le \ell$ for some $J \subset I$.
We collect all possible $\calX$-monomials generated from $X_{2,5}$, and then enumerate them by
\begin{equation} \label{eq: generated monomials from X25 in G2}
	\dots < m_v < \dots < m_0 = X_{2,5},
\end{equation}
where $<$ is a total order compatible with $\lN$ at $q = 1$.
Let $\tm_v$ be a $\calX_q$-monomial determined inductively from $F_{J,q}(\tm_u) \in \bigcap_{i \in J} \frakK_{i,q}(\sfg)$ for some $u < v$ and $J \subset I$, where $F_{J,q}(\tm_u)$ contains $\tm_u$ as a unique $J$-dominant monomial.
Note that $\tm_v$ is uniquely determined up to a coefficient in $q^{\frac{1}{2}\Z}$.
For this reason, we fix an order \eqref{eq: order for monomials} on spectral parameters to write them uniquely.

Put $J \subsetneq I$.
We denote by
	$\left( s(m_v)(q) \right)_{v \in \Z_{\ge 0}}$ and
	$\left( s_J(m_v)(q) \right)_{v \in \Z_{\ge 0}}$
sequences in $\Z[q^{\pm \frac{1}{2}}]^{\Z_{\ge 0}}$ defined inductively as follows:
\begin{equation} \label{eq: coloring and coefficient in q-algorithm}
\begin{split}
	s_J(m_v)(q) &= \sum_{u < v} \left( s(m_u)(q) - s_J(m_u) \right) c_{J}(q)(m_v), \\
	s(m_v)(q) &=
	\begin{cases}
		s_J(m_v)(q) & \text{if $m_v$ is not $J$-dominant,} \\
		0 & \text{if $m_v$ is dominant,}
	\end{cases}
\end{split}
\end{equation}
where
$s(m_0)(q) = 1$, $s_J(m_0)(q) = 0$ and
$c_J(q)(m_v)$ is a $\Z[q^{\pm \frac{1}{2}}]$-coefficient of $\tm_v$ in $F_{J,q}(\tm_u)$.
Here we assume that $F_{J,q}(\tm_u) = 0$ if $m_u$ is not $J$-dominant, so $c_J(q)(m_u) = 0$ in this case.
Note that the sequences $\left( s(m_v)(q) \right)_{v \in \Z_{\ge 0}}$ and $\left( s_J(m_v)(q) \right)_{v \in \Z_{\ge 0}}$ are well-defined, and $s(m_v)(q)$ does not depend on the choice of $J \subsetneq I$ (e.g.~see Remark \ref{rem: Ft vs Fq and independence of q-algorithm to J}).
%
Finally, we verify that the sum of all $s(\tm_v)(q)\,\tm_v$ for $v \ge 0$
is in $\ker (S_{i,q})$ for all $i \in I$ (e.g.~see \eqref{eq: Fq25 = sum of F1q's}). In other words,
	$F_q(\tX_{2,5}) = \sum_{v \ge 0} s(\tm_v)(q) \,\tm_v \in \frakK_q(\bfg)$.

We emphasize that the monomials in \eqref{eq: generated monomials from X25 in G2} might be infinitely countable, but non-zero $\Z[q^{\pm \frac{1}{2}}]$-coefficients of them should be finite (see Proposition \ref{prop:folding lemma}). Thus, the formula on the right-hand side makes sense, and it is actually a finite sum.
For $n\in \Z \setminus \{ 0 \}$, we use $\tX_{i,p}^n$ to denote $\tX_{i,p}^{*n}$ for simplicity.
In fact, one may observe that $F_q(\tX_{2,5}) \in {\rm ker}(S_{i,q})$ for $i = 1, 2$   (see the complete computation in \cite{JLO1}).  
For example, $F_q(\tX_{2,5})$ is written as
\begin{equation} \label{eq: Fq25 = sum of F1q's}
\begin{split}
	F_q(\tX_{2,5}) &=
	F_{1,q}(\tX_{2,5}) +
	q^3 F_{1,q}(\tX_{1,6}^3 * \tX_{2,7}^{-1}) +
	(q^{-3} + q^{-1} + q) F_{1,q}(\tX_{1,6} * \tX_{1,8} * \tX_{2,9}^{-1}) \\
	& \quad + (q^{-3} + q^3) F_{1,q}(\tX_{2,7} * \tX_{2,9}^{-1}) +
	q^6 F_{1,q}(\tX_{1,8}^3 * \tX_{2,9}^{-2}) +
	q^{-3} F_{1,q}(\tX_{2,11}^{-1})
\end{split}
\end{equation}
which is clearly in ${\rm ker}(S_{1,q})$ (recall Proposition \ref{prop: Fiq}).
The case of $i = 2$ is similar.
In general, for a dominant $\calX_q$-monomial $\tm$,
the $q$-algorithm allows us to write $F_q(\tm)$ as a linear combination of $F_{i,q}(\,\cdot\,)$'s over $\Z[q^{\pm \frac{1}{2}}]$. 
This fact plays a crucial role in the construction of $F_q(\tm)$. 
%
We present that the directed colored graphs $\Gamma(X_{2,5})$ and $\Gamma(\tX_{2,5})$ of $F(X_{2,5}) \in \frakK(\sfg)$ and $F_q(\tX_{2,5}) \in \frakK_q(\sfg)$ are given as follows, respectively:
  \begin{equation} \label{ex: folding}
  \begin{split}
  \qquad \quad
  \raisebox{13pc}{
      \scalebox{0.5}{$
      \begin{tikzpicture}[baseline=(current  bounding  box.center), every node/.style={scale=1.31}]
      \node (node_2) at (6, 0) {$X_{2,{5}}$};
      \node (node_3) at (6, -2) {$X_{1,{6}}^{3} X_{2,{7}}^{-1}$};
      \node (node_4) at (6, -4) {$3X_{1,{6}}^{2} X_{1,{8}}^{-1}$};
      \node (node_5) at (6, -6) {$3X_{1,{6}} X_{2,{7}} X_{1,{8}}^{-2}$};
      \node (node_6) at (3, -8) {$X_{2,{7}}^{2} X_{1,{8}}^{-3}$};
      \node (node_7) at (9, -8) {$3X_{1,{6}} X_{1,{8}} X_{2,{9}}^{-1}$};
      \node (node_8) at (3, -10) {$2X_{2,{7}} X_{2,{9}}^{-1}$};
      \node (node_10) at (9, -10) {$3X_{1,{6}} X_{1,{10}}^{-1}$};
      \node (node_9) at (3, -12) {$X_{1,{8}}^{3} X_{2,{9}}^{-2}$};
      \node (node_1) at (9, -12) {$3X_{2,{7}} X_{1,{8}}^{-1} X_{1,{10}}^{-1}$};
      \node (node_0) at (6, -14) {$3X_{1,{8}}^{2}X_{2,{9}}^{-1}X_{1,{10}}^{-1}$};
      \node (node_11) at (6, -16) {$3X_{1,{8}} X_{1,{10}}^{-2}$};
      \node (node_12) at (6, -18) {$X_{2,{9}} X_{1,{10}}^{-3}$};
      \node (node_13) at (6, -20) {$X_{2,{11}}^{-1}$};
      \draw[-{>[scale=1.5]}] (node_8) -- node[fill=white, anchor=center]  {\tiny 2, 8} (node_9);
      \draw[-{>[scale=1.5]}] (node_10) -- node[fill=white, anchor=center] {\tiny 1, 7} (node_1);
      \draw[-{>[scale=1.5]}] (node_7) -- node[fill=white, anchor=center]  {\tiny 1, 9} (node_10);
      \draw[-{>[scale=1.5]}] (node_5) -- node[fill=white, anchor=center]  {\tiny 2, 8} (node_7);
      \draw[-{>[scale=1.5]}] (node_5) -- node[fill=white, anchor=center]  {\tiny 1, 7} (node_6);
      \draw[-{>[scale=1.5]}] (node_4) -- node[fill=white, anchor=center]  {\tiny 1, 7} (node_5);
      \draw[-{>[scale=1.5]}] (node_11) -- node[fill=white, anchor=center] {\tiny 1, 9} (node_12);
      \draw[-{>[scale=1.5]}] (node_1) -- node[fill=white, anchor=center]  {\tiny 2, 8} (node_0);
      \draw[-{>[scale=1.5]}] (node_0) -- node[fill=white, anchor=center]  {\tiny 1, 9} (node_11);
      \draw[-{>[scale=1.5]}] (node_12) -- node[fill=white, anchor=center] {\tiny 2, 10} (node_13);
      \draw[-{>[scale=1.5]}] (node_3) -- node[fill=white, anchor=center]  {\tiny 1, 7} (node_4);
      \draw[-{>[scale=1.5]}] (node_6) -- node[fill=white, anchor=center]  {\tiny 2, 8} (node_8);
      \draw[-{>[scale=1.5]}] (node_9) -- node[fill=white, anchor=center]  {\tiny 1, 9} (node_0);
      \draw[-{>[scale=1.5]}] (node_2) -- node[fill=white, anchor=center]  {\tiny 2, 6} (node_3);
        \end{tikzpicture}$}}
      & \quad \quad \quad
      \raisebox{13pc}{
      \scalebox{0.5}{$
      \begin{tikzpicture}[baseline=(current  bounding  box.center), every node/.style={scale=1.31}]
      \node (node_2) at (6, 0) {$\tX_{2,{5}}$};
      \node (node_3) at (6, -2) {$q^3 \tX_{1,{6}}^{3} *\tX_{2,{7}}^{-1}$};
      \node (node_4) at (6, -4) {$(q^{-2} + 1 + q^2)\tX_{1,{6}}^{2} *\tX_{1,{8}}^{-1}$};
      \node (node_5) at (6, -6) {$(1 + q^2 + q^4)\tX_{1,{6}} *\tX_{2,{7}}* \tX_{1,{8}}^{-2}$};
      \node (node_6) at (3, -8) {$q^9  \tX_{2,{7}}^{2} *\tX_{1,{8}}^{-3}$};
      \node (node_7) at (11, -8) {$(q^{-3} + q^{-1} + q)\tX_{1,{6}} *\tX_{1,{8}} *\tX_{2,{9}}^{-1}$};
      \node (node_8) at (3, -10) {$(q^{-3} + q^3)\tX_{2,{7}} *\tX_{2,{9}}^{-1}$};
      \node (node_10) at (11, -10) {$(q^{-4} + q^{-2} + 1)\tX_{1,{6}}* \tX_{1,{10}}^{-1}$};
      \node (node_9) at (3, -12) {$q^6  \tX_{1,{8}}^{3}* \tX_{2,{9}}^{-2}$};
      \node (node_1) at (11, -12) {$(q^{-2} + 1 + q^2)\tX_{2,{7}}* \tX_{1,{8}}^{-1} *\tX_{1,{10}}^{-1}$};
      \node (node_0) at (6, -14) {$(q^{-2} + 1 + q^2)\tX_{1,{8}}^{2} *\tX_{2,{9}}^{-1} *\tX_{1,{10}}^{-1}$};
      \node (node_11) at (6, -16) {$(q^{-3} + q^{-1} + 1)\tX_{1,{8}} *\tX_{1,{10}}^{-2}$};
      \node (node_12) at (6, -18) {$q^3  \tX_{2,{9}} *\tX_{1,{10}}^{-3}$};
      \node (node_13) at (6, -20) {$q^{-3}  \tX_{2,{11}}^{-1}$};
      \draw[-{>[scale=1.5]}] (node_8) -- node[fill=white, anchor=center]  {\tiny 2, 8} (node_9);
      \draw[-{>[scale=1.5]}] (node_10) -- node[fill=white, anchor=center] {\tiny 1, 7} (node_1);
      \draw[-{>[scale=1.5]}] (node_7) -- node[fill=white, anchor=center]  {\tiny 1, 9} (node_10);
      \draw[-{>[scale=1.5]}] (node_5) -- node[fill=white, anchor=center]  {\tiny 2, 8} (node_7);
      \draw[-{>[scale=1.5]}] (node_5) -- node[fill=white, anchor=center]  {\tiny 1, 7} (node_6);
      \draw[-{>[scale=1.5]}] (node_4) -- node[fill=white, anchor=center]  {\tiny 1, 7} (node_5);
      \draw[-{>[scale=1.5]}] (node_11) -- node[fill=white, anchor=center] {\tiny 1, 9} (node_12);
      \draw[-{>[scale=1.5]}] (node_1) -- node[fill=white, anchor=center]  {\tiny 2, 8} (node_0);
      \draw[-{>[scale=1.5]}] (node_0) -- node[fill=white, anchor=center]  {\tiny 1, 9} (node_11);
      \draw[-{>[scale=1.5]}] (node_12) -- node[fill=white, anchor=center] {\tiny 2, 10} (node_13);
      \draw[-{>[scale=1.5]}] (node_3) -- node[fill=white, anchor=center]  {\tiny 1, 7} (node_4);
      \draw[-{>[scale=1.5]}] (node_6) -- node[fill=white, anchor=center]  {\tiny 2, 8} (node_8);
      \draw[-{>[scale=1.5]}] (node_9) -- node[fill=white, anchor=center]  {\tiny 1, 9} (node_0);
      \draw[-{>[scale=1.5]}] (node_2) -- node[fill=white, anchor=center]  {\tiny 2, 6} (node_3);
       \end{tikzpicture}$}}
    \end{split}
\end{equation}
\noindent
Here $F(X_{2,5})$ is obtained from ${\rm ev}_{t=1}(\sfT_5(F_t(\tY_{2,0})))$\,\,(see~Example \ref{ex: F20 in D4} for $F_t(\tY_{2,0})$) by folding the $\calY$-monomials (recall Remark \ref{rem: F vs Ft in ADE}).
\end{example}

\begin{remark} \label{rem: Ft vs Fq and independence of q-algorithm to J} 
Each $\Z[q^{\pm \frac{1}{2}}]$-coefficient determined by the $q$-algorithm does not depend on the choice of $J$, which can be proved as in \cite[Lemma 5.20]{H04}.
Let us recall that $F(X_{\oi,p})$ is obtained from $F(Y_{\im,p})$ by folding the monomials of $F(Y_{\im,p})$ via \eqref{eq: osigma} (see~Corollary~\ref{cor: some F to F}). However, we would like to emphasize that we do not know yet whether $F_q(\tX_{\oi,p})$ could be obtained directly from $F_t(\tY_{\im,p})$ by folding $\calY_t$-monomials with some modification of coefficients in $\Z[q^{\pm \frac{1}{2}}]$.
\end{remark}

\begin{definition} \label{def: monomials of Xq}
For $f \in \calX_q$, we set
\begin{align*}
\calM(f) & \seteq \{   \ev_{q=1}(\tm)  \ | \ \tm \text{ is a monomial in $f$}  \},  && \;\; \ucalM(f) \seteq \{ \um  \ | \  m \in \calM(f) \}, \\
\calM_+(f) &\seteq \{  \ev_{q=1}(\tm)  \ | \ \tm \text{ is a dominant monomial in $f$} \}, &&
 \ucalM_+(f) \seteq \{ \um  \ | \ m \in \ucalM_+(f) \}.
\end{align*}
\end{definition}

For $P \in \frakK_q(\sfg)$, a monomial $m$ in $P$ is called {\it   maximal   monomial} (resp. {\it   minimal   monomial}) if its weight is not lower (resp. not higher) than any other monomial in $P$ with respect to $\leN$.

\begin{lemma}{\em (cf. \cite[Lemma 5.6]{FM01})} \label{lem:maximal or minimal monomials}
For $P \in \frakK_q(\sfg)$, any maximal $($resp. minimal$)$ monomial in $P$ is dominant $($resp. anti-dominant$)$.	
\end{lemma}

\begin{proof}
Let us first consider a maximal monomial in $P$, denoted by $\tm$. Take $i \in I$.
By Definition \ref{def: virtual quantum G ring} and Proposition \ref{prop: YA com}, we have
\begin{equation*}
	P \in \frakK_{i,q}(\sfg) = \Z[q^{\pm \frac{1}{2}}][\tX_{j,l}^{\pm 1}]_{(j,l)\in \widetilde{\sDynkinf}_0,\, j\neq i} \otimes \Z[q^{\pm \frac{1}{2}}][\tX_{i,l} + q_i^{-1}\tX_{i,l} * \tB_{i,l+1}^{-1}]_{(i,l)\in \widetilde{\sDynkinf}_0}.
\end{equation*}
Hence, the element $P$ can be written in the following form:
\begin{equation*}
	P = \sum \tm_{(1)} * \tp_{(2)},
\end{equation*}
where $\tm_{(1)} \in \Z[q^{\pm \frac{1}{2}}][\tX_{j,l}]_{(j,l)\in\widetilde{\sDynkinf}_0,\, j\neq i}$ are monomials and,
$\tp_{(2)} \in \Z[q^{\pm \frac{1}{2}}][\tX_{i,l}+q_i^{-1}\tX_{i,l} *\tB_{i,l+1}^{-1}]_{(i,l)\in\widetilde{\sDynkinf}_0}$ are of the form
\begin{align*}
\tp_{(2)} = n \; c(q) \st^{\to}_{  \substack{ (i,l) \in \widetilde{\sDynkinf}_0,  \\ \text{ finite } }}  (\tX_{i,l}+q_i^{-1}\tX_{i,l} * \tB_{i,l+1}^{-1})^{n_{i,l}}
\end{align*}
for some $n_{i,l} \in \Z_{\ge 1}$, $n \in \Z$ and $c(q) \in q^{\frac{1}{2}\Z}$. 
\noindent
In particular, the maximal monomial $\tm$ is a monomial in $\tm_{(1)} * \tp_{(2)}$.
Since $ X_{i,l} B_{i,l+1}^{-1}   \lN   X_{i,l}$
the monomial $\tm$ should be  obtained from $\tm_{(1)}$ and $ \tX_{i,l}$'s, .
Otherwise, it contradicts the assumption that $m$ is a maximal monomial.
Since $i \in I$ is arbitrary and $P \in \frakK_q(\sfg)$, the maximal monomial $m$ should be dominant.
\noindent
In the case of minimal monomials, the proof is almost identical because of another characterization of $\frakK_{i,q}(\sfg)$ in Remark \eqref{rem: some remark for screenings}.
\end{proof}

For an interval $[a,b]$, $i \in I$, $(i,t) \in \tDynkin_0$ and $k \in \Z_{\ge 1}$, we define
\begin{equation}  \label{eq: dominant monomial of KR type}
\begin{aligned}
m^{(i)}[a,b]  \seteq \displaystyle\prod_{\substack{(i,p) \in \widetilde{\sDynkinf}_0 \\ p \in [a,b]}} X_{i,p}\quad \text{ and }\quad
m^{(i)}_{k,t} \seteq \displaystyle \prod_{s=0}^{k-1} X_{i,t+2s}.
\end{aligned}
\end{equation}
We define $m^{(i)}(a,b]$, $m^{(i)}[a,b)$, and $m^{(i)}(a,b)$ in a similar way.

\smallskip 

The following proposition plays a crucial role in proving fundamental results established on $\frakK_{q}(\sfg)$.

\begin{proposition} \label{prop:folding lemma}
For $p < s$,
let $m^{(i)}[p,s]$ be given such that $\osigma(\bfm^{(\im)}[p,s]) = m^{(i)}[p,s]$ $($i.e.~$\overline{\im} = i$$)$.
\ben
\item \label{folding lemma:1} For each $m' \in \calM(F_q(\um^{(i)}[p,s]))$,
there exists  $\bfm' \in \calM(F_t(\ubfm^{(\im)}[p,s]))$ such that $\osigma(\bfm') =  m'$.
%
\item \label{folding lemma:2} We have $F_q(\um^{(i)}[p,s]) \in \frakK_q(\sfg)$.
\ee
\end{proposition}

\begin{proof}
We will give a proof of Proposition \ref{prop:folding lemma} in Section \ref{subsec:proof of folding lemma}.
\end{proof}

\begin{definition} \label{def: KR-type}
We call an element of the form $F_q(\um^{(i)}[p,s])$ a \emph{KR-polynomial}. In particular,
we call $F_q(X_{i,p})$ a \emph{fundamental polynomial}. We also call a monomial of the form $m^{(i)}[p,s]$ a \emph{KR-monomial}.
\end{definition}

\begin{corollary} \label{cor: folding monomials of Ft = monomials of Fq}
For $p < s$,
let $m^{(i)}[p,s]$ be such that $\osigma(\bfm^{(\im)}[p,s]) = m^{(i)}[p,s]$ $($i.e.~$\overline{\im} = i$$)$.
Then we have
\begin{equation*}
	\osigma\left( \calM(F_t(\ubfm^{(\im)}[p,s])) \right) = \calM(F_q(\um^{(i)}[p,s])).
\end{equation*}
\end{corollary}
\begin{proof}
The inclusion $\supset$ follows from Proposition \ref{prop:folding lemma}  \ref{folding lemma:1}.
Let us prove the opposite inclusion $\subset$.
Let $m \in \osigma\left( \calM(F_t(\ubfm^{(\im)}[p,s])) \right)$ be an $\calX$-monomial, where we write $m = \osigma(\bfm)$ for some $\calY$-monomial $\bfm \in \calM(F_t(\ubfm^{(\im)}[p,s]))$.
We have seen
\begin{equation} \label{eq: evaluation KR types}
	{\rm ev}_{t=1} \left( F_t(\ubfm^{(\im)}[p,s]) \right) = F(\bfm^{(\im)}[p,s])
\end{equation}
(see Remark \ref{rem: F vs Ft in ADE}), and then
the quantum positivity for $F_t(\ubfm^{(\im)}[p,s])$ in Theorem \ref{thm: Nak positivity} with \eqref{eq: evaluation KR types} implies that all the coefficients of $F(\bfm^{(\im)}[p,s])$ should be positive. In particular, the coefficient of $\bfm$ in $F(\bfm^{(\im)}[p,s])$ is positive.
Since it follows from Corollary \ref{cor: some F to F} and Theorem \ref{thm:monomials for KR modules in types ADE} that
\begin{equation} \label{eq: folding KR types}
	\osigma\left( F(\bfm^{(\im)}[p,s]) \right) = F(m^{(i)}[p,s]),
\end{equation}
the $\calX$-monomial $m$ appears in $F(m^{(i)}[p,s])$ with a positive coefficient. But, we have
\begin{equation*}
	{\rm ev}_{q=1} \left( F_q(\um^{(i)}[p,s]) \right) = F(m^{(i)}[p,s])
\end{equation*}
(see Corollary \ref{cor: evaluation folded ver}), which implies that there exists a term $f(q) \um$ in $F_q(\um^{(i)}[p,s])$ such that ${\rm ev}_{q=1}(f(q)\um) = f(1)m$ is a term in $F(m^{(i)}[p,s])$ with $f(1) > 0$. 
\end{proof}

\begin{proposition} \label{prop:properties of Fq}
For each $(i,p) \in \tbDynkinf_0$, we have
\bna
\item \label{it: anti sfg}  $F_q(\sfX_{i,p}) = F_q(\sfX^{-1}_{i^*,p+\sfh})$ contains only one anti-dominant monomial $\sfX^{-1}_{i^*,p+\sfh}$.
\item \label{it: btw sfg} All $\calX_q$-monomials of $F_q(\sfX_{i,p}) - \sfX_{i,p}$ $ -\sfX_{i,p+\sfh}^{-1}$ are products of $\tX_{j,u}^{\pm1}$ with $p < u$ {$  <p+\sfh$}.
\item \label{it: commute sfg} $F_q(\tX_{i,p})$ and $F_q(\tX_{j,p})$  $((i,p),(j,p) \in \tbDynkinf_0, \; j \ne i)$ commute.
\ee
\end{proposition}
\begin{proof}
Since $F_q(\sfX_{i,p})$ is an element in $\frakK_q(\sfg)$, it contains an anti-dominant monomial by Lemma~\ref{lem:maximal or minimal monomials}. Then Theorem~\ref{thm:monomials for KR modules in types ADE} and  Proposition~\ref{prop:folding lemma}
tell that $F_q(\sfX_{i,p})$ has the unique antidominant monomial $\sfX_{i^*,p+\sfh}^{-1}$. Thus \ref{it: anti sfg} follows.
By~\eqref{eq: preserving Nakajima ordering}, \ref{it: btw sfg} follows from~\ref{it: anti sfg}, Theorem~\ref{thm:monomials for KR modules in types ADE}
and Proposition~\ref{prop:folding lemma}. Finally,~\ref{it: commute sfg} follows from the same argument as in \cite[Lemma 5.12 (iv)]{H04}.
\end{proof}

\begin{example} \label{ex: F110 in G2}
As in Example \ref{ex: F25 in G2},
one may compute the formula of $F_q(\tX_{1,10})$ given by
\begin{equation*}
\tX_{1,10} +
q^{2} \tX_{2,11} \tX_{1,12}^{-1} +
q^{2} \tX_{1,12}^{2} \tX_{2,13}^{-1} +
(q^{-1} + q) \tX_{1,12} \tX_{1,14}^{-1} + q^{3} \tX_{2,13} \tX_{1,14}^{-2} + \tX_{1,14} \tX_{2,15}^{-1} +
q^{-1} \tX_{1,16}^{-1}.
\end{equation*}
Then $F_q(\sfX_{1,10}) = (\sfX_{1,10} * \tX_{1,10}^{-1}) F_q(\tX_{1,10}) = q^{\frac{1}{2}}F_q(\tX_{1,10}) \in \mathfrak{K}_q(\sfg)$ is bar-invariant.
Note that $\tX_{1,10} * \tX_{2,10} = \tX_{2,10} * \tX_{1,10}$ 
and there is no dominant $\calX_q$-monomial in $\calM_+(F_q(\tX_{1,10}) * F_q(\tX_{2,10}))$ except for $\tX_{1,10}*\tX_{2,10}$ (cf.~Example \ref{ex: F25 in G2}).
Hence we have
$F_q(\tX_{1,10}) * F_q(\tX_{2,10}) = F_q(\tX_{2,10}) * F_q(\tX_{1,10})$.
\end{example}

\begin{example} \label{ex: fundamentals in B3}
By the $q$-algorithm starting from $\tX_{i,p}$ as in Example \ref{ex: F25 in G2}, 
one can compute the explicit formulas of $F_q(\tX_{i,p})$ for $1 \le i \le 3$ of the finite type $B_3$ as follows:
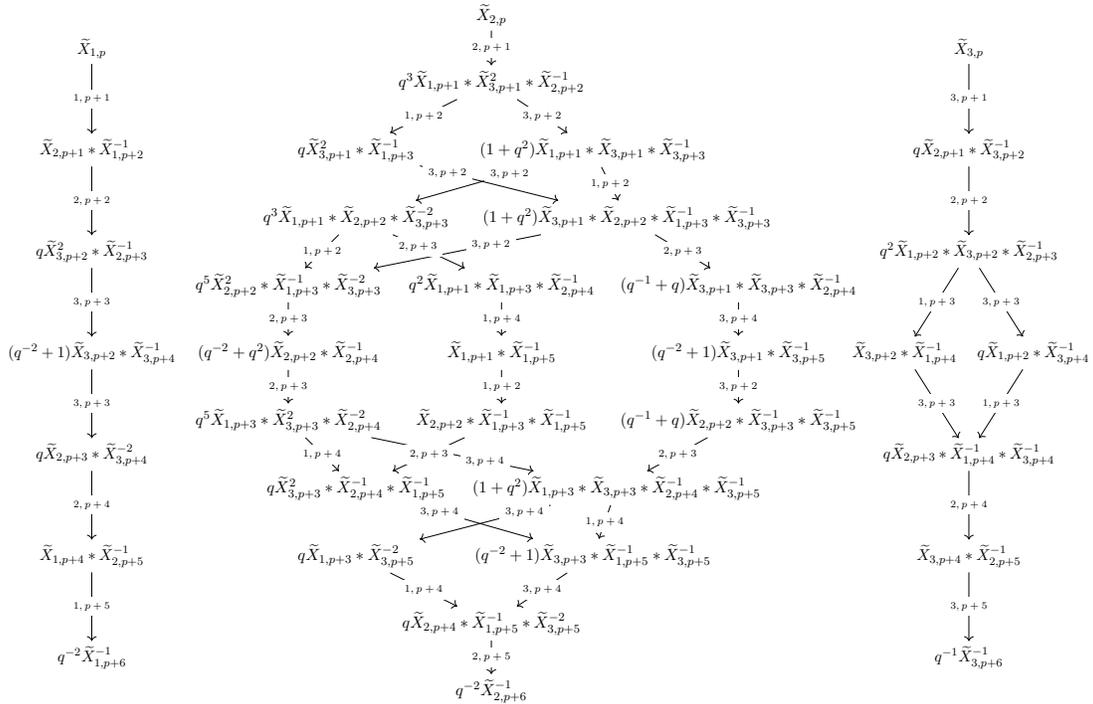
\begin{equation*}
\begin{split}
\raisebox{10pc}{
	\scalebox{0.45}{$
	\begin{tikzpicture}[baseline=(current  bounding  box.center), every node/.style={scale=1.2}]
	\node (node_0) at (2, 0) {$\tX_{1,p}$};
	\node (node_1) at (2, -3) {$\tX_{2,p+1} * \tX_{1,p+2}^{-1}$};
	\node (node_2) at (2, -6) {$q \tX_{3,p+2}^{2} * \tX_{2,p+3}^{-1}$};
	\node (node_3) at (2, -9) {$(q^{-2} + 1) \tX_{3,p+2} * \tX_{3,p+4}^{-1}$};
	\node (node_4) at (2, -12) {$q \tX_{2,p+3} * \tX_{3,p+4}^{-2}$};
	\node (node_5) at (2, -15) {$\tX_{1,p+4} * \tX_{2,p+5}^{-1}$};
	\node (node_6) at (2, -18) {$q^{-2} \tX_{1,p+6}^{-1}$};
	\draw[-{>[scale=1.5]}] (node_3) -- node[fill=white, anchor=center, pos=0.5] {\tiny $3, p+3$} (node_4);
	\draw[-{>[scale=1.5]}] (node_0) -- node[fill=white, anchor=center, pos=0.5] {\tiny $1, p+1$} (node_1);
	\draw[-{>[scale=1.5]}] (node_1) -- node[fill=white, anchor=center, pos=0.5] {\tiny $2, p+2$} (node_2);
	\draw[-{>[scale=1.5]}] (node_5) -- node[fill=white, anchor=center, pos=0.5] {\tiny $1, p+5$} (node_6);
	\draw[-{>[scale=1.5]}] (node_4) -- node[fill=white, anchor=center, pos=0.5] {\tiny $2, p+4$} (node_5);
	\draw[-{>[scale=1.5]}] (node_2) -- node[fill=white, anchor=center, pos=0.5] {\tiny $3, p+3$} (node_3);
	\end{tikzpicture}$}}
	\raisebox{10pc}{
	\scalebox{0.45}{$
	\begin{tikzpicture}[baseline=(current  bounding  box.center), every node/.style={scale=1.25}]
\node (node_0) at (11, 0) {$\tX_{2,p}$};
\node (node_1) at (11, -2) {$q^{3} \tX_{1,p+1}*\tX_{3,p+1}^{2}*\tX_{2,p+2}^{-1}$};
\node (node_2) at (7, -4) {$q \tX_{3,p+1}^{2}*\tX_{1,p+3}^{-1}$};
\node (node_3) at (14, -4) {$(1 + q^{2}) \tX_{1,p+1}*\tX_{3,p+1}*\tX_{3,p+3}^{-1}$};
\node (node_4) at (7, -6) {$q^{3} \tX_{1,p+1}*\tX_{2,p+2}*\tX_{3,p+3}^{-2}$};
\node (node_5) at (15, -6) {$(1 + q^{2}) \tX_{3,p+1}*\tX_{2,p+2}*\tX_{1,p+3}^{-1}*\tX_{3,p+3}^{-1}$};
\node (node_6) at (5, -8) {$q^{5} \tX_{2,p+2}^{2}*\tX_{1,p+3}^{-1}*\tX_{3,p+3}^{-2}$};
\node (node_7) at (11.3, -8) {$q^{2} \tX_{1,p+1}*\tX_{1,p+3}*\tX_{2,p+4}^{-1}$};
\node (node_8) at (18.3, -8) {$(q^{-1} + q) \tX_{3,p+1}*\tX_{3,p+3}*\tX_{2,p+4}^{-1}$};
\node (node_9) at (5, -10) {$(q^{-2} + q^{2}) \tX_{2,p+2}*\tX_{2,p+4}^{-1}$};
\node (node_10) at (5, -12) {$q^{5} \tX_{1,p+3}*\tX_{3,p+3}^{2}*\tX_{2,p+4}^{-2}$};
\node (node_11) at (11.3, -10) {$\tX_{1,p+1}*\tX_{1,p+5}^{-1}$};
\node (node_12) at (11.3, -12) {$\tX_{2,p+2}*\tX_{1,p+3}^{-1}*\tX_{1,p+5}^{-1}$};
\node (node_13) at (18.3, -10) {$(q^{-2} + 1) \tX_{3,p+1}*\tX_{3,p+5}^{-1}$};
\node (node_14) at (18.3, -12) {$(q^{-1} + q) \tX_{2,p+2}*\tX_{3,p+3}^{-1}*\tX_{3,p+5}^{-1}$};
\node (node_15) at (7, -14) {$q \tX_{3,p+3}^{2}*\tX_{2,p+4}^{-1}*\tX_{1,p+5}^{-1}$};
\node (node_16) at (14.7, -14) {$(1 + q^{2}) \tX_{1,p+3}*\tX_{3,p+3}*\tX_{2,p+4}^{-1}*\tX_{3,p+5}^{-1}$};
\node (node_17) at (7, -16) {$q \tX_{1,p+3}*\tX_{3,p+5}^{-2}$};
\node (node_18) at (14, -16) {$(q^{-2} + 1) \tX_{3,p+3}*\tX_{1,p+5}^{-1}*\tX_{3,p+5}^{-1}$};
\node (node_19) at (11, -18) {$q \tX_{2,p+4}*\tX_{1,p+5}^{-1}*\tX_{3,p+5}^{-2}$};
\node (node_20) at (11, -20) {$q^{-2} \tX_{2,p+6}^{-1}$};
\draw[-{>[scale=1.5]}] (node_16) -- node[fill=white, anchor=center, pos=0.5] {\tiny $1, p+4$} (node_18);
\draw[-{>[scale=1.5]}] (node_16) -- node[fill=white, anchor=center, pos=0.2] {\tiny $3, p+4$} (node_17);
\draw[-{>[scale=1.5]}] (node_5) -- node[fill=white, anchor=center,  pos=0.5] {\tiny $2, p+3$} (node_8);
\draw[-{>[scale=1.5]}] (node_5) -- node[fill=white, anchor=center,  pos=0.3] {\tiny $3, p+2$} (node_6);
\draw[-{>[scale=1.5]}] (node_18) -- node[fill=white, anchor=center, pos=0.5] {\tiny $3, p+4$} (node_19);
\draw[-{>[scale=1.5]}] (node_3) -- node[fill=white, anchor=center,  pos=0.5] {\tiny $1, p+2$} (node_5);
\draw[-{>[scale=1.5]}] (node_3) -- node[fill=white, anchor=center,  pos=0.2] {\tiny $3, p+2$} (node_4);
\draw[-{>[scale=1.5]}] (node_9) -- node[fill=white, anchor=center,  pos=0.5] {\tiny $2, p+3$} (node_10);
\draw[-{>[scale=1.5]}] (node_14) -- node[fill=white, anchor=center, pos=0.5] {\tiny $2, p+3$} (node_16);
\draw[-{>[scale=1.5]}] (node_8) -- node[fill=white, anchor=center,  pos=0.5] {\tiny $3, p+4$} (node_13);
\draw[-{>[scale=1.5]}] (node_13) -- node[fill=white, anchor=center, pos=0.5] {\tiny $3, p+2$} (node_14);
\draw[-{>[scale=1.5]}] (node_10) -- node[fill=white, anchor=center, pos=0.7] {\tiny $3, p+4$} (node_16);
\draw[-{>[scale=1.5]}] (node_10) -- node[fill=white, anchor=center, pos=0.5] {\tiny $1, p+4$} (node_15);
\draw[-{>[scale=1.5]}] (node_17) -- node[fill=white, anchor=center, pos=0.5] {\tiny $1, p+4$} (node_19);
\draw[-{>[scale=1.5]}] (node_12) -- node[fill=white, anchor=center, pos=0.5] {\tiny $2, p+3$} (node_15);
\draw[-{>[scale=1.5]}] (node_6) -- node[fill=white, anchor=center,  pos=0.5] {\tiny $2, p+3$} (node_9);
\draw[-{>[scale=1.5]}] (node_2) -- node[fill=white, anchor=center,  pos=0.2] {\tiny $3, p+2$} (node_5);
\draw[-{>[scale=1.5]}] (node_19) -- node[fill=white, anchor=center, pos=0.5] {\tiny $2, p+5$} (node_20);
\draw[-{>[scale=1.5]}] (node_15) -- node[fill=white, anchor=center, pos=0.2] {\tiny $3, p+4$} (node_18);
\draw[-{>[scale=1.5]}] (node_7) -- node[fill=white, anchor=center,  pos=0.5] {\tiny $1, p+4$} (node_11);
\draw[-{>[scale=1.5]}] (node_11) -- node[fill=white, anchor=center, pos=0.5] {\tiny $1, p+2$} (node_12);
\draw[-{>[scale=1.5]}] (node_4) -- node[fill=white, anchor=center,  pos=0.5] {\tiny $1, p+2$} (node_6);
\draw[-{>[scale=1.5]}] (node_4) -- node[fill=white, anchor=center,  pos=0.35] {\tiny $2, p+3$} (node_7);
\draw[-{>[scale=1.5]}] (node_1) -- node[fill=white, anchor=center,  pos=0.5] {\tiny $3, p+2$} (node_3);
\draw[-{>[scale=1.5]}] (node_1) -- node[fill=white, anchor=center,  pos=0.5] {\tiny $1, p+2$} (node_2);
\draw[-{>[scale=1.5]}] (node_0) -- node[fill=white, anchor=center,  pos=0.5] {\tiny $2, p+1$} (node_1);
\end{tikzpicture}$}}
	\!\!\!\!\!
	\raisebox{10pc}{
	\scalebox{0.45}{$
\begin{tikzpicture}[baseline=(current  bounding  box.center), every node/.style={scale=1.2}]
\node (node_0) at (6, 0) {$\tX_{3,p}$};
\node (node_1) at (6, -3) {$q \tX_{2,p+1} * \tX_{3,p+2}^{-1}$};
\node (node_2) at (6, -6) {$q^2 \tX_{1,p+2} * \tX_{3,p+2} * \tX_{2,p+3}^{-1} $};
\node (node_3) at (4.1, -9) {$\tX_{3,p+2} * \tX_{1,p+4}^{-1}$};
\node (node_4) at (7.9, -9) {$q \tX_{1,p+2} * \tX_{3,p+4}^{-1}$};
\node (node_5) at (6, -12) {$q \tX_{2,p+3} * \tX_{1,p+4}^{-1} * \tX_{3,p+4}^{-1}$};
\node (node_6) at (6, -15) {$\tX_{3,p+4} * \tX_{2,p+5}^{-1}$};
\node (node_7) at (6, -18) {$q^{-1} \tX_{3,p+6}^{-1}$};
\draw[-{>[scale=1.5]}] (node_2) -- node[fill=white, anchor=center, pos=0.5] {\tiny $3, p+3$} (node_4);
\draw[-{>[scale=1.5]}] (node_2) -- node[fill=white, anchor=center, pos=0.5] {\tiny $1, p+3$} (node_3);
\draw[-{>[scale=1.5]}] (node_4) -- node[fill=white, anchor=center, pos=0.5] {\tiny $1, p+3$} (node_5);
\draw[-{>[scale=1.5]}] (node_5) -- node[fill=white, anchor=center, pos=0.5] {\tiny $2, p+4$} (node_6);
\draw[-{>[scale=1.5]}] (node_3) -- node[fill=white, anchor=center, pos=0.5] {\tiny $3, p+3$} (node_5);
\draw[-{>[scale=1.5]}] (node_6) -- node[fill=white, anchor=center, pos=0.5] {\tiny $3, p+5$} (node_7);
\draw[-{>[scale=1.5]}] (node_1) -- node[fill=white, anchor=center, pos=0.5] {\tiny $2, p+2$} (node_2);
\draw[-{>[scale=1.5]}] (node_0) -- node[fill=white, anchor=center, pos=0.5] {\tiny $3, p+1$} (node_1);
\end{tikzpicture}$}}
\end{split}
\end{equation*}
\end{example}

For a dominant monomial $m \in \calM_+^\sfg$, we set
\begin{equation} \label{eq: Eqm}
\begin{split}
E_q(\um) \seteq q^b
\left(  \st_{p\in\Z}^{\to} \left(  \st_{i\in I; (i,p)\in\sDynkinf_0} F_q(\sfX_{i,p} )^{u_{i,p}(m)}  \right) \right) \in \frakK_q(\sfg),
\end{split}
\end{equation}
where $b$ is an element in $\frac{1}{2}\Z$ such that $\um$ appears in $E_q(\um) $ with the coefficient $1$.
By Proposition \ref{prop:folding lemma},
we have
\begin{equation} \label{eq: finiteness of Eqm}
	E_q(\um) \in \frakK_q(\sfg)
\end{equation}
and there are finitely many dominant monomials in $E_q(\um)$.
As we regard $E_q(\um)$ as an element of $\frakK_q^{\infty,\,f}(\sfg)$ (recall Proposition \ref{prop: basis Fq in Kqf_infty}),
we obtain a uni-triangular transition map as
in~\eqref{eq: uni t1} between $\{ E_q(\um) \}$ and $\{ F_q(\um) \}$ in $\frakK_q^{\infty,\,f}(\sfg)$ by Proposition \ref{prop: basis Fq in Kqf_infty}:
\begin{equation}\label{eq: uni t2 sfg}
\begin{split}
E_q(\um) = F_q(\um) + \sum_{m' \lN m} C_{m,m'} F_q(\um') \quad \text{ in \,\,$\frakK_q^{\infty,\,f}(\sfg)$, }
\end{split}
\end{equation}
where $C_{m,m'} \in \Z[q^{\pm \frac{1}{2}}]$.
Note that the summation in \eqref{eq: uni t2 sfg} is finite since $E_q(\um)$ has finitely many dominant monomials.
Hence, \eqref{eq: uni t2 sfg} implies that $F_q(\um)$ can be written as a linear combination of $E_q(\um')$ for $\um' \leN \um$, so $F_q(\um) \in \frakK_q(\sfg)$ by \eqref{eq: Eqm} and \eqref{eq: finiteness of Eqm}.
Until now, we have proved the following.

\begin{proposition} \label{prop: bases Eq and Fq in Kqg}
The sets
\begin{equation*}
\begin{split}
& \sfE_q \seteq \{ E_q(\um) \ |  \ m \in \calM^\sfg_+\} \quad \text{ and } \quad
\sfF_q \seteq \{   F_q(\um) \ | \ m \in \calM^\sfg_+ \}
\end{split}
\end{equation*}
are $\Z[q^{\pm \frac{1}{2}}]$-bases of $\frakK_q(\sfg)$, respectively.
\end{proposition}
In particular, we call $\sfE_q$ the \emph{standard basis} of $\frakK_q(\sfg)$.
Now let us further investigate the basis $\sfF_q$ of $\frakK_q(\sfg)$, which is characterized as follows:

\begin{theorem} \label{thm: F_q}
Let $\tm \in \calX_q$ be a dominant $($resp.~anti-dominant$)$ monomial.
\bna
	\item \label{it:unique dominant monomial} The Laurent $($non-commutative$)$ polynomial $F_q(\tm)$ is the
unique element in $\frakK_q(\sfg)$ such that $\tm$ is the unique
dominant $($resp.~anti-dominant$)$ monomial occurring in $F_q(\tm)$.
	\vskip 1mm

	\item \label{it:right negativeness} Every monomial in $F_q(\tm) - \tm$ is strictly less $($resp.~greater$)$ than $\tm$ with respect to $\lN$.
	\vskip 1mm
	
	\item \label{it:bar-invariant basis} The set
$\sfF_q$ forms a  bar-invariant
$\Z[q^{\pm \frac{1}{2}}]$-basis of $\frakK_q(\sfg)$.
\ee
\end{theorem}

\begin{proof}
We prove only the dominant case because the proof for the anti-dominant case is almost identical.
Let us first prove \ref{it:unique dominant monomial}.
Thanks to \eqref{eq: uni t2 sfg},
$F_q(\tm)$ can be written as a linear combination of $E_q(\tm')$ for $\tm' \leN \tm$, where the sum is finite due to Proposition \ref{prop:folding lemma}(\ref{folding lemma:2}).
Hence, $F_q(\tm) \in \frakK_q(\sfg)$.
Note that $F_q(\tm)$ has the unique dominant monomial $\tm$ by its construction through the $q$-algorithm (see \eqref{eq: coloring and coefficient in q-algorithm}).
Let $G_q(\tm)$ be another element in $\frakK_q(\sfg)$ such that $\tm$ is the unique dominant monomial occurring in $G_q(\tm)$.
Then $F_q(\tm) - G_q(\tm) \in \frakK_q(\sfg)$ has no dominant monomial and should be $0$ by Lemma \ref{lem:maximal or minimal monomials}, otherwise it has a dominant monomial different from $\tm$. Thus, the uniqueness of $F_q(\tm)$ in $\frakK_q(\sfg)$ follows.
Second, the \ref{it:right negativeness} is a direct consequence of the $q$-algorithm.

Finally, let us prove the \ref{it:bar-invariant basis}.
The linear independence follows from the uniqueness of the dominant $\calX_q$-monomial of $F_q(\um)$.
Take an element $\chi \in \frakK_q(\sfg)$.
We enumerate $\calM_+(\chi)$ by $m_0,\, m_1,\, \dots,\, m_L$.
Let us write $\la_k \in \Z[q^{\pm \frac{1}{2}}]$ be the coefficients of $\um_k$ in $\chi$ for $k=0,1,\ldots,L$.
Then, the element $\chi - \sum_{k=0}^L \lambda_k F_q(\um_k) \in \frakK_q(\sfg)$ has no dominant $\calX_q$-monomial.
If it is non-zero, then it has at least one dominant $\calX_q$-monomial by Lemma \ref{lem:maximal or minimal monomials}, which yields a contradiction.
Hence, we conclude that the set $\{   F_q(\um) \ | \ m \in \calM^\sfg_+ \}$ generates $\frakK_q(\sfg)$.
\end{proof}

\begin{corollary} \label{cor: evaluation folded ver}
Let $\tm \in \calX_q$ be a dominant monomial. Then we have
\begin{equation*} 
	{\rm ev}_{q=1} (F_q(\tm)) = F({\rm ev}_{q=1}(\tm)).
\end{equation*}
\end{corollary}
\begin{proof}
It follows from \eqref{eq: commutative diagram for finite types} that ${\rm ev}_{q=1} (F_q(\tm)) \in \frakK(\sfg)$, where ${\rm ev}_{q=1} (F_q(\tm))$ has the unique dominant monomial $\ev_{q=1}(\tm) \in \calX$ by Theorem \ref{thm: F_q}~\ref{it:unique dominant monomial}. 
Thus our assertion is proved from Theorem \ref{thm: FHR}(2).
\end{proof}

As in the simply-laced cases \eqref{eq: translation by r}, we have
\begin{equation} \label{eq: shift of spectral parameters}
\sfT_r( F_q( \um^{(i)}[p,s] ) ) =  F_q( \um^{(i)}[p+r,s+r] ) \quad \text{ for any $r \in 2\Z$,}
\end{equation}
where $r \in 2\Z$ and $\sfT_{r}$ is the $\Z[q^{\pm \frac{1}{2}}]$-algebra automorphism of $\calX_q$ sending $\tX_{i,p}$ to $\tX_{i,p+r}$.

\begin{proposition} \label{prop: range of folded KR}
For $(i,p),(i,s) \in \tbDynkinf_0$ with $p<s$, the element $F_q( \um^{(i)}[p,s] )$ is of the form
\begin{equation} \label{eq:right negativeness of KR}
F_q( \um^{(i)}[p,s] )   = \um^{(i)}[p,s] * ( 1+ \tB^{-1}_{i,s+1} * \chi)
\end{equation}
where $\chi$ is a $($non-commutative$)$ $\Z[q^{\pm \frac{1}{2}}]$-polynomial in $\tB^{-1}_{j,k+1}$, $(j,k) \in \tbDynkinf_0$. In particular, we have
\begin{align} \label{eq: Fq m[ps] = Fq m-[ph sh]}
F_q( \um^{(i)}[p,s] ) = F_q( \um^{(i^*)}_-[p+\sfh,s+\sfh] )
\end{align}
and
\ben
\item \label{it: n KR1} $F_q( \um^{(i)}[p,s] )$ contains the unique dominant monomial $\um^{(i)}[p,s]$,
\item \label{it: n KR2} $F_q( \um^{(i)}[p,s] )$ contains the unique anti-dominant monomial $\um^{(i^*)}_-[p+\sfh,s+\sfh]$,
\item \label{it: right-negative of Fqmps} each $\calX_q$-monomial of $F_q(\um^{(i)}[p,s]) - \um^{(i)}[p,s]-\um^{(i^*)}_-[p+\sfh,s+\sfh] $ is a product of $\tX_{j,u}^{\pm1}$ with $p \le u \le s+\sfh$, having at least one of it factors from $p < u < s+\sfh$, 
and is right-negative. In particular, each monomial of $F_q(\underline{\tX_{\im,p}}) - \underline{\tX_{\im,p}}-\underline{\tX_{\im^*,p+\sfh}}$ is a product of $\tX_{j,u}^{\pm1}$ with $p < u < p+\sfh$. 
\ee
\end{proposition}
\begin{proof}
\ref{it: n KR1} follows from Theorem \ref{thm: F_q}~\ref{it:unique dominant monomial}.
~\ref{it: n KR2} and~\eqref{eq: Fq m[ps] = Fq m-[ph sh]} follow from the reversed version of the $q$-algorithm (see Remark~\ref{rem: reversed q-algorithm}) and~\ref{it: n KR1}.
Finally, ~\eqref{eq:right negativeness of KR} and ~\ref{it: right-negative of Fqmps}  are the direct consequences of Theorem \ref{thm:monomials for KR modules in types ADE} and Proposition \ref{prop:folding lemma}.
\end{proof}

%

\begin{conjecture} \label{conjecture 1}
For $(i,p),(i,s) \in \tbDynkinf_0$ with $p<s$,  every monomial in $F_q( \um^{(i)}[p,s] )$ has a quantum positive coefficient; that means, each coefficient of a monomial in $F_q( \um^{(i)}[p,s] )$ is contained in $\Z_{\ge0}[q^{\pm \frac{1}{2}}]$.
\end{conjecture}

\begin{remark}
In the proof of Corollary \ref{cor: folding monomials of Ft = monomials of Fq},
we have seen that the coefficients of monomials in $F(\um^{(i)}[p,s])$ are positive.
In Section \ref{sec: qCA on Kqxi}, we will provide a quantum cluster algebra theoretic algorithm for computing $F_q( \um^{(i)}[p,s] )$, which starts from
an initial quantum cluster variable $\um^{(i)}[p,s]$ (see Proposition~\ref{prop: m times mu} below).
This may be viewed as an evidence of Conjecture \ref{conjecture 1}, which is compatible with the quantum positivity conjecture of quantum cluster algebras (\cite{BZ05}). 
\end{remark}

By the following theorem, we have the third basis, denoted by 
\begin{equation*}
\sfL_q \seteq \{ L_q(\um) \ | \ m \in \calM_+^\sfg\},
\end{equation*}
and called the \emph{canonical basis} of $\frakK_q(\sfg)$.
We remark that the reason why we call it the canonical basis is further explained in \cite{JLO2}.

\begin{theorem} \label{thm: bar L_q(m)}
For $m \in \calM^\sfg_+$, there exists a unique element $L_q(\um)$ in $\frakK_q(\sfg)$ such that
\bna
\item \label{it: L p1} $\overline{L_q(\um)}=L_q(\um)$,
\item \label{it: unit2}
$E_q(\um) = L_q(\um) + \displaystyle\sum_{m' \lN m }  P_{m,m'}(q) L_q(\um')$  with $P_{m,m'}(q) \in q\Z[q]$.
\ee
\end{theorem}

\begin{proof}
For $\um \in \calM^\sfg_+$, we will construct $L_q(\um)$ inductively using some dominant $\calX_q$-monomials below, which are all less than $\um$ with respect to the Nakajima order $\leN$.
\smallskip

\noindent
{\it Step 1}. Let us first collect all dominant $\calX_q$-monomials obtained from $\um$ in an inductive way.
Let $\ucalM_1 \seteq \ucalM_+\bl E_q(\um) \br = \left\{ \um_{1,1},\, \um_{1,2},\, \dots ,\, \um_{1, \ell_1} = \um \right\}$. Then we define
\begin{equation*}
\begin{split}
	\ucalM_n &\seteq \displaystyle\bigcup_{1 \le k \le \ell_{n-1}} \ucalM_+\bl E_q(\um_{n-1,k}) \br,  
\end{split}
\end{equation*}
where $\ucalM_{n-1} = \left\{\, \um_{n-1,1},\, \um_{n-1,2},\, \dots,\, \um_{n-1,\ell_{n-1}} \,\right\}$ for $n \ge 2$.
Note that
\begin{equation*}
	\ucalM_+\bl E_q(\um) \br = \ucalM_1 \subset \ucalM_2 \subset \ucalM_3 \subset \cdots
\end{equation*}
The above chain has finite length, that is, there exists $N$ such that $\ucalM_n = \ucalM_{n+1}$ for $n \ge N$ 
because we can apply the same argument as in the proof of \cite[Lemma 3.13 and Lemma 3.14]{H04}.
For simplicity, let us relabel the dominant $\calX_q$-monomials in $\ucalM_N$
as follows:
\begin{align} \label{eq: dominant monomials}
\usfm_1 \boldsymbol{<}  \usfm_2 \boldsymbol{<} \dots \boldsymbol{<} \usfm_M = \um.
\end{align}
where $\boldsymbol{<}$ is also a total order compatible with $\leN$. In particular, $E_q(\usfm_1)$ has no dominant $\calX_q$-monomial other than $\usfm_1$ by construction.

\vskip 1mm

\noindent
{\it Step 2.}
We construct $L_q(\um)$ by inductive argument on \eqref{eq: dominant monomials} as follows.
Since $E_q(\usfm_1)$ has the unique dominant $\calX_q$-monomial $\usfm_1$ by  construction,   
we have $E_q(\usfm_1) = F_q(\usfm_1)$.
If we set $L_q(\usfm_1) = E_q(\usfm_1)$, then the initial step is done because $\overline{E_q(\usfm_1)} = \overline{F_q(\usfm_1)} = F_q(\usfm_1) = E_q(\usfm_1)$.
\smallskip

Suppose that $L_q(\usfm_k)$ is well-defined and uniquely determined for $1 \le k \le M-1$.
By the property \ref{it: unit2}, one can write
\begin{equation*}
L_q(\usfm_k) = E_q(\usfm_k) + \sum_{\usfm_l \lN \usfm_k} Q_{\usfm_l, \usfm_k}(q) E_q(\usfm_l).
\end{equation*}
By \eqref{eq: uni t2 sfg}, $L_q(\usfm_k)$ can be written as a linear combination of $F_q(\usfm_l)$ for $1\le l \le k$.
In particular, the coefficient of $F_q(\usfm_k)$ is $1$ due to the property \ref{it: L p1}.
Hence, the finiteness described in~\eqref{eq: dominant monomials} implies that
\begin{align}\label{eq: uni T}
\text{$F_q(\usfm_k)$ can be written as a linear combination of $L_q(\usfm_l)$ for $1\le l \le k$. }
\end{align}
By replacing $F_q(\um')$ in \eqref{eq: uni t2 sfg} with~\eqref{eq: uni T}, we have
\begin{equation} \label{eq: Equm = Fqum + lower sum}
\begin{split}
    E_q(\um) = F_q(\um) + \sum_{1 \le l \le M-1} \alpha_{l}(q) L_q(\usfm_l).
\end{split}
\end{equation}

Let us take $\beta_l(q) \in \Z[q^{\pm 1}]$ such that $\beta_l(q)$ is symmetric in $q$ and $q^{-1}$, and $\alpha_l(q) - \beta_l(q) \in q\Z[q]$ for all $1 \le l \le L-1$. This is possible by the following way.
Let us write $\alpha_l(q)$ by $\alpha_l^+(q) + \alpha_l^0(q) + \alpha_l^-(q)$, where $\alpha_l^\pm(q) \in q^{\pm 1} \Z[q^{\pm 1}]$ and $\alpha_l^0(q) \in \Z$. Then we define $\beta_l(q) = \beta_l^+(q) + \beta_l^0(q) +  \beta_l^-(q)$ by setting $\beta_l^+(q) = \alpha_l^-(q^{-1})$, $\beta_l^-(q)=\alpha_l^-(q)$ and $\beta_l^0(q) = \alpha_l^0(q)$.
Now, we define
\begin{equation*}
    L_q(\um) = F_q(\um) + \sum_{1 \le l \le M-1} \beta_{l}(q) L_q(\usfm_l) \in \frakK_q(\sfg).
\end{equation*}
Then, $L_q(\um)$ satisfies the properties \ref{it: L p1} and \ref{it: unit2} due to the our choice of $\beta_l(q)$, which is the desired element of $\frakK_q(\sfg)$.
Note that it follows from Proposition \ref{prop: bases Eq and Fq in Kqg} and \ref{it: unit2} that $\sfL_q$ is a $\Z[q^{\pm \frac{1}{2}}]$-basis of $\frakK_q(\sfg)$.
\vskip 1mm

\noindent
{\it Step 3}. Let us prove the uniqueness of $L_q(\um)$.
Assume that $L_q'(\um) \in \frakK_q(\sfg)$ satisfies \ref{it: L p1} and \ref{it: unit2}.
By \eqref{eq: dominant monomials} and \ref{it: unit2}, we have
\begin{equation*}
	L_q'(\usfm_1) = E_q(\usfm_1) = L_q(\usfm_1).
\end{equation*} 
By induction on \eqref{eq: dominant monomials},
we suppose that $L_q(\usfm_k) = L_q'(\usfm_k)$ for $1 \le k \le M-1$.
By \ref{it: unit2} and induction hypothesis, $E_q(\um)$ is written as 
\begin{equation*}
E_q(\um) = L_q(\um) + \displaystyle\sum_{1 \le k \le M-1}  P_{m,\sfm_k}(q) L_q(\usfm_k) = L_q'(\um) + \displaystyle\sum_{1 \le k \le M-1}  P_{m,\sfm_k}'(q) L_q(\usfm_k),
\end{equation*}
where $P_{m,\sfm_k}(q), P_{m,\sfm_k}'(q) \in q\Z[q]$.
Hence we have
\begin{equation} \label{eq: Lq = L'q + R}
\begin{split}
	L_q(\um) = L_q'(\um) + \sum_{1 \le k \le M-1} \left( P_{m,\sfm_k}'(q) - P_{m,\sfm_k}(q) \right) L_q(\usfm_k).
\end{split}
\end{equation}
By taking the bar involution on both sides of \eqref{eq: Lq = L'q + R},
it follows from \ref{it: L p1} that for $1 \le k \le M-1$,
\begin{equation*}
	\overline{P_{m,\sfm_k}'(q) - P_{m,\sfm_k}(q)} = P_{m,\sfm_k}(q) - P_{m,\sfm_k}'(q) \in q\Z[q] \cap q^{-1} \Z[q^{-1}] = \left\{ \,0\, \right\}.
\end{equation*}
This implies that $L_q'(\um) = L_q(\um)$ by \eqref{eq: Lq = L'q + R}.
\end{proof}

\begin{remark} \label{rem: new KL-type polynomials}
In the viewpoint of Kazhdan--Lusztig theory (explained briefly in Remark \ref{rem: KL theory}),
we regard the polynomials $P_{m,m'}(q)$'s as new KL-type polynomials, which generalize Nakajima's KL-type polynomials, since the $t$-quantized Cartan  matrices for types ADE are equal to the quantum Cartan matrices and the basis in Theorem \ref{thm: bar L_q(m)} essentially coincides with Nakajima's as explained in \cite{H04, H05}.
It would be very interesting to find a geometric or representation theoretic interpretation behind $P_{m,m'}(q)$ in the spirit of Kazhdan--Lusztig theory.
\end{remark}

\begin{remark} \label{rem: difference between Lt and Lq}
We emphasize that the basis $\sfL_q = \{ L_q(\um) \ | \ m \in \calM_+^\sfg\}$ of $\frakK_q(\sfg)$ is quite different from the $\bfL_\ttt$ of $\frakK_t(\bfg) \simeq \calK_\ttt(\scrC_\bfg^0)$, that is, $L_q(\um)$ cannot be obtained from $L_t(\ubfm)$ by folding $\calY_t$-monomials with some modification of coefficients in $\Z[t^{\pm \frac{1}{2}}]$, where $m = \osigma(\bfm)$. 
We give an example to illustrate this phenomenon.
Let us consider $L_t(\tY_{1,1})$ and $L_t(\tY_{4,-2})$ of the finite type $A_5$.
One may observe that $L_t(\tY_{1,1})$ $q$-commutes with $L_t(\tY_{4,-2})$, which implies that $L_t(\tY_{1,1}*\tY_{4,-2})$ coincides with $L_t(\tY_{1,1})*L_t(\tY_{4,-2})$ up to $q^\Z$ \cite[Corollary 5.5]{HL15}.
On the other hand, for type $C_3$, $L_q(\tX_{1,1})$ does not $q$-commute with $L_q(\tX_{2,-2})$.
This implies that $L_q(\tX_{1,1}*\tX_{2,-2})$ is not equal to $L_q(\tX_{1,1}) * L_q(\tX_{2,-2})$ up to $q^\Z$.
In fact, $L_q(\tX_{1,1}*\tX_{2,-2})$ has two dominant $\calX_q$-monomials, while $L_t(\tY_{1,1}*\tY_{4,-2})$ has only one dominant $\calY_t$-monomial.
\end{remark}

\begin{conjecture} \label{conjecture 2}
For $(i,p),(i,s) \in \tbDynkinf_0$ with $p<s$,  we have
\begin{equation*}
L_q( \um^{(i)}[p,s] ) = F_q( \um^{(i)}[p,s] ),
\end{equation*}
where $\um^{(i)}[p,s] \seteq \underline{m^{(i)}[p,s]}$ denotes the bar-invariant $\calX_q$-monomial corresponding to $m^{(i)}[p,s]$ \eqref{eq: dominant monomial of KR type} as in {\rm Remark \ref{rem: bar-invariant}}.
\end{conjecture}

\begin{example} \label{ex: X25 * X110 in G2}
Let us illustrate Theorem \ref{thm: bar L_q(m)} in the case of $L_q(\underline{X_{2,5} X_{1,10}})$ for type $G_2$. 

\noindent
{\it Step 1}. By \eqref{eq: Eqm}, we have
\begin{equation*}
	E_q(\underline{X_{2,5} X_{1,10}})
	=
	q^{\frac{3}{2}} F_q(\sfX_{2,5}) * F_q(\sfX_{1,10}).
\end{equation*}
Let us recall the formulas of $F_q(\sfX_{2,5})$ and $F_q(\sfX_{1,10})$ in Examples \ref{ex: F25 in G2} and \ref{ex: F110 in G2}, respectively.
Then we observe that there exist two bar-invariant dominant $\calX_q$-monomials with $\Z[q^{\pm\frac{1}{2}}]$-coefficients in $E_q(\underline{X_{2,5} X_{1,10}})$, namely, $\underline{X_{2,5} X_{1,10}}$ and $\left( q^{-1} + q + q^{3} \right) \underline{X_{1,6}}$.
\smallskip

\noindent
{\it Step 2}. By {\it Step 1}, we have
\begin{equation*}
	E_q(\underline{X_{2,5} X_{1,10}})
	= F_q(\underline{X_{2,5} X_{1,10}}) + \left( q^{-1} + q + q^{3} \right) F_q(\underline{X_{1,6}}),
\end{equation*}
which corresponds to \eqref{eq: Equm = Fqum + lower sum} in this case, that is, $M = 2$ and $\alpha_1(q) = q^{-1} + q + q^3$.
Set $\beta_1(q) = q + q^{-1}$ by construction in the proof of Theorem \ref{thm: bar L_q(m)}. Then we have
\begin{equation*}
	L_q(\underline{X_{2,5} X_{1,10}})
	=
	F_q(\underline{X_{2,5} X_{1,10}}) + \left( q^{-1} + q \right) F_q(\underline{X_{1,6}}),
\end{equation*}
which is bar-invariant.
Note that $L_q(\underline{X_{2,5} X_{1,10}})$ has two dominant $\calX_q$-monomials $\underline{X_{2,5} X_{1,10}}$ and $\underline{X_{1,6}}$.
Moreover, we verify
\begin{equation*}
	E_q(\underline{X_{2,5} X_{1,10}})
	=
	L_q(\underline{X_{2,5} X_{1,10}}) + P_{X_{2,5}X_{1,10}, X_{1,6}}(q) L_q(\underline{X_{1,6}}),
\end{equation*}
where $P_{X_{2,5}X_{1,10}, X_{1,6}}(q) = q^3 \in q\Z_{\ge 0}[q]$, that is, $L_q(\underline{X_{2,5} X_{1,10}})$ is the unique element in $\frakK_q(\sfg)$ satisfying the properties \ref{it: L p1} and \ref{it: unit2} in Theorem \ref{thm: bar L_q(m)}.
\end{example}

\subsection{Proof of Proposition \ref{prop:folding lemma}} \label{subsec:proof of folding lemma}
To prove Proposition \ref{prop:folding lemma},
we utilize some analogues of the results in \cite{H04}, where we will skip some proof of them when they can be obtained from the corresponding arguments in \cite{H04}.
%

For $J^\bfg \subset I^\bfg$, we set $\overline{J^\bfg} = \{\, \overline{\im}  \, | \,\, \im \in J^\bfg \,\} \subset I^\sfg$.
Let $J \subset I$ be given such that $J = \overline{J^\bfg}$ for some $J^\bfg \subset I^\bfg$.
Let us define $\frakK_J(\sfg) \subset \calX$ as follows:
\begin{align*}
 \frakK_J(\sfg) = \bigcap_{j \in J} \left(\Z[X^{\pm 1}_{k,l} \ |  \  (k,l) \in \tDynkinf^\sfg_0,  j \ne k \in I ]
\otimes \Z[X_{j,l}(1+B_{j,l+1}^{-1})  \ |  \  (j,l)  \in \tDynkinf^\sfg_0] \right).
\end{align*}
Note that $\frakK_I(\sfg) = \frakK(\sfg)$.
We also define $\frakK_{J,q}(\sfg) \subset \calX_q$ as above by replacing the letters $X$ and $B$ with $\widetilde{X}$ and $\widetilde{B}$, respectively.

\begin{proposition} \label{prop: Fiq}
Let $J \subset I$ with $|J| \le 2$.
For a $J$-dominant monomial $m$, there exists a unique $F_{J,q}(\um) \in \frakK_{J,q}(\sfg)$ such that $\um$ is the unique $J$-dominant $\calX_q$-monomial of $F_{J,q}(\um)$.
Moreover, $$\text{$\{ F_{J,q}(\um) \, | \, \text{$m$ is $J$-dominant} \, \}$ is a $\Z[q^{\pm \frac{1}{2}}]$-basis of $\frakK_{J,q}(\sfg)$.}$$
\end{proposition}

For $m \in \calM_+^J $, we define
\begin{equation} \label{eq: EJm}
\begin{split}
	E_J(m) = \prod_{j \in J; (j,p)\in\sDynkinf_0} F_J(X_{j,p})^{u_{j,p}(m)} \in \frakK_J(\sfg),
\end{split}
\end{equation}
where $F_J(X_{j,p}) := {\rm ev}_{q=1}(F_{J,q}(\sfX_{i,p}))$ is a unique element in $\frakK_J(\sfg)$ such that $X_{j,p}$ is the unique dominant monomial of $F_J(X_{j,p})$ (cf.~Remark \ref{rem: evaluation of Kqg} and Remark \ref{rem: F vs Ft in ADE}).  
Let $\frakK_{i,q}^\infty(\sfg)$ be the completion of $\frakK_{i,q}(\sfg)$  
given by the method in \cite[Section 5.2.2]{H04}.
Put $\frakK_{J,q}^\infty(\sfg) = \bigcap_{j \in J} \frakK_{j,q}^\infty(\sfg)$.

\begin{lemma} \label{lem: characterization of KJq}
\hfill
\ben
\item \label{it: at least one} A non-zero element of $\frakK_{J,q}^\infty$ has at least one $J$-dominant $\calX_q$-monomial.
\item \label{it: intersection inf}
We have
\begin{align*}
\frakK_{J,q}(\sfg) = \frakK_{J,q}^\infty(\sfg) \,\scalebox{0.9}{$\bigcap$}\, \calX_q.
\end{align*}
\ee
\end{lemma}

For $i \in I^\sfg$, take $\im \in I^\bfg$ such that $\overline{\im} = i$ and put
\begin{itemize}
	\item $D_{\bfm^{(\im)}[p,s]}^\bfg = (\bfm^{(k)})_{k \ge 0}$ : the countable set as in \cite[Section 5.2.3]{H04} associated with $\bfm^{(\im)}[p,s]$,
	\item $D_{m^{(i)}[p,s]}^\sfg = (m^{(k)})_{k \ge 0}$ : the analogue of the above one for $m^{(i)}[p,s]$ in terms of \eqref{eq: EJm}.
\end{itemize}

\begin{remark}
The set $D_{\bfm^{(\im)}[p,s]}^\bfg$ may be an infinitely countable set.
If we enumerate the monomials in the countable set as follows:
\begin{align*}
\dots < \bfm^{(2)} < \bfm^{(1)} < \bfm^{(0)} = \bfm^{(\im)}[p,s].
\end{align*}
Then the $t$-algorithm determines $\Z[t^{\pm \frac{1}{2}}]$-coefficients of the monomials $\ubfm^{(k)}$'s.
Let $(\mathsf{c}^\bfg(\bfm^{(r)}))_{r \ge 0}$ be the sequence of $\Z[t^{\pm \frac{1}{2}}]$-coefficients for $\ubfm^{(r)}$'s determined by the $t$-algorithm starting from $\ubfm^{(\im)}[p,s]$.
It was known in \cite{H05} that the sequence $(\mathsf{c}^\bfg(\bfm_k))_{k \ge 0}$ should have finitely many non-zero coefficients, that is,
$F_t(\ubfm^{(\im)}[p,s]) \in \frakK_t(\bfg)$.
Note that $\ucalM(F_t(\ubfm^{(\im)}[p,s])) \subset \{\, \ubfm^{(k)} \, \mid \, k \ge 0 \,\}$.
\end{remark}
\smallskip

Let us enumerate the finite set $\calM(F_t(\ubfm^{(\im)}[p,s]))$ as follows:
\begin{align*}
\bfm_N < \dots < \bfm_2 < \bfm_1 < \bfm_0 = \bfm^{(\im)}[p,s],
\end{align*}
where $<$ is a total order compatible with $\lN$. In particular, $\bfm_N$ is an anti-dominant $\calY$-monomial, i.e. $\bfm_N = \bfm_-^{(\im^*)}[p+\sfh,s+\sfh]$ by Theorem \ref{thm:monomials for KR modules in types ADE}.
It follows from Corollary \ref{cor: some F to F} and Theorem \ref{thm:monomials for KR modules in types ADE} that
\begin{equation*}
\mathsf{M} := \left\{ \, \osigma(\bfm_k) \, | \, 1 \le k \le N \, \right\} \subset D_{m^{(i)}[p,s]}^\sfg.
\end{equation*}
 Then we enumerate the $\calX$-monomials in $\mathsf{M}$ by
\begin{align} \label{eq: enumeration of folded monomials from FYip}
	m_-^{(i)}[p+\sfh,s+\sfh] = \mathsf{m}_{N'} <' \dots <' \mathsf{m}_1 <' \mathsf{m}_0 = m^{(i)}[p,s],
\end{align}
where $<'$ is a total order compatible with $\lN$.

\begin{definition} \label{def: restricted q-algorithm}
Set $\widetilde{\mathsf{c}}^\sfg (m^{(i)}[p,s]) = 1$ and $\widetilde{\mathsf{c}}^\sfg_J(m^{(i)}[p,s]) = 0$.
For $J \subset I$ with $|J| \le 2$ and $\mathsf{m} \in \mathsf{M}$ such that $\mathsf{m} \neq m^{(i)}[p,s]$, we define
\begin{equation*} 
\begin{split}
	\widetilde{\mathsf{c}}_J^{\sfg}(\mathsf{m}) &= \sum_{\substack{\mathsf{m} \in \mathsf{M} \\ \mathsf{m}\, <' \, \mathsf{m}'}} \left( \widetilde{\mathsf{c}}^{\sfg}(\mathsf{m}') - \widetilde{\mathsf{c}}_J^{\sfg}(\mathsf{m}') \right) \left[ F_{J,q}(\underline{\mathsf{m}'}) \right]_{\underline{\mathsf{m}}}, \\
	\widetilde{\mathsf{c}}^{\sfg}(\mathsf{m}) &=
	\begin{cases}
		\widetilde{\mathsf{c}}_J(\mathsf{m}) & \text{if $\mathsf{m}$ is not $J$-dominant,} \\
		0 & \text{if $\mathsf{m}$ is dominant,}
	\end{cases}
\end{split}
\end{equation*}
where $\left[ F_{J,q}(\underline{\mathsf{m}'}) \right]_{\underline{\mathsf{m}}}$ is a $\Z[q^{\pm \frac{1}{2}}]$-coefficient of $\underline{\mathsf{m}}$ in $F_{J,q}(\underline{\mathsf{m}}')$.
Here $F_{J,q}(\underline{\mathsf{m}}')$ is assumed to be $0$ when $\mathsf{m}'$ is not $J$-dominant.
\end{definition}

\noindent
Since the proof of the following lemma is similar to \cite{H03}, so we omit it, but the complete proof can be found in \cite{JLO1}.

\begin{lemma} \label{lem: restricted sequences}
The sequences $(\,\widetilde{\mathsf{c}}^\sfg_J(\mathsf{m})\,)_{\mathsf{m} \in \mathsf{M}}$ and $(\,\widetilde{\mathsf{c}}^\sfg(\mathsf{m})\,)_{\mathsf{m} \in \mathsf{M}}$ are well-defined, and $(\,\widetilde{\mathsf{c}}^\sfg(\mathsf{m})\,)_{\mathsf{m} \in \mathsf{M}}$ is not depend on the choice of $J$ with $|J| \le 2$.
\end{lemma}

By Proposition \ref{prop: Fiq} and Lemma \ref{lem: restricted sequences}, we set
$\chi \seteq \sum_{\mathsf{m} \in \mathsf{M}} \widetilde{\mathsf{c}}^\sfg(\mathsf{m}) \underline{\mathsf{m}} \in \calX_q$, and 
\begin{equation*}
    \chi_i \seteq \sum_{\mathsf{m} \in \mathsf{M}}\mu_i(\mathsf{m}) F_{i,q}(\underline{\mathsf{m}}) \in \frakK_{i,q}(\sfg),
\end{equation*}
where $\mu_i(\mathsf{m}) = \widetilde{\mathsf{c}}^\sfg(\mathsf{m}) - \widetilde{\mathsf{c}}^\sfg_i(\mathsf{m})$.
%
%
Now, we are ready to prove Proposition \ref{prop:folding lemma}.

\begin{proof} [Proof of Proposition \ref{prop:folding lemma}.]
Let us compute the coefficient of $\underline{\mathsf{m}}'$ in $\chi - \chi_i$ for $\mathsf{m}' \in \mathsf{M}$.
\smallskip

\noindent
{\it Case 1}. $\mathsf{m}'$ is not $i$-dominant. By definition of $\widetilde{\mathsf{c}}^\sfg(\mathsf{m}')$,
we have
\begin{align*}
\text{(coefficient of $\underline{\mathsf{m}'}$ in $\chi - \chi_i$)}	&= \widetilde{\mathsf{c}}^\sfg(\mathsf{m}') - \sum_{\substack{\mathsf{m} \in \mathsf{M} \\ \mathsf{m}'\, \le'\, \mathsf{m}}} \mu_i(\mathsf{m}) \left[ F_{i,q}(\underline{\mathsf{m}}) \right]_{\underline{\mathsf{m}'}} \\
&= (\widetilde{\mathsf{c}}^\sfg(\mathsf{m}') - \widetilde{\mathsf{c}}_i^\sfg(\mathsf{m}')) \left[ F_{i,q}(\underline{\mathsf{m}'}) \right]_{\underline{\mathsf{m}'}} = 0,
\end{align*}
where $F_{i,q}(\underline{\mathsf{m}'}) = 0$ since $\mathsf{m}'$ is not $i$-dominant.
\smallskip

\noindent
{\it Case 2}. $\mathsf{m}'$ is $i$-dominant. By uniqueness of $i$-dominant $\calX_q$-monomial for $F_{i,q}(\underline{\mathsf{m}})$ with $\mathsf{m}' \le' \mathsf{m}$, we have $\widetilde{\mathsf{c}}_i^\sfg(\mathsf{m}') = 0$, and the coefficient of $\underline{\mathsf{m}'}$ in $\chi_i$ is $\mu_i(\mathsf{m}') = \widetilde{\mathsf{c}}^\sfg(\mathsf{m}') - \widetilde{\mathsf{c}}_i^\sfg(\mathsf{m}') = \widetilde{\mathsf{c}}^\sfg(\mathsf{m}')$.
This implies that the coefficient of $\underline{\mathsf{m}}'$ in $\chi - \chi_i$ is $0$ in this case.
\smallskip

By {\it Case 1} and {\it Case 2}, we have $\chi = \chi_i \in \frakK_{i,q}(\sfg)$ and then $\chi \in \frakK_q(\sfg)$.
Note that $\chi$ has unique dominant $\calX_q$-monomial $\um^{(i)}[p,s]$ by Definition \ref{def: restricted q-algorithm} (or our choice of $\mathsf{M}$).
Since $F_q(\um^{(i)}[p,s]) - \chi \in \frakK_q^\infty(\sfg)$ has no dominant $\calX_q$-monomial,
we conclude $F_q(\um^{(i)}[p,s]) = \chi \in \frakK_q(\sfg)$ by Lemma \ref{lem: characterization of KJq}.
\end{proof}

\section{Subrings of \texorpdfstring{$\frakK_q(\g)$ } aand  the quantum folded \texorpdfstring{$T$}--systems} \label{sec:subrings and T-systems}

\noindent
 In this section, we prove the quantum folded  $T$-systems, which play a crucial role in this paper.
To do this, we consider a subring $\frakK_{q,\xi}(\g)$ of $\frakK_q(\g)$ for a height function $\xi$. 
We mainly employ the framework in~\cite{HL15,HL16} (see also~\cite{B21}).

\subsection{Subring} Let $\sfS$ be a convex set of $\tDynkinf_0$ (recall Definition \ref{def: convex subset}\,\ref{def: convexity}).
We denote by ${}^{\sfS}\calX$ the subring of $\calX$ generated by $X_{i,p}^{\pm 1}$ for $(i, p) \in \sfS$.
Let $\lScalMp$ be the set all dominant monomials in the variables $X_{i,p}$'s for $(i,p) \in \sfS$.
We define the $\Zq$-module $\frakK_{q,\sfS}(\g)$ as the $\Zq$-submodule of $\frakK_q(\g)$ given by
\begin{align} \label{eq: KqsfS}
\frakK_{q,\sfS}(\g)  \seteq \soplus_{ m \in  \lScalMp } \Zq F_q(\um).
\end{align}

\begin{lemma}      [{cf.\ \cite[Lemma 5.6]{FHOO}}]  \label{lem: convexity}
The set $\lScalMp$  is an ideal of the partially ordered set $(\calM_+,\leN)$; i.e., it is closed under taking smaller elements in $\calM_+$ with respect to $\leN$.
\end{lemma}

\begin{proof}
Let $m \in \lScalMp$ and $m M \in \calM_+$ where $M \in \Bq^{-k}$ for some $k \in \Z_{\ge 1}$.
For a  factor $B_{i,p}^{-1}$ of $M$,
the monomial $m$ should have factors $X_{i, p-1}$ and $X_{i, p+1}$ due to \eqref{eq: Bip}.
Thus we have an oriented path from $(i,p+1)$ to a vertex in $\sfS$ and another oriented path from a vertex in $\sfS$ to $(i,p-1)$ (these paths are possibly of length zero)
in $\hDynkinf_0$.
Hence we have an oriented path whose end points are in $\sfS$
factoring through both $(i,p-1)$ and $(i,p+1)$.
By convexity of $\sfS$ and the definition of $B_{i,p}$ \eqref{eq: Bip}, $M \in {}^\sfS\calX$ and $m M \in \lScalMp$ as we desired.
\end{proof}

\begin{proposition} \label{prop: KsfS is subalg}
For a convex subset $\sfS$ in $\tDynkinf_0$, the $\Zq$-module $\frakK_{q,\sfS}(\g)$ is a $\Zq$-subalgebra of $\frakK_q(\g)$. Moreover, we have
\begin{align}\label{eq: other presentations}
\frakK_{q,\sfS}(\g) =  \soplus_{ m \in  \lScalMp } \Zq E_q(\um) = \soplus_{ m \in  \lScalMp } \Zq L_q(\um).
\end{align}
\end{proposition}

\begin{proof}
Let $m_1,m_2 \in \lScalMp$.
By Theorem~\ref{thm: F_q} and Proposition~\ref{prop: range of folded KR}, $F_q(\um_1)*  F_q(\um_2) \in \frakK_q(\g)$ is written as shown below.
\begin{align} \label{eq: linear expansion of Fqm1 Fqm2}
F_q(\um_1) *  F_q(\um_2) =
\sum_{\substack{m \in \calM_+ \\ m \,\leN\, m_1 m_2}} c_{\um} F_q(\um),
\end{align}
where $c_{\um} \in \Z[q^{\pm \frac{1}{2}}] \setminus \{\, 0 \,\}$.
Then it follows from Lemma~\ref{lem: convexity} that $m \in \lScalMp$ for a monomial $m\leN m_1 m_2$ above.
Hence, we conclude that $\frakK_{q,\sfS}(\g)$ is a $\Zq$-subalgebra of $\frakK_q(\g)$ by definition \eqref{eq: KqsfS} of $\frakK_{q,\sfS}(\g)$.

Since $\frakK_{q,\sfS}(\g)$ is given by \eqref{eq: KqsfS}, ~\eqref{eq: other presentations} follows from   $\frakK_q(\g)$-analogue  of~\eqref{eq: uni t1} and ~\ref{it: unit2} in Theorem \ref{thm: bar L_q(m)}.
\end{proof}

\subsection{Truncation} Let $\xi$ be a height function of $\Dynkinf$.
For a (non-commutative) Laurent polynomial $x\in \calX_q$, we denote
by $x_{\le \xi}$ the element of $\lxi\calX_q$ obtained from $x$ by
discarding all the monomials containing $\tX_{i,p}^{\pm 1}$ with
$(i,p) \in \tDynkinf_0 \setminus \lxi\tDynkinf_0$.

The map
\begin{align*}
( \cdot)_{\le \xi} :  \calX_q \longrightarrow \lxi\calX_q \quad \text{ given by } \quad x \longmapsto x_{\le \xi}
\end{align*}
is a $\Zq$-linear map, which is not $\Zq$-algebra homomorphism.
For $m \in \calM_+$, we denote by $F_q(\um)_{\le \xi}$ the image of $F_q(\um)$ under the map $( \cdot)_{\le \xi}$.

Let us recall Definition \ref{def: convex subset} 
and~\eqref{eq: KqsfS}. We set
\begin{align} \label{eq: Kqxi}
\frakK_{q,\xi}(\g) \seteq  \frakK_{q,\lxi\widetilde{\sDynkinf}_0}(\g).
\end{align}

\begin{proposition} \label{prop: inj}
For a height function $\xi$ on $\Dynkinf$, the map $( \cdot)_{\le \xi}$ restricts to the injective $\Zq$-algebra homomorphism
\begin{align*}
( \cdot)_{\le \xi} \,:\,
 \frakK_{q,\xi}(\g)    \hookrightarrow   \lxi\calX_q.
\end{align*}
\end{proposition}

\begin{proof}
The injectivity follows from Theorem \ref{thm: F_q}.
Let us take $m_1,m_2 \in \lxi\calM_+$.
We consider a linear expansion of $F_q(\um_1) * F_q(\um_2)$ as in \eqref{eq: linear expansion of Fqm1 Fqm2}.
Then we claim that
\begin{align} \label{eq: truncated linear exp}
F_q(\um_1)_{\le \xi}  * F_q(\um_2)_{\le \xi} =
\sum_{\substack{m \in \lScalMp \\ m \,\leN\, m_1 m_2}} c_{\um} F_q(\um)_{\le \xi}
\quad (c_{\um} \neq 0).
\end{align}

Take a $\calX_q$-monomial $\tm'$ (resp. $\tm''$) appearing in $F_q(\um_1)_{\le \xi}$ (resp. $F_q(\um_2)_{\le \xi}$).
If $\ev_{q=1}(\tm' \tm'') \in \calM_+$, then $\ev_{q=1}(\tm' \tm')' \in {}^{\sfS}\calM_+$ by Lemma~\ref{lem: convexity}.
Furthermore, by Theorem~\ref{thm: F_q} and definition of ${}^{\xi}\calX_q$, $F_q(\um_1)_{\le \xi} * F_q(\um_2)_{\le \xi}$ is written as a linear combination of $\{   F_q(\um)_{\le \xi} \ | \ m \in {}^{\sfS}\calM_+ \}$.
Thus, $F_q(\tm' \tm'')_{\le \xi}$ appears in the right-hand side of \eqref{eq: truncated linear exp} up to $\Z[q^{\pm \frac{1}{2}}]$. This proves the above claim.

Finally, we have
\begin{align*}
(\cdot)_{\le \xi}(F_q(\um_1) *  F_q(\um_2)) = \sum_{\substack{m \in {}^{\sfS}\calM_+ \\ m \,\leN\, m_1 m_2}} c_{\um} F_q(\um)_{\le \xi} = F_q(\um_1)_{\le \xi} *  F_q(\um_2)_{\le \xi}.
\end{align*}
by Proposition \ref{prop: KsfS is subalg} and \eqref{eq: truncated linear exp}, which completes the proof.
\end{proof}

\begin{definition} 
For $m \in \calM_+$, we say $L_q(m)$ (resp.~$F_q(m)$) \emph{real} if, for any $k \in \Z_{\ge1}$, we have
$(L_q(m))^k=q^t L_q(m^k)$ (resp.~$(F_t(m))^k=q^t F_q(m^k)$) for some $t \in \Z$. 
\end{definition}

\begin{corollary}
For each KR-monomial $\um^{(i)}[p,s]$, $F_q(\um^{(i)}[p,s])$ is real. 
\end{corollary}

\begin{proof}
Let $\xi$ be a height function with $\xi_i=s$. Then we have
\begin{align*}
(F_q(\um^{(i)}[p,s]))_{\le \xi} = \um^{(i)}[p,s],
\end{align*}
by~\eqref{eq:right negativeness of KR} in Proposition~\ref{prop: range of folded KR}.
Since
\begin{align*}
\ev_{q=1}\left( \bl F_q(\um^{(i)}[p,s])^{* n} \br_{\le \xi} \right)=  (m^{(i)}[p,s])^n  =  \ev_{q=1}\left( \bl F_q(\um^{(i)}[p,s]^{* n}) \br_{\le \xi} \right),
\end{align*}
our assertion follows from Proposition~\ref{prop: inj}.
\end{proof}

\begin{conjecture} 
For $m \in \calM_+$, if   $L_q(m)$ is real, then 
$L_q(m)$ has a quantum positive coefficient.  
\end{conjecture}

\subsection{Quantum folded $T$-system}
For $f,g \in \calX_q$, we say that $f$ and $g$ \emph{$q$-commute or are $q$-commutative} if $fg = q^k gf$ for some $k \in \frac{1}{2}\Z$.
In this subsection, we shall prove the functional equations among KR-polynomials $F_q(m^{(i)}[p,s])$'s, called the \emph{quantum folded $T$-system}.
For simply-laced finite type, the quantum folded $T$-system is nothing but the quantum $T$-system, investigated in \cite{HL15} (see also \cite{HO19,FHOO}).

\begin{lemma} \label{lem: m(i;p,s)}
For $(i,p),(i,s) \in \tDynkin_0$ with $p<s$, let $j , j' \in \Dynkinf_0$ such that $d(i,j)=d(i,j')=1$. Then we have
$$
\FMq{ \um^{(j)}(p,s) } *  \FMq{ \um^{(j')}(p,s) } = \FMq{ \um^{(j')}(p,s) } *  \FMq{ \um^{(j)}(p,s) }.
$$
\end{lemma}

\begin{proof}
Note that one can take a height function on $\Dynkinf$ such that $\xi_{j}=\xi_{j'}=s-1$ and $\xi_j = \max\{\,\xi_i \, | \, i \in I\,\}$.
By~\eqref{eq:right negativeness of KR},
\begin{align*}
\FMq{ \um^{(j)}(p,s) }_{\le \xi} = \um^{(j)}(p,s) \quad \text{ and } \quad
\FMq{ \um^{(j')}(p,s) }_{\le \xi}= \um^{(j')}(p,s).
\end{align*}
By Proposition~\ref{prop: inj}, we have
\begin{align*}
\FMq{ \um^{(j)}(p,s) } *  \FMq{ \um^{(j')}(p,s) } = q^\be \FMq{ \um^{(j')}(p,s) } *  \FMq{ \um^{(j)}(p,s) }
\end{align*}
for some $\be \in \frac{1}{2}\Z$.

Now, let us prove that $\be =0$ by induction on $k=(p-s)/2$.
When $k=1$, we have $\um^{(j)}(p,s)=\tX_{j,p+1}$. In this case, $\be =0$ by \eqref{eq: commute}.
Suppose that $k > 1$.
By the induction hypothesis, we have
\begin{align*}
\um^{(j)}(p,s-2)  *  \um^{(j')}(p,s-2)  =  \um^{(j')}(p,s-2)  *  \um^{(j)}(p,s-2).
\end{align*}
Then we have
\begin{align*}
\um^{(j)}(p,s) * \um^{(j')}(p,s)  =   q^{\ucalN(\tX_{j,s-1},m^{(j')}(p,s)) + \ucalN(m^{(j)}(p,s),\tX_{j',s-1}) )  } \um^{(j)}(p,s)  * \um^{(j')}(p,s).
\end{align*}
Since
\begin{align*}
\ucalN(\tX_{j,s-1},m^{(j')}(p,s)) & = \sum_{i=0}^{(p-s)/2-1} \ucalN(j,s-1;j',p+1+2i) \allowdisplaybreaks\\
 &  = \sum_{i=0}^{(p-s)/2-1}\tfb_{j,j'}(s-p-2i-3)-\tfb_{j,j'}(s-p-2i-1), \allowdisplaybreaks\\
\ucalN(m^{(j)}(p,s),\tX_{j',s-1}) & = -\ucalN(\tX_{j',s-1},m^{(j)}(p,s))  = -\sum_{i=0}^{(p-s)/2-1} \ucalN(j',s-1;j,p+1+2i) \allowdisplaybreaks\\
&= \sum_{i=0}^{(p-s)/2-1}-\tfb_{j',j}(s-p-2i-3)+\tfb_{j',j}(s-p-2i-1),
\end{align*}
our assertion follows from the fact that $\tfb_{j,j'}(u)=\tfb_{j',j}(u)$ for all $u \in \Z$ (cf. \cite[Section 4]{KO23}).
\end{proof}

\begin{lemma}  \label{lem: m[p,s] m(p,s)}
For $(i,p),(i,s) \in \tDynkinf_0$ with $p<s$, we have
$$
\FMq{\um^{(i)}[p,s]}   * \FMq{\um^{(i)}(p,s)} = \FMq{\um^{(i)}(p,s)} *   \FMq{\um^{(i)}[p,s]}.
$$
\end{lemma}

\begin{proof}
Let us first show that $F_q(\um^{(i)}[p,s])$ and $F_q(\um^{(i)}(p,s))$ are $q$-commutative.
Since $[p,s] \neq (p,s)$, we cannot apply the same argument as in the proof of Lemma \ref{lem: m(i;p,s)}.
Instead, to show their $q$-commutativity,
we shall apply the $\mathfrak{sl}_2$-reduction argument as in \cite[Remark 9.10]{HO19} (see also \cite[Proposition 6.10]{FHOO}).
We should remark that the $\mathfrak{sl}_2$-reduction argument in \cite[Remark 9.10]{HO19} is based on \cite[Proposition 5.3, Lemma 5.6]{H06} (cf.~\cite[Lemma 9.9]{HO19}).
By~\eqref{eq: F to F KR},  one can prove the $\frakK(\sfg)$-analogues of \cite[Proposition 5.3, Lemma 5.6]{H06}.
Although we do not know the positivity of $\FMq{\um^{(i)}[p,s]}$ yet, it is enough to know $\calX_q$-monomials (except for their coefficients in $\Z[q^{\pm \frac{1}{2}}]$) appearing in $F_q(\um)$ for applying the $\mathfrak{sl}_2$-reduction argument to our first claim. This is done by Proposition \ref{prop:folding lemma} (see also Section \ref{subsec:proof of folding lemma}).

Now, let us prove our first claim.
The multiplicities of bar-invariant dominant monomials in the polynomials $F_q(\um^{(i)}[p,s]) * F_q(\um^{(i)}(p,s))$ and $F_q(\um^{(i)}(p,s)) *  F_q(\um^{(i)}[p,s])$
are the same as those in the corresponding polynomials for the $\mathfrak{sl}_2$-case up to overall power of $q^{1/2}$.
This follows from the fact that the $q$-commutation relations between $\tB_{i,s}^{-1}$, $\tB_{i,s'}^{-1}$
and between $\tX_{i,s}$, $\tB_{i,s'}^{-1}$ are the same as in the $\mathfrak{sl}_2$-case by Proposition~\ref{prop: YA com}. Thus, as the corresponding (non-commutative) polynomials in the $\mathfrak{sl}_2$-case commute up to a power of $q$, we obtain the $q$-commutativity.

Finally, to complete our assertion, it suffices to show that
\begin{align*}
\um^{(i)}[p,s] *  \um^{(i)}(p,s) = \um^{(i)}(p,s) *  \um^{(i)}[p,s].
\end{align*}
By an induction on $(p-s)/2$, we have
\begin{align*}
\ucalN(\um^{(i)}[p,s], \um^{(i)}(p,s))  & =  \ucalN(\tX_{i,s}, \um^{(i)}(p,s-2]) +  \ucalN( \um^{(i)}[p,s-2),\tX_{i,s-2})  \\
& \overset{\star}{=}\ucalN(\tX_{i,s}, \um^{(i)}(p,s]) +  \ucalN( \um^{(i)}[p,s),\tX_{i,s-2})  \\
&  \overset{\dagger}{=} \ucalN(\tX_{i,s}, \um^{(i)}(p,s]) +  \ucalN( \um^{(i)}(p,s],\tX_{i,s}) =0,
\end{align*}
where $\overset{\star} {=}$ follows from $\ucalN(\tX_{i,t},\tX_{i,t})=0$
and $\overset{\dagger}{=}$ follows from  $\ucalN(\tX_{i,t},\tX_{i,t'})=\ucalN(\tX_{i,t \pm 2},\tX_{i,t' \pm 2})$.
\end{proof}

For $(i,p),(i,s) \in \tDynkin_0$ with $p<s$, we set
$m(i; p,s) \seteq \prod_{ j; \; d(i,j)=1} m^{(j)}(p,s)^{-\sfc_{j,i}}$,
where $m^{(j)}(p,s)$ is given as in \eqref{eq: dominant monomial of KR type}.

\begin{lemma} \label{lem: last term in T-system}
For $(i,p),(i,s) \in \tDynkin_0$ with $p<s$,  we have
\begin{align*}
F_q( \um(i;p,s) ) =   \prod_{j; \; d(i,j)=1}   F_q(\um^{(j)}(p,s))^{-\sfc_{j,i}},
\end{align*}
where the order of the product does not matter.
\end{lemma}

\begin{proof} By Lemma~\ref{lem: m(i;p,s)}, $\prod_{j; \; d(i,j)=1}   F_q(\um^{(j)}(p,s))^{-\sfc_{j,i}}$ is well-defined.
Let $\xi$ be a height function on $\Dynkinf$ such that $\xi_i =s$ and $\xi_{j}=s-1$ for $j \in \Dynkinf_0$ with $d(i,j)=1$. Then we have
\begin{align*}
 \left( \prod_{j; \; d(i,j)=1}   F_q(\um^{(j)}(p,s))^{-\sfc_{j,i}} \right)_{\le \xi}  =  \um(i;p,s),
\end{align*}
which implies the assertion.
\end{proof}

Now, we are in a position to state and prove the quantum folded $T$-system (cf.~Theorem \ref{thm: Quantum T-system}).

\begin{theorem} [Quantum folded $T$-system] \label{thm: quantum folded}   
For $(i,p),(i,s) \in \tDynkinf_0$ with $p<s$ and $k = (s-p)/2 \in \Z_{\ge 1}$, we have
\begin{align*}
\FMq{ \um^{(i)}[p,s) } *  \FMq{\um^{(i)}(p,s]} = q^{\al(i,k)}  \FMq{\um^{(i)}(p,s)}  * \FMq{\um^{(i)}[p,s]} + q^{\ga(i,k)}  \hspace{-1ex} \prod_{j; \; d(i,j)=1}   \FMq{\um^{(j)}(p,s)}^{-\sfc_{j,i}},
\end{align*}
where $\ga(i,k) = \dfrac{1}{2}\left( \tfb_{i,i}(2k-1)+\tfb_{i,i}(2k+1) \right)$ and
$\al(i,k)=\ga(i,k)-d_i$.
\end{theorem}

\begin{proof}
First, we claim that
\begin{align*}
\FMq{ \um^{(i)}[p,s) } * \FMq{\um^{(i)}(p,s]} = q^{\al} \FMq{\um^{(i)}[p,s]} \cdot  \FMq{\um^{(i)}(p,s)} + q^{\ga}   \FMq{  \um(i; p,s) }
\end{align*}
for some $\al,\ga \in \dfrac{1}{2}\Z$.
By using the $q$-algorithm and the argument in \cite[Lemma 5.6]{H06} (or \cite[Theorem 9.6, Lemma 9.9]{HO19}), the product of
 $\FMq{ \um^{(i)}[p,s) }$ and $\FMq{\um^{(i)}(p,s]}$ has exactly distinct $k$ dominant monomials
 \begin{align*}
M_1,\, M_2,\, \ldots,\, M_k,
 \end{align*}
 where
$\ev_{q=1}(M_1) = m^{(i)}[p,s)   m^{(i)}(p,s]$.
Moreover, $M_1,\ldots,M_{k-1}$
exhaust the dominant monomials occurring in  $F_q(\um^{(i)}[p,s])  F_q(\um^{(i)}(p,s))$  and
\begin{align*}
\ev_{q=1}(M_k) = \left(  m^{(i)}[p,s)   B^{-1}_{i,s-1}   B^{-1}_{i,s-3} \cdots  B^{-1}_{i,p+1} \right) m^{(i)}(p,s]  =  m(i; p,s).
\end{align*}
Hence, our claim follows from Theorem~\ref{thm: F_q} and Lemma \ref{lem: last term in T-system}.

\smallskip

Second, we compute $\al = \al(i,k)$ and $\ga = \ga(i,k)$ explicitly.
By Theorem~\ref{thm: F_q}, Lemma~\ref{lem: m(i;p,s)} implies that
\begin{align*}
 \FMq{  \um(i; p,s) } = \prod_{j; \; d(i,j)=1}  \hspace{-1.5ex}  \FMq{\um^{(j)}(p,s)}^{-\sfc_{j,i}}.
\end{align*}
Also, by Lemma~\ref{lem: m[p,s] m(p,s)}, we also have
\begin{align*}
\FMq{ \um^{(i)}[p,s] } *   \FMq{\um^{(i)}(p,s)} = \FMq{\um^{(i)}(p,s)} * \FMq{ \um^{(i)}[p,s] }.
\end{align*}
Thus it suffices to compute
$\al,\ga$ such that
\begin{align*}
 \um^{(i)}[p,s)   * \um^{(i)}(p,s] = q^\al    \um^{(i)}[p,s]   *  \um^{(i)}(p,s)  = q^\al  \um^{(i)}(p,s) * \um^{(i)}[p,s]
\end{align*}
and
\begin{align*}
\left( \um^{(i)}[p,s) \cdot   \tB^{-1}_{i,s-1} \cdot  \tB^{-1}_{i,s-3} \cdots \tB^{-1}_{i,p+1} \right)  *  \um^{(i)}(p,s]  = q^\ga \;  \um(i; p,s).
\end{align*}
The coefficient $\al$ can be computed as follows:
\begin{align*}
\al &= \sum_{a=1}^{k-1} \ucalN(i,p; i,p+2a)  + \dfrac{1}{2} \ucalN(i,p;i,p+2k)\allowdisplaybreaks \\
&= \sum_{a=1}^{k-1} \left( \tfb_{i,i}(2a+1) -\tfb_{i,i}(2a-1)  \right)  + \dfrac{1}{2} \left( \tfb_{i,i}(2k+1) -\tfb_{i,i}(2k-1)  \right) \allowdisplaybreaks\\
&= -\tfb_{i,i}(1) + \dfrac{1}{2} \left( \tfb_{i,i}(2k+1) + \tfb_{i,i}(2k-1)  \right) = -d_i +  \dfrac{1}{2} \left( \teta_{i,i}(2k+1) + \teta_{i,i}(2k-1)  \right).
\end{align*}

\noindent
Note that $\um \seteq \left( \um^{(i)}[p,s) \cdot   \tB^{-1}_{i,s-1} \cdot  \tB^{-1}_{i,s-3} \cdots \tB^{-1}_{i,p+1} \right) $ is contained in $F_q( \um^{(i)}[p,s))$ with coefficient $1$, and
$\um \cdot \um^{(i)}(p,s]   =\displaystyle \prod_{j; \; d(i,j)=1} \um^{(j)}(p,s)^{-\sfc_{j,i}}$. Thus we have
\begin{align*}
\um * \um^{(i)}(p,s] = \left( \bl \um^{(i)}(p,s] \br^{-1} \cdot \prod_{j; \; d(i,j)=1} \um^{(j)}(p,s)^{-\sfc_{j,i}}   \right) * \um^{(i)}(p,s] = q^{\ga}  \prod_{j; \; d(i,j)=1} \um^{(j)}(p,s)^{-\sfc_{j,i}},
\end{align*}
where
\begin{align*}
\ga & = \dfrac{1}{2} \sum_{j; \; d(i,j)=1} \hspace{-2ex} -\sfc_{j,i}  \sum_{a=1}^k \sum_{b=1}^k    \ucalN(j,p+2a-1; i,p+2b)\allowdisplaybreaks \\
 & = \dfrac{1}{2} \sum_{j; \; d(i,j)=1} \hspace{-2ex} -\sfc_{j,i} \sum_{a=1}^k \sum_{b=1}^k    \left(
\tfb_{j,i}(2(a-b)-2)- \tfb_{j,i}(2(a-b)) - \tfb_{j,i}(2(b-a)) +\tfb_{j,i}(2(b-a)+2)
\right)\allowdisplaybreaks\\
 & = \dfrac{1}{2} \sum_{j; \; d(i,j)=1} \hspace{-2ex} -\sfc_{j,i} \sum_{a=1}^k \left(
\tfb_{j,i}(2(a-k)-2)- \tfb_{j,i}(2(a-1)) - \tfb_{j,i}(2(1-a)) +\tfb_{j,i}(2(k-a)+2)
\right) \allowdisplaybreaks\\
 & = \dfrac{1}{2} \sum_{j; \; d(i,j)=1} \hspace{-2ex} -\sfc_{j,i} \sum_{a=1}^k \left(
- \tfb_{j,i}(2(a-1))  +\tfb_{j,i}(2(k-a)+2)
\right) =  \dfrac{1}{2} \sum_{j; \; d(i,j)=1} \hspace{-2ex} -\sfc_{j,i}  \tfb_{j,i}(2k) \allowdisplaybreaks\\
& =    \dfrac{1}{2} \sum_{j; \; d(i,j)=1} \hspace{-2ex} -\sfc_{j,i}  \teta_{j,i}(2k) =    \dfrac{1}{2} \sum_{ j; \; d(i,j)=1} \hspace{-2ex} -\sfc_{j,i}  \teta_{i,j}(2k).
\end{align*}
Then our proof is completed by Lemma~\ref{lem: teta}.
\end{proof}

\begin{example} \label{ex: quantum folded T-system of Fq25 Fq27 in G2}
Let us recall the formula of $F_q(\tX_{2,5})$ in \eqref{ex: folding}.
Also, $F_q(\sfX_{2,5}) = q^{\frac{3}{2}}F_q(\tX_{2,5})  \in \mathfrak{K}_q(\sfg)$ and it is bar-invariant with respect to \eqref{eq: bar involution}.
Note that $F_q(\tX_{2,7}) = \sfT_2(F_q(\tX_{2,5}))$ and $F_q(\sfX_{2,7}) = q^{\frac{3}{2}}F_q(\tX_{2,7})$.
Clearly, these computations implies that $F_q(\sfX_{2,5}) * F_q(\sfX_{2,7})$ has two dominant $\calX_q$-monomials, namely, $\underline{X_{2,5} X_{2,7}}$ and $\underline{X_{1,6}^3}$.
By Theorem \ref{thm: F_q}, we should have
\begin{equation} \label{eq: Fq25 = q al Fq2527 + q gam Fq16 3}
	F_q(\sfX_{2,5}) * F_q(\sfX_{2,7})
	=
	q^{\frac{3}{2}} F_q(\underline{X_{2,5} X_{2,7}})
	+
	q^{\frac{9}{2}} F_q(\underline{X_{1,6}})^3.
\end{equation}
On the other hand, we obtain
\begin{equation*}
\begin{split}
       d_2=3, \quad
	\gamma(2, 1) = \frac{1}{2} \left( \tfb_{2,2}(1)+\tfb_{2,2}(3) \right) = \frac{9}{2}, \quad
	\alpha(2, 1) = \gamma(2,1) - d_2 = \frac{3}{2},
	\quad
	-\sfc_{1,2} = 3,
\end{split}
\end{equation*}
where $\tfb_{2,2}(1)=3$ and $\tfb_{2,2}(3)=6$ from \eqref{eq: inverse of tBt}.
Hence \eqref{eq: Fq25 = q al Fq2527 + q gam Fq16 3} illustrates Theorem \ref{thm: quantum folded}.
\end{example}

\section{Quantum cluster algebra} \label{sec:QCA}

\noindent
In this section we recall the definition of skew-symmetrizable quantum cluster algebras of infinite rank, following \cite{BZ05}, \cite[\S 8]{GLS13}, \cite{HL16} and \cite{KKOP2}.

\subsection{Quantum seed}
Let $\sfK$ be an index set described in Section~\ref{subsec:
Bi-deco}. Let $L=(\la_{i,j})_{i,j \in \sfK}$ be a skew symmetric
integer-valued $\sfK \times \sfK$-matrix. Let $q$ be an indeterminate.

\begin{definition} \label{def: quantum torus}
We define $(\scrP(L) ,\star) $ as the $\Z[q^{\pm \frac{1}{2}}]$-algebra, called the \emph{quantum torus associated to $L$}, generated by a family of elements
$\{ Z_{i} \}_{i \in \sfK}$ with the defining relations
\begin{align*}
Z_i  \star Z_j = q^{\la_{i,j}} Z_j \star  Z_i \qquad (i,j \in \sfK).
\end{align*}
We denote by $\frakF(L)$ the skew field of fractions of $\scrP(L)$.
\end{definition}

For $\bfa = (a_i)_{i \in \sfK} \in \Z^{\oplus \sfK}$, we define the element $Z^{\bfa}$ of $\scrF(L)$ as
\begin{align} \label{eq: ordered prod on Z}
  Z^{\bfa}  \seteq q^{\frac{1}{2} \sum_{i > j} a_ia_j \la_{i,j}}   \sta_{i \in \sfK}^\to Z_i^{a_i}
\end{align}
(cf.~\eqref{eq: inv mono1}).  Here we take a total order $<$ on the set $\sfK$.  Note that $Z^\bfa$ does not depend on the choice of a total order on $\sfK$.
We have
$$
Z^{\bfa} \star  Z^{\bfb} = q^{\frac{1}{2} \sum_{i,j \in \sfK}a_ib_j \la_{i,j}}Z^{\bfa+\bfb}.
$$

Let $(\calA,\star)$ be a $\Z[q^{\pm \frac{1}{2} }]$-algebra. We say that a family
$\{z_i\}_{i\in\K}$ of elements of $\calA$ is {\em $L$-commuting} if
it satisfies $z_i \star z_j=q^{\la_{i,j}}z_j \star z_i$ for any $i,j\in\K$.
In that case we can define $z^\mathbf{a}$ for any $\mathbf{a}\in\Z_{\ge0}^{ \oplus \sfK}$  as in \eqref{eq: ordered prod on Z}.
We say that an $L$-commuting family $\{z_i\}_{i\in\sfK}$ is {\em algebraically independent} if the algebra map
$\mathscr P(L)\to \calA$ given by $Z_i\mapsto z_i$ is injective.

Let $\tB = (b_{i,j})_{i \in \sfK, j \in \Kex}$ be an integer-valued $\sfK \times \Kex$-exchange matrix satisfying \eqref{eq: ex matrix}.
We say that the pair $(L,\tB)$  is \emph{compatible with a diagonal
matrix} ${\rm diag}( \sfd_i \in \Z_{\ge 1} \ | \  i \in \sfK)$, if we have
\begin{equation} \label{eq: comp D}
\begin{split}
\sum_{k \in \sfK} b_{ki} \lambda_{kj} = \delta_{i,j} \sfd_i, \quad \text{ equivalently, }  \ \
(L\tB)_{ji} = -\delta_{i,j} \sfd_i,
\end{split}
\end{equation}
for any $i \in \sfK_\ex$ and $j \in \sfK$.
We also call the pair $(L, \tB)$ a \emph{compatible pair} for short.

Let $(L,  \widetilde B)$ be a compatible pair and $\calA$ a $\Z[q^{\pm1/2}]$-algebra.
We say that $\seed = (\{z_i\}_{i\in\K},L, \widetilde B)$ is
a {\em quantum seed} in $\calA$ if $\{z_i\}_{i\in\K}$ is
an algebraically independent $L$-commuting family of elements of $\calA$.
The set $\{z_i\}_{i\in\K}$ is called the {\em quantum cluster} of $\seed$ and
its elements  the {\em quantum cluster variables}.
The quantum cluster variables $z_i$ ($i\in\Kfr$) are called  the {\em frozen variables}.
The elements $z^{\bf a}$ (${\bf a}\in \Z_{\ge0}^{\oplus \K}$) are called  the {\em quantum cluster monomials}.


\subsection{Mutation}
For $k\in\Kex$, we define a $\K\times \K$-matrix  $E=(e_{i,j})_{i,j\in\K}$ and a
$\Kex\times \Kex$-matrix $F=(f_{i,j})_{i,j\in\Kex}$ as follows:
\eqn
e_{i,j}=
\begin{cases}
  \delta_{i,j} & \text{if} \ j \neq k, \\
  -1 & \text{if} \ i= j = k, \\
  \max(0, -b_{i,k}) & \text{if} \ i \neq  j = k,
\end{cases}
\hs{10ex}
f_{i,j}=
\begin{cases}
  \delta_{i,j} & \text{if} \ i \neq k, \\
  -1 & \text{if} \ i= j = k, \\
  \max(0, b_{k,j}) & \text{if} \ i = k \neq j.
\end{cases}
\eneqn
The {\em mutation $\mu_k(L,\widetilde B)\seteq(\mu_k(L),\mu_k(\widetilde B))$ of a compatible pair $(L,\widetilde B)$ in direction $k$} is given by
\eqn
\mu_k(L)\seteq(E^T) \, L \, E, \quad \mu_k(\widetilde B)\seteq E \, \widetilde B \, F.
\eneqn

\noindent
We define
\eq \label{eq: a' and a''}
&&a_i'=
\begin{cases}
  -1 & \text{if} \ i=k, \\
 \max(0,b_{i,k}) & \text{if} \ i\neq k,
\end{cases} \qquad
a_i''=
\begin{cases}
  -1 & \text{if} \ i=k, \\
 \max(0,-b_{i,k}) & \text{if} \ i\neq k.
\end{cases}
\label{eq:aa}
\eneq
and set ${\bf a}'\seteq(a_i')$ and ${\bf a}''\seteq(a_i'') \in \Z^{\oplus \sfK}$.

Let $\calA$ be a $\Z[q^{\pm1/2}]$-algebra contained in a skew field $K$.
Let $\scrS=(\{z_i\}_{i\in\K},L, \widetilde B)$ be a quantum seed in $\calA$.
Define the elements $\mu_k(z)_i$ of $K$ by
\eq \label{eq:quantum mutation}
\mu_k(z)_i\seteq
\begin{cases}
  z^{{\bf a}'}  +   z^{{\bf a}''} & \text{if} \ i=k, \\
 z_i & \text{if} \ i\neq k.
\end{cases}
\eneq
Then $\{\mu_k(z)_i\}$ is an algebraically independent
$\mu_k(L)$-commuting family in $K$.
We call
\eqn
\mu_k(\scrS)\seteq\bl\{\mu_k(z)_i\}_{i\in\K},\mu_k(L),\mu_k(\widetilde B)\br
\eneqn
the {\em mutation
of $\seed$ in direction $k$}.
It becomes a new quantum seed in $K$; that means,
\ben
\item $\bl \mu_k(L),\mu_k(\widetilde B) \br$ is compatible with the diagonal matrix of $(L,\tB)$,
\item $\{\mu_k(z)_i\}_{i\in\K}$ is $\mu_k(L)$-commuting.
\ee

\begin{definition} \label{def: mutation equi}
Let $\seed=( \{ z_i \}_{i \in \K}, L,\tB )$ and $\seed'=( \{ z_i'\}_{i \in \K'}, L',\tB' )$ be quantum seeds in a $\Z[q^{\pm1/2}]$-algebra $\calA$.

\bnum
\item We say that \emph{$\seed'$ is mutated from $\seed$}
if the following condition is satisfied: For any finite subset $\mathsf{J}$ of $\K'$,
there exist
\bna
\item a finite sequence $(k_1,k_2,\ldots,k_r)$ in $\Kex$,
\item an injective map $\upsigma\colon \sfJ \to \K$, depending on the choice of $\sfJ$, such that
\ben
\item $ \upsigma(\mathsf{J}_\ex) \subset \K_\ex$, where $\mathsf{J}_{\ex} \seteq \mathsf{J} \cap (\K')_\ex$,

\item $z'_{j}=\upmu(z)_{\upsigma(j)}$  for all $j \in \sfJ$,
\item $ (\tB')_{(i,j)}=\upmu(\tB)_{\upsigma(i),\upsigma(j)}$
for any $(i,j)\in \sfJ \times {\sfJ}^\ex$,
\ee
\quad where $\upmu \seteq  \mu_{k_r} \circ \cdots \circ \mu_{k_1}$.
\ee
\item  We say that the quantum seeds $\seed$ and $\seed'$ are \emph{mutation equivalent} if $\seed'$ is mutated from $\seed$ and $\seed$ is also mutated from $\seed'$.
In this case, we write $\seed \simeq \seed'$.
\ee
\end{definition}

\subsection{Mutation of valued quiver} Recall that we can associate the valued quiver $\calQ_{\tB}$ to an exchange matrix $\tB$. Here we describe the algorithm transforming a valued quiver $\calQ$
into a new valued quiver $\mu_k(\calQ)$ $(k \in \Kex)$, which
corresponds to $\mu_k(\tB)$.

\begin{algorithm} \label{Alg. mutation}
For $k \in \sfK_\ex$, the \emph{valued quiver mutation $\mu_k$} transforms $\calQ$ into a new valued quiver $\mu_k(\calQ)$ via the following rules, where we assume {\rm (i)} $ac>0$ or $bd>0$, and {\rm (ii)} we do not perform {\rm ($\mathcal{NC}$)} and {\rm ($\mathcal{C}$)} below, if $i$ and $j$ are frozen at the same time$:$
\ben
\item[]\hspace{-0.69cm}{\rm ($\mathcal{NC}$)} For each full-subquiver   $\xymatrix@!C=7mm@R=1mm{ i \ar[r]_{\ulcorner  a,b \lrcorner }  \ar@/^0.7pc/[rr]^{\ulcorner  e,f \lrcorner} & k  \ar[r]_{\ulcorner  c,d \lrcorner} & j }$ in $\calQ$,
we change the value of the arrow from $i$ to $j$ into $\ulcorner e+ac,f-bd \lrcorner:$
\begin{align*}
\xymatrix@!C=20mm@R=1mm{ i  \ar[r]_{\ulcorner e+ac,f-bd \lrcorner} & j}.
\end{align*}
\item[{\rm ($\mathcal{C}$)}]  For each  full-subquiver   $\xymatrix@!C=7mm@R=1mm{ i \ar[r]_{\ulcorner  a,b \lrcorner }  & k  \ar[r]_{\ulcorner  c,d \lrcorner} & j  \ar@/_0.7pc/[ll]_{\ulcorner  e,f \lrcorner}  }$ with $(e,f) \ne (0,0)$ in $\calQ$, we change the  valued arrow between $i$ and $j$ as follows:
\begin{align*}
\bc
\xymatrix@!C=20mm@R=7mm{ i  \ar@{<-}[r]_{\ulcorner e-bd,f+ac \lrcorner} & j}   & \text{ if }  f+ac \le  0 \le e-bd,\\
\xymatrix@!C=20mm@R=7mm{ i  \ar[r]_{\ulcorner f+ac,e-bd \lrcorner} & j}  & \text{ if }   f+ac \ge  0 \ge  e-bd.
\ec
\end{align*}

\item[{\rm ($\mathcal{R}$)}] Reverse the direction of each arrow incident to the vertex $k$ and change the value $\ulcorner  a,b \lrcorner $  of each arrow into $\ulcorner  -b,-a \lrcorner$.
\ee
Here if there is no arrow between $i$ and $j$ in {\rm ($\mathcal{NC}$)} and {\rm ($\mathcal{C}$)}, then put $e = f = 0$ and follow the same rule.
\end{algorithm}

\begin{example} \label{ex: mutations in B3}
Consider the following $9 \times 6$ integer-valued matrix:
\begin{align} \label{eq: B3}
\tB= \tiny{
\left(\begin{array}{rrrrrr}
0 & -1 & 0 & 1 & 0 & 0 \\
1 & 0 & -1 & -1 & 1 & 0 \\
0 & 2 & 0 & 0 & -2 & 1 \\
-1 & 1 & 0 & 0 & -1 & 0 \\
0 & -1 & 1 & 1 & 0 & -1 \\
0 & 0 & -1 & 0 & 2 & 0 \\
0 & 0 & 0 & -1 & 1 & 0 \\
0 & 0 & 0 & 0 & -1 & 1 \\
0 & 0 & 0 & 0 & 0 & -1
\end{array}\right)}
\end{align}
By taking $\sfK_\ex=\{1,2,3,4,5,6\}$ and  $\sfK_\fr=\{7,8,9\}$, one can see that
its principal part is skew-symmetrizable with $\sfS={\rm diag}(2,2,1,2,2,1)$.

Using Convention~\ref{conv: bi-deco}, the valued quiver $\calQ$ associated to $\tB$ in~\eqref{eq: B3} can be drawn as
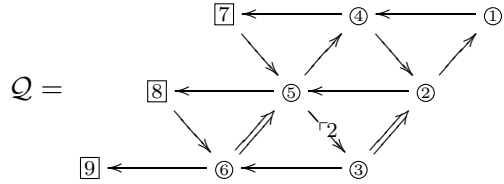
\begin{align*}
\calQ =
\raisebox{11mm}{\xymatrix@!C=0.5mm@R=5mm{
&& \reced{7} \ar@{->}[dr]   && \circled{4} \ar@{->}[dr]\ar[ll] &&\circled{1}  \ar[ll] \\
&  \reced{8}\ar@{->}[dr]    && \circled{5} \ar@{->}[dr]|{ \; \ulcorner \hspace{-.2ex} 2}\ar@{->}[ur]\ar[ll]&& \circled{2}  \ar@{->}[ur]\ar[ll] \\
  \reced{9}  &&\circled{6}  \ar@{=>}[ur] \ar[ll]&& \circled{3} \ar@{=>}[ur]\ar[ll]  }}
\end{align*}
Here $\reced{k}$ denotes  $k \in \sfK_\fr$.
Then $\mu_2(\calQ)$, $\mu_5(\calQ)$, 
are depicted as follows:
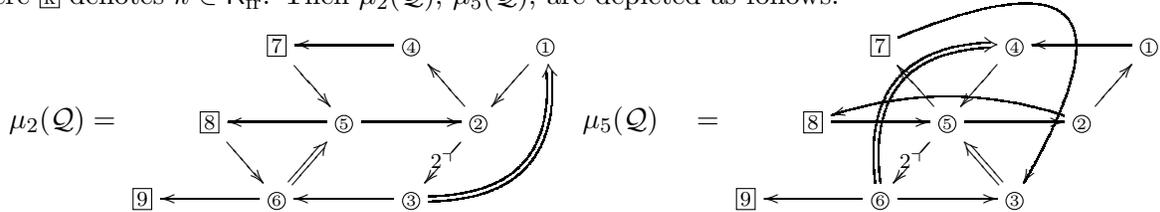
\begin{align*}
 \mu_2(\calQ) &=
\raisebox{11mm}{\xymatrix@!C=0.5mm@R=5mm{
&& \reced{7} \ar@{->}[dr]   && \circled{4}  \ar[ll] &&\circled{1}    \ar[dl]\\
&  \reced{8}\ar@{->}[dr]    && \circled{5}   \ar[ll]\ar[rr]&& \circled{2}   \ar[ul] \ar[dl]|{ 2 \urcorner \; }\\
  \reced{9}  &&\circled{6}  \ar@{=>}[ur] \ar[ll]&& \circled{3} \ar[ll]  \ar@/_2pc/@{=>}[uurr]  }}  \ \
\mu_5(\calQ)& =
\raisebox{11mm}{\xymatrix@!C=0.5mm@R=5mm{
&& \reced{7}  \ar@/^5pc/@{->}[ddrr]      && \circled{4} \ar@{->}[dl]  &&\circled{1}  \ar[ll] \\
&  \reced{8}\ar[rr]     && \circled{5} \ar@{<=}[dr] \ar@{->}[ul] \ar@{->}[dl]|{2 \urcorner  \; } \ar[rr] && \circled{2}
\ar@/_0.8pc/[llll]   \ar@{->}[ur]  \\
  \reced{9}  &&\circled{6}   \ar@/^2pc/@{=>}[uurr]    \ar[rr] \ar[ll]&& \circled{3}  }}
\end{align*}
\end{example}

\subsection{Quantum cluster algebra}
Let $\seed=(\{z_i\}_{i\in\K},L, \widetilde B)$ be a quantum seed in a $\Z[q^{\pm1/2}]$-algebra $\calA$.
   The {\em quantum cluster algebra $\mathscr A_{q^{1/2}}(\seed)$} associated to the quantum seed $\seed$ is the $\Z[q^{\pm \frac{1}{2}}]$-subalgebra of the skew field $K$ generated by all the quantum cluster variables in the quantum seeds obtained from $\seed$ by any \emph{finite} sequence of mutations.
Here we call $\seed$ the {\em initial quantum seed}
of the quantum cluster algebra $\mathscr A_{q^{1/2}}(\seed)$.

\begin{lemma} \label{lem: mutation equiv on QCA}
Let $\seed$ and $\seed'$ be quantum seeds in $\calA$. If $\seed'$ is mutated from $\seed$, then
$ \scrA_{q^{1/2}}(\seed')$ is isomorphic to $\Z[q^{\pm \frac{1}{2}}]$-subalgebra of $\scrA_{q^{1/2}}(\seed).$
Furthermore, if $\seed$ and $\seed'$   are mutation equivalent to each other, then we have
$$   \scrA_{q^{1/2}}(\seed')\simeq \scrA_{q^{1/2}}(\seed).$$
\end{lemma}

\begin{proof}
This assertion follows from Definition~\ref{def: mutation equi}.
\end{proof}

\smallskip

\begin{definition} \label{def: quantum cluster algebra structure}
A {\it  quantum cluster algebra structure} associated with a quantum seed $\seed$ in  a $\Z[q^{\pm 1/2}]$-algebra $\calA$, contained in a skew field $K$,  is a family $\scrF$ of
quantum seeds in $\calA$ satisfying the following conditions:
\bna
\item For any quantum seed $\seed$ in $\scrF$, the quantum cluster algebra $\mathscr A_{q^{1/2}}(\seed)$ is isomorphic to $\calA$ as a $\Z[q^{\pm 1/2}]$-algebra.
\item Any mutation of a quantum seed in $\scrF$ is in $\scrF$.
\item For any pair $\seed$, $\seed'$ of quantum seeds in $\scrF$, we have  $\seed' \simeq \seed$.
\ee
\end{definition}

\section{Quantum cluster algebra structure on $\frakK_{q,\xi}(\g)$} \label{sec: qCA on Kqxi}

\noindent
In this section, we will prove that the ring $\frakK_{q,\xi}(\g)$ has a quantum cluster algebra structure  based on the recent work \cite{KO23} by Kashiwara--Oh. 
As applications, we obtain
\begin{enumerate}[$\bullet$]
\item a quantum cluster algebra algorithm to compute the KR-polynomials $F_q(\um^{(i)}[a,b])$ for KR-monomials $m^{(i)}[a,b]$,
\vskip 1mm
\item a $q$-commutativity for KR-polynomials $F_q(\um^{(i)}_{k,r})$ and $F_q(\um^{(j)}_{l,t})$ satisfying certain conditions on the pair of their KR-monomials $(m^{(i)}_{k,r},\,m^{(j)}_{l,t})$.
\end{enumerate}

In this section, we shall employ the framework in \cite{HL16,B21} for our goal.

\subsection{Compatible pair}
Let $\sfS$ be a convex subset of $\tDynkinf_0$ with an upper bound (recall Definition \ref{def: convex subset}).
For each $j \in \Dynkinf_0$, we set
$$ \upxi_j  \seteq \max( s \ | \  (j,s) \in \sfS ).$$

Recall the exchange matrices in Definition \ref{def: tB associated with Dynkin} and Definition~\ref{def: convex subset}.

\begin{theorem} \cite[Theorem 7.1]{KO23} $($see also \cite{FHOO2}$)$ \label{thm: compatible pairs}
Define
\begin{align*}
\La_{(i,p),(j,s)} = \ucalN(  m^{(i)}[p,\upxi_i], m^{(j)}[s,\upxi_j]  ) \qquad (i,p),(j,s) \in \sfS.
\end{align*}
Then the pair $ (    (\La_{(i,p),(j,s)})_{(i,p),(j,s) \in \sfS},  {}^\sfS\tB )$ is compatible with  ${\rm diag}(  2d_{i,p} \seteq 2d_i \ |  \ (i,p) \in  \sfS )$.
\end{theorem}

Recall that the subset $\lxi\tDynkinf_0$ is convex without frozen
indices. Thus the pair
  $(\lxi L, \lxi\tB)$ is compatible with  ${\rm diag}(  2d_{i,p} \seteq 2d_i \ |  \ (i,p) \in   \lxi\tDynkinf_0 )$,
where
\begin{align} \label{eq: matrix xi L}
\lxi L=\bl   \La_{(i,p),(j,s)} \br_{(i,p),(j,s) \in \lxi\widetilde{\sDynkinf}_0} \quad\text{ and }\quad
\La_{(i,p),(j,s)}  = \ucalN \bl m^{(i)}[p,\xi_i], m^{(j)}[s,\xi_j]  \br.
\end{align}

\subsection{Sequence of mutations} Let us consider the valued quiver $\tDynkinxi$ associated to the height function $\xi$
of $Q$.  Note that, for a source $i$ of $Q$,
\eq &&
\parbox{85ex}{
\bnum
\item \label{it: boundary} the vertex $(i,\xi_i)$ is located at the boundary   of ${}^{\xi}\tDynkin$ determined by $\xi$,   and vertically sink and horizontally source,
\item $s_i\xi$ is a height function defined as in~\eqref{eq: si xi}.
\ee
}\label{eq:  i xi}
\eneq

For a source $i$ of $Q$, we set a sequence of mutations
\begin{equation} \label{bi_xi upmu}
   {}^{\bii}_{\xi}\upmu  \seteq    \cdots \circ \mu_{(i,\xi_i-4)} \circ \mu_{(i,\xi_i-2)} \circ \mu_{(i,\xi_i)}
\end{equation}
and call it the \emph{forward shift} at $i$ (see \cite{HL16} for $\calK_\ttt(\scrC_\g^0)$-cases).

\begin{proposition} \label{prop: one step}
For a Dynkin quiver $Q=(\Dynkin,\xi)$ and a source $i$, we have
$$
 {}^{\bii}_{\xi}\upmu(  {}^{\xi}\tDynkin ) \simeq  {}^{s_{i}\xi}\tDynkin.
$$
\end{proposition}

\begin{proof}
We shall prove our assertion by an inductive argument on the sequence ${}^{\bii}_{\xi}\upmu$.
For this, we observe first two steps $\mu_{(i,\xi_i)}$ and $\mu_{(i,\xi_i-2)} \circ \mu_{(i,\xi_i)}$.

\smallskip

\noindent
{\it Step 1}.
Let us consider  $ \mu_{(i,\xi_i)}(\tDynkinxi)$. In this case, the vertex $(i,\xi)$ in $\tDynkinxi$ (marked with $*$ below) is vertically sink and horizontally source in $\tDynkinxi$ by~\eqref{eq: ex tDynkin} and ~\eqref{eq:  i xi}~\ref{it: boundary} as follows:
\begin{align*}
\scalebox{0.8}{\xymatrix@!C=10mm@R=6mm{
 (j,\xi_j-4)\ar@{->}[dr]|{ \ulcorner -\sfc_{j,i},\sfc_{i,j} \lrcorner}
&& (j,\xi_j-2)\ar@{->}[dr]|{ \ulcorner -\sfc_{j,i},\sfc_{i,j} \lrcorner} \ar[ll] &&  (j,\xi_j) \ar@{->}[dr]|{ \ulcorner -\sfc_{j,i},\sfc_{i,j} \lrcorner} \ar[ll] \\
& (i,\xi_i-4) \ar@{->}[dr]|{ \ulcorner -\sfc_{ij'},\sfc_{j'i} \lrcorner}\ar@{->}[ur]|{ \ulcorner -\sfc_{i,j},\sfc_{j,i} \lrcorner}\ar[l]
&& (i,\xi_i-2) \ar@{->}[dr]|{ \ulcorner -\sfc_{ij'},\sfc_{j'i} \lrcorner}\ar@{->}[ur]|{ \ulcorner -\sfc_{i,j},\sfc_{j,i} \lrcorner}\ar[ll]
&& (i,\xi_i)^* \ar[ll] \\
 (j',\xi_{j'}-4)  \ar@{->}[ur]|{ \ulcorner -\sfc_{j'i},\sfc_{ij'} \lrcorner}
&& (j',\xi_{j'}-2)  \ar@{->}[ur]|{ \ulcorner -\sfc_{j'i},\sfc_{ij'} \lrcorner} \ar[ll]&& (j',\xi_{j'}) \ar@{->}[ur]|{ \ulcorner -\sfc_{j'i},\sfc_{ij'} \lrcorner}\ar[ll]  }}
\end{align*}
Here $j$ and $j'$ are indices in $\Dynkin_0$ such that $d(i,j)=d(i,j')=1$.
Note that, in order to observe the behavior with respect to $\mu_{(i,\xi_i)}$, it suffices to consider the full-subquiver described as above.

Applying Algorithm~\ref{Alg. mutation},  $ \mu_{(i,\xi_i)}(\tDynkinxi)$ can be depicted as follows:
\begin{align*}
\scalebox{0.8}{\xymatrix@!C=10mm@R=6mm{
 (j,\xi_j-4)\ar@{->}[dr]|{ \ulcorner -\sfc_{j,i},\sfc_{i,j} \lrcorner}
&& (j,\xi_j-2)\ar@{->}[dr]|{ \ulcorner -\sfc_{j,i},\sfc_{i,j} \lrcorner} \ar[ll] &&  (j,\xi_j) \ar@{<-}[dr]|{ \ulcorner -\sfc_{i,j},\sfc_{j,i} \lrcorner} \ar[ll] \\
& (i,\xi_i-4) \ar@{->}[dr]|{ \ulcorner -\sfc_{ij'},\sfc_{j'i} \lrcorner}\ar@{->}[ur]|{ \ulcorner -\sfc_{i,j},\sfc_{j,i} \lrcorner}\ar[l]
&& (i,\xi_i-2)^*  \ar[ll] \ar[rr]
&& (i,\xi_i)   \\
 (j',\xi_{j'}-4)  \ar@{->}[ur]|{ \ulcorner -\sfc_{j'i},\sfc_{ij'} \lrcorner}
&& (j',\xi_{j'}-2)  \ar@{->}[ur]|{ \ulcorner -\sfc_{j'i},\sfc_{ij'} \lrcorner} \ar[ll]&& (j',\xi_{j'}) \ar@{<-}[ur]|{ \ulcorner -\sfc_{ij'},\sfc_{j'i} \lrcorner}\ar[ll]  }}
\end{align*}
in which the vertex $(i,\xi_i-2)$ (marked with $*$ above) becomes vertically sink and horizontally source.
\smallskip

\noindent
{\it Step 2}.
Let us consider $ (\mu_{(i,\xi_i-2)} \circ \mu_{(i,\xi_i)} ) (\tDynkinxi)$.  Applying Algorithm~\ref{Alg. mutation} again,
$ (\mu_{(i,\xi_i-2)} \circ \mu_{(i,\xi_i)} ) (\tDynkinxi)$ becomes
\begin{align*}
\scalebox{0.8}{\xymatrix@!C=10mm@R=6mm{
 (j,\xi_j-4)\ar@{->}[dr]|{ \ulcorner -\sfc_{j,i},\sfc_{i,j} \lrcorner}
&& (j,\xi_j-2)\ar@{<-}[dr]|{ \ulcorner -\sfc_{i,j} ,\sfc_{j,i}\lrcorner} \ar[ll] \ar@{->}[drrr]|{ \ulcorner -\sfc_{j,i},\sfc_{i,j} \lrcorner}
 &&  (j,\xi_j) \ar@{<-}[dr]|{ \ulcorner -\sfc_{i,j},\sfc_{j,i} \lrcorner} \ar[ll] \\
& (i,\xi_i-4)^* \ar[rr] \ar[l]   && (i,\xi_i-2)
&& (i,\xi_i)   \ar[ll]   \\
 (j',\xi_{j'}-4)  \ar@{->}[ur]|{ \ulcorner -\sfc_{j'i},\sfc_{ij'} \lrcorner}
&& (j',\xi_{j'}-2)  \ar@{->}[ur]|{ \ulcorner -\sfc_{j'i},\sfc_{ij'} \lrcorner} \ar[ll]  \ar@{->}[urrr]|{ \ulcorner -\sfc_{j'i},\sfc_{ij'} \lrcorner}
&&  (j',\xi_{j'}) \ar@{<-}[ur]|{ \ulcorner -\sfc_{ij'},\sfc_{j'i} \lrcorner}\ar[ll]  }}
\end{align*}
which is isomorphic to
\begin{align} \label{eq: 2-steps}
\scalebox{0.8}{\raisebox{14mm}{\xymatrix@!C=10mm@R=6mm{
&(j,\xi_j-4) \ar@{->}[dr]|{ \ulcorner -\sfc_{j,i},\sfc_{i,j} \lrcorner}\ar[l]
&&&&(j,\xi_j-2) \ar[llll] \ar[dr]|{\ulcorner  -\sfc_{j,i},\sfc_{i,j}  \lrcorner}
&& (j,\xi_j) \ar[ll]\\
(i,\xi_i-6)  \ar[ur]|{\ulcorner  -\sfc_{i,j},\sfc_{j,i}  \lrcorner} \ar[dr]|{\ulcorner  -\sfc_{ij'},\sfc_{j'i}  \lrcorner}
&&(i,\xi_i-4)^* \ar[rr] \ar[ll]
&&(i,\xi_i-2) \ar[ur]|{\ulcorner  -\sfc_{i,j},\sfc_{j,i}  \lrcorner} \ar[dr]|{\ulcorner  -\sfc_{ij'},\sfc_{j'i}  \lrcorner}
&& (i,\xi_i)\ar[ur]|{\ulcorner  -\sfc_{i,j},\sfc_{j,i}  \lrcorner} \ar[dr]|{\ulcorner  -\sfc_{ij'},\sfc_{j'i}  \lrcorner}\ar[ll] \\
&(j',\xi_{j'}-4) \ar@{->}[ur]|{ \ulcorner -\sfc_{j'i},\sfc_{ij'} \lrcorner}\ar[l]
&&&&(j',\xi_{j'}-2) \ar[ur]|{\ulcorner  -\sfc_{j'i},\sfc_{ij'}  \lrcorner}  \ar[llll]
&& (j',\xi_{j'}) \ar[ll]
}}}
\end{align}
Here the vertex $(i,\xi_i-4)$ (marked with $*$ in \eqref{eq: 2-steps}) becomes also vertically sink and
horizontally source.
\smallskip

By {\it Step 1} and {\it Step 2}, we observe that the full-subquiver consisting of the
rightmost 6-vertices in~\eqref{eq: 2-steps} are isomorphic to the
rightmost 6-vertices of ${}^{s_i\xi}\tDynkin$.
Furthermore, since the local circumstance of $(i,\xi_i-4)$ in  $ (\mu_{(i,\xi_i-2)} \circ \mu_{(i,\xi_i)} ) (\tDynkinxi)$
is the same as the one of $(i,\xi_i-2)$ in  $   \mu_{(i,\xi_i)}(\tDynkinxi)$, we can apply an induction on $k$ for the valued quiver
\begin{align*}
( \mu_{(i,\xi_i-2k)} \cdots \circ \mu_{(i,\xi_i-4)} \circ \mu_{(i,\xi_i-2)} \circ \mu_{(i,\xi_i)} ) (\tDynkinxi)  \qquad \text{ for all  } k \in \Z_{\ge 1}.
\end{align*}
Finally, our assertion comes from the definition of ${}^{s_i\xi}\tDynkin$.
\end{proof}

The following proposition is a direct consequence of Proposition~\ref{prop: one step} and the definition of
$\tDynkinxi$.

\begin{proposition} \label{prop: seq of mutations on quivers}
Let $i,j$ be sources of $Q=(\Dynkin,\xi)$. Then we have
\begin{align*}
( {}^{\bij}_{s_i\xi}\upmu \circ  {}^{\bii}_{\xi}\upmu ) ({}^{\xi}\tDynkin) \simeq    ({}^{\bii}_{s_j\xi}\upmu \circ  {}^{\bij}_{\xi}\upmu ) ({}^{\xi}\tDynkin).
\end{align*}
Thus, for any $Q$-adapted reduced expression $s_{i_1}\cdots s_{i_n}$ of the Coxeter element $\tau_Q$,
\begin{align}\label{eq: Coxeter mutation}
{}^\biQ \upmu \seteq  {}^{\bii_n}_{s_{i_n}\cdots s_{i_1}\xi}   \upmu    \circ         \cdots  \circ {}^{\bii_2}_{s_{i_1}\xi}  \upmu \circ
{}^{\bii_1}_{\xi} \upmu \text{ is well-defined.}
\end{align}
\end{proposition}

\begin{theorem} \label{thm:main1}
For Dynkin quivers $Q=(\Dynkin,\xi)$ and $Q=(\Dynkin,\xi')$, there exists a sequence of mutations $\upmu$ such that
\begin{align*}
\upmu(\tDynkinxi) \simeq   {}^{\xi'}\tDynkin \qquad  \text{as valued quivers}.
\end{align*}
In particular, we have
\begin{align*}
{}^\biQ \upmu ({}^{\xi}\tDynkin) \simeq {}^{\xi}\tDynkin \qquad  \text{as valued quivers}.
\end{align*}
\end{theorem}

\begin{proof}
This assertion follows from~\eqref{eq: differ int}, ~\eqref{eq: -2} and Proposition~\ref{prop: one step}.
\end{proof}

\subsection{Quantum cluster algebra structure on $\frakK_{q,\xi}(\g)$}
For each $s \in \Z$, we denote by ${}^{(s)}\xi$ the height function such that $\lhbs\xi_i \in \{ s,s-1\}$ for all $i \in \Dynkinf_0$ and $\lhbs Q = (\Dynkinf,\lhbs\xi)$.
For a height function $\xi$, we set
\begin{align*}
{}^{\xi}\frakK_q(\g)  \seteq \soplus_{ m \in  \lxi \calM_+ } \Zq (F_q(\um))_{\le \xi} \subset {}^{\xi}\calX_q.
\end{align*}
Note that ${}^{\xi}\frakK_q(\g) = (\cdot)_{\le \xi}\left( \frakK_{q,\xi}(\g) \right)  \simeq \frakK_{q,\xi}(\g)$.

\smallskip

For simplicity of notations, we set
\bna
\item ${}^{s}\tDynkinf \seteq  {}^{\lhbs\xi} \tDynkinf$, ${}^{s}\calX_q  \seteq  {}^{\lhbs\xi} \calX_q$, $ (-)_{\le s} \seteq (-)_{\le {}^{(s)}\xi}$,
\item $\lhs L \seteq {}^{\lhbs\xi}   L$,  $\lhs B \seteq {}^{\lhbs\xi} \tB$,
\item $\lhs\upmu \seteq {}^{\lhbs \biQ} \upmu$,  $\lhs\frakK_{q}(\g) \seteq {}^{{}^{(s)}\xi}\frakK_{q}(\g)$
and  $ \frakK_{q,s}(\g) \seteq \frakK_{q,{}^{(s)}\xi}(\g).$
\ee

From now on,
we fix $s \in \Z$ and $\tDynkinf$.
Let us denote by $\lhs \calA_q$ the quantum cluster algebra whose initial seed is
\begin{align} \label{eq: initial seed}
\lhs \calS   \seteq ( \{ v_{i,p} \seteq\ul{m^{(i)}[p,s]}   \}_{(i,p) \in {}^{s}\widetilde{\sDynkinf}_0 } , \lhs L, \lhs B ).
\end{align}

For  $n \ge 0$, let $\vip^{(n)}$ be the quantum cluster variable obtained at vertex $(i,p)$ after applying the sequence of mutations $\lhs\upmu$ $n$-times.
Then we give a quantum cluster algebra algorithm to compute $F_q(\um^{(i)}[a,b])$ for KR-monomials $m^{(i)}[a,b]$.

\begin{proposition}  {\em (cf.~\cite[Theorem 3.1]{HL16}, \cite[Proposition 6.3.1]{B21})}  \label{prop: m times mu}
For each $(i,p) \in \lhs\tDynkinf_0$ and $n \ge 0$,
\begin{align}\label{eq: trun in process}
\vip^{(n)} = \lhs F_q(\um^{(i)}[p-2n, s-2n] ) \seteq (F_q(\um^{(i)}[p-2n, s-2n] ))_{\le s}.
\end{align}
In particular, if $2n \ge \sfh$, we have
\begin{equation*}
\vip^{(n)} = F_q(\um^{(i)}[p-2n, s-2n]).
\end{equation*}
\end{proposition}

\begin{proof}
Let us apply induction on $n$ for this assertion. For $n=0$, it follows from~\eqref{eq:right negativeness of KR} in Proposition~\ref{prop: range of folded KR}.
Let $n \ge 0$ and $(i,p) \in \lhs\tDynkinf_0$. Suppose we have applied $\lhs\upmu$ $n$-times on $\lhs \calS$,  and $(n+1)$-times
on all vertices preceding $(i,p)$ in the sequence $\lhs\upmu$, and that all those previous vertices satisfy~\eqref{eq: trun in process}.

Thanks to Theorem \ref{thm:main1}, the corresponding valued quivers coincide up to a shift of spectral parameters in labeling of vertices.
Then, the argument in the proof of Proposition~\ref{prop: one step} tells us that the vertex $(i,p)$ is vertically sink or horizontally source, that is, one of the following configurations:
\begin{align} \label{eq: 2-steps p}
\scalebox{0.7}{\raisebox{14mm}{\xymatrix@!C=5mm@R=6mm{
&(j, p+\lhbs\xi_{j,i}) \ar[dr]|{ \ulcorner -\sfc_{j,i},\sfc_{i,j} \lrcorner} \\
(i,p-2) && (i,p) \ar[ll]\ar[rr]&& (i,p+2) \\
&(j', p+\lhbs\xi_{j'i})  \ar[ur]|{ \ulcorner -\sfc_{j'i},\sfc_{ij'} \lrcorner}
}}}
\quad
\scalebox{0.7}{\raisebox{14mm}{\xymatrix@!C=5mm@R=6mm{
&(j, p+\lhbs\xi_{j,i}) \ar[dr]|{ \ulcorner -\sfc_{j,i},\sfc_{i,j} \lrcorner} \\
(i,p-2) && (i,p) \ar[ll]\ar[rr]&& (i,p+2) \\
&
}}}
\quad
\scalebox{0.7}{\raisebox{12mm}{\xymatrix@!C=5mm@R=6mm{
& \\
(i,p-2) && (i,p) \ar[ll]\ar[rr]&& (i,p+2) \\
&(j', p+\lhbs\xi_{j'i})  \ar[ur]|{ \ulcorner -\sfc_{j'i},\sfc_{ij'} \lrcorner}
}}}
\end{align}
where $\lhbs\xi_{k,i} \seteq (-1)^{\delta\bl\lhbs\xi_k < \lhbs\xi_i \br}$ for $k \in \Dynkinf_0$ with $d(i,k)=1$.
In this proof, we only consider the first one in \eqref{eq: 2-steps p} since the computation below is almost identical for the other cases.

By the definition of $\lhbs\xi$, we have
$\lhbs\xi_{j,i}= \lhbs\xi_{j',i}$ for all $j,j' \in \Dynkinf_0$ with  $d(j,i)=d(j',i)=1$.
Now let us assume that $i$ is a source of $\lhbs\xi$ since the proof for the cases when $i$ is a sink of $\lhbs\xi$ is similar.
Then the quantum exchange relation has the form
\begin{equation} \label{eq:quantum exchange relation}
\begin{split}
\vip^{(n+1)}  * \vip^{(n)} = q^\al  v_{i,p+2}^{(n+1)}  \cdot   v_{i,p-2}^{(n)} + q^\be \prod_{j; \; d(j,i)=1} \bl v_{i,p-1}^{(n)} \br^{-\sfc_{j,i}}
\end{split}
\end{equation}
for some $\al,\be \in \frac{1}{2}\Z$, where
\begin{align}\label{eq: bar-inv}
q^\al  \left(v_{i,p+2}^{(n+1)} \cdot  v_{i,p-2}^{(n)} \right) * \bl \vip^{(n)} \br^{-1}  \ \  \text{and} \ \ q^\be \left( \prod_{j; \; d(j,i)=1} \bl v_{i,p-1}^{(n)} \br^{-\sfc_{j,i}} \right) * \bl \vip^{(n)} \br^{-1} \text{ are bar-invariant.}
\end{align}
Here the dot product $\cdot$ is given in \eqref{eq: dot product}.

The rest of this proof is devoted to show that the above quantum exchange relation coincides with the truncated image of the quantum folded $T$-system in Theorem~\ref{thm: quantum folded}.
For this, it suffices to assume that $s=0$ and hence $p \in \Z_{\le 0}$.
For each $(i,p) \in {}^{0}\tDynkinf_0$, we set $k  \seteq \max(u \ | \ p+2u \le 0)$.
By the induction hypothesis, we have
\begin{align*}
\uip^{(n+1)} *    \lhz F_q(\um^{(i)}_{k ,p-2n})   & = q^\al \left( \lhz  F_q(\um^{(i)}_{k -1,p-2n}) \cdot \lhz  F_q(\um^{(i)}_{k +1,p-2n-2})  \right) +q^\ga \hspace{-2ex} \prod_{j; \; d(j,i)=1}  \hspace{-2ex} \lhz  F_q(\um^{(j)}_{k,p-1-2n})^{-\sfc_{j,i}}
\end{align*}

On the other hand, the corresponding truncated image of the quantum folded $T$-system in Theorem~\ref{thm: quantum folded} is
\begin{equation}\label{eq: bar-inv2}
\begin{aligned}
\lhz F_q(\um^{(i)}_{k ,p-2n-2})    *  \lhz F_q(m^{(i)}_{k ,p-2n})  & = q^{\al'}    \left( \lhz  F_q(\um^{(i)}_{k -1,p-2n})  \cdot    \lhz  F_q(\um^{(i)}_{k +1,p-2n-2})  \right) \\
& \hspace{17ex} +  q^{\ga'}  \prod_{j; \; d(j,i)=1}   \lhz F_q(\um^{(j)}_{k,p-1-2n})^{-\sfc_{j,i}},
\end{aligned}
\end{equation}
where
$$
\ga' = \dfrac{1}{2}  \left(\tfb_{i,i}(2k-1)+ \tfb_{i,i}(2k+1)\right) \quad \text{ and }\quad \al' = \ga'+d_i.
$$

By using the dominant monomials in~\eqref{eq: bar-inv2} and bar-invariance in~\eqref{eq: bar-inv},
\begin{align*}
q^\al \bl \um^{(i)}_{k -1,p-2n} \cdot \um^{(i)}_{k +1,p-2n-2} \br * \bl \um^{(i)}_{k ,p-2n} \br^{-1} \,\,\text{ and }\,\,
q^{\ga}  \left( \prod_{j; \; d(j,i)=1} \bl \um^{(j)}_{k,p-1-2n} \br ^{-\sfc_{j,i}} \right) * \bl \um^{(i)}_{k ,p-2n} \br^{-1}
\end{align*}
are bar-invariant. Thus we have
\begin{align*}
\al &  = \frac{1}{2} \sum_{a=0}^{k-1} \left(   \sum_{b=0}^{k-2} \bl \teta_{i,i}(2(a-b)+1) - \teta_{i,i}(2(a-b)-1)  \br + \sum_{b=0}^{k} \bl \teta_{i,i}(2(a-b)+3) - \teta_{i,i}(2(a-b)+1)    \br  \right) \\
 &  = \frac{1}{2} \sum_{a=0}^{k-1} \left( \teta_{i,i}(2a+1)  - \teta_{i,i}(2a-2k+3)  + \teta_{i,i}(2a+3)  - \teta_{i,i}(2a-2k+1) \right) \\
 &  = \frac{1}{2}  \left(     \teta_{i,i}(2k+1)  + \teta_{i,i}(2k-1)      \right) + \teta_{i,i}(1) =  \al'
\end{align*}
and
\begin{align*}
\ga & = \frac{1}{2}\sum_{j; d(i,j)=1} -\sfc_{j,i} \left( \sum_{a=0}^{k-1} \left(   \sum_{b=0}^{k-1}   \teta_{i,i}(2(a-b)+2) - \teta_{i,i}(2(a-b))   \right) \right) \\
 &= \frac{1}{2} \sum_{j; d(i,j)=1} -\sfc_{j,i} \left( \sum_{a=0}^{k-1} \left(  \teta_{i,i}(2a+2) - \teta_{i,i}(2a-2k+2) \right) \right) = \frac{1}{2} \sum_{j; d(i,j)=1} -\sfc_{j,i}\teta_{i,i}(2k) \\
& \overset{\dagger}{=}  \frac{1}{2}  \left(     \teta_{i,i}(2k+1)  + \teta_{i,i}(2k-1)      \right) =\ga'.
\end{align*}
Here $\overset{\dagger}{=} $ holds by Lemma~\ref{lem: teta}.

Since $\lhz F_q(m^{(i)}_{k ,p-2n})$
is invertible in the skew-field of fractions $\lhz\frakF_q$
of the quantum torus $\lhz\calX_q$, we conclude that
\begin{align*}
\vip^{(n+1)} =  \lhz F_q(\um^{(i)}_{k ,p-2n-2}),
\end{align*}
as desired.   The second assertion follows from  Proposition~\ref{prop: range of folded KR}.
\end{proof}

Let $\lhbs\calT_q$ be the quantum torus associated with ${}^{s}L$  generated by $\vip$ for $(i,p) \in \lhs \tDynkinf_0$.
Then, $\lhbs\calT_q$ is isomorphic to ${}^{s}\calX_{q}$.
Thus, $\lhs\calA_q$ can be understood as a $\Zq$-subalgebra in $\lhbs\calT_q$.
\smallskip

By following the argument in the proof of \cite[Lemma 6.4.1]{B21}, we have the following lemma:

\begin{lemma} \label{lem: Upomega}
The assignment
\begin{align*}
\Upomega:  \vip \mapsto F_q(\um^{(i)}[p,s])
\end{align*}
extends to a well-defined injective $\Z[q^{\pm \frac{1}{2}}]$-algebra homomorphism
\begin{align*}
\Upomega:  \lhbs\calT_q \mapsto \calX_q.
\end{align*}

Moreover, the restriction of $\Upomega$ to the quantum cluster algebra $\lhs \calA_q$ has its image in the
quantum torus $\calX_q$ and the $\Z[q^{\pm\frac{1}{2}}]$-algebra homomorphisms $\Upomega$ and $(\cdot)_{\le s}$ satisfy the following
commutative diagram:
\begin{align} 
\raisebox{8mm}{\xymatrix@!C=20mm@R=6mm{
 \lhs \calA_q \ar[dr]_{\lhs\Upomega} \ar[r]^{\Upomega}& \calX_q \\
&  \lhs\calX_q \ar@{<-}[u]_{ (\cdot)_{\le s}},
}}
\end{align}
where $\lhs\Upomega$ is the map induced from the assignment  $\vip \to \um^{(i)}[p,s]$.
\end{lemma}

Let $\lhs R_q$ be the image of the quantum cluster algebra $\lhs \calA_q$ under the map $\Upomega$:
$$   R_{q,s} \seteq \Upomega(\lhs \calA_q).$$

\begin{proposition}{\em (cf.~\cite[Theorem 5.1]{HL16}, \cite[Proposition 6.4.2]{B21})} \label{prop: R=A}
We have
$$   R_{q,s} = \frakK_{q,s}(\g).$$
\end{proposition}

\begin{proof}
Let us recall $v_{i,p} \seteq\ul{m^{(i)}[p,s]}$ and $\lhbs\xi_i \in \{ s-1,  s\}$.
By Proposition \ref{prop: m times mu} and Lemma \ref{lem: Upomega}, we have
\begin{align*}
\Omega\left(v_{i,\lhbs\xi_i}^{(n)} \right) = F_q\left(  \sfX_{i,\lhbs\xi_i-2n} \right) \quad \text{for $i \in \Dynkinf_0$ and $n \in \Z_{\ge 0}$}.
\end{align*}
Since $\frakK_{q,s}(\g)$ is generated by $F_q\left( \sfX_{i,p} \right)$ for all $(i,p) \in {}^{s}\tDynkinf$ as a $\Z[q^{\pm \frac{1}{2}}]$-algebra by Theorem \ref{thm: F_q} (see also \eqref{eq: Eqm} below), we have the following inclusion:
\begin{align*}
\frakK_{q,s}(\g) \subset  R_{q,s}.
\end{align*}

Next, let us prove the reverse inclusion.  As we see in Section \ref{subsec:quantization}, there exist $\Z[q^{\pm \frac{1}{2}}]$-derivations $S_{i,q}: \calX_q \to \calX_{i,q}$
such that
\begin{equation} \label{eq: Kqg = cap Siq}
\begin{split}
\bigcap_{i \in \widetilde{\sDynkinf}_0} {\rm Ker}(S_{i,q}) = \frakK_q(\g).
\end{split}
\end{equation}

Let us prove by induction that all cluster variables $Z$ in $\lhs \calA_q$ satisfy $\Upomega(Z) \in \frakK_{q,s}(\g)$.
Let $Z$ be a quantum cluster variable in $\lhs \calA_q$. If $Z$ belongs to the initial cluster variables,  it is done by definition of $\Upomega$.
Let us assume that $Z$ does not belong to the initial cluster variables.
Then $Z$ is obtained from a finite sequence of mutations. Then we have
\begin{align*}
Z  Z_1 = q^\al M_1 + q^\be M_2,
\end{align*}
where $Z_1$,
 $M_1$ and $M_2$ are quantum cluster monomials of $\lhs\calA_q$.
 By the induction hypothesis,
\begin{align} \label{eq: induction on Upomega images}
\Upomega(Z_1),\,\,\, \Upomega(M_1),\,\,\,  \Upomega(M_2) \in \frakK_{q,s}(\g).
\end{align}
Note that $\Upomega(Z_1) \neq 0$.
By Lemma \ref{lem: Upomega}, we have
\begin{align*}
\Upomega(Z)  * \Upomega(Z_1) = q^\al \Upomega(M_1) + q^\be \Upomega(M_2).
\end{align*}

\noindent
Since $S_{i,q}$ ($i \in \Dynkinf_0$) is a $\Z[q^{\pm \frac{1}{2}}]$-linear derivation (Proposition \ref{prop: Siq is derivation}),
\begin{align*}
S_{i,q}(\Upomega(Z)  *\Upomega(Z_1) )&
= \Upomega(Z) \bcdot  S_{i,q}( \Upomega(Z_1) )  + S_{i,q}(\Upomega(Z) )  \bcdot \Upomega(Z_1)\\
&=  q^\al S_{i,q}(\Upomega(M_1)) + q^\be S_{i,q}(\Upomega(M_2)).
\end{align*}

\noindent
By the induction hypothesis and \eqref{eq: Kqg = cap Siq}, we have
\begin{align*}
S_{i,q}( \Upomega(Z_1) ) = S_{i,q}(\Upomega(M_1)) = S_{i,q}(\Upomega(M_2))=0.
\end{align*}

\noindent
Then  Lemma \ref{lem:Xiq is free} tells us that
$S_{i,q}(\Upomega(Z)) =0$,  that is, $\Upomega(Z) \in \frakK_{i,q}(\g)$ for all $i \in \Dynkinf_0$. Hence,
$\Upomega(Z)  \in \frakK_{q,s}(\g)$ due to \eqref{eq: Kqg = cap Siq} and \eqref{eq: induction on Upomega images}, as we desired.
\end{proof}

Now, we present the main result in this section.

\begin{theorem} \label{them:main2}
For each heigh function $\xi$ on $\Dynkinf$,
$\frakK_{q,\xi}(\g)$ has a quantum cluster algebra structure
whose initial quantum seed is
\begin{align}    \label{eq: scr seed S}
\seed_\xi = \bl \{  F_q(\um^{(i)}[p,\xi_i]) \}_{(i,p) \in \lxi\widetilde{\sDynkinf}_0}, \lxi L, \lxi\tB  \br.  
\end{align}
\end{theorem}

\begin{proof}
Our assertion for $\lhbs \xi$ already holds   by Proposition~\ref{prop: R=A}.
Let $j \in \lhs \tDynkinf_0$ be a source of $\lhbs \xi$. Then we have
\begin{align} \label{eq: connectedness on seeds}
{}_{\lhbs \xi}^{\hspace{0.37cm} \bij}\upmu \left(  \seed_{\lhbs \xi}   \right)
= \left( \{  F_q(\um^{(i)}[p-2\delta_{i,j},\lhbs\xi_i-2\delta_{i,j}]) \}_{(i,p) \in {}^{s}\widetilde{\sDynkinf}_0}, {}^{s_j \lhbs \xi} \hspace{-.2ex} L , {}^{s_j \lhbs \xi} \hspace{-.2ex} \tB   \right)
=  \seed_{s_j \lhbs \xi},
\end{align}
by Proposition~\ref{prop: one step} and Proposition~\ref{prop: m times mu}.
Let $Q$ (resp. $\lhbs Q$) be the Dynkin quiver of $\Dynkinf$ corresponding to $\xi$ (resp. $\lhbs \xi$).
Since
any Dynkin quivers of $\Dynkinf$ are connected by a finite sequence of reflections (up to constant on their height functions), so are $Q$ and $\lhbs Q$.
Then the quantum seed $\seed_{\lhbs \xi}$ is mutation equivalent
to $\seed_{\xi}$ by \eqref{eq: connectedness on seeds} and $\sfT_{r}$  $(r \in 2\Z)$.
Hence, it follows from Lemma \ref{lem: mutation equiv on QCA} and Proposition~\ref{prop: R=A}  that
$\frakK_{q,s}(\g) \simeq \mathscr A_{q^{1/2}}(\seed_{\lhbs \xi}) \simeq \mathscr A_{q^{1/2}}(\seed_{\xi}) \simeq \frakK_{q,\xi}(\g)
$, so $\frakK_{q,\xi}(\g)$ has a quantum cluster algebra structure.
\end{proof}

As an application of Theorem \ref{them:main2}, we obtain $q$-commutativities of $F_q(\um^{(i)}_{k,r})$ and $F_q(\um^{(j)}_{l,t})$ satisfying certain conditions as follows.

\begin{theorem} \label{thm: folded KR-commuting condition}
For a pair $(m^{(i)}_{k,r},\, m^{(j)}_{l,t} )$,
$\bl F_q(\um^{(i)}_{k,r}), F_q(\um^{(j)}_{l,t}) \br$ is a $q$-commuting pair if
\bna
\item \label{it: equation1} $r-d(i,j)  \le t \le t+2(l-1) \le  r+2(k-1)+d(i,j)$ or
\item \label{it: equation2} $t-d(i,j)  \le r  \le r+2(k-1) \le  t+2(l-1)+d(i,j)$.
\ee
In particular, $F_q(m^{(i)}_{k,r})$ q-commutes with $F_q(X_{j,p})$ if $$
 r-d(i,j)  \le p \le r+2(k-1)+d(i,j).
$$
\end{theorem}

\begin{proof}
Under the conditions ~\ref{it: equation1} and ~\ref{it: equation2}, there exists a height function $\xi$ on $\Dynkinf$ such that $\xi_i =r+2(k-1)$ and $\xi_j =t+2(l-1)$.
Then we have
\begin{align*}
F_q(\um^{(i)}_{k,r}) =F_q(\um^{(i)}[\xi_i-2(k-1),\xi_i]) \,\,\, \text{and} \,\,\,
F_q(\um^{(j)}_{l,t}) =F_q(\um^{(i)}[\xi_j-2(l-1),\xi_j])
\end{align*}
which can be viewed as initial quantum cluster variables in  $\seed_\xi$.
Thus our assertion follows from Theorem \ref{them:main2}.
\end{proof}

The conjecture below is proved in \cite{OS19} when $\g$ is of finite $AD$-type.

\begin{conjecture} \label{conj: condition for q-comm}  
For a pair $(m^{(i)}_{k,r},\, m^{(j)}_{l,t} )$,
$F_q(\um^{(i)}_{k,r})$ and $F_q(\um^{(j)}_{l,t})$ $q$-commute unless there exist $1 \le u \le \sfh$ and $0 \le s \le \min(k,l)-1$ satisfying
\begin{equation} \label{eq: commuting condition}
	|k+r-l-t| = u+ |k-l| +2s  \quad \text{ and }  \quad \tfb_{i,j}(u-1) \ne 0.
\end{equation}
\end{conjecture}

\section{Extension to $\frakK_q(\g)$} \label{sec:extension}

\noindent
In this section,
we will extend Theorem \ref{them:main2} to $\frakK_q(\g)$, that is, the quantum virtual Grothendieck ring $\frakK_q(\g)$ has also a quantum cluster algebra structure (of skew-symmetrizable type) isomorphic to its subalgebra $\frakK_{q,\xi}(\g)$ for each height function $\xi$ on $\tDynkinf$.

\subsection{Sink-source quiver}
For an integer $s \in \Z$, recall the height function  ${}^{(s)}\xi$ on $\Dynkinf$.
Now let us consider a new valued quiver $\lhs \overset{\gets}{\para}$ whose set of vertices is ${}^s \tDynkinf_0$ and the exchange matrix ${}^{s} \frakB$ is given as follows:
\begin{align}\label{eq: ssq b}
\mathfrak{b}_{(i,p),(j,t)} =
\bc
-\sfc_{i,j}  &\text{ if either (a) $t-p=1$,  $i \ne j$ and $p \equiv_4 \xi_i =s-1$},\\
&  \text{ $\hspace{4.9ex}$ or (b)  $p-t=1$,  $i \ne j$ and $p \not \equiv_4 \xi_i =s$,} \\
\sfc_{i,j}  &\text{ if either (a$'$) $p-t=1$,  $i \ne j$ and $p \equiv_4 \xi_i =s$}, \\
&  \text{ $\hspace{4.9ex}$  or (b$'$)  $t-p=1$,  $i \ne j$ and $p \not \equiv_4 \xi_i =s-1$,} \\
1 & \text{ if either (A)  $|p-t|=2$, $i=j$ and $p \equiv_4 \xi_i =s$,  }  \\
&  \text{ $\hspace{4.9ex}$ or (B)  $|p-t|=2$,  $i = j$ and $p \not \equiv_4 \xi_i =s-1$,} \\
-1 & \text{ if either (A$'$)  $|p-t|=2$, $i=j$ and $p \not\equiv_4 \xi_i =s$,  }  \\
&  \text{ $\hspace{4.9ex}$ or (B$'$)  $|p-t|=2$,  $i = j$ and $p  \equiv_4 \xi_i =s-1$,} \\
0 & \text{otherwise}.
\ec
\end{align}

\noindent
Note that $\lhs\frakB$ satisfies ~\eqref{eq: ex matrix} with the sequence
 $S = ( s_{i,p} \; | \;  s_{i,p} =d_i )$ and without frozen vertices.

\begin{example} \label{ex: SS quiver BCFG}
Here are a couple of examples of $\ssq{s}$  for non-simply-laced types: \ben
\item For $s=0$ and $\g$ of type $B_3$, $\ssq{0}$ can be depicted as follows:
$$
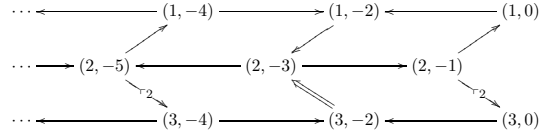
\scalebox{0.6}{
\xymatrix@!C=10mm@R=6mm{
\cdots&&(1,-4)  \ar[ll]\ar[rr] &&(1,-2) \ar[dl]  &&(1,0)  \ar[ll] \\
\cdots\ar[r]&(2,-5) \ar[ur]\ar[dr]|{\ulcorner 2}  &&(2,-3) \ar[rr] \ar[ll]&&(2,-1) \ar[ur]\ar[dr]|{\ulcorner 2} \\
\cdots&&(3,-4)  \ar[rr]\ar[ll]   &&(3,-2) \ar@{=>}[ul]&&(3,0) \ar[ll]
}
}
$$

\item For $s=0$ and $\g$ of type $G_2$, $\ssq{0}$ can be depicted as follows:
$$
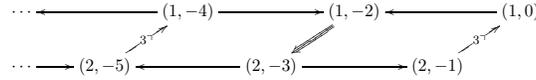
\scalebox{0.6}{
\xymatrix@!C=10mm@R=6mm{
\cdots&&(1,-4)  \ar[ll]\ar[rr] &&(1,-2) \ar@{=>}[dl]  &&(1,0)  \ar[ll] \\
\cdots\ar[r]&(2,-5) \ar[ur]|{ 3 \urcorner}  &&(2,-3) \ar@{-}[ur]\ar[rr] \ar[ll]&&(2,-1) \ar[ur]|{ 3 \urcorner}
}
}
$$
\ee
\end{example}

\begin{remark} \label{rmk: sink source}
Note that every vertex $(i,p)$ in $ \ssq{s}_0$ is either
\bnum
\item vertically sink and horizontally source, or
\item vertically source and horizontally sink.
\ee
 More precisely,  when
\bnump
\item  $\xi_i=s$ and $p \equiv_4 s$,  or  $\xi_i=(s-1)$ and $p \not\equiv_4 s-1$, $(i,p)$ satisfies (i),
\item  $\xi_i=s$ and $p \not\equiv_4 s$,  or $\xi_i=(s-1)$ and $p \equiv_4 s-1$, $(i,p)$ satisfies (ii).
\ee
\end{remark}

For each $(i,p) \in \ssq{s}_0 = {}^s\tDynkinf_0$, we assign
${}^s\fraku_{i,p} \in \frakK_q(\g)$ at $(i,p)$, which is defined by
\begin{equation} \label{eq: def of uip}
	{}^{s}\fraku_{i,p} \seteq
F_q\left( \um^{(i)} [ \lhs o_{i,p},   \lhs o_{i,p}+2  \; \lhs l_{i,p}] \right),
\end{equation}
where
\begin{equation} \label{eq: slip and soip}
\lhs l_{i,p} \seteq \left({}^{(s)}\xi_i - p \right)/2  \in \Z_{\ge 0}   \quad \text{ and } \quad   \lhs o_{i,p} \seteq  \lh{(s)}\xi_i - 2 \times \left\lfloor \frac{ \lhs l_{i,p} +\delta\bl {}^{(s)}\xi_ i = s  \br}{2} \right\rfloor.
\end{equation}

\begin{example} \label{ex: F_q quiver}
By replacing vertices $(i,p)$'s in $\ssq{s}_0$ with $\fraku_{i,p}$ in Example~\ref{ex: SS quiver BCFG}, we can obtain the following pictures:
\ben
\item \label{it: B_3 F_q} For $s=0$ and $\g$ of type $B_3$, we have
\begin{equation*}
\scalebox{0.6}{
\xymatrix@!C=12mm@R=6mm{
\cdots&&F_q(m^{(1)}[-2,2] )   \ar[ll]\ar[rr] &&F_q(m^{(1)}[-2,0] ) \ar[dl]  &&F_q(m^{(1)}[0,0] )   \ar[ll] \\
\cdots\ar[r]&F_q(m^{(2)}[-3,1] ) \ar[ur]\ar[dr]|{\ulcorner 2}  &&F_q(m^{(2)}[-1,1] ) \ar[rr] \ar[ll]&& F_q(m^{(2)}[-1,-1] )  \ar[ur]\ar[dr]|{\ulcorner 2} \\
\cdots&&F_q(m^{(3)}[-2,2] )  \ar[rr]\ar[ll]   &&F_q(m^{(3)}[-2,0] )  \ar@{=>}[ul]&&F_q(m^{(3)}[0,0] )  \ar[ll]
}
}
\end{equation*}
\item For $s=0$ and $\g$ of type $G_2$, we have
\begin{equation*}
\scalebox{0.6}{
\xymatrix@!C=12mm@R=6mm{
\cdots&&F_q(m^{(1)}[-2,2] )  \ar[ll]\ar[rr] &&F_q(m^{(1)}[-2,0] ) \ar@{=>}[dl]  &&F_q(m^{(1)}[0,0] )  \ar[ll] \\
\cdots\ar[r]&F_q(m^{(2)}[-3,1] ) \ar[ur]|{ 3 \urcorner}  &&F_q(m^{(2)}[-1,1] ) \ar@{-}[ur]\ar[rr] \ar[ll]&&F_q(m^{(2)}[-1,-1] ) \ar[ur]|{ 3 \urcorner}
}
}
\end{equation*}
\ee
\end{example}

Let us define a matrix $\lhs\Uplambda = (
\lhs\Uplambda_{(i,p),(j,t)})_{(i,p),(j,t) \in {}^s
\widetilde{\sDynkinf}_0}$ such that
$$\lhs\Uplambda_{(i,p),(j,t)} = \ucalN(  m^{(i)} [ \lhs o_{i,p},  \lhs  o_{i,p}+2   \; \lhs l_{i,p}],m^{(j)} [ \lhs o_{j,t},   \lhs o_{j,s}+2 \;  \lhs  l_{j,t}] ) .$$

\begin{theorem} \label{thm: compatible pair of whole ring}
The pair $(\lhs\Uplambda,\lhs \frakB)$ is compatible  with ${\rm diag}(2d_{i,p} \seteq 2d_i \, \mid \, (i,p) \in \ssq{s}_0)$.
\end{theorem}

\begin{proof}
Let $(i,p),(j,t) \in \ssq{s}_0$.
In this proof, we only consider the case of $\xi_j=s$ and $t \equiv_4 \xi_j$, since the other cases are similar.
Set $a_1=\lhs o_{i,p}$, $a_2=a_1+2\lhs l_{i,p}$, $b_1=\lhs o_{j,t}$ and $b_2=b_1+2\lhs l_{l,t}$.
By~\eqref{eq: ssq b}, we have
\begin{equation*}  
\begin{split}
-(\lhs\Uplambda   \lhs\frakB )_{(i,p),(j,t)} &= \delta(t \ne \xi_j) \; \lhs\Uplambda_{(i,p),(j,t+2)}+\lhs\Uplambda_{(i,p),(j,t-2)}
+ \sum_{k; \; d(j,k)=1} \sfc_{k,j} \lhs\Uplambda_{(i,p),(k,t-1)}  \allowdisplaybreaks\\
& = \delta(t \ne \xi_j) {}^{s} \ucalN(  m^{(i)}[a_1,a_2], m^{(j)}[b_1,b_2-2] ) + \ucalN(  m^{(i)}[a_1,a_2], m^{(j)}[b_1-2,b_2] ) \allowdisplaybreaks\\
& \hspace{2.5ex} + \sum_{k; \; d(j,k)=1} \sfc_{k,j}  \ucalN(  m^{(i)}[a_1,a_2], m^{(k)}[b_1-1,b_2-1] ) \\
& \overset{*}{=} \ucalN(  m^{(i)}[a_1,a_2],\,  B_{j,b_1-1} B_{j,b_1+3} \cdots  B_{j,b_2-1})
\end{split}
\end{equation*}
where $\overset{*}{=}$ holds by \eqref{eq: Bip} and \eqref{eq: ucalN for two monomials}.
Then
it follows from ~\eqref{eq;beta} in Proposition~\ref{prop: YA com} that
\begin{align*}
-(\lhs\Uplambda   \lhs\frakB  )_{(i,p),(j,t)}
&= \ucalN(  m^{(i)}[a_1,a_2],  B_{j,b_1-1} B_{j,b_1+3}  \cdots  B_{j,b_2-1} )\allowdisplaybreaks \\
&= \sum_{ x =0}^{\frac{a_2-a_1}{2}} \sum_{ y =0}^{\frac{b_2-b_1}{2}} \delta_{i,j} ( -\delta(a_1+2x- b_1-2y= -2)  + \delta(a_1+2x- b_1-2y =0)     ) 2d_i \allowdisplaybreaks\\
&=  \delta_{i,j} \sum_{ x =0}^{\frac{a_2-a_1}{2}}   ( -\delta(a_1+2x- b_1 =-2) +  \delta(a_1+2x- b_2 = 0)   ) 2d_i.
\end{align*}
If $i=j$, we have the following:
\ben
\item $[a_1,a_2]$ and $[b_1,b_2]$ are inclusive, that is, either
$[a_1,a_2] \subset [b_1,b_2] \text{ or } [b_1,b_2] \subset [a_1,a_2]$;
\item if $a_k=b_k$, then $b_l-a_l =2$ or $0$ for $\{ k,l \}=\{1,2\}$.
\ee
Thus we can conclude that
\begin{align*}
-(\lhs\Uplambda   \lhs\frakB  )_{(i,p),(j,t)}   = \delta((i,p)=(j,t)) 2d_i,
\end{align*}
as we desired.
\end{proof}

\begin{lemma}
The set $\{ {}^s\fraku_{i,p} \}_{ \in {}^s \widetilde{\sDynkinf}_0}$ forms a $q$-commuting family in $\frakK_q(\g)$.  
\end{lemma}

\begin{proof}
From Theorem~\ref{thm: folded KR-commuting condition},  our assertion easily follows.
\end{proof}

\begin{theorem} \label{thm :extension}
The family of quantum seeds
\begin{align}\label{eq: qseed frak}
\frakS_s  = (\{ {}^s\fraku_{i,p} \}_{ (i,p) \in {}^s
\widetilde{\sDynkinf}_0}, {}^{s}\Uplambda,  \lhs \frakB   )
\quad
\text{for $s \in \Z$,}
\end{align}
 gives a
quantum cluster algebra structure on $\frakK_q(\g)$.
\end{theorem}

The rest of this paper will be devoted to proving Theorem \ref{thm :extension}.
Let $\lhs\scrA_q(\g)$ be the quantum cluster algebra generated by the quantum seed $\frakS_s$.
To prove Theorem \ref{thm :extension}, we need to show that
\begin{equation} \label{eq: sAqg = Kqg}
	\lhs\scrA_q(\g) = \frakK_q(\g).
\end{equation}
Then the proof of \eqref{eq: sAqg = Kqg} is separated into two steps as follows:
\smallskip

{\it Step 1}.
For the inclusion $\lhs\scrA_q(\g) \subset \frakK_q(\g)$, we will prove the following proposition in Section \ref{subsec: step 1 in proof of extension}.
\begin{proposition} \label{prop: sub 1}
For any finite sequence $\mu$ of mutations, a cluster variable in $\mu\bl \frakS_s \br$ is contained in $\frakK_q(\g)$.
\end{proposition}

\noindent
The key observation for proving Proposition \ref{prop: sub 1} is that the mutated variables from $\frakS_s$ are understood as the ones from $\seed_{s'}$ for some $s' \in \Z$, which implies $\lhs\scrA_q(\g) \subset \frakK_q(\g)$.
\smallskip

{\it Step 2}. The opposite inclusion will be proved as the following proposition is shown in Section \ref{subsec: step 2 in proof of extension}.
\begin{proposition}\label{prop: sub 2}
For $(i,p) \in \tDynkinf_0$, there exists a finite sequence $\mu$ of mutations such that
$\mu\bl \frakS_s \br$ 
contains $F_q(\sfX_{i,p})$ as its cluster variable.
\end{proposition}
\smallskip

\noindent
Since $\frakK_q(\g)$ is generated by $F_q(\sfX_{i,p})$ for $(i,p) \in \tDynkinf_0$ by Theorem \ref{eq: Eqm} (see also \eqref{eq: Eqm} below),
the opposite inclusion for proving \eqref{eq: sAqg = Kqg} follows from Proposition \ref{prop: sub 2}.

\subsection{Proof of Theorem \ref{thm :extension}: Step 1} \label{subsec: step 1 in proof of extension}
For $k \le s$, we set
\begin{align*}
\ang{k} \seteq  \{  (i,k)  \in \lhs\tDynkinf_0    \}   \quad \text{ and } \quad  \ang{k,s} \seteq  \{  (i,p)  \in \lhs\tDynkinf_0  \ | \  k \le p \le s  \}.
\end{align*}
We understand $\ang{k,s} =\emptyset$ for $k >s$.

\begin{lemma} \label{lem: mu equiv start}
For the valued quiver $\lhs\tDynkinf$,
we have
\begin{align*}
\mu_{(i_1,s)} \circ \mu_{(i_2,s)} \circ \cdots  \circ \mu_{(i_r,s)} (\lhs\tDynkinf ) \simeq   \mu_{(j_1,s)} \circ \mu_{(j_2,s)} \circ \cdots  \circ \mu_{(j_r,s)} (\lhs\tDynkinf ),
\end{align*}
where 
$ \{  (i_t, s) \}_{1 \le t \le r} = \{  (j_t, s) \}_{1 \le t \le r}  = \ang{s}$.
Thus, $\mu_{\ang{s}} $ is well-defined  on $\lhs\tDynkinf$, that is, $\mu_{\ang{s}}(\lhs\tDynkinf)$ is uniquely determined.
\end{lemma}

\begin{proof}
Note that   {\rm (a)} each $(i_k,s) \in \ang{s}$ is vertically sink and horizontally source,   {\rm (b)} all the length 2 paths passing through $(i_k,s)$
start from $(i',s-1)$ and end at $(i_k,s-2)$ where $d(i',i_k)=1$, and
{\rm (c)} there is no arrow between $(i_k,s)$ and $(i_{k'},s)$ for $i_k \ne i_{k'}$.
\begin{align*}
\lhs\tDynkinf =
\scalebox{0.8}{\raisebox{3.5em}{    \xymatrix@!C=15mm@R=6mm{
 \cdots&(i_{k},s-4)\ar[l]\ar@{->}[dr]|{ \ulcorner -\sfc_{i_k,i'},\sfc_{i',i_k} \lrcorner} &&  (i_{k},s-2)\ar@{->}[dr]|{ \ulcorner -\sfc_{i_k,i'},\sfc_{i',i_k} \lrcorner} \ar[ll] &&  (i_{k},s) \ar[ll]  \\
\cdots&& (i',s-3) \ar@{->}[dr]|{ \ulcorner -\sfc_{i',i_{k'}},\sfc_{i_{k'},i'} \lrcorner}\ar@{->}[ur]|{ \ulcorner -\sfc_{i',i_k},\sfc_{i_k,i'} \lrcorner}\ar[ll]&& (i',s-1) \ar@{->}[dr]|{ \ulcorner -\sfc_{i',i_{k'}},\sfc_{i_{k'},i'} \lrcorner}\ar@{->}[ur]|{ \ulcorner -\sfc_{i',i_k},\sfc_{i_k,i'} \lrcorner}\ar[ll] \\
\cdots &(i_{k'},s-4)\ar[l] \ar@{->}[ur]|{ \ulcorner -\sfc_{i_{k'},i'},\sfc_{i',i_{k'}} \lrcorner}  && (i_{k'},s-2)  \ar@{->}[ur]|{ \ulcorner -\sfc_{i_{k'},i'},\sfc_{i',i_{k'}} \lrcorner} \ar[ll]&& (i_{k'},s)\ar[ll]   }}}
\end{align*}
Thus the mutation $\mu_{(i_k,s)}$ of  $\lhs\tDynkinf$ at $(i_k,s)$ does not affect the local circumstance of $(i_{k'},s)$ and the arrows between $(i_k,s-2)$ and $(i',s-1)$ for $d(i_k,i')=1$ are canceled out by the mutation $\mu_{(i_k,s)}$.
\begin{align} \label{eq: mu ang s}
\begin{split}
& \overset{\mu_{(i_k,s)}}{\longrightarrow }
\scalebox{0.8}{\raisebox{3.5em}{    \xymatrix@!C=15mm@R=6mm{
 \cdots&(i_{k},s-4)\ar[l]\ar@{->}[dr]|{ \ulcorner -\sfc_{i_k,i'},\sfc_{i',i_k} \lrcorner} &&  (i_{k},s-2)  \ar[ll] &&  (i_{k},s) \ar@{<-}[ll]  \\
\cdots&& (i',s-3) \ar@{->}[dr]|{ \ulcorner -\sfc_{i',i_{k'}},\sfc_{i_{k'},i'} \lrcorner}\ar@{->}[ur]|{ \ulcorner -\sfc_{i',i_k},\sfc_{i_k,i'} \lrcorner}\ar[ll]&& (i',s-1) \ar@{->}[dr]|{ \ulcorner -\sfc_{i',i_{k'}},\sfc_{i_{k'},i'} \lrcorner}\ar@{<-}[ur]|{ \ulcorner -\sfc_{i_k,i'},\sfc_{i',i_k} \lrcorner}\ar[ll] \\
\cdots &(i_{k'},s-4)\ar[l] \ar@{->}[ur]|{ \ulcorner -\sfc_{i_{k'},i'},\sfc_{i',i_{k'}} \lrcorner}  && (i_{k'},s-2)  \ar@{->}[ur]|{ \ulcorner -\sfc_{i_{k'},i'},\sfc_{i',i_{k'}} \lrcorner} \ar[ll]&& (i_{k'},s)\ar[ll]   }}} \\ 
& \hspace{15ex} \overset{\mu_{(i_{k'},s)}}{\longrightarrow }
\scalebox{0.8}{\raisebox{3.5em}{    \xymatrix@!C=15mm@R=6mm{
 \cdots&(i_{k},s-4)\ar[l]\ar@{->}[dr]|{ \ulcorner -\sfc_{i_k,i'},\sfc_{i',i_k} \lrcorner} &&  (i_{k},s-2)  \ar[ll] &&  (i_{k},s) \ar@{<-}[ll]  \\
\cdots&& (i',s-3) \ar@{->}[dr]|{ \ulcorner -\sfc_{i',i_{k'}},\sfc_{i_{k'},i'} \lrcorner}\ar@{->}[ur]|{ \ulcorner -\sfc_{i',i_k},\sfc_{i_k,i'} \lrcorner}\ar[ll]&& (i',s-1) \ar@{<-}[dr]|{ \ulcorner -\sfc_{i_{k'},i'},\sfc_{i',i_{k'}} \lrcorner}\ar@{<-}[ur]|{ \ulcorner -\sfc_{i_k,i'},\sfc_{i',i_k} \lrcorner}\ar[ll] \\
\cdots &(i_{k'},s-4)\ar[l] \ar@{->}[ur]|{ \ulcorner -\sfc_{i_{k'},i'},\sfc_{i',i_{k'}} \lrcorner}  && (i_{k'},s-2)  \ar[ll]&& (i_{k'},s)\ar@{<-}[ll]   }}}  
\end{split}
\end{align}
Hence the assertions follow.
\end{proof}

\begin{lemma} \label{lem: mu equiv process}
For the valued quiver $\lhs\tDynkinf$ and $k \le s$,   the valued quiver
\begin{align}\label{eq: k,s mu}
 \mu_{\ang{k,s}} (\lhs\tDynkinf)   \seteq  \mu_{\ang{k}} \circ  \mu_{\ang{k+1}} \circ \cdots \circ   \mu_{\ang{s}} (\lhs\tDynkinf)  \text{ is uniquely determined}.
\end{align}
Thus $\mu_{\ang{k}} $ is well-defined  on $\mu_{\ang{ k+1,s }}(\lhs\tDynkinf)$ and hence $\mu_{\lan k,s \ran} $  is well-defined  on $\lhs\tDynkinf$.
\end{lemma}

\begin{proof}
The assertion for $k=s$ holds by the previous lemma. As we can observe in~\eqref{eq: mu ang s},  {\rm (a)} each $(i',s-1) \in \ang{s-1}$ is vertically sink and horizontally source,
{\rm (b)}  all the length 2 paths passing through $(i',s-1)$
start from $(i,s)$ and end at $(i',s-3)$ where $d(i',i)=1$ and
{\rm (c)} there is no path between $(i',s-1)$ and $(i'',s-1)$. Thus   $\mu_{(i',s-1)} \circ \mu_{(i'',s-1)}= \mu_{(i'',s-1)} \circ \mu_{(i',s-1)}$  on $\mu_{\ang{s}}(\lhs\tDynkinf)$. Thus the assertion holds for $k=s-1$, and
$\mu_{\ang{s-1}}$ yields arrows from $(i,s)$ to $(i',s-3)$, and  hence $\mu_{\ang{s-1,s}} (\lhs\tDynkinf) $ can be depicted as follows:
\begin{align} \label{eq: mu ang s-1,s}
\scalebox{0.8}{\raisebox{4em}{    \xymatrix@!C=15mm@R=8mm{
 \cdots&(i_{k},s-4)\ar[l]\ar@{->}[dr]|{ \ulcorner -\sfc_{i_k,i'},\sfc_{i',i_k} \lrcorner} &&  (i_{k},s-2)  \ar[ll] &&  (i_{k},s) \ar@{<-}[ll]  \ar[dlll]|{ \ulcorner -\sfc_{i_k,i'},\sfc_{i',i_k} \lrcorner} \\
(i',s-5) \ar@{->}[dr]|{ \ulcorner -\sfc_{i',i_{k'}},\sfc_{i_{k'},i'} \lrcorner}\ar@{->}[ur]|{ \ulcorner -\sfc_{i',i_k},\sfc_{i_k,i'} \lrcorner}&& (i',s-3) \ar@{->}[dr]|{ \ulcorner -\sfc_{i',i_{k'}},\sfc_{i_{k'},i'} \lrcorner}\ar@{->}[ur]|{ \ulcorner -\sfc_{i',i_k},\sfc_{i_k,i'} \lrcorner}\ar[ll]&& (i',s-1) \ar[dr]|{ \ulcorner -\sfc_{i',i_{k'}},\sfc_{i_{k'},i'} \lrcorner}\ar[ur]|{ \ulcorner -\sfc_{i',i_k},\sfc_{i_k,i'}  \lrcorner}\ar@{<-}[ll] \\
\cdots &(i_{k'},s-4)\ar[l] \ar@{->}[ur]|{ \ulcorner -\sfc_{i_{k'},i'},\sfc_{i',i_{k'}} \lrcorner}  && (i_{k'},s-2)  \ar[ll]&& (i_{k'},s)\ar@{<-}[ll] \ar[ulll]|{ \ulcorner -\sfc_{i_{k'},i'},\sfc_{i',i_{k'}} \lrcorner}  }}}
\end{align}
By the same reasons for $\mu_{\ang{s}}$ and $\mu_{\ang{s-1}}$,  the sequence of mutations $\mu_{\ang{s-2}}$ is well-defined. Furthermore, by the mutation rules,  the arrows between $(i_k,(s -2)\pm 2)$ and $(i',s-3)$ for $d(i_k,i')=1$ are canceled out by the mutation $\mu_{\ang{s-1}}$. Thus $\mu_{\ang{s-2,s}}(\lhs\tDynkinf)$ can be depicted as follows:
\begin{align} \label{eq: mu ang s-2,s}
\scalebox{0.77}{\raisebox{4em}{    \xymatrix@!C=13mm@R=8mm{
\cdots&(i_{k},s-6) \ar[l]\ar@{->}[dr]|{ \ulcorner -\sfc_{i_k,i'},\sfc_{i',i_k} \lrcorner} &&(i_{k},s-4)\ar[ll] &&  (i_{k},s-2)  \ar@{<-}[ll] &&  (i_{k},s) \ar[ll]    \\
(i',s-7)\ar@{->}[dr]|{ \ulcorner -\sfc_{i',i_{k'}},\sfc_{i_{k'},i'} \lrcorner}\ar@{->}[ur]|{ \ulcorner -\sfc_{i',i_k},\sfc_{i_k,i'} \lrcorner}&&(i',s-5) \ar[ll]\ar@{->}[dr]|{ \ulcorner -\sfc_{i',i_{k'}},\sfc_{i_{k'},i'} \lrcorner}\ar@{->}[ur]|{ \ulcorner -\sfc_{i',i_k},\sfc_{i_k,i'} \lrcorner}&& (i',s-3) \ar@{<-}[dr]|{ \ulcorner -\sfc_{i_{k'},i'},\sfc_{i',i_{k'}} \lrcorner}\ar@{<-}[ur]|{ \ulcorner -\sfc_{i_k,i'},\sfc_{i',i_k} \lrcorner}\ar[ll]&& (i',s-1) \ar[dr]|{ \ulcorner -\sfc_{i',i_{k'}},\sfc_{i_{k'},i'} \lrcorner}\ar[ur]|{ \ulcorner -\sfc_{i',i_k},\sfc_{i_k,i'}  \lrcorner}\ar@{<-}[ll] \\
\cdots&(i_{k'},s-6) \ar[l] \ar@{->}[ur]|{ \ulcorner -\sfc_{i_{k'},i'},\sfc_{i',i_{k'}} \lrcorner}   &&(i_{k'},s-4)\ar[ll]   && (i_{k'},s-2)  \ar@{<-}[ll]&& (i_{k'},s)\ar[ll]    }}}
\end{align}
As in the previous cases, $\mu_{\ang{s-3}}$ is well-defined,
$\mu_{\ang{s-3}}$ yields arrows from $(i,s-2)$ to $(i',s-3 \pm 2)$ as $\mu_{\ang{s-1}}$ did, and hence   $\mu_{\ang{s-3,s}}(\lhs\tDynkinf)$ can be depicted as follows:
\begin{align} \label{eq: mu ang s-3,s}
\scalebox{0.8}{\raisebox{4em}{    \xymatrix@!C=15mm@R=8mm{
\cdots&(i_{k},s-6) \ar[l]\ar@{->}[dr]|{ \ulcorner -\sfc_{i_k,i'},\sfc_{i',i_k} \lrcorner} &&(i_{k},s-4)\ar[ll] &&  (i_{k},s-2)  \ar[dr]|{ \ulcorner -\sfc_{i_k,i'},\sfc_{i',i_k} \lrcorner}  \ar[dlll]|{ \ulcorner -\sfc_{i_k,i'},\sfc_{i',i_k} \lrcorner}  \ar@{<-}[ll] &&  (i_{k},s) \ar[ll]    \\
\cdots&&(i',s-5) \ar[ll]\ar@{->}[dr]|{ \ulcorner -\sfc_{i',i_{k'}},\sfc_{i_{k'},i'} \lrcorner}\ar@{->}[ur]|{ \ulcorner -\sfc_{i',i_k},\sfc_{i_k,i'} \lrcorner}&
& (i',s-3) \ar[dr]|{ \ulcorner -\sfc_{i',i_{k'}},\sfc_{i_{k'},i'} \lrcorner}\ar[ur]|{ \ulcorner -\sfc_{i',i_k},\sfc_{i_k,i'} \lrcorner}\ar@{<-}[ll]&& (i',s-1) \ar[dr]|{ \ulcorner -\sfc_{i',i_{k'}},\sfc_{i_{k'},i'} \lrcorner}\ar[ur]|{ \ulcorner -\sfc_{i',i_k},\sfc_{i_k,i'}  \lrcorner}\ar[ll] \\
\cdots&(i_{k'},s-6) \ar[l] \ar@{->}[ur]|{ \ulcorner -\sfc_{i_{k'},i'},\sfc_{i',i_{k'}} \lrcorner}   &&(i_{k'},s-4)\ar[ll]   && (i_{k'},s-2)  \ar[ur]|{ \ulcorner -\sfc_{i_{k'},i'},\sfc_{i',i_{k'}} \lrcorner}  \ar[ulll]|{ \ulcorner -\sfc_{i_{k'},i'},\sfc_{i',i_{k'}} \lrcorner}  \ar@{<-}[ll]&& (i_{k'},s)\ar[ll]    }}}
\end{align}
Then one can see that
\bnum
\item the full-subquiver of $\mu_{\ang{s-2,s}}(\lhs\tDynkinf)$ obtained by excluding vertices in $\ang{s}$ is isomorphic to the valued quiver $\mu_{\ang{s}}(\lhs\tDynkin)$ in~\eqref{eq: mu ang s},
\item the full-subquiver of $\mu_{\ang{s-3,s}}(\lhs\tDynkinf)$ obtained by excluding vertices in $\ang{s-1,s}$ is isomorphic to the valued quiver $\mu_{\ang{s-1,s}}(\lhs\tDynkin)$ in~\eqref{eq: mu ang s-1,s}.
\ee
Thus the induction works.
\end{proof}

\begin{remark} \label{rmk: pattern}
In the previous lemmas, we observe the following:
\ben
\item \label{it: mu=T}
Each $\mu_{(i,p)}$ in $\mu_{\ang{k,s}}$ happens when $(i,p)$ is vertically sink and horizontally source, and the arrows adjacent to $(i,p)$ are given as follows: for any $j$ with $d(i,j)=1$,
$$
\bc
 \scalebox{0.6}{\raisebox{2.5em}{     \xymatrix@!C=13mm@R=8mm{    (i,p-2) &&  (i,p) \ar[ll]\ar[rr]  && (i,p+2)   \\ & (j,p-1) \ar[ur]|{ \ulcorner -\sfc_{j,i},\sfc_{i,j} \lrcorner}   }}}  & \text{ if } p \equiv_2 s, \\
  \scalebox{0.6}{\raisebox{2.5em}{     \xymatrix@!C=13mm@R=8mm{    (i,p-2) &&  (i,p) \ar[ll]\ar[rr]  && (i,p+2)   \\ &&& (j,p+1) \ar[ul]|{ \ulcorner -\sfc_{j,i},\sfc_{i,j} \lrcorner}   }}} & \text{ if } p \not\equiv_2 s.
\ec
$$
\item \label{it: far enough} Each $\mu_{(i,p)}$ in $\mu_{\ang{k,s}}$ does not affect on the local circumstance of the vertex $(j,s)$ for $|s-p| > 2$ in the valued quiver obtained by applying the preceding mutations on $\lhs\tDynkinf$.
\ee
\end{remark}

\begin{example} \label{ex: illustration of Lemma 9.9}
By applying $\mu_{\ang{s-4}}$  on the valued quiver $\mu_{\ang{s-3,s}}(\lhs\tDynkinf)$ in  ~\eqref{eq: mu ang s-3,s}, we observe that the local circumstance of vertices in $\ang{s-1,s}$ are preserved as explained in Remark~\ref{rmk: pattern}~\ref{it: far enough}:
\begin{align} \label{eq: mu ang s-4,s}
\mu_{\ang{s-4,s}}(\lhs\tDynkinf) =
\scalebox{0.7}{\raisebox{4em}{    \xymatrix@!C=12mm@R=8mm{
\cdots&(i_{k},s-6) \ar[l] &&(i_{k},s-4)\ar@{<-}[ll] &&  (i_{k},s-2)  \ar[dr]|{ \ulcorner -\sfc_{i_k,i'},\sfc_{i',i_k} \lrcorner}  \ar[ll] &&  (i_{k},s) \ar[ll]    \\
(i',s-7) \ar@{->}[dr]|{ \ulcorner -\sfc_{i',i_{k'}},\sfc_{i_{k'},i'} \lrcorner}\ar@{->}[ur]|{ \ulcorner -\sfc_{i',i_k},\sfc_{i_k,i'} \lrcorner} &&(i',s-5) \ar[ll]\ar@{<-}[dr]|{ \ulcorner -\sfc_{i_{k'},i'},\sfc_{i',i_{k'}} \lrcorner}\ar@{<-}[ur]|{ \ulcorner -\sfc_{i_k,i'},\sfc_{i',i_k} \lrcorner}&
& (i',s-3) \ar[dr]|{ \ulcorner -\sfc_{i',i_{k'}},\sfc_{i_{k'},i'} \lrcorner}\ar[ur]|{ \ulcorner -\sfc_{i',i_k},\sfc_{i_k,i'} \lrcorner}\ar@{<-}[ll]&& (i',s-1) \ar[dr]|{ \ulcorner -\sfc_{i',i_{k'}},\sfc_{i_{k'},i'} \lrcorner}\ar[ur]|{ \ulcorner -\sfc_{i',i_k},\sfc_{i_k,i'}  \lrcorner}\ar[ll] \\
\cdots&(i_{k'},s-6) \ar[l]   &&(i_{k'},s-4)\ar@{<-}[ll]   && (i_{k'},s-2)  \ar[ur]|{ \ulcorner -\sfc_{i_{k'},i'},\sfc_{i',i_{k'}} \lrcorner}    \ar[ll]&& (i_{k'},s)\ar[ll]    }}}
\end{align}
\end{example}

For notational simplicity, let us keep the following notations:
\begin{itemize}
	\item $\Upsilon_s(\ang{k,s}) \seteq \mu_{\ang{k,s}} (\lhs\tDynkinf)$ (in~\eqref{eq: k,s mu}), \quad $\Upsilon_s \seteq \lhs\tDynkinf$,  \quad $\Theta_s \seteq \lhs\tpara$, \quad $\Omega_s \seteq\lhs\opara$,
	\item for a valued quiver $\Upgamma$, a quiver ${}^X\Upgamma$ denotes the full-subquiver of $\Upgamma$ whose vertices are in $X \subseteq \Upgamma_0$,
\end{itemize}
where $\lhs\opara$ is the quiver obtained from $\lhs\tpara$ by reversing the orientation of arrows in $\lhs\tpara$.
\noindent
By Remark~\ref{rmk: pattern}~\ref{it: far enough}, we have
\begin{align} \label{eq: far left}
{}^{\ang{-\infty,k-3}} \Upsilon_s(\ang{k,s}) \simeq  {}^{\ang{-\infty,k-3}} \Upsilon_s,
\end{align}
for any $k \le s$.
The lemma below concerns ${}^{\ang{k-3,s}} \Upsilon_s(\ang{k,s})$.

\begin{lemma} \label{lem: right part}
For $r \in \Z_{\ge 0}$, as a finite quiver, 
\bna
\item \label{it: odd far right} ${}^{\ang{s-2r+1,s}}\Upsilon_s(\ang{s-2r+1,s}) \simeq {}^{\ang{s-2r+1,s}}\Upsilon_s$.
\item \label{it: 0 equiv 2 right} ${}^{\ang{s-2r+1,s}}\Upsilon_s(\ang{s-2r,s}) \simeq {}^{\ang{s-2r+1,s}}\Upsilon_s$ and
$${}^{\ang{s-2r-3,s-2r+2}}\Upsilon_s(\ang{s-2r,s}) \simeq
\bc
 {}^{\ang{s-3,s}}\Omega_s   &\text{ if } r=0, \\
 {}^{\ang{s-5,s}}\Theta_s   &\text{ otherwise.}
\ec$$
\ee
\end{lemma}

\begin{proof}
\ref{it: odd far right}
Recall $\ang{s-2r+1,s} = \emptyset$ if $r=0$, so this case trivially holds.
The cases of $r=1$ and $r=2$ are already verified in ~\eqref{eq: mu ang s-1,s} and ~\eqref{eq: mu ang s-3,s}, respectively. 
One observes that in the general case (i.e.~$r \ge 3$), the mutation patterns in the intermediate steps are identical with \eqref{eq: mu ang s} and \eqref{eq: mu ang s-2,s} up to the shift of the second parameters.
This completes the proof of \ref{it: odd far right}.
\smallskip

\noindent
\ref{it: 0 equiv 2 right}
Let us consider the cases of $0 \le r \le 2$ precisely as follows:

\smallskip

\noindent
{\it Case 1}. $r=0$.
By~\eqref{eq: mu ang s},  ${}^{\ang{s-3,s}}\Upsilon_{s}(\ang{s})$ and   ${}^{\ang{s+1,s}}\Upsilon_{s}(\ang{s})$ are
\begin{align}  \label{eq: dec of s}
\scalebox{0.7}{\raisebox{3.5em}{    \xymatrix@!C=8mm@R=8mm{
 &&  (i_{k},s-2)   &&  (i_{k},s) \ar@{<-}[ll]  \\
 & (i',s-3) \ar@{->}[dr]|{ \ulcorner -\sfc_{i',i_{k'}},\sfc_{i_{k'},i'} \lrcorner}\ar@{->}[ur]|{ \ulcorner -\sfc_{i',i_k},\sfc_{i_k,i'} \lrcorner} && (i',s-1) \ar@{<-}[dr]|{ \ulcorner -\sfc_{i_{k'},i'},\sfc_{i',i_{k'}} \lrcorner}\ar@{<-}[ur]|{ \ulcorner -\sfc_{i_k,i'},\sfc_{i',i_k} \lrcorner}\ar[ll] \\
&& (i_{k'},s-2)  && (i_{k'},s)\ar@{<-}[ll]   }}}  \quad \text{ and } \quad \emptyset,
\end{align}
\smallskip

\noindent
{\it Case 2}. $r=1$.
By~\eqref{eq: mu ang s-2,s}, ${}^{\ang{s-5,s}}\Upsilon_{s}(\ang{s-2,s})$ and ${}^{\ang{s-1,s}}\Upsilon_{s}(\ang{s-2,s})$  are
\begin{align*}
\scalebox{0.77}{\raisebox{4em}{    \xymatrix@!C=8mm@R=8mm{
&(i_{k},s-4)  &&  (i_{k},s-2)  \ar@{<-}[ll] &&  (i_{k},s) \ar[ll]    \\
(i',s-5) \ar@{->}[dr]|{ \ulcorner -\sfc_{i',i_{k'}},\sfc_{i_{k'},i'} \lrcorner}\ar@{->}[ur]|{ \ulcorner -\sfc_{i',i_k},\sfc_{i_k,i'} \lrcorner}&& (i',s-3) \ar@{<-}[dr]|{ \ulcorner -\sfc_{i_{k'},i'},\sfc_{i',i_{k'}} \lrcorner}\ar@{<-}[ur]|{ \ulcorner -\sfc_{i_k,i'},\sfc_{i',i_k} \lrcorner}\ar[ll]&& (i',s-1) \ar[dr]|{ \ulcorner -\sfc_{i',i_{k'}},\sfc_{i_{k'},i'} \lrcorner}\ar[ur]|{ \ulcorner -\sfc_{i',i_k},\sfc_{i_k,i'}  \lrcorner}\ar@{<-}[ll] \\
&(i_{k'},s-4)   && (i_{k'},s-2)  \ar@{<-}[ll]&& (i_{k'},s)\ar[ll]    }}}
 \quad \text{ and } \quad
\scalebox{0.77}{\raisebox{4em}{    \xymatrix@!C=12mm@R=8mm{
&  (i_{k},s)     \\
 (i',s-1) \ar[dr]|{ \ulcorner -\sfc_{i',i_{k'}},\sfc_{i_{k'},i'} \lrcorner}\ar[ur]|{ \ulcorner -\sfc_{i',i_k},\sfc_{i_k,i'}  \lrcorner}  \\
& (i_{k'},s)    }}}
\end{align*}
\smallskip

\noindent
{\it Case 3}. $r=2$.
By~\eqref{eq: mu ang s-4,s}, ${}^{\ang{s-7,s-2}}\Upsilon_{s}(\ang{s-4,s})$ and  ${}^{\ang{s-3,s}}\Upsilon_{s}(\ang{s-4,s})$  are
\begin{align}  \label{eq: dec of s-4,s}
\scalebox{0.7}{\raisebox{4em}{    \xymatrix@!C=8mm@R=8mm{
 &(i_{k},s-6)  &&(i_{k},s-4)\ar@{<-}[ll] &&  (i_{k},s-2) \ar[ll]   \\
(i',s-7) \ar@{->}[dr]|{ \ulcorner -\sfc_{i',i_{k'}},\sfc_{i_{k'},i'} \lrcorner}\ar@{->}[ur]|{ \ulcorner -\sfc_{i',i_k},\sfc_{i_k,i'} \lrcorner} &&(i',s-5) \ar[ll]\ar@{<-}[dr]|{ \ulcorner -\sfc_{i_{k'},i'},\sfc_{i',i_{k'}} \lrcorner}\ar@{<-}[ur]|{ \ulcorner -\sfc_{i_k,i'},\sfc_{i',i_k} \lrcorner}\ar[rr]& & (i',s-3)
\ar[dr]|{ \ulcorner -\sfc_{i',i_{k'}},\sfc_{i_{k'},i'} \lrcorner}\ar[ur]|{ \ulcorner -\sfc_{i',i_k},\sfc_{i_k,i'} \lrcorner} \\
 &(i_{k'},s-6)   &&(i_{k'},s-4)\ar@{<-}[ll]   && (i_{k'},s-2)   \ar[ll]   }}}
 \quad \text{ and } \quad
\scalebox{0.7}{\raisebox{4em}{    \xymatrix@!C=12mm@R=8mm{
 &  (i_{k},s-2)  \ar[dr]|{ \ulcorner -\sfc_{i_k,i'},\sfc_{i',i_k} \lrcorner}  &&  (i_{k},s) \ar[ll]    \\
 (i',s-3) \ar[dr]|{ \ulcorner -\sfc_{i',i_{k'}},\sfc_{i_{k'},i'} \lrcorner}\ar[ur]|{ \ulcorner -\sfc_{i',i_k},\sfc_{i_k,i'} \lrcorner} && (i',s-1) \ar[dr]|{ \ulcorner -\sfc_{i',i_{k'}},\sfc_{i_{k'},i'} \lrcorner}\ar[ur]|{ \ulcorner -\sfc_{i',i_k},\sfc_{i_k,i'}  \lrcorner}\ar[ll] \\
 & (i_{k'},s-2)  \ar[ur]|{ \ulcorner -\sfc_{i_{k'},i'},\sfc_{i',i_{k'}} \lrcorner}    && (i_{k'},s)\ar[ll]    }}}
\end{align}
\noindent

One may further observe from {\it Case 1}--{\it Case 3} that
\begin{itemize}
	\item ${}^{\ang{s-2r+1,s}}\Upsilon_s(\ang{s-2r,s}) \simeq {}^{\ang{s-2r+1,s}}\Upsilon_s$ for $r \ge 1$ (by similar argument as in \ref{it: odd far right}),
	\item ${}^{\ang{s-2r-3,s-2r+2}}\Upsilon_s(\ang{s-2r,s})$ stabilizes for $r \ge 1$ up to the shift of the second parameters, which is isomorphic to ${}^{\ang{s-5,s}}\Theta_s$ as a finite quiver, where ${}^{\ang{s-3,s}}\Upsilon_s(\ang{s}) \simeq {}^{\ang{s-3,s}}\Omega_s$.
\end{itemize}
Hence we complete the proof of \ref{it: 0 equiv 2 right}.
\end{proof}

For $r \in \Z_{\ge 1}$, we define
\begin{equation*}
\mu_{\rang{s-2r,s}}  \seteq \bc
\mu_{\ang{s}}  \circ\mu_{\ang{s-4,s}}  \circ  \cdots   \circ \mu_{\ang{s-2r+4,s}}  \circ  \mu_{\ang{s-2r,s}}  & \text{ if } r \equiv_2 0 , \\
\mu_{\ang{s-2,s}}  \circ\mu_{\ang{s-6,s}}  \circ  \cdots   \circ \mu_{\ang{s-2r+4,s}}   \circ \mu_{\ang{s-2r,s}}  & \text{ if } r \equiv_2 1.
\ec
\end{equation*}
By Lemma~\ref{lem: right part}~\ref{it: 0 equiv 2 right},
${}^{\ang{s-2r+1,s}}\Upsilon_s(\ang{s-2r,s}) \simeq {}^{\ang{s-2r+1,s}}\Upsilon_s$.
By Lemma~\ref{lem: mu equiv process} and Remark~\ref{rmk: pattern}~\ref{it: far enough},
$\mu_{\ang{s-2r+4,s}}$ is well-defined on ${}^{\ang{s-2r+1,s}}\Upsilon_s(\ang{s-2r,s})$. 
Thus it makes sense to define
\begin{equation*}
\Upsilon_s(\rang{s-2r,s}) \seteq \mu_{\rang{s-2r,s}}(\lhs\tDynkinf).
\end{equation*}
Then we have a generalization of Lemma \ref{lem: right part}.

\begin{proposition} \label{prop: right mut}
For $r \in \Z_{\ge 0}$,
we have
\begin{equation*}
{}^{\ang{s-2r-3,s}}\Upsilon_s(\rang{s-2r,s})
\simeq
\bc
{}^{\ang{s-2r-3,s}}\Omega_s & \text{ if } r \equiv_2 0, \\
{}^{\ang{s-2r-3,s}}\Theta_s & \text{ if } r \equiv_2 1.
\ec
\end{equation*}
\end{proposition}

\begin{proof}
We first consider the case $0 \le r \le 2$, and then the general case $r\ge 3$.
\smallskip

\noindent
{\it Case 1}. $0 \le r \le 2$.
The assertion for $r=0$ and $r=1$ is shown by~\eqref{eq: mu ang s} and \eqref{eq: mu ang s-2,s}, respectively.
Let us consider the case  $r=2$.
By \eqref{eq: mu ang s-4,s},
we may consider
${}^{\ang{s-7,s}} \Upsilon_s(\ang{s-4,s})$
by separating it into two parts ${}^{\ang{s-7,s-2}} \Upsilon_s(\ang{s-4,s})$ and ${}^{\ang{s-2,s}} \Upsilon_s(\ang{s-4,s})$,
where each one is shown in \eqref{eq: dec of s-4,s}.
Then ${}^{\ang{s-7,s}} \Upsilon_s(\rang{s-4,s})$ is understood as a concatenation of ${}^{\ang{s-7,s-2}} \Upsilon_s(\ang{s-4,s})$
and $\mu_{\ang{s}} \bl  {}^{\ang{s-2,s}} \Upsilon_s(\ang{s-4,s})\br$ due to Remark~\ref{rmk: pattern}~\ref{it: far enough}, where the common vertices are overlapped.
Since  $\mu_{\ang{s}} \bl  {}^{\ang{s-2,s}} \Upsilon_s(\ang{s-4,s})\br$ is isomorphic to the valued quiver in~\eqref{eq: dec of s}, the assertion for $r=2$ is proved.

\smallskip

\noindent
{\it Case 2}. $r \ge 3$.
The proof idea in this case is identical with {\it Case 1}, that is,
by using the same argument as in {\it Case 1}, we observe that the finite valued quiver
$$\Upgamma_1 \seteq {}^{\ang{s-2r-3,s-2r+6}} \left( \mu_{\ang{s-2r+4,s}}  \circ  \mu_{\ang{s-2r,s}} (\lhs\tDynkinf) \right) $$ is a concatenation of
\begin{align*}
\scalebox{0.6}{\raisebox{4em}{    \xymatrix@!C=8mm@R=8mm{
 &(i_{k},s-2r-2)  &&(i_{k},s-2r)\ar@{<-}[ll] &&  (i_{k},s-2r+2) \ar[ll]   \\
(i',s-2r-3) \ar@{->}[dr]|{ \ulcorner -\sfc_{i',i_{k'}},\sfc_{i_{k'},i'} \lrcorner}\ar@{->}[ur]|{ \ulcorner -\sfc_{i',i_k},\sfc_{i_k,i'} \lrcorner} &&(i',s-2r-1) \ar[ll]\ar@{<-}[dr]|{ \ulcorner -\sfc_{i_{k'},i'},\sfc_{i',i_{k'}} \lrcorner}\ar@{<-}[ur]|{ \ulcorner -\sfc_{i_k,i'},\sfc_{i',i_k} \lrcorner}\ar[rr]& & (i',s-2r+1)
\ar[dr]|{ \ulcorner -\sfc_{i',i_{k'}},\sfc_{i_{k'},i'} \lrcorner}\ar[ur]|{ \ulcorner -\sfc_{i',i_k},\sfc_{i_k,i'} \lrcorner} \\
 &(i_{k'},s-2r-2)   &&(i_{k'},s-2r)\ar@{<-}[ll]   && (i_{k'},s-2r+2)   \ar[ll]   }}}
\text{ and  } \
\scalebox{0.6}{\raisebox{4em}{    \xymatrix@!C=8mm@R=8mm{
 &(i_{k},s-2r+2)  &&(i_{k},s-2r+4)\ar@{<-}[ll] &&  (i_{k},s-2r+6) \ar[ll]   \\
(i',s-2r+1) \ar@{->}[dr]|{ \ulcorner -\sfc_{i',i_{k'}},\sfc_{i_{k'},i'} \lrcorner}\ar@{->}[ur]|{ \ulcorner -\sfc_{i',i_k},\sfc_{i_k,i'} \lrcorner} &&(i',s-2r+3) \ar[ll]\ar@{<-}[dr]|{ \ulcorner -\sfc_{i_{k'},i'},\sfc_{i',i_{k'}} \lrcorner}\ar@{<-}[ur]|{ \ulcorner -\sfc_{i_k,i'},\sfc_{i',i_k} \lrcorner}\ar[rr]& & (i',s-2r+5)
\ar[dr]|{ \ulcorner -\sfc_{i',i_{k'}},\sfc_{i_{k'},i'} \lrcorner}\ar[ur]|{ \ulcorner -\sfc_{i',i_k},\sfc_{i_k,i'} \lrcorner} \\
 &(i_{k'},s-2r+2)   &&(i_{k'},s-2r+4)\ar@{<-}[ll]   && (i_{k'},s-2r+6)   \ar[ll]   }}}
\end{align*}
where we regard the common vertices to be overlapped in the concatenation.
Since
\begin{equation*}
{}^{\ang{s-2r+5,s}} \left( \mu_{\ang{s-2r+4,s}}  \circ  \mu_{\ang{s-2r,s}} (\lhs\tDynkinf) \right)  \simeq    {}^{\ang{s-2r+5,s}}\Upsilon_s \quad \text{by Lemma~\ref{lem: right part}~\ref{it: 0 equiv 2 right}},
\end{equation*}
and $\mu_{\ang{s-2r+8,s}}$ does not contribute to $\Upgamma_1$, we complete the proof
by applying the same argument to ${}^{\ang{s-2r+5,s}}\Upsilon_s$ as in {\it Case 1}. 
\end{proof}

\smallskip

Let us write $\mu$ in Proposition~\ref{prop: sub 1} as
\begin{align} \label{eq: finite mu}
\mu = \mu_{(i_l,p_l)} \circ \mu_{(i_{l-1},p_{l-1})} \circ \cdots \circ \mu_{(i_1,p_1)}.
\end{align}
Take $t \in \Z$ such that $t \ll \min( p_k  \ | \  1 \le k \le l )$ and  $s-t \equiv_4 2$.
By our choice of $t$, it follows from Proposition~\ref{prop: right mut} that
\begin{align*}
 {}^{\ang{t-3,s}}\Upsilon_{s}(\rang{t,s}) \simeq    {}^{\ang{t-3,s}}\Theta_s \quad \text{  as a valued quiver,}
\end{align*}
where 
\begin{align} \label{eq: u}
 s-t = 4u + 2 \qquad \text{ for some $u \in \Z_{\ge 0}$.}
\end{align}
Recall the quantum seeds
\begin{align}
& \seed_s = \bl \{ {}^s\frakv_{i,p} \seteq  F_q(\um^{(i)}[p, {}^{(s)}\xi_i]) \}_{(i,p) \in \lhs\widetilde{\sDynkinf}_0}, \lhs L, \lhs\tB   \br \text{ associated to ${}^{(s)}\xi$ in~\eqref{eq: scr seed S}} \label{eq: seed_s}, \\
& \frakS_s  = \left(\left\{ {}^s\fraku_{i,p} = F_q\bl \um^{(i)} [ \lhs o_{i,p},   \lhs o_{i,p}+2   \; \lhs l_{i,p}] \br \right\}_{ (i,p) \in {}^s
\widetilde{\sDynkinf}_0}, {}^{s}\Uplambda,{}^{s}\frakB \right) \text{ in ~\eqref{eq: qseed frak}.} \nonumber
\end{align}
 
\begin{proposition} \label{prop: mu=T-sys}
Every mutation $\mu_{(i,p)}$ in $\mu_{\rang{t,s}}$ on the cluster  $ \{ {}^s\frakv_{i,p} \}$ corresponds to the quantum folded T-system in Theorem~\ref{thm: quantum folded}.  Furthermore, each mutation $\mu_{(i,p)}$ of the quantum cluster variable sitting at $(i,p)$ corresponds to $\sfT_{-2}$.
\end{proposition}

\begin{proof}
First, let us consider a mutation $\mu_{(i,p)}$ in $\mu_{\ang{t,s}}$.
When $(i,p)=(i,s)$ (i.e.~one of the vertices located in the right-most of $\lhs\tDynkinf_0$), the local circumstance of $(i,s)$ described in Remark ~\ref{rmk: pattern}~\ref{it: mu=T} tells us that the quantum exchange relation is given by
\begin{align*}
 \mu_{(i,s)}\big(F_q(\sfX_{i,s}) \big) *F_q(\sfX_{i,s}) &= q^{\al(i,1)} F_q(\um^{(i)}[s-2,s]) +  q^{\ga(i,1)} \prod_{j; \; d_{i,j}=1} F_q(\sfX_{j,s-1})^{-\sfc_{j,i}},
\end{align*}
where $q^{\al(i,1)}$ and $q^{\ga(i,1)}$ are determined to be bar-invariant as in the sense of \eqref{eq:quantum exchange relation}.  
Consequently, it corresponds to the quantum folded T-system in Theorem~\ref{thm: quantum folded} and hence $\mu_{(i,s)}(F_q(\sfX_{i,s})) = F_q(\sfX_{i,s-2})$ as we desired.
Note that another mutation at $(i',s)$ does not affect the mutation at $(i,s)$ as shown in Lemma~\ref{lem: mu equiv start}.
\smallskip

Second, let us consider a mutation at $(j,s-1)$, which appears later than any $(i,s)$ in $\mu_{\ang{t,s}}$.
Let us keep in mind that the cluster variable located at $(i',s)$ is already mutated by former mutations, which is $F_q(\sfX_{i',s-2})$.
Then the quantum exchange relation is given as follows (recall Remark ~\ref{rmk: pattern}~\ref{it: mu=T}):
\begin{align*}
 \mu_{(j,s-1)}\big(F_q(\sfX_{j,s-1})\big)  * F_q(\sfX_{j,s-1})  &= q^{\al(j,1)} F_q(\um^{(j)}[s-3,s-1]) +  q^{\ga(j,1)} \prod_{i; \; d_{j,i}=1} F_q(\sfX_{i,s-2})^{-\sfc_{i,j}},
\end{align*}
which coincides with the quantum folded T-system in Theorem~\ref{thm: quantum folded}.
Hence $\mu_{(j,s-1)}(F_q(\sfX_{j,s-1}))  = F_q(\sfX_{j,s-3})$, as we desired.
\smallskip

Finally, by using this argument and the local circumstance of $(k,p)$ in the order for applying $\mu_{(k,p)}$, described in Remark ~\ref{rmk: pattern}~\ref{it: mu=T}, one can conclude that each mutation $\mu_{(i,p)}$ in $\mu_{\ang{t,s}}$ corresponds to shifting the second parameters of cluster variables by $-2$.
The assertion for mutations in $\mu_{\ang{t+4r,s}}$ $(r \ge 1)$ follows from Lemma~\ref{lem: right part}~\ref{it: 0 equiv 2 right}, Remark~\ref{rmk: pattern}~\ref{it: far enough} and the argument for mutations in $\mu_{\ang{t,s}}$.
\end{proof}

Recall $u \in \Z_{\ge0}$ in~\eqref{eq: u} depending on $\ang{t,s}$. For $(j,a) \in \lhs\tDynkinf_0$ with $t \le a \le s$, we remark that
\begin{itemize}
	\item[(A)]  there exists $0 \le e \le u$   
such that  $ s-4e-2   \le  a< \min(s+1,s-4e+2 )$,  equivalently
$$ a \in \{ s-4e-2,s-4e-1,s-4e,s-4e+1\},$$
	\item[(B)] $\lh{(s)}\xi_j=s$ if $a= s-4e-2$ or $s-4e$, and $\lh{(s)}\xi_j=s-1$, otherwise,
	\item[(C)] since $\mu_{(j,a)}$ appears $(u+1-e)$-times  in $\mu_{\rang{t,s}}$ and ${}^s\frakv_{j,a} =   F_q(  \um^{(j)}[ a,s -\delta( \lh{(s)} \xi_j \ne s)])$, it follows from Proposition~\ref{prop: mu=T-sys} that
\begin{equation*} 
\bl \mu_{  \rang{t,s} } (  \{ {}^s\frakv_{k,p} \} ) \br_{ (j,a) } =   F_q(  \um^{(j)}[ a'+2e,s'+2e -\delta(a \not\equiv_2 s )]).
\end{equation*}
\end{itemize}

\begin{proposition} \label{prop: cluster mut}
For $(j,a) \in \lhs\tDynkinf_0$ with $t \le a \le s$,
\begin{equation*} 
\bl \mu_{  \rang{t,s} } (  \{ {}^s\frakv_{k,p} \} ) \br_{ (j,a) } =
F_q\bl \um^{(j)} [ \lh{s'} o_{j,a'},   \lh{s'}  o_{j,a'}+2 \lh{s'}  l_{j,a'}] ),
\end{equation*}
where $s'=s-2(u+1)$ and $a'=a-2(u+1)$ for $u \in \Z$ in {\rm (C)}.
\end{proposition}
\begin{proof}
We have
\begin{align*}
\bl \lh{s'} o_{j,a'} , \lh{s'} l_{j,a'} \br =\bc
  \bl s'-2( e+1) ,2e+1 \br  & \text{ if } s-a =4e+2, \\
  \bl s'-2e ,2e\br  & \text{ if } s-a =4e, \\
  \bl  s'-2e ,2e\br  & \text{ if }( s-1)-a =4e, \\
  \bl s'-2(e-1),2e-1\br & \text{ if } ( s-1)-a =4e-2,
\ec
\end{align*}
where the integers on the left-hand side are defined in \eqref{eq: slip and soip}.
Then one can easily check that
\begin{align*}
\lh{s'} o_{j,a'} = a'+2e \quad \text{ and } \quad  \lh{s'} o_{j,a'} +2 \; \lh{s'} l_{j,a'} = s'+2e -\delta(a \not\equiv_2 s ),
\end{align*}
which implies our assertion.
\end{proof}

Now, we are ready to prove Proposition~\ref{prop: sub 1}.

\begin{proof}[Proof of Proposition~\ref{prop: sub 1}]
Write $\mu$ in Proposition~\ref{prop: sub 1} as in~\eqref{eq: finite mu}.
Let us set
\begin{align*}
Z \seteq  (\mu( \{   \lhs\fraku_{k,p}  \} ) )_{(i_l,p_l)}.
\end{align*} 
By Proposition~\ref{prop: right mut}  and Proposition~\ref{prop: cluster mut},
we have
\begin{align*}
(\{   \lh{s'}\frakv_{k',p'}  \})_{\ang{t',s'}} =(\{   \lh{s}\fraku_{k,p}  \})_{\ang{t,s}},
\end{align*}
that is,
$Z$ can be understood as a mutated variable from $\{   \lh{s'}\frakv_{k',p'}  \}$ as follows:
\begin{align*}
Z = \bl  \mu \circ \mu_{\rang{t',s'}} (  \{   \lh{s'}\frakv_{k',p'}  \}  ) \br_{ (i_l, p_l')},
\end{align*}
Here  $t'=t+2(u+1)$, $s'=s+2(u+1)$ and $p_l'=p_l+2(u+1)$.  
Since
%
\begin{align*}
\seed_{s'} = \bl \{  F_q(\um^{(i)}[p,\lh{(s')}\xi_i]) \}_{(i,p) \in \lh{s'}\widetilde{\sDynkinf}_0}, \lh{s'} L, \lh{s'}\tB  \br
\end{align*}
is an initial quantum seed of the quantum cluster algebra $\frakK_{q,s'}(\g) \subset \frakK_q(\g)$,
$Z$ is contained in $\frakK_q(\g)$, which completes the proof.
\end{proof}

\subsection{Proof of Theorem \ref{thm :extension}: Step 2} \label{subsec: step 2 in proof of extension}
For $k \le s$, we set
\begin{align*}
\ang{k}^- &\seteq  \{  (i,k)  \in \ssq{s}_0 \ |  \  \text{ $(i,k)$ is vertically sink and horizontally source in $\ssq{s}_0$ }     \} , \\
\ang{k}^+ &\seteq  \{  (i,k)  \in \ssq{s}_0 \ |  \  \text{ $(i,k)$ is vertically source and horizontally sink in $\ssq{s}_0$ }     \}.
\end{align*}
For $k \in  \Z_{\le s} \sqcup \{ -\infty \}$,
\begin{align*}
\ang{k,s}^- &\seteq  \bigsqcup_{k  \le t \le s} \ang{t}^-   \quad \text{and} \quad \ang{k ,s}^+ \seteq  \bigsqcup_{k  \le t \le s} \ang{t}^+.
\end{align*}
If $k>s$, then we understand those sets as empty set.   Note that there is no arrows between vertices in $\ang{k,s}^\pm$ for any $k \in  \Z_{\le s} \sqcup \{ -\infty \}$.

\begin{lemma} \label{lem: mu equiv start pm T} 
For $ \{  (i_t, p_t) \}_{1 \le t \le r} = \{  (j_t, q_t) \}_{1 \le t \le r}  = \ang{k,s}^\pm$, as a valued quiver, 
\begin{align*}
\mu_{(i_1,p_1)} \circ \mu_{(i_2,p_2)} \circ \cdots  \circ \mu_{(i_r,p_r)} (\ssq{s}) \simeq   \mu_{(j_1,q_1)} \circ \mu_{(j_2,q_2)} \circ \cdots  \circ \mu_{(j_r,q_r)} (\ssq{s}),
\end{align*}
that is, $\mu_{\ang{k,s}^\pm}(\ssq{s})$ is uniquely determined. 
\end{lemma}

We remark that an analog of Lemma \ref{lem: mu equiv start pm T} by replacing $\ssq{s}$ with $\rssq{s}$ also holds.

\begin{proof} 
In this proof, we only consider the case of $\ang{k,s}^+$ since the proof of $\ang{k,s}^-$ is similar.
Let $(i,p), (j,s) \in \ang{k,s}^+$ such that $(i,p) \ne (j,s)$. 
The neighborhood of $(i,p)$ on the valued quiver $\ssq{s}$ is depicted as follows:
\begin{align*}
\ssq{s} = \scalebox{0.6}{\raisebox{4em}{    \xymatrix@!C=8mm@R=8mm{
\cdots&&(j,p-3)\ar[rr]\ar[ll] &&(j,p-1)\ar[dl]|{ \ulcorner -\sfc_{j,i},\sfc_{i,j}}    &&(j,p+1) \ar[rr]\ar[ll] &&  (j,p+3)   \ar[dl]|{ \ulcorner -\sfc_{j,i},\sfc_{i,j}} &\cdots\ar[l] \\
\cdots\ar[r]&(i,p-4)\ar[dr]|{ \ulcorner -\sfc_{i,j'},\sfc_{j,i'} \lrcorner}\ar[ur]|{ \ulcorner -\sfc_{i,j},\sfc_{j,i} \lrcorner}    &&(i,p-2) \ar[ll]\ar[rr]  & & (i,p) \ar[dr]|{ \ulcorner -\sfc_{i,j'},\sfc_{j,i'} \lrcorner}\ar[ur]|{ \ulcorner -\sfc_{i,j},\sfc_{j,i} \lrcorner}  & & (i,p+2) \ar[ll]  \ar[rr]
&&\cdots \\
\cdots  &&(j',p-3) \ar[rr]\ar[ll]   &&(j',p-1) \ar[ul]|{ \ulcorner -\sfc_{j',i},\sfc_{i,j'}}   &&(j',p+1) \ar[rr]  \ar[ll]    && (j',p+3)   \ar[ul]|{ \ulcorner -\sfc_{j',i},\sfc_{i,j'}}   &\cdots \ar[l]  }}}
\end{align*}
By Algorithm \ref{Alg. mutation}, we have
\begin{align*}
{\scriptsize \mu_{(i,p)}(\ssq{s})}&=
\scalebox{0.6}{\raisebox{4em}{    \xymatrix@!C=9mm@R=8mm{
\cdots&&(j,p-3)\ar[rr]\ar[ll] &&(j,p-1)\ar[dl]|{ \ulcorner -\sfc_{j,i},\sfc_{i,j}\lrcorner}    &&(j,p+1) \ar[rr]\ar[ll] &&  (j,p+3)   \ar[dl]|{ \ulcorner -\sfc_{j,i},\sfc_{i,j}\lrcorner} &\cdots\ar[l] \\
\cdots\ar[r]&(i,p-4)\ar[dr]|{ \ulcorner -\sfc_{i,j'},\sfc_{j,i'} \lrcorner}\ar[ur]|{ \ulcorner -\sfc_{i,j},\sfc_{j,i} \lrcorner}
&&(i,p-2) \ar[ll]  \ar[urrr]|{ \ulcorner -\sfc_{i,j},\sfc_{j,i} \lrcorner} \ar[drrr]|{ \ulcorner -\sfc_{i,j'},\sfc_{j',i} \lrcorner}
& & (i,p) \ar[ll]\ar[rr]  \ar@{<-}[dr]|{ \ulcorner -\sfc_{j,i'},\sfc_{i,j'} \lrcorner}\ar@{<-}[ur]|{ \ulcorner -\sfc_{j,i},\sfc_{i,j} \lrcorner}
& & (i,p+2)  \ar[rr] \ar[ul]|{ \ulcorner -\sfc_{i,j},\sfc_{j,i} \lrcorner}\ar[dl]|{ \ulcorner -\sfc_{i,j'},\sfc_{j',i} \lrcorner}
&&\cdots \\
\cdots  &&(j',p-3) \ar[rr]\ar[ll]   &&(j',p-1) \ar[ul]|{ \ulcorner -\sfc_{j',i},\sfc_{i,j'} \lrcorner}   &&(j',p+1) \ar[rr]  \ar[ll]    && (j',p+3)   \ar[ul]|{ \ulcorner -\sfc_{j',i},\sfc_{i,j'}\lrcorner}   &\cdots \ar[l]  }}}
\end{align*}
Here one can observe that
\begin{itemize}
	\item $\mu_{(i,p)}(\ssq{s})$ has arrows between $(i,p  \pm 2)$ and $(j,p+1)$ for $d(i,j)=1$, where $(i,p  \pm 2) ,(j,p+1) \in \ang{k,s}^-$,
	\item the arrows adjacent to $(j,p-1)$ and $(j,p+3)$ are not changed by $\mu_{(i,p)}$.
\end{itemize}
Hence, for $(x,y) \in \{ (j,p-1), (j,p+3) \ | \ d(i,j)=1\}$, the mutation $\mu_{(x,y)}\large(  \mu_{(i,p)}(\ssq{s})\large)$
yields arrows between $(x,y \pm 2)$ and $(k,y+1)$ for $d(x,k)=1$, one of which disappears due to an arrow from $\mu_{(i,p)}(\ssq{s})$.
For instance, 
\begin{align*}
{\scriptsize \mu_{(j',p-1)}( \mu_{(i,p)}(\ssq{s}))}&=
\scalebox{0.6}{\raisebox{4em}{    \xymatrix@!C=9mm@R=8mm{
\cdots&&(j,p-3)\ar[rr]\ar[ll] &&(j,p-1)\ar[dl]|{ \ulcorner -\sfc_{j,i},\sfc_{i,j}\lrcorner}    &&(j,p+1) \ar[rr]\ar[ll] &&  (j,p+3)   \ar[dl]|{ \ulcorner -\sfc_{j,i},\sfc_{i,j} \lrcorner} &\cdots\ar[l] \\
\cdots\ar[r]&(i,p-4)\ar[dr]|{ \ulcorner -\sfc_{i,j'},\sfc_{j,i'} \lrcorner}\ar[ur]|{ \ulcorner -\sfc_{i,j},\sfc_{j,i} \lrcorner}
&&(i,p-2) \ar[ll]  \ar[urrr]|{ \ulcorner -\sfc_{i,j},\sfc_{j,i} \lrcorner}
& & (i,p) \ar[ll]\ar[rr]  \ar@{<-}[dr]|{ \ulcorner -\sfc_{j,i'},\sfc_{i,j'} \lrcorner}\ar@{<-}[ur]|{ \ulcorner -\sfc_{j,i},\sfc_{i,j} \lrcorner}
& & (i,p+2)  \ar[rr] \ar[ul]|{ \ulcorner -\sfc_{i,j},\sfc_{j,i} \lrcorner}\ar[dl]|{ \ulcorner -\sfc_{i,j'},\sfc_{j',i} \lrcorner}
&&\cdots \\
\cdots  &&(j',p-3) \ar[ll]   \ar[ur]|{ \ulcorner -\sfc_{j',i},\sfc_{i,j'} \lrcorner}    &&(j',p-1) \ar[ll]\ar[rr] \ar@{<-}[ul]|{ \ulcorner -\sfc_{i,j'},\sfc_{j',i}  \lrcorner }   &&(j',p+1) \ar[rr]      && (j',p+3)   \ar[ul]|{ \ulcorner -\sfc_{j',i},\sfc_{i,j'} \lrcorner}   &\cdots \ar[l]  }}}
\end{align*}
Here the arrow from $(i,p-2)$ to $(j',p+1)$ on $\mu_{(i,p)}(\ssq{s})$ disappeared by the new arrow from $(j',p+1)$ to $(i,p-2)$ generated when we apply the mutation $\mu_{(j',p-1)}$ to $( \mu_{(i,p)}(\ssq{s}))$. In fact, one may observe that
\begin{align*}
\mu_{(i,p)} \circ \mu_{(j',p-1)}  =   \mu_{(j',p-1)}\circ \mu_{(i,p)} \,\, \text{  on  } \,\, \ssq{s},  
\end{align*}
and the arrows among $(i,p)$, $(i,p-2)$, $(j',p+1)$ and $(j',p-1)$  in  $\mu_{(j',p-1)}\circ \mu_{(i,p)} (\ssq{s})$  are reversed.  
Furthermore, one may generalize the above as follows:
\begin{align*}
\mu_{(i,p)} \circ \mu_{(j,s)}  =   \mu_{(j,s)}\circ \mu_{(i,p)} \,\, \text{ on } \,\, \mu_{(i_k,p_k)} \circ  \cdots  \circ \mu_{(i_r,p_r)} (\ssq{s})
\end{align*}
for $(i,p),(j,s) \in \ang{k,s}^+ \setminus \{ (i_k,p_k),(i_{k+1},p_{k+1}) , \ldots, (i_r,p_r) \} $, which proves our assertion.
\end{proof}

For $s \in \Z$ and $k \in  \Z_{\le s} \sqcup \{ -\infty \}$, put
\begin{equation*}
	\Theta_s(\ang{k,s}^\pm) \seteq \mu_{\ang{k,s}^\pm} (\ssq{s}), \quad
	\Omega_s(\ang{k,s}^\pm) \seteq \mu_{\ang{k,s}^\pm} (\rssq{s}).
\end{equation*}

\begin{lemma} \label{lem: reverse}
We have  
\begin{equation*}
\Omega_s(\ang{-\infty,s}^\pm)  \simeq  \Theta_s  \,\,\, \text{ and } \,\,\, \Theta_s(\ang{-\infty,s}^\pm)  \simeq  \Omega_s \quad \text{as valued quivers}.
\end{equation*}
\end{lemma}

\begin{proof}
We only prove the second isomorphism for $\ang{k,s}^+$ since the proof of the other cases is almost identical.
In the proof of Lemma~\ref{lem: mu equiv start pm T},
we have seen that a mutation $\mu_{(i,p)}$ for $(i,p) \in \ang{k,s}^+$
generates arrows between vertices in $\ang{k,s}^-$ and then they disappear in the course of the mutations $\mu_{(j,p')}$'s for $(j,p') \in \ang{k,s}^+$ located near $(i,p)$.
Moreover, the arrows adjacent to $(i,p)$ are reversed during the mutations.
Hence we have $\Theta_s(\ang{-\infty,s}^+)  \simeq  \Omega_s$.
\end{proof}

\begin{proposition} \label{prop: mu=T-sys T}
Every mutation $\mu_{(i,p)}$ in $\mu_{\ang{-\infty,s}^\pm}$ on the cluster  $ \{ {}^s\fraku_{i,p} \}$ corresponds to the quantum folded T-system in Theorem~\ref{thm: quantum folded}.
Furthermore, each mutation $\mu_{(i,p)}$ on ${}^s\fraku_{i,p}$ in $\mu_{\ang{-\infty,s}^\pm}$  corresponds to $\sfT_{\pm 2}$.  
\end{proposition}

\begin{proof} For $(i,p) \ne (j,t) \in \ang{-\infty,s}^\pm$,
recall that the mutation $\mu_{(i,p)}$ does not affect the arrows adjacent to $(j,s)$. Thus it suffices to consider $(i,p)$ and vertices connected to $(i,p)$ by arrows. Assume first that $(i,p) \in \ang{-\infty,s}^+$. Then by replacing vertices in $\ssq{s}$ with $\lhs\fraku_{k,q}$'s, we have the following:
\begin{align*}
 \scalebox{0.6}{\raisebox{2.5em}{     \xymatrix@!C=13mm@R=8mm{
 & F_q(\um^{(k)}[a+1,b+1]) & && && \\
 F_q(\um^{(i)}[a,b+2] ) \ar[rr] && F_q(\um^{(i)}[a,b] )  \ar@{->}[dl]|{ \ulcorner -\sfc_{i,j},\sfc_{j,i} \lrcorner} \ar@{->}[ul]|{ \ulcorner -\sfc_{i,k},\sfc_{k,i} \lrcorner} && \ar[ll] F_q(\um^{(i)}[a+2,b])      \\ & F_q(\um^{(j)}[a+1,b+1])   }}}
\quad \raisebox{-1pc}{\text{for $j,k$ with $d(i,j),d(i,k)\le1$}}
\end{align*}
where $\lhs\fraku_{i,p} = F_q(\um^{(i)}[a,b] )$. Note that $\lhs\fraku_{k,q}$ for $(k,q) \in \ang{-\infty,s}^-$ never mutate by $\mu_{\ang{-\infty,s}^+}$. Hence the mutation rule for cluster variables can be expressed as
\begin{align*}
F_q(\um^{(i)}[a,b]) *  \mu_{(i,p)} (F_q(\um^{(i)}[a,b] )) &= q^{\al(i, (b+2-a)/2 )}  F_q(\um^{(i)}[a+2,b]) \cdot F_q(\um^{(i)}[a,b+2])  \\ &  \qquad \qquad +  q^{\ga(i, (b+2-a)/2 )} \prod_{j; d(i,j)=1} F_q(\um^{(j)}[a+1,b+1])^{-\sfc_{j,i}}.
\end{align*}
Here $q^{\al(i, (b+2-a)/2 )}$ and $q^{\ga(i, (b+2-a)/2 )} $ are computed by bar-invariance. Hence,  as in Proposition~\ref{prop: m times mu}, and the above equation coincides  
 with the formula in Theorem~\ref{thm: quantum folded}. Thus we have
\begin{align*}
\mu_{(i,p)} (F_q(\um^{(i)}[a,b] )) = F_q(\um^{(i)}[a+2,b+2]).
\end{align*}
Thus the assertion for $\ang{-\infty,s}^+$ follows.

Similarly,  the arrows adjacent to $(i,p)$ for $(i,p) \in \ang{-\infty,s}^-$ can be depicted as follows:
$$
 \scalebox{0.6}{\raisebox{2.5em}{     \xymatrix@!C=13mm@R=8mm{
  && & F_q(\um^{(k)}[a-1,b-1]) & && \\
     F_q(\um^{(i)}[a-2,b] )  && \ar[ll] F_q(\um^{(i)}[a,b] ) \ar@{<-}[ur]|{ \ulcorner -\sfc_{k,i},\sfc_{i,k} \lrcorner}  \ar[rr] && F_q(\um^{(i)}[a,b-2])      \\ &&& F_q(\um^{(j)}[a-1,b-1]) \ar@{->}[ul]|{ \ulcorner -\sfc_{j,i},\sfc_{i,j} \lrcorner}   }}}
\quad \raisebox{-1pc}{\text{for $j,k$ with $d(i,j),d(i,k)\le1$.}}
$$
Then as in $\ang{-\infty,s}^+$, we can conclude that
\begin{align*}
\mu_{(i,p)} (F_q(\um^{(i)}[a,b] )) = F_q(\um^{(i)}[a-2,b-2]),
\end{align*}
which proves our assertion.
\end{proof}

\begin{example} \label{ex: mutated quiver by mu}
By replacing vertices $(i,p)$ in $ \mu_{\ang{-\infty,0}^\pm}(\ssq{0})$ with $\mu_{\ang{-\infty,0}^\pm}(\lh{0}\fraku_{i,p})$ obtained from Example~\ref{ex: F_q quiver}~\ref{it: B_3 F_q}, we have the following
by Lemma~\ref{lem: reverse} and Proposition~\ref{prop: mu=T-sys T}:
\smallskip

\noindent
{\it Case 1}. $\mu_{\ang{-\infty,0}^+}(\ssq{0})$.
\begin{align*}
\raisebox{2em}{ \scalebox{0.5}{
\xymatrix@!C=10mm@R=6mm{
\cdots\ar[rr] &&F_q(m^{(1)}[-2,2] )  \ar[dl]    &&F_q(m^{(1)}[0,2] ) \ar[ll]\ar[rr] &&F_q(m^{(1)}[0,0] ) \ar[dl]    \\
\cdots&F_q(m^{(2)}[-1,3] ) \ar[l]    \ar[rr]  &&F_q(m^{(2)}[-1,1] )\ar[dr]|{2 \lrcorner} \ar[ur] && F_q(m^{(2)}[1,1] ) \ar[ll]  \\
\cdots\ar[rr] &&F_q(m^{(3)}[-2,2] )   \ar@{=>}[ul]    &&F_q(m^{(3)}[0,2] ) \ar[ll]\ar[rr]&&F_q(m^{(3)}[0,0] )   \ar@{=>}[ul]
}}}  \hspace{-2ex} \simeq
\raisebox{2em}{ \scalebox{0.5}{
\xymatrix@!C=10mm@R=6mm{
\cdots\ar[r] &F_q(m^{(1)}[-2,2] )  \ar[dr]    &&F_q(m^{(1)}[0,2] ) \ar[ll]\ar[rr] &&F_q(m^{(1)}[0,0] ) \ar[dr]    \\
\cdots&&F_q(m^{(2)}[-1,3] ) \ar[ll]    \ar[rr]  &&F_q(m^{(2)}[-1,1] )\ar[dl]|{\llcorner 2 } \ar[ul] && F_q(m^{(2)}[1,1] ) \ar[ll]  \\
\cdots\ar[r] &F_q(m^{(3)}[-2,2] )   \ar@{=>}[ur]    &&F_q(m^{(3)}[0,2] ) \ar[ll]\ar[rr]&&F_q(m^{(3)}[0,0] )   \ar@{=>}[ur]
}}}
\end{align*}
where the parameters of quantum cluster variables located at vertices that are vertically sink and horizontally source are shifted by $2$, and the orientation of all arrows is reversed.

\noindent
{\it Case 2}. $\mu_{\ang{-\infty,0}^-}(\ssq{0})$.

\begin{align*}
\raisebox{2em}{ \scalebox{0.5}{
\xymatrix@!C=10mm@R=6mm{
\cdots\ar[rr] &&F_q(m^{(1)}[-4,0] )  \ar[dl]    &&F_q(m^{(1)}[-2,0] ) \ar[ll]\ar[rr] &&F_q(m^{(1)}[-2,-2] ) \ar[dl]    \\
\cdots&F_q(m^{(2)}[-3,1] ) \ar[l]    \ar[rr]  &&F_q(m^{(2)}[-3,-1] )\ar[dr]|{2 \lrcorner} \ar[ur] && F_q(m^{(2)}[-1,-1] ) \ar[ll]  \\
\cdots\ar[rr] &&F_q(m^{(3)}[-4,0] )   \ar@{=>}[ul]    &&F_q(m^{(3)}[-2,0] ) \ar[ll]\ar[rr]&&F_q(m^{(3)}[-2,-2] )   \ar@{=>}[ul]
}}} \hspace{-2ex} \simeq
\raisebox{2em}{ \scalebox{0.5}{
\xymatrix@!C=10mm@R=6mm{
\cdots\ar[r] &F_q(m^{(1)}[-4,0] )  \ar[dr]    &&F_q(m^{(1)}[-2,0] ) \ar[ll]\ar[rr] &&F_q(m^{(1)}[-2,-2] ) \ar[dr]    \\
\cdots&&F_q(m^{(2)}[-3,1] ) \ar[ll]    \ar[rr]  &&F_q(m^{(2)}[-3,-1] )\ar[dl]|{ \llcorner 2} \ar[ul] && F_q(m^{(2)}[-1,-1] ) \ar[ll]  \\
\cdots\ar[r] &F_q(m^{(3)}[-4,0] )   \ar@{=>}[ur]    &&F_q(m^{(3)}[-2,0] ) \ar[ll]\ar[rr]&&F_q(m^{(3)}[-2,-2] )   \ar@{=>}[ur]
}}}
\end{align*}
where the parameters of quantum cluster variables located at vertices that are vertically sink and horizontally source are shifted by $-2$, and the orientation of all arrows is reversed.
\smallskip

\noindent
Thus we can conclude that
\begin{equation*}
 \mu_{\ang{-\infty,0}^+}(\frakS_0) \simeq  \frakS_1 \qquad \text{ and } \qquad
 \mu_{\ang{-\infty,0}^-}(\frakS_0) \simeq    \frakS_{-1}.
\end{equation*}
 \end{example}

Following Example \ref{ex: mutated quiver by mu},
it is straightforward to check the following proposition.

\begin{proposition} \label{prop: mu Ss = Ss+1}
For $s \in \Z$, we have
$$   \mu_{\ang{-\infty,s}^{\pm}}(\frakS_s) \simeq  \frakS_{s\pm1}.$$
\end{proposition}

For $s \in \Z$, put
\begin{equation*}
	\frakU \seteq \mu_{\ang{-\infty,s}^\pm}(\{ {}^s\fraku_{i,p} \}).
\end{equation*}

\begin{proposition} \label{prop: mu=T-sys O}
Every mutation $\mu_{(i,p)}$ in $\mu_{\ang{-\infty,s}^\pm}$ on the cluster $\frakU$  
 corresponds to the quantum folded T-system in Theorem~\ref{thm: quantum folded}.
Furthermore, each mutation $\mu_{(i,p)}$ on the quantum cluster variable at $(i,p)$ in $\mu_{\ang{-\infty,s}^\pm}$  corresponds to $\sfT_{\mp 2}$.
\end{proposition}

\begin{proof}
Set
$$ \{ {}^s\frakz^+ _{i,p} \} \seteq   \mu_{\ang{-\infty,s}^+}(\{ {}^s\fraku_{i,p} \})
\quad \text{ and }\quad
 \{ {}^s\frakz^-_{i,p} \} \seteq   \mu_{\ang{-\infty,s}^-}(\{ {}^s\fraku_{i,p} \}).
$$
In this proof, we only consider the case of $\{ {}^s\frakz^+ _{i,p} \}$ since the proof of $\{ {}^s\frakz^- _{i,p} \}$ is parallel.
Let $(i,p) \in \ang{-\infty,s}^\pm$.
\smallskip

\noindent
{\it Case 1}.
$(i,p) \in \ang{-\infty,s}^+$. By replacing vertices in $\rssq{s}$ with $\lhs\frakz_{k,q}^+$'s, we have the following:
$$
 \scalebox{0.6}{\raisebox{2.5em}{     \xymatrix@!C=13mm@R=8mm{
 & F_q(\um^{(k)}[a+1,b+1]) \ar@{->}[dr]|{ \ulcorner -\sfc_{k,i},\sfc_{i,k} \lrcorner} & && && \\
 F_q(\um^{(i)}[a,b+2] ) && F_q(\um^{(i)}[a+2,b+2] )\ar[rr] \ar[ll]    && F_q(\um^{(i)}[a+2,b])      \\ & F_q(\um^{(j)}[a+1,b+1]) \ar@{->}[ur]|{ \ulcorner -\sfc_{j,i},\sfc_{i,j} \lrcorner}   }}}
\quad \raisebox{-1pc}{\text{for $j,k$ with $d(i,j),d(i,k)\le1$}}
$$
where $\lhs\frakz_{i,p}^+ = F_q(\um^{(i)}[a+2,b+2] )$.  Hence the mutation rule for quantum cluster variables can be expressed as
\begin{align*}
 \mu_{(i,p)} (F_q(\um^{(i)}[a+2,b+2] ))  * F_q(\um^{(i)}[a+2,b+2])  &= q^{\al(i, (b+2-a)/2 )} F_q(\um^{(i)}[a+2,b]) \cdot  F_q(\um^{(i)}[a,b+2])  \\ &  \hspace{-2ex} +  q^{\ga(i, (b+2-a)/2 )} \prod_{j; d(i,j)=1} F_q(\um^{(j)}[a+1,b+1])^{-\sfc_{j,i}}.
\end{align*}
Thus we have
$$
\mu_{(i,p)} (F_q(\um^{(i)}[a+2,b+2] )) = F_q(\um^{(i)}[a,b]).
$$
\smallskip

\noindent
{\it Case 2}.
 $(i,p) \in \ang{-\infty,s}^-$. The arrows adjacent to $(i,p)$ for $(i,p) \in \ang{-\infty,s}^-$ are depicted as follows:
$$
 \scalebox{0.6}{\raisebox{2.5em}{     \xymatrix@!C=13mm@R=8mm{
  && & F_q(\um^{(k)}[a+1,b+1]) & && \\
     F_q(\um^{(i)}[a,b+2] )  \ar[rr] && F_q(\um^{(i)}[a,b] ) \ar@{->}[ur]|{ \ulcorner -\sfc_{i,k},\sfc_{k,i} \lrcorner}  && F_q(\um^{(i)}[a+2,b])  \ar[ll]    \\ &&& F_q(\um^{(j)}[a+1,b+1]) \ar@{<-}[ul]|{ \ulcorner -\sfc_{i,j},\sfc_{j,i} \lrcorner}   }}}
\quad \raisebox{-1pc}{\text{for $j,k$ with $d(i,j),d(i,k)\le1$}}
$$
where $\lhs\frakz_{i,p}^+ = F_q(\um^{(i)}[a,b] )$.
Then as in {\it Case 1}, we have
\[
\mu_{(i,p)} (F_q(\um^{(i)}[a,b] )) = F_q(\um^{(i)}[a+2,b+2]). \qedhere
\]
\end{proof}

Now, we are ready to prove Proposition~\ref{prop: sub 2}.

\begin{proof} [Proof of Proposition~\ref{prop: sub 2}]
Let us define
\begin{align*}
\bmup \seteq \mu_{\ang{-\infty,s+1}^-} \circ \mu_{\ang{-\infty,s}^+}  \quad \text{ and } \quad   \bmum \seteq \mu_{\ang{-\infty,s-1}^+} \circ \mu_{\ang{-\infty,s}^-}.
\end{align*}
It follows from Propositions~\ref{prop: mu=T-sys T}, \ref{prop: mu Ss = Ss+1} and \ref{prop: mu=T-sys O} that
\begin{align*}
\bmup (\frakS_s) \simeq \frakS_{s+2} \qquad \text{ and } \qquad  \bmum (\frakS_s) \simeq \frakS_{s-2}.
\end{align*}
By applying $\bmup$ repeatedly,
we obtain $F_q(\sfX_{i,p})$ for $(i,p) \in \tDynkinf_0$ with $p \ge s$ as a cluster variable of $\lhs\scrA_q(\g)$.
Similarly, we obtain every $F_q(\sfX_{i,p})$ for $(i,p) \in \tDynkinf_0$ with $p \le s$ as a cluster variable of $\lhs\scrA_q(\g)$  by using the repetition of $\bmum$.
Thus the cluster algebra $\lhs\scrA_q(\g)$ contains every $F_q(\sfX_{i,p})$ associated to $\tDynkinf_0$ as its cluster variables.
\end{proof}

\begin{conjecture} \label{conjecture 3}
Let $s$ be an arbitrary integer.
If $F_q(\um^{(j)}[a,b]) \in \frakK_q(\sfg)$ $q$-commutes with $\lhs\fraku_{i,p}$ for all $\lhs\fraku_{i,p} \in \frakS_s$, then there exists $(j,l) \in \tbDynkinf_0$ such that
\begin{equation*}
\lhs\fraku_{j,l} = F_q(\um^{(j)}[a,b]).
\end{equation*}
\end{conjecture}

\bibliographystyle{amsplain}
\bibliography{ref}{}

\end{document}